\documentclass[a4paper,11pt]{article}
\usepackage{amsmath,amssymb,amsthm}
\usepackage{mathrsfs}
\usepackage{a4wide}
\usepackage{xfp}

\usepackage{graphicx,xcolor}
\usepackage{tikz}
\usetikzlibrary{arrows.meta}

\usepackage[hang,small,bf]{caption}
\usepackage[subrefformat=parens]{subcaption}
\usepackage{here}
\usepackage{url}
\usepackage{array,booktabs}
\usepackage{cite}
\numberwithin{equation}{section}

\usepackage[pdfencoding=auto, psdextra, backref=page]{hyperref}
\hypersetup{
    plainpages=false,
    unicode=true,
    pdftoolbar=true,
    pdfmenubar=true,
    pdffitwindow=true,
    pdftitle={Learning the Willmore flow with a lightweight neural operator},
    pdfauthor={Elie Bretin, Roland Denis, Tokuhiro Eto, Simon Masnou},
    pdfsubject={Numerical approximation of the Willmore flow by phase field neural operators},
    pdfkeywords={Willmore flow; phase field; neural operator; rotationally invariant kernels; non-oriented interfaces; surface reconstruction},
    pdfnewwindow=true,      
    colorlinks=true,       
    linkcolor=blue,          
    citecolor=blue,        
    filecolor=blue,      
    urlcolor=magenta,          
    urlbordercolor=0 1 1
}

\newtheorem{rem}{Remark}[section]

\newcommand{\Section}[1]{Section~\ref{#1}}
\newcommand{\Figure}[1]{Figure~\ref{#1}}

\newcommand{\Table}[1]{Table~\ref{#1}}

\DeclareMathOperator*{\dist}{dist}

\renewcommand{\epsilon}{\varepsilon}

\definecolor{Blue}{rgb}{0.,0.,1}

\def\R{\mathbb{R}}

\def\bv#1{\mbox{\boldmath{$#1$}}}    

\newcommand{\dL}[1]{dx}
\newcommand{\dH}[1]{d\mathscr{H}^{#1}}

\newcommand{\mH}[1]{\mathscr{H}^{#1}}

\newcommand{\wnn}[2]{W^{NN}_{#1,\fpeval{2 * #2}}[q]}
\newcommand{\wnnSimple}[2]{W^{NN}_{#1,\fpeval{2 * #2}}}
\newcommand{\wpnn}[2]{W^{NN}_{iso,#1,\fpeval{2 * #2}}[q]}
\newcommand{\wpnnSimple}[2]{W^{NN}_{iso,#1,\fpeval{2 * #2}}}
\newcommand{\wnnprim}[2]{W^{NN}_{#1,\fpeval{2 * #2}}[q']}
\newcommand{\wpnnprim}[2]{W^{NN}_{iso,#1,\fpeval{2 * #2}}[q']}
\newcommand{\snn}[2]{S^{NN}_{#1,#2}[q]}
\newcommand{\snnSimple}[2]{S^{NN}_{#1,#2}}

\newcommand{\spnnSimple}[2]{{S}^{NN}_{iso,#1,#2}}

\newcommand{\snnprim}[2]{S^{NN}_{#1,#2}[q']}
\newcommand{\ts}[0]{\tau}
\newcommand{\tts}[0]{\widetilde{\tau}}

\newcommand{\wenergy}[1]{\mathcal{W}_{#1}}
\newcommand{\cenergy}[1]{\mathcal{P}_{#1}}
\newcommand{\loss}[1]{\mathcal{L}_{#1}}

\newcommand{\reaction}[0]{\mathcal{R}}
\newcommand{\diffusion}[0]{\mathcal{D}}

\newcommand{\beltrami}[0]{\Delta_{\Gamma}}
\newcommand{\operatorprim}[0]{\Psi}
\newcommand{\ft}[2]{\mathcal{F}_{#1}\left[#2\right]}
\newcommand{\ift}[2]{\mathcal{F}^{-1}_{#1}\left[#2\right]}
\newcommand{\radial}[1]{\phi_{#1}}
\newcommand{\relu}[0]{\mathrm{ReLU}}
\newcommand{\Krad}[0]{K}
\newcommand{\sdist}[1]{d(\cdot,{#1})}

\definecolor{bleuperso}{rgb}{0.1,0.1,0.9}

\newif\ifdraftnotes
\draftnotestrue          

\title{Learning the Willmore flow with a lightweight neural operator}
\author{
    Elie Bretin\thanks{INSA Lyon, CNRS, Ecole Centrale de Lyon, Universit\'e Claude Bernard Lyon 1, Universit\'e Jean Monnet, ICJ UMR5208, 69622 Villeurbanne, France. Email: elie.bretin@insa-lyon.fr}
    \and
    Roland Denis\thanks{CNRS, Universit\'{e} Claude Bernard Lyon 1, Centrale Lyon, INSA Lyon, Universit\'{e} Jean Monnet, ICJ UMR5208, 69622 Villeurbanne, France. Email: denis@math.univ-lyon1.fr}
    \and
    Tokuhiro Eto\thanks{Universit\'{e} Claude Bernard Lyon 1, CNRS, Centrale Lyon, INSA Lyon, Universit\'{e} Jean Monnet, ICJ UMR5208, 69622 Villeurbanne, France. Email: eto@math.univ-lyon1.fr}
    \and
    Simon Masnou\thanks{Universit\'{e} Claude Bernard Lyon 1, CNRS, Centrale Lyon, INSA Lyon, Universit\'{e} Jean Monnet, ICJ UMR5208, 69622 Villeurbanne, France. Email: masnou@math.univ-lyon1.fr}
    }
\date{}

\begin{document}

\maketitle

\begin{abstract}

We extend to the Willmore flow the phase-field neural operator approach proposed in~\cite{BretinDenisMasnouTerii2022} for mean curvature motion in two and three dimensions.
The resulting lightweight, single-time-step architecture is motivated by a Bence--Merriman--Osher-type expansion of the Allen--Cahn semigroup,
and is trained using trajectories generated by a hybrid minimizing movement scheme. In our experiments with non-oriented interfaces,
unconstrained convolution kernels produce anisotropic dynamics,
whereas radial kernels enhance the evolution process. We also construct three-dimensional kernels from trained two-dimensional kernels via their radial Fourier profiles,
which provides an effective initialization for the non-oriented 3D model.
Finally, we show applications to curve and surface reconstruction from unoriented point clouds, including reconstruction using the transferred model without the need for additional 3D training.
\end{abstract}



\section{Introduction}

The Willmore energy measures the bending of an interface and arises in
mechanics, biology, image processing and geometry processing. For a smooth
bounded domain $\Omega\subset\R^d\,(d = 2,\,3)$, it is defined by
\begin{equation}\label{WillmoreFunctional}
    \mathcal W(\Omega)=\frac12\int_{\partial\Omega}H^2\, \dH{d-1},
\end{equation}
where $H=\operatorname{div}_{\partial\Omega}\nu$ is the scalar mean curvature with
respect to the outer unit normal $\nu$ and $\mH{d-1}$ is the \((d-1)\)-dimensional Hausdorff measure.
Thus, spheres have positive mean curvature. The associated $L^2$-gradient flow has outer normal velocity
\begin{equation}\label{eq:willmore}
    V=\beltrami H-\frac12 H^3+H\|A\|^2,
\end{equation}
where $\beltrami$ is the Laplace--Beltrami operator and $A$ is the second
fundamental form, with $\|A\|^2=\sum_{i=1}^{d-1}k_i^2$ in terms of the
principal curvatures $k_i\,(1\leq i\leq d-1)$.
The analysis of this fourth-order geometric flow is
challenging even for smooth interfaces
\cite{dziuk2002evolution,kuwert2002gradient,kuwert2004removability,kuwert2012willmore}.
Its numerical approximation must also resolve curvature variations while
handling possible changes of topology.

\subsection{Phase-field formulation}\label{sec:intro_pf}

Phase-field methods are widely used in applied mathematics and physics to model interfaces and phase transitions. They represent an interface implicitly by a smooth transition of width \(\varepsilon>0\) between two phases. For a scalar field \(u\), the Modica--Mortola energy is given by
\begin{equation}
    \label{eq:CHenergy-et-potential}
    \cenergy{\varepsilon}(u)=\int_{\R^d}
    \left(\frac{\varepsilon}{2}|\nabla u|^2+
    \frac1\varepsilon W(u)\right)\,\dL{d},
    \qquad W(u)=\frac12u^2(1-u)^2.
\end{equation}
This energy approximates
\(c_W\) times the perimeter, with
\(
    c_W=\int_0^1\sqrt{2W(s)}\,ds,
\)
in the sense of the Modica--Mortola \(\Gamma\)-convergence
\cite{Modica1977,Modica1977a}. Its rescaled \(L^2\)-gradient flow is described by the Allen--Cahn equation
\[
    \partial_tu
    =\Delta u-\frac1{\varepsilon^2}W'(u),
\]
whose sharp-interface limit is motion by mean curvature
\cite{Chen1992,DeMottoniSchatzman1995,bellettini1995quasi}.

Assuming well-prepared initial data, let \(\Omega(t)\) denote the limiting region corresponding to the phase
\(u\to1\) as \(\varepsilon\to0\), and let
\(
    \Gamma(t)=\partial\Omega(t)
\)
be the corresponding sharp interface. Assuming that \(\Gamma(t)\) is smooth, the leading-order inner profile of the phase field is
\(
    u(x,t)\simeq
    q\!\left(\frac{\sdist{\Omega(t)}(x)}{\varepsilon}\right),
\)
where \(\sdist{\Omega(t)}\) denotes the signed distance to \(\Gamma(t)\), taken negative inside \(\Omega(t)\). The one-dimensional profile \(q\) satisfies
\[
    q''-W'(q)=0,\qquad q(0)=\tfrac12,\qquad
    q(-\infty)=1,\quad q(+\infty)=0.
\]
For the potential \(W\) in \eqref{eq:CHenergy-et-potential}, we have
\[
    q(z)=\frac{1}{1+e^z},
    \qquad
    q'(z)=-\frac{1}{4\cosh^2(z/2)}.
\]

The classical diffuse Willmore energy is given by
\begin{equation}\label{eq:classical-willmore-energy}
    \wenergy{\varepsilon}(u)=\frac1{2\varepsilon}\int_{\R^d}
    \left(\varepsilon\Delta u-\frac{W'(u)}{\varepsilon}\right)^2\,dx.
\end{equation}
In dimensions two and three, the corresponding $\Gamma$-convergence results, under
bounded diffuse perimeter, relate this functional to $c_W\mathcal W$ on
smooth limiting interfaces~\cite{belpao93,Roger_schatzle_2006,naga-tone}; see also 
\cite{MR3383330} for phase-field approximations of the energy and flow.
Its rescaled $L^2$-gradient flow is governed by
\begin{equation}\label{flow:classic}
\begin{cases}
    \varepsilon^2\partial_tu=\Delta\mu-\frac1{\varepsilon^{2}}W''(u)\mu,\\
    \mu=W'(u)-\varepsilon^2\Delta u.
\end{cases}
\end{equation}
Convergence to the smooth Willmore flow has been studied through formal
asymptotic and rigorous analysis
\cite{MR1757511,MR2377282,MR3383330,rigorousWillmoreconvergence}.
Numerical methods include finite elements, finite differences and Fourier
spectral schemes~\cite{du2004phase,du2006convergence,esedoglu2012colliding,MR3383330}.
Resolving a thin interface together with fourth-order dynamics remains
computationally demanding.

\begin{figure}[!htbp]
    \centering
    \begin{subfigure}[b]{0.40\textwidth}
        \centering
        \includegraphics[width=0.69\textwidth,height=0.60\textwidth]{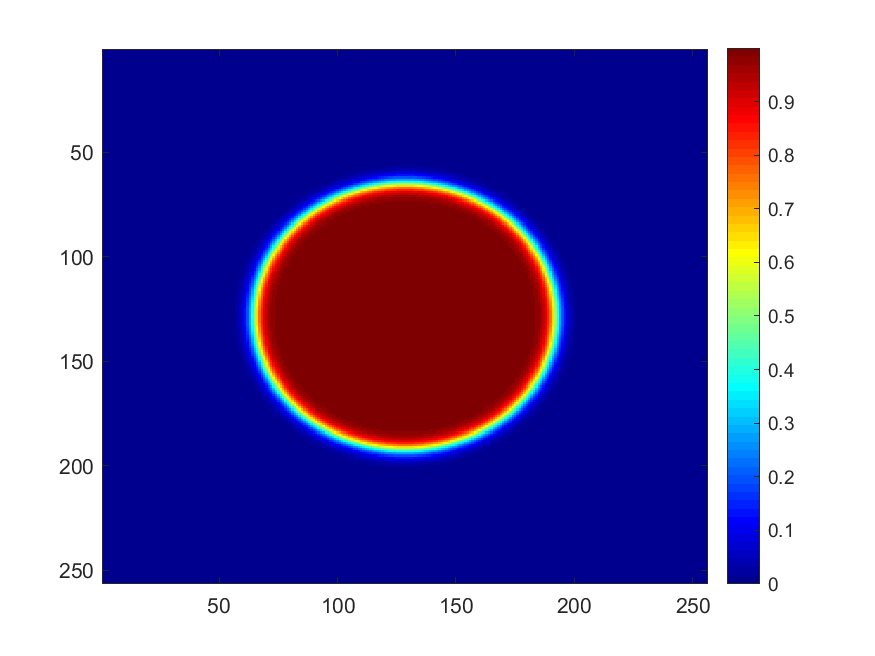}
        \caption{Oriented: $q(\sdist{\Omega}/\varepsilon)$.}
    \end{subfigure}
    \hspace{0.04\textwidth}
    \begin{subfigure}[b]{0.40\textwidth}
        \centering
        \includegraphics[width=0.69\textwidth,height=0.60\linewidth]{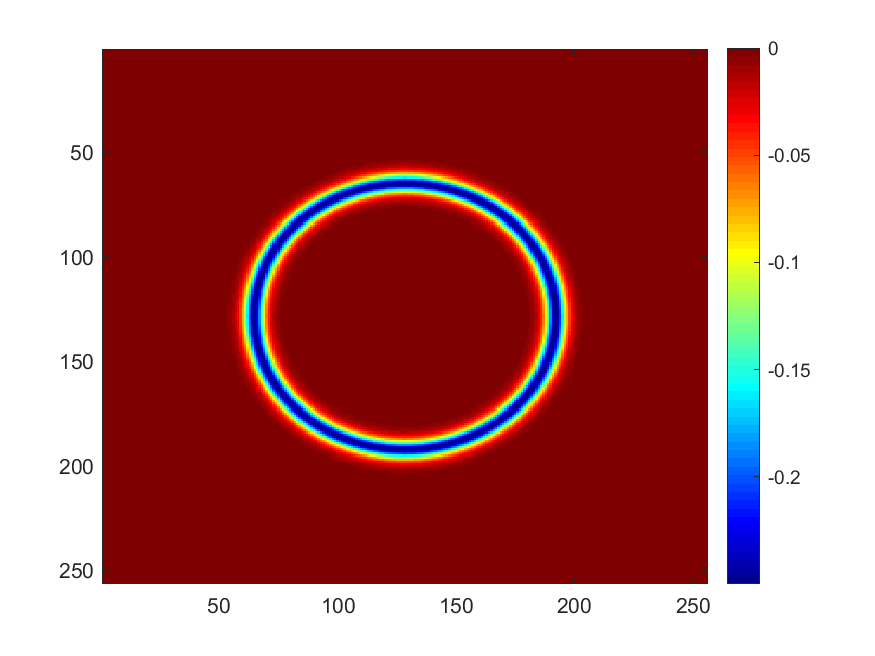}
        \caption{Non-oriented: $q'(\dist(\cdot,\Gamma)/\varepsilon)$.}
    \end{subfigure}
    \caption{Phase-field representations of a circle
    $\Gamma=\partial\Omega$. Left: the profile $q$ is close to $1$ inside
    the disk and $0$ outside. Right: the profile $q'$ is minimal on the
    interface, where it equals $-1/4$, and tends to $0$ on both sides.}
    \label{fig_phase_field_representation}
\end{figure}

Alongside this oriented description, we use the non-oriented encoding
$u=q'(\dist(\cdot,\Gamma)/\varepsilon)$ of an interface $\Gamma$, with $\dist(\cdot,\Gamma)$ the distance function to $\Gamma$.
Since $q'$ is even, negative, and minimal on the interface, this encoding
requires no distinction between an inside and an outside.
It has been used for neural approximations of
mean curvature motion~\cite{BretinDenisMasnouTerii2022}, it has motivated a
variational phase-field model~\cite{bretin2025cahnhilliardwillmorephasefieldmodel},
and been applied to Plateau problems~\cite{bonnivard2025phasefieldapproximationplateaus}.
In our learning framework, the corresponding training data are obtained
by post-processing oriented trajectories; see
\Section{sec:no_willmore_nn_2d}. This representation allows us to explore
configurations with junctions and surfaces with boundary, as well as
reconstruction from unoriented point clouds.
\Figure{fig_phase_field_representation} illustrates the difference between
the two representations, oriented and non-oriented, for a circular interface.

\subsection{Related numerical and learning approaches}\label{sec:related_learning}
Many recent studies have explored the use of neural networks for the
numerical approximation of phase-field solutions. Without attempting
an exhaustive review, we distinguish three broad, partly overlapping
families of approaches.

Physics-informed approaches (PINNs) do not require precomputed solutions for training.
In their standard formulation, the network represents the
solution itself, and the PDE
residual together with the initial and boundary conditions defines
the training loss
\cite{raissi2017physicsinformeddeeplearning,MR3881695,MR3874585};
see \cite{cuomo2022scientific,deryck2024numerical} for recent overviews.
For Allen--Cahn- or Cahn--Hilliard-type phase-field equations,
minimizing this loss is challenging because of the stiffness associated
with the double-well potential and the thinness of the diffuse
interface. Several variants have been proposed to address these
difficulties: adaptive or sequential training strategies
\cite{wight2020solving,mattey2022novel}, losses based on the
energy-dissipation structure of the flow rather than on the pointwise
PDE residual in the spirit of the Deep Ritz method
\cite{yu2018deep,kutuk2024energy,geng2024deep}, and
physics-informed operator variants that combine this idea with
DeepONet- or finite-element-based architectures
\cite{li2023phasefield,yamazaki2025finite,zhang2024energydissipative}.

Operator learning focuses on approximating the solution operator
of the PDE, so that evaluating the learned operator for
a new initial condition requires only a forward pass.
DeepONet~\cite{Lu_2021} and the Fourier neural
operator~\cite{li2021fourier} are two prominent architectures
of this kind. These architectures are designed for broad classes
of PDEs rather than tailored to a specific evolution equation, they can require large numbers of trainable parameters and substantial amounts of training data.

The approach we follow here belongs to a third, intermediate family,
in which the architecture is derived from a discretization scheme for the specific equation considered. The resulting network approximates a solution
operator, as in the second family, but its structure incorporates
information that a generic architecture would otherwise have to learn
from data. In our previous work \cite{BretinDenisMasnouTerii2022},
we explored this strategy for mean curvature flow by training a lightweight
network to approximate its evolution over a single time step.

For the Willmore flow, minimizing movement schemes have also been
developed in both parametric \cite{Perl_2013} and phase-field
\cite{FrRuWi11} settings. These schemes characterize each time
step through a minimization problem. Within this framework,
the Rumpf--Sassen--Smoch (RSS) scheme
\cite{Rumpf_Sassen_Smoch25} reduces the computational burden
by incorporating a learned mean curvature operator into the
minimization procedure, while still requiring the solution of an optimization problem at each time step.

\subsection{Contributions and outline}

This work extends our previous approach for learning mean curvature flow
\cite{BretinDenisMasnouTerii2022} to the Willmore flow. We train a lightweight
neural operator on trajectories generated by the Rumpf--Sassen--Smoch (RSS) scheme~\cite{Rumpf_Sassen_Smoch25}. Once trained,
it replaces the optimization required at each RSS time step by a single
network evaluation. We consider both oriented and non-oriented interfaces
in two and three dimensions.

Our contributions address three aspects of this extension. First, we
construct a compact architecture with parallel branches motivated by an
expansion of the Allen--Cahn semigroup. The training data include shapes
with nonconstant curvature, so that the network encounters the variations
entering $\beltrami H$ in \eqref{eq:willmore}.
Second, we show that small training errors can coexist with anisotropic
dynamics for non-oriented interfaces, and that imposing radial convolution
kernels remedies this behavior in our experiments.
Third, we transfer the radial Fourier profiles of trained two-dimensional
kernels to initialize the three-dimensional model. This enables training
in our non-oriented 3D experiments, where random initialization was
unsuccessful. The transferred model also provides a useful regularizer
for surface reconstruction from unoriented point clouds, without
additional 3D training.

The paper is organized as follows.
\Section{sec:one-step-NN} presents the two-dimensional architecture and
the RSS-based procedure used to generate training data.
\Section{sec:no_willmore_nn_2d} addresses non-oriented curves, introduces
radial convolution kernels and discusses curve reconstruction from
point clouds. \Section{sec:o_willmore_nn_3d} extends the approach to
oriented surfaces in three dimensions.
Finally, \Section{sec:nonoriented_3d} considers non-oriented surfaces,
including the transfer of trained two-dimensional models and applications
to surface reconstruction.

\section{A one-step neural operator for Willmore flow in 2D}\label{sec:one-step-NN}

\subsection{Mean curvature network building blocks}\label{sec:intro_nn}

We recall the diffusion reaction operators of \cite{BretinDenisMasnouTerii2022}
that form the building blocks of our architecture. The networks
$\snnSimple{\theta}{1}$ and $\snnSimple{\theta}{2}$ are inspired by first-
and second-order splittings of the Allen--Cahn equation, respectively.
Their parameters consist of convolution kernels, pointwise reaction weights
and channel-mixing coefficients.

At the continuous level, a diffusion layer acts on the phase field $u$ by
convolution with an integrable kernel $\mathcal D:\R^d\to\R$:
\begin{equation}\label{eq:convolution_continuous}
    (\mathcal D*u)(x) := \int_{\R^d}\mathcal D(x-y)u(y)\,dy,
    \qquad x\in\R^d.
\end{equation}
On a periodic computational domain, $u$ is understood to be extended
periodically. This continuous formulation is the starting point for the
radial kernel construction in \Section{sec:IsotropicNN_MCF}.

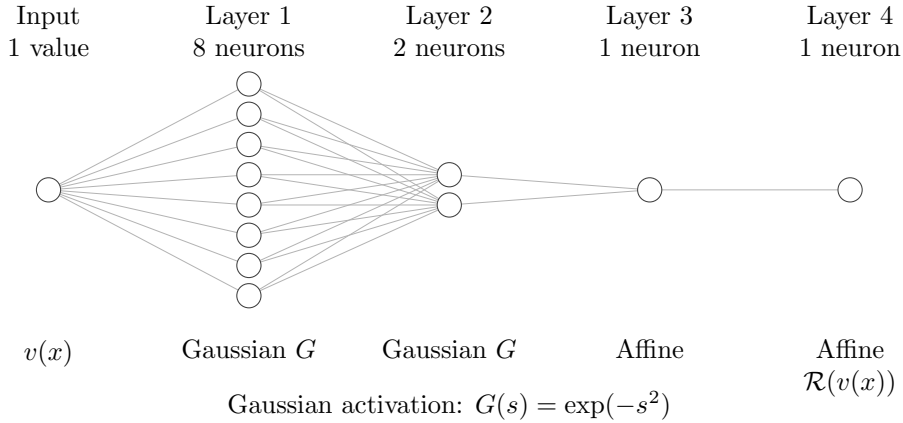
\begin{figure}[H]
    \centering
    \begin{tikzpicture}[font=\small,
        neuron/.style={circle,draw=black!75,fill=white,minimum size=3.2mm,inner sep=0pt},
        connection/.style={draw=black!30,line width=0.35pt}]
        \coordinate (p0-0) at (0.000,0.000);
        \coordinate (p1-0) at (2.650,1.400);
        \coordinate (p1-1) at (2.650,1.000);
        \coordinate (p1-2) at (2.650,0.600);
        \coordinate (p1-3) at (2.650,0.200);
        \coordinate (p1-4) at (2.650,-0.200);
        \coordinate (p1-5) at (2.650,-0.600);
        \coordinate (p1-6) at (2.650,-1.000);
        \coordinate (p1-7) at (2.650,-1.400);
        \coordinate (p2-0) at (5.300,0.200);
        \coordinate (p2-1) at (5.300,-0.200);
        \coordinate (p3-0) at (7.950,0.000);
        \coordinate (p4-0) at (10.600,0.000);
        \draw[connection] (p0-0) -- (p1-0);
        \draw[connection] (p0-0) -- (p1-1);
        \draw[connection] (p0-0) -- (p1-2);
        \draw[connection] (p0-0) -- (p1-3);
        \draw[connection] (p0-0) -- (p1-4);
        \draw[connection] (p0-0) -- (p1-5);
        \draw[connection] (p0-0) -- (p1-6);
        \draw[connection] (p0-0) -- (p1-7);
        \draw[connection] (p1-0) -- (p2-0);
        \draw[connection] (p1-0) -- (p2-1);
        \draw[connection] (p1-1) -- (p2-0);
        \draw[connection] (p1-1) -- (p2-1);
        \draw[connection] (p1-2) -- (p2-0);
        \draw[connection] (p1-2) -- (p2-1);
        \draw[connection] (p1-3) -- (p2-0);
        \draw[connection] (p1-3) -- (p2-1);
        \draw[connection] (p1-4) -- (p2-0);
        \draw[connection] (p1-4) -- (p2-1);
        \draw[connection] (p1-5) -- (p2-0);
        \draw[connection] (p1-5) -- (p2-1);
        \draw[connection] (p1-6) -- (p2-0);
        \draw[connection] (p1-6) -- (p2-1);
        \draw[connection] (p1-7) -- (p2-0);
        \draw[connection] (p1-7) -- (p2-1);
        \draw[connection] (p2-0) -- (p3-0);
        \draw[connection] (p2-1) -- (p3-0);
        \draw[connection] (p3-0) -- (p4-0);
        \node[neuron] at (p0-0) {};
        \node[align=center] at (0.000,2.10) {Input\\1 value};
        \node[align=center,anchor=north] at (0.000,-1.85) {$v(x)$};
        \node[neuron] at (p1-0) {};
        \node[neuron] at (p1-1) {};
        \node[neuron] at (p1-2) {};
        \node[neuron] at (p1-3) {};
        \node[neuron] at (p1-4) {};
        \node[neuron] at (p1-5) {};
        \node[neuron] at (p1-6) {};
        \node[neuron] at (p1-7) {};
        \node[align=center] at (2.650,2.10) {Layer 1\\8 neurons};
        \node[align=center,anchor=north] at (2.650,-1.85) {Gaussian $G$};
        \node[neuron] at (p2-0) {};
        \node[neuron] at (p2-1) {};
        \node[align=center] at (5.300,2.10) {Layer 2\\2 neurons};
        \node[align=center,anchor=north] at (5.300,-1.85) {Gaussian $G$};
        \node[neuron] at (p3-0) {};
        \node[align=center] at (7.950,2.10) {Layer 3\\1 neuron};
        \node[align=center,anchor=north] at (7.950,-1.85) {Affine};
        \node[neuron] at (p4-0) {};
        \node[align=center] at (10.600,2.10) {Layer 4\\1 neuron};
        \node[align=center,anchor=north] at (10.600,-1.85) {Affine\\$\reaction(v(x))$};
        \node at (5.300,-2.85) {Gaussian activation: $G(s)=\exp(-s^2)$};
    \end{tikzpicture}
    \caption{Architecture of the pointwise reaction: four affine layers with
    Gaussian activation after the first two layers ($39$ learnable parameters).
    The same network is applied independently to the scalar value at every
    spatial point $x$, with shared parameters.}
    \label{fig:pointwise_reaction}
\end{figure}

For the numerical implementation in two dimensions, the computational
domain is discretized by a uniform Cartesian
grid with $N$ points in each direction and periodic boundary conditions, and
the discrete phase field function is represented by the array
$u = (u_{i,j})_{i,j=1}^N$. Using the same notation for the discrete stencil,
with quadrature weights absorbed into its coefficients, we introduce
$\mathcal D = ({\mathcal D}_{i,j})_{i,j=-\Krad}^{\Krad}$,
where $\Krad\geq 1$ is an integer. The stencil thus has $2\Krad+1$ points
per coordinate direction. In practice, we use $\Krad=8$ in both two and
three dimensions, giving $17\times17$ and $17\times17\times17$ kernels,
respectively. The associated discrete convolution is
\begin{equation}\label{eq:convolution_2d}
    (\mathcal D*u)_{i,j} := \sum_{i'=-\Krad}^{\Krad} \sum_{j'=-\Krad}^{\Krad} \mathcal D_{i',j'} u_{i-i', j-j'} \qquad\mbox{for } i,j\in\mathbb{Z},
\end{equation}
where, consistently with the periodic boundary conditions, the indices
$i-i'$ and $j-j'$ are understood modulo $N$; in other words all convolutions are
performed with circular padding, and no artificial boundary is introduced.
The reaction $\reaction:\R\to\R$ is represented by a small multilayer
perceptron (MLP). Its layer widths and Gaussian activations are shown in
\Figure{fig:pointwise_reaction}. It acts pointwise: at each spatial
point $x$, the same scalar map is applied to the value $v(x)$, with shared
parameters at every grid point. Spatial coupling is provided by the
convolution layers.

\begin{figure}[H]
    \centering
    \begin{subfigure}[b]{0.45\textwidth}
        \centering
        \includegraphics[width=\linewidth]{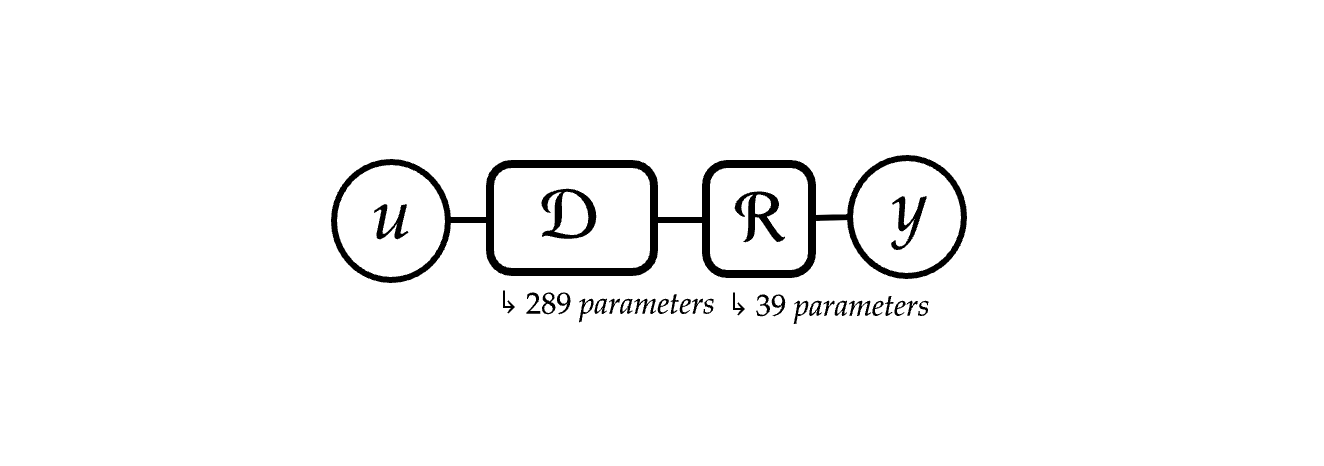}
        \caption{$\snnSimple{\theta}{1}$ ($328$ learnable parameters)}
    \end{subfigure}
    \begin{subfigure}[b]{0.45\textwidth}
        \centering
        \includegraphics[width=\linewidth]{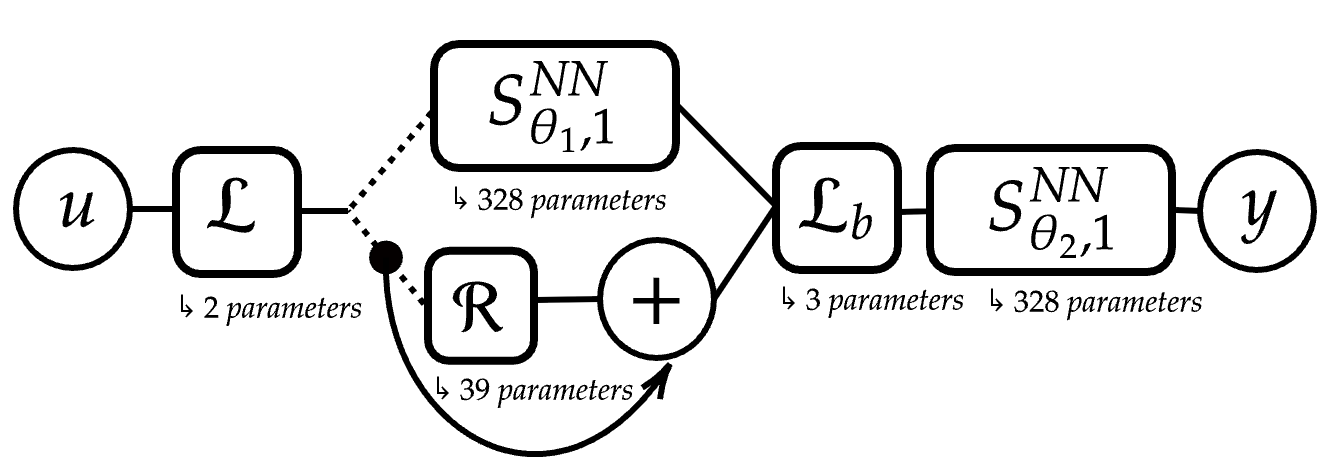}
        \caption{$\snnSimple{\theta}{2}$ ($700$ learnable parameters)}
    \end{subfigure}
    \caption{
        Two network structures obtained from first (left) and second (right) order discretization schemes of the Allen--Cahn equation.
        The number of the learnable parameters in each component is displayed below the corresponding block.
        }
    \label{fig:phase_field_neural_network_structure}
\end{figure}

The structures of the networks
$\snnSimple{\theta}{1}$ and $\snnSimple{\theta}{2}$ are shown in~\Figure{fig:phase_field_neural_network_structure}, where the arrow entering the $+$ symbol represents the residual (skip) connection of $\reaction$,
that is, the map $x\mapsto \reaction(x) + x$, $x\in\mathbb{R}$, applied
pointwise.  The symbols $\mathcal{L}_b$ and $\mathcal{L}$ denote channel-mixing
layers with and without bias, respectively; they act pointwise in space and
only recombine the channels.

The corresponding explicit expressions for the two networks are:

\begin{align*}
    \snnSimple{\theta}{1}(u) &:= \mathcal{R}\circ(\mathcal{D}*u),\\
    \snnSimple{\theta}{2}(u) &:= \snnSimple{\theta_2}{1}\left(\mathcal L_b\circ\begin{pmatrix}
        \snnSimple{\theta_1}{1}(\mathcal L u)_1\\
        (\mathcal L u)_2 + \mathcal R\circ(\mathcal L u)_2\\
    \end{pmatrix}\right),
\end{align*}
where $\theta$, $\theta_1$ and $\theta_2$ collect all the learnable parameters,
namely the entries of the kernel $\mathcal{D}$ and the weights and biases of
$\mathcal{R}$, $\mathcal{L}$ and $\mathcal{L}_b$. The composition operator \(\circ\) is understood entrywise.  The layer $\mathcal{L}$ duplicates
the single input channel into two channels with learnable weights, that is,
$\mathcal{L}u = (a\,u,\, b\,u)^\top$ with $a,b\in\R$, whereas $\mathcal{L}_b$
recombines two channels into one, $\mathcal{L}_b(v_1,v_2)^\top = a'v_1 + b'v_2 + c'$
with $a',b',c'\in\R$.

These operators use fixed local stencils and can therefore act on arrays
of different sizes. Applying a trained model on another grid requires
preserving the ratio $\varepsilon/h$ between interface width and mesh size;
this is the setting used for reconstruction in
\Section{sec:applications_3d}. This compatibility condition alone does not
establish convergence under mesh refinement.

Throughout the remainder of the paper, we write $\snn{\theta}{1}$ and $\snn{\theta}{2}$ when the networks are
trained on the oriented profile $q$, and $\snnprim{\theta}{1}$ and
$\snnprim{\theta}{2}$ when they are trained on the non-oriented profile
$q'$ introduced in \Section{sec:intro_pf}.

\subsection{Architecture motivated by the Allen--Cahn semigroup}

Inspired by the Bence--Merriman--Osher (BMO) approach \cite{merriman1994motion}, and its extension to high order geometric motions by Esedo\=glu, Ruuth and Tsai \cite{Esedoglu_Ruuth_Tsai08}, one can show that the
Willmore flow is recovered by combining several Allen--Cahn flows run over different times.  We recall here the argument of Terii
\cite[Chapter 11]{Terii22}, which motivates the architecture proposed below.

Assume that $u$ is a classical solution to the Allen--Cahn equation, and let
$\mathcal{AC}(t)$ denote the solution operator of the semigroup associated with
that equation, namely $[\mathcal{AC}(t)]u(s, \cdot) = u(s+t, \cdot)$ for every
$s,\, t \geq 0$.  Write
\[
F(u) := \Delta u - \frac{1}{\epsilon^2}W'(u),
\qquad\mbox{so that}\quad \partial_t u = F(u).
\]
The $L^2$-gradients of the two energies introduced in \Section{sec:intro_pf}
are
\begin{equation}\label{eq:gradients}
\nabla\cenergy{\varepsilon}(u) = -\varepsilon\Delta u + \frac{1}{\varepsilon}W'(u)
 = -\varepsilon F(u),
\qquad
\nabla\wenergy{\varepsilon}(u) = \varepsilon\,F'(u)\bigl[F(u)\bigr],
\end{equation}
the second identity following from
$\wenergy{\varepsilon}(u) = \frac{\varepsilon}{2}\int_{\R^d}F(u)^2\,dx$ and
from the self-adjointness of the linearized operator
$F'(u) = \Delta - \frac{1}{\varepsilon^2}W''(u)\operatorname{Id}$.  The Taylor expansion of
$\mathcal{AC}(t)u$ at $t=0$ reads
\[
\mathcal{AC}(t)u = u + t\,\partial_t u + \frac{t^2}{2}\,\partial_{tt}^2 u + \mathcal{O}(t^3),
\]
and, since
\[
\partial_{tt}^2 u = F'(u)\bigl[F(u)\bigr] =
\Delta\left(\Delta u - \frac{1}{\epsilon^2}W'(u)\right)
-\frac{1}{\epsilon^2}W''(u)\left(\Delta u - \frac{1}{\epsilon^2}W'(u)\right),
\]
we deduce from \eqref{eq:gradients} that $\partial_t u = -\frac{1}{\varepsilon}\nabla\cenergy{\varepsilon}(u)$
and $\partial^2_{tt} u = \frac{1}{\varepsilon}\nabla\wenergy{\varepsilon}(u)$, whence
\begin{equation}\label{eq:expansion_ac}
\mathcal{AC}(t)u =  u - \frac{t}{\varepsilon} \nabla\cenergy{\varepsilon}(u)
\;\mathbf{+}\;\frac{t^2}{2\varepsilon} \nabla\wenergy{\varepsilon}(u)
+\mathcal{O}(t^3),
\end{equation}
where $\nabla\cenergy{\varepsilon}(u)$ and $\nabla\wenergy{\varepsilon}(u)$
denote the $L^2$-gradients of the Cahn--Hilliard energy and of the phase field
Willmore energy, respectively.

Formula \eqref{eq:expansion_ac} shows in particular that the second-order term
of $\mathcal{AC}(t)$ can be isolated by a combination of the form
\begin{equation}\label{eq:two_term_bmo}
\left[
2 \mathcal{AC}\left(\sqrt{t}\right)
- \mathcal{AC}\left(2\sqrt{t}\right)
\right]u
= u - \frac{t}{\varepsilon}\nabla\wenergy{\varepsilon}(u)
+\mathcal{O}(t^{3/2}),
\end{equation}
in which the Cahn--Hilliard term, of order $\sqrt{t}$, cancels while the
Willmore term, of order $t$, survives with the sign of a gradient descent.

To obtain more accurate schemes, one possibility is to consider additional
combinations.  For instance, a higher-order cancellation is obtained by
choosing the pairs $(\alpha_i,\beta_i)_{1\leq i\leq 4}$, with the $\alpha_i>0$
pairwise distinct, so that
\begin{equation*}
\sum_{i=1}^{4} \beta_i \mathcal{AC}(\alpha_i t)u
= u - \frac{t^2}{\varepsilon}\nabla\wenergy{\varepsilon}(u)
+\mathcal{O}(t^4).
\end{equation*}
Indeed, inserting \eqref{eq:expansion_ac} gives
\[
\sum_{i=1}^{4} \beta_i \mathcal{AC}(\alpha_i t)u
= \Bigl(\sum_i\beta_i\Bigr) u
 - \frac{t}{\varepsilon}\Bigl(\sum_i\beta_i\alpha_i\Bigr)\nabla\cenergy{\varepsilon}(u)
 + \frac{t^2}{2\varepsilon}\Bigl(\sum_i\beta_i\alpha_i^2\Bigr)\nabla\wenergy{\varepsilon}(u)
 + \mathcal{O}(t^3),
\]
so that the announced identity holds as soon as $\sum_i\beta_i = 1$,
$\sum_i\beta_i\alpha_i = 0$, $\sum_i\beta_i\alpha_i^2 = -2$ and, to cancel the
next order, $\sum_i\beta_i\alpha_i^3 = 0$.

It suffices to enforce the following Vandermonde system:
\[
\begin{pmatrix}
1 & 1 & 1 & 1 \\
\alpha_1 & \alpha_2 & \alpha_3 & \alpha_4 \\
\alpha_1^2 & \alpha_2^2 & \alpha_3^2 & \alpha_4^2 \\
\alpha_1^3 & \alpha_2^3 & \alpha_3^3 & \alpha_4^3
\end{pmatrix}
\begin{pmatrix}
\beta_1\\
\beta_2\\
\beta_3\\
\beta_4
\end{pmatrix}
=
\begin{pmatrix}
1\\
0\\
-2\\
0
\end{pmatrix}.
\]

This matrix equation admits a unique solution if and only if the $\alpha_i$ are pairwise distinct. In this case, the coefficients $\beta_i$ are uniquely determined and satisfy
\[
\beta_i =
\frac{\prod_{j\neq i}\alpha_j}{\prod_{j\neq i}(\alpha_j-\alpha_i)}
-2\frac{\sum_{j\neq i}\alpha_j}{\prod_{j\neq i}(\alpha_j-\alpha_i)},
\qquad i=1,\ldots,4.
\]

Thus, choosing the $\alpha_i$ appropriately yields a higher-order scheme.  The
advantage of such an approach is that it is significantly easier to simulate
Allen--Cahn flows, and combinations thereof, than to simulate directly the
gradient flow of the energy $\wenergy{\varepsilon}$.  However, the parameters
$\alpha_i$, and consequently the $\beta_i$, must be chosen carefully in order to
ensure the stability of the scheme, because the cancellation conditions alone do not ensure stability.

This suggests an architecture combining in parallel $2$ or $4$ copies of the
network adapted to the Allen--Cahn flow, recombined by a linear layer.  We
denote the resulting networks by $\wnnSimple{\theta}{1}$ and
$\wnnSimple{\theta}{2}$, respectively.  Their main advantage is that the
parameters $\alpha_i$, which are encoded implicitly in the diffusion kernels of each
branch, the weights $\beta_i$, and the coefficients of the linear
recombination are not prescribed but optimized during training.

\begin{figure}[H]
    \centering
    \begin{subfigure}[b]{0.48\textwidth}
        \centering
        \includegraphics[width=\linewidth]{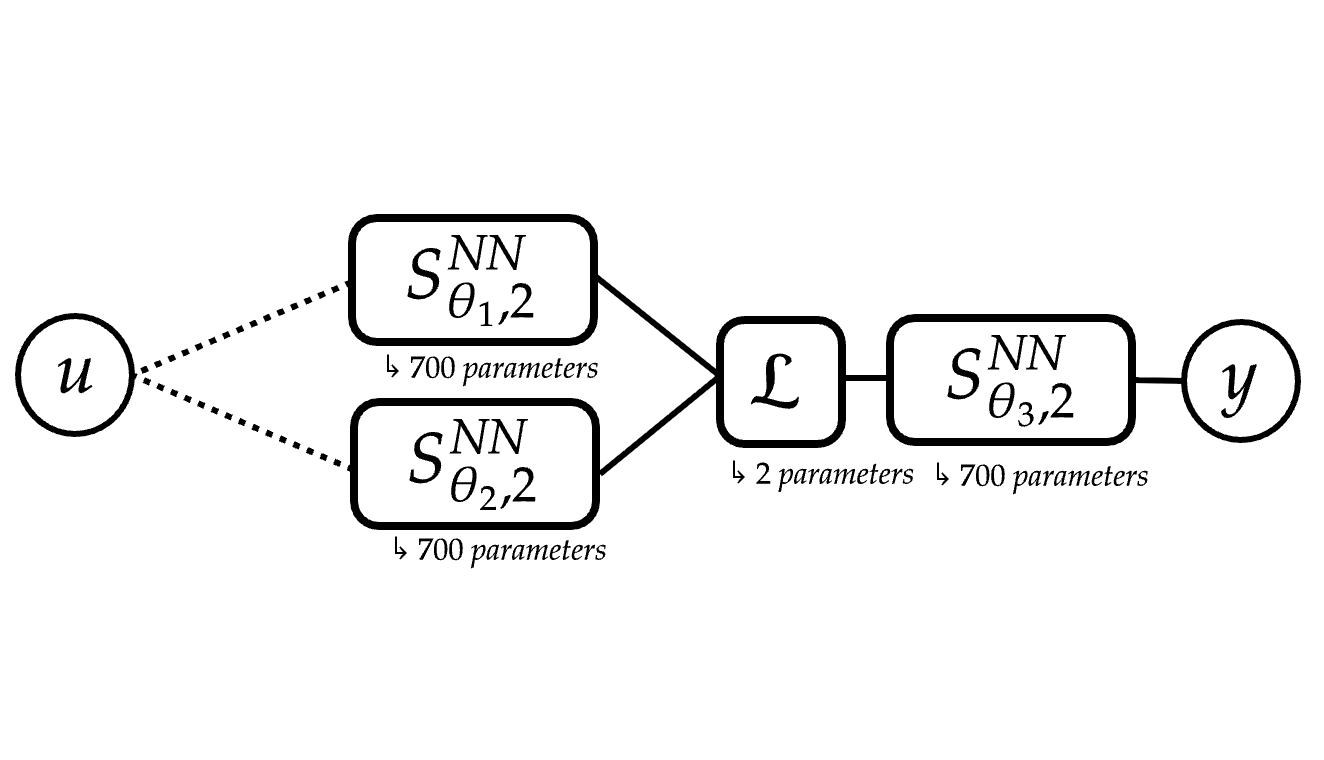}
        \caption{$\wnnSimple{\theta}{1}$ ($2102$ learnable parameters)}
    \end{subfigure}
    \hfill
    \begin{subfigure}[b]{0.48\textwidth}
        \centering
        \includegraphics[width=\linewidth]{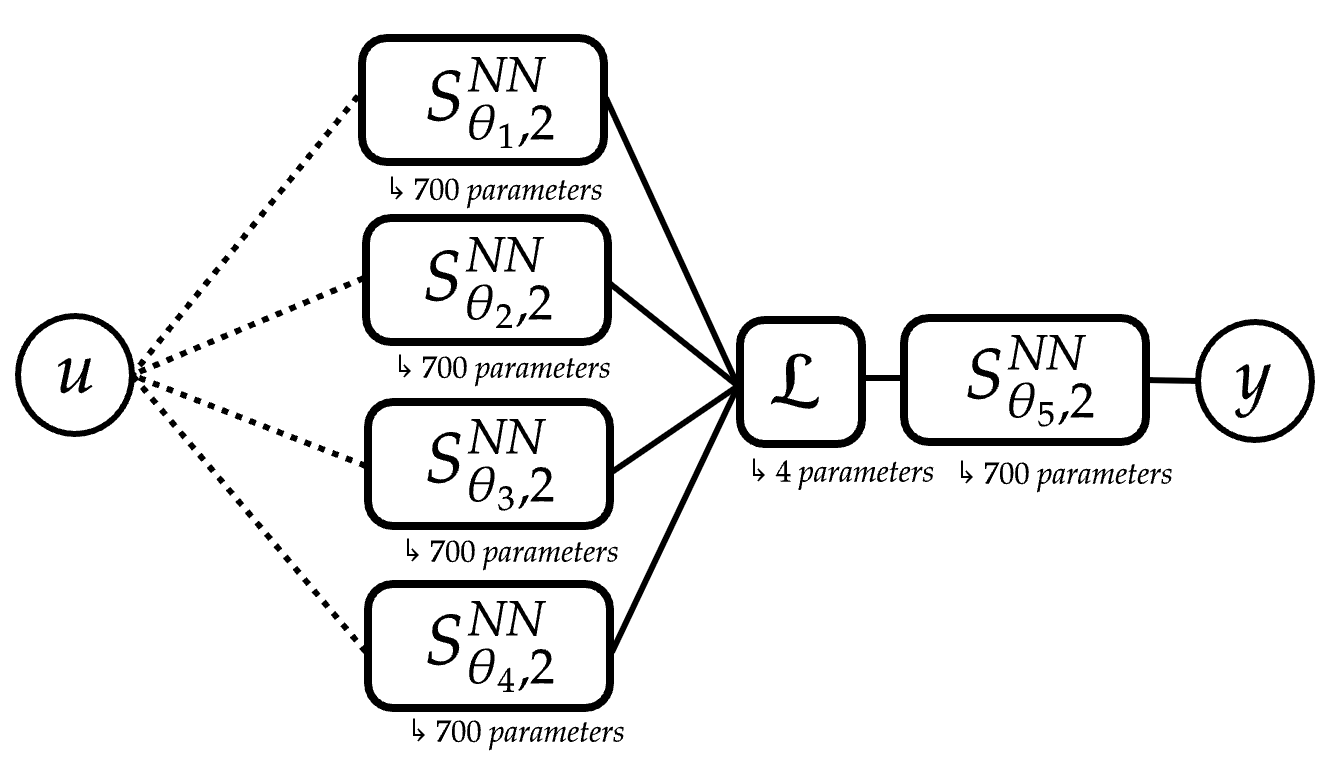}
        \caption{$\wnnSimple{\theta}{2}$ ($3504$ learnable parameters)}
    \end{subfigure}
    \caption{
        The network structures (two (left) and four (right) parallel models) inspired by the BMO scheme.
        The number of the learnable parameters in each component is displayed below the corresponding block.
        }
    \label{fig:nn_structure}
\end{figure}

Precisely, we define $\wnnSimple{\theta}{1}$ and $\wnnSimple{\theta}{2}$ as follows, see also \Figure{fig:nn_structure}:

\begin{equation}
    \label{eq:formula_WillmoreNN}
    \wnnSimple{\theta}{1}(u) := \snnSimple{\theta_3}{2}\left(\mathcal{L}\begin{pmatrix}
        \snnSimple{\theta_1}{2}(u)\\
        \snnSimple{\theta_2}{2}(u)
    \end{pmatrix}\right),\qquad
    \wnnSimple{\theta}{2}(u) := \snnSimple{\theta_5}{2}\left(\mathcal{L}\begin{pmatrix}
        \snnSimple{\theta_1}{2}(u)\\
        \snnSimple{\theta_2}{2}(u)\\
        \snnSimple{\theta_3}{2}(u)\\
        \snnSimple{\theta_4}{2}(u)\\
    \end{pmatrix}\right).
\end{equation}
We note that the input of $\wnnSimple{\theta}{1}$ (resp.\ $\wnnSimple{\theta}{2}$)
is copied into $2$ channels (resp.\ $4$ channels).
The linear layers of both models do not include any bias term.
The diffusion kernels have size $(2\Krad+1)\times(2\Krad+1)$, with
$\Krad=8$ in practice, as specified in \Section{sec:intro_nn}.

As mentioned in the introduction, training our networks to approximate the
Willmore flow requires a substantially richer dataset than the evolution of
circles in 2D, or of spheres in 3D. This motivates the use of a dedicated numerical scheme to produce reliable approximations of the Willmore flow. In this study, we adopt the approach proposed by Rumpf, Sassen, and Smoch \cite{Rumpf_Sassen_Smoch25}, which is based on a minimizing movement scheme. A detailed description of this method is provided in the following subsection.

\subsection{Willmore flow approximation using the hybrid RSS approach}\label{sec:hybrid_rss}
We now describe the hybrid scheme of Rumpf, Sassen and Smoch (RSS)~\cite{Rumpf_Sassen_Smoch25}, which we use throughout the paper to generate
training data.  It combines a minimizing movement scheme for the Willmore
energy with the mean curvature networks of \cite{BretinDenisMasnouTerii2022}:
given the state $V$ reached at the previous time step, one minimizes
\begin{equation}\label{eq:energy}
    E_{\varepsilon,\ts,\tts}(U, V) := \varepsilon\|U - V\|^2_{L^2(\Omega)}
    + \frac{\varepsilon\tts}{\ts^2}\|\mathcal{N}(U) - U\|^2_{L^2(\Omega)},
\end{equation}
where $\mathcal{N}$ is either of the networks $\snn{\theta}{1}$,
$\snn{\theta}{2}$ trained in \cite{BretinDenisMasnouTerii2022} to approximate
the mean curvature flow.  The first term is the usual minimizing movement
penalization; the second is a discrete surrogate for
$\|\nabla\cenergy{\varepsilon}(U)\|^2_{L^2}$, the network supplying at
negligible cost the one-step mean curvature evolution against which $U$ is
compared.  The two time steps should not be confused: $\ts$ is the step over
which $\mathcal{N}$ was trained, and is thus a fixed characteristic of the
network, whereas $\tts$ is the step of the minimizing movement scheme, which is
ours to choose.

The state $U_n$ is obtained from $U_{n-1}$ by $J$ iterations of the ADAM
optimizer applied to $U \mapsto E_{\varepsilon,\ts,\tts}(U, U_{n-1})$,
\begin{equation}\label{eq:adam}
    U_{n,0} := U_{n-1}, \qquad
    U_{n,k} := U_{n,k-1} - \eta_k\, G_{n,k}, \quad k = 1,\dots, J,
    \qquad U_{n} := U_{n,J},
\end{equation}
where $\eta_k$ and $G_{n,k}$ are the step size and the descent direction
produced by ADAM from the gradients of $E_{\varepsilon,\ts,\tts}(\cdot\,,
U_{n-1})$.
We depart here from \cite{Rumpf_Sassen_Smoch25}, where this inner
minimization is carried out by a more efficient Newton-type method.  Throughout,
$J$ denotes the number of these inner iterations, $\Krad$ being reserved for
the half-width of the convolution stencils in \eqref{eq:convolution_2d}. \\

We now show an example of an evolution computed with the RSS scheme
involving singularities.  We take $\Omega := (0,1)^2$, discretized as a Cartesian grid $\mathcal{M}$ with $N = 2^7$ points in each direction, and set $\varepsilon = 1/N = 2^{-7}$ and
$\ts = \varepsilon^2 = 2^{-14}$.

\begin{figure}[H]
    \centering
    \begin{subfigure}[b]{0.23\textwidth}
        \centering
        \includegraphics[width=\linewidth]{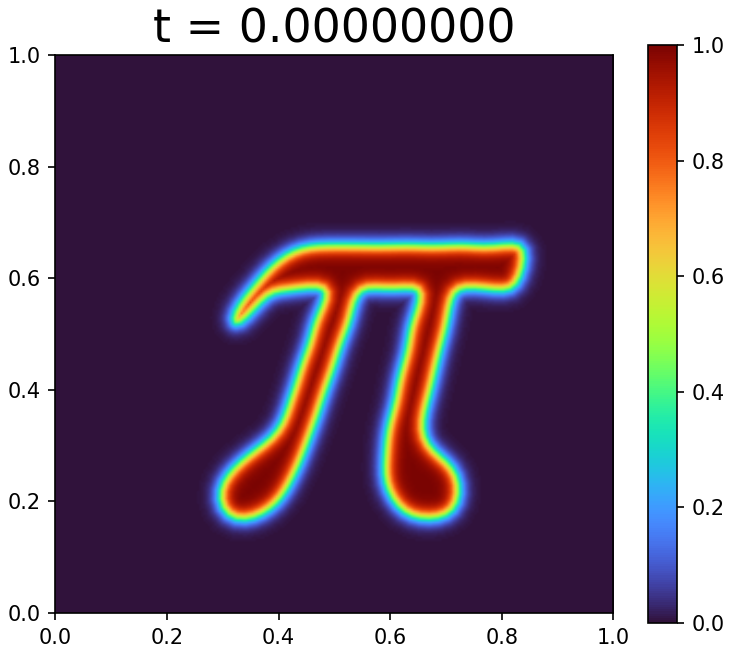}
    \end{subfigure}
    \hfill
    \begin{subfigure}[b]{0.23\textwidth}
        \centering
        \includegraphics[width=\linewidth]{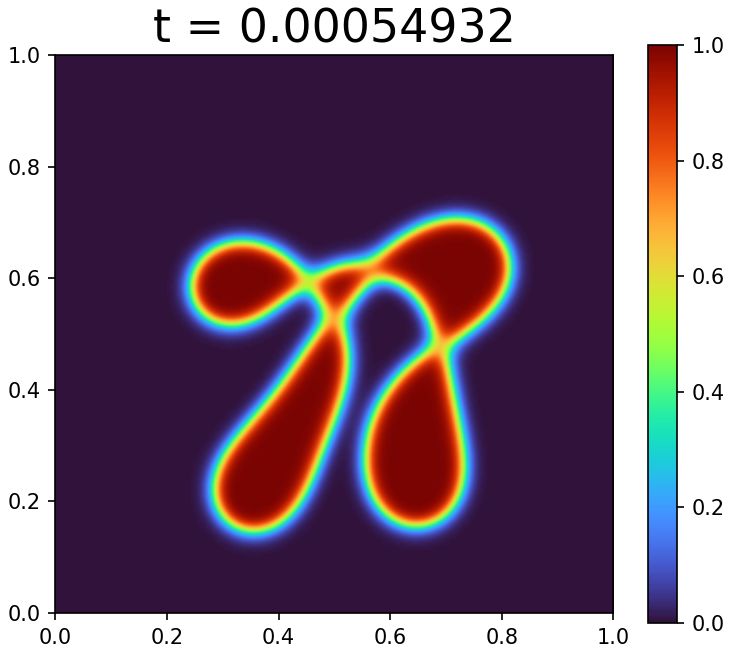}
    \end{subfigure}
    \hfill
    \begin{subfigure}[b]{0.23\textwidth}
        \centering
        \includegraphics[width=\linewidth]{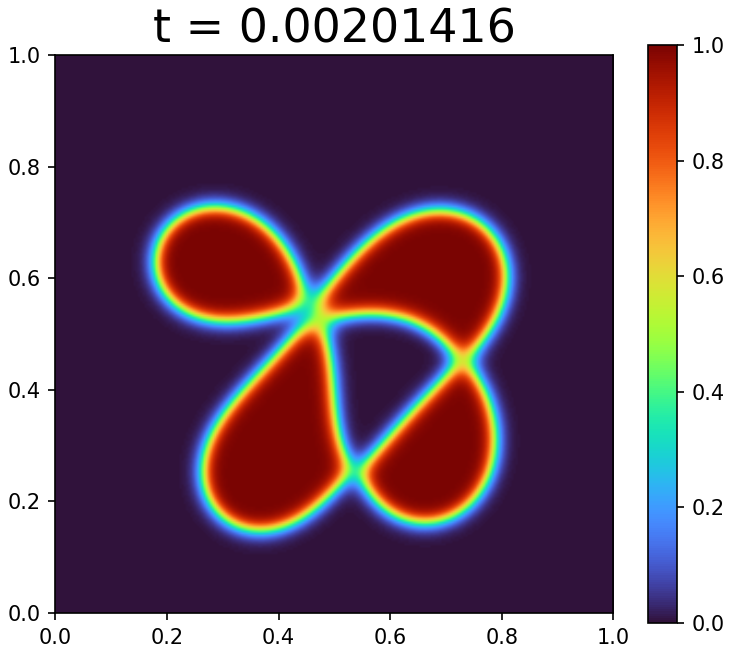}
    \end{subfigure}
    \hfill
    \begin{subfigure}[b]{0.23\textwidth}
        \centering
        \includegraphics[width=\linewidth]{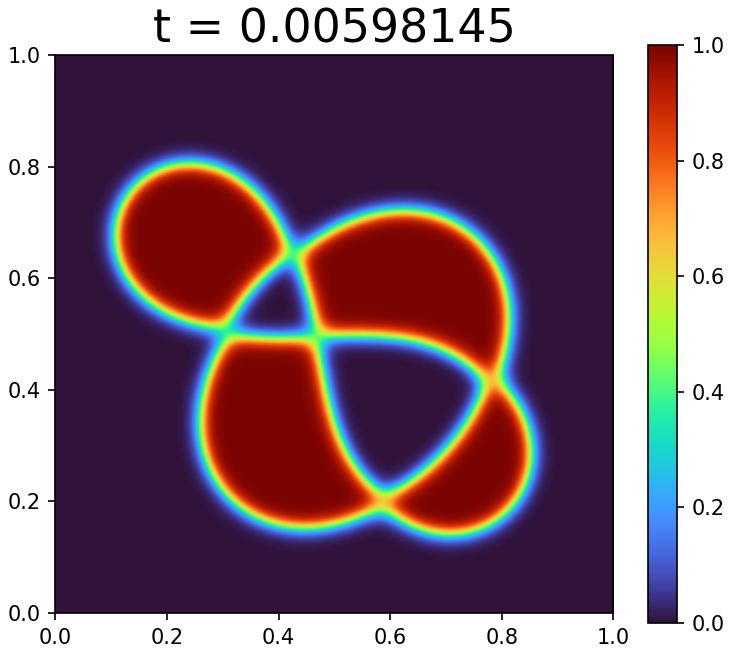}
    \end{subfigure}
    \caption{Evolution of a $\pi$-shaped curve under the RSS scheme with $\tts = \ts$, $J = 500$, and $T = 0.006$.}
    \label{fig:willmore-pi}
\end{figure}

The evolution starts from a $\pi$-shaped curve
(\Figure{fig:willmore-pi}).
Its initial coordinates are taken from \textit{WolframAlpha}\footnote{\url{https://ja.wolframalpha.com/}} and are scaled by the factor $1/3000$.

In  \Figure{fig:willmore-pi},
the evolution
deforms the initial curves into more circular shapes over time while preserving the crossing singularities.  We note that the scheme remains well-behaved for
the comparatively large time step $\tts = \ts$; this is consistent with
\cite[Proposition 2.4]{Rumpf_Sassen_Smoch25}, where the existence of a
minimizer of \eqref{eq:energy} is established even though the functional
involves the neural network $\mathcal{N}$.

The scheme does, however, leave two parameters essentially unconstrained: the
number $J$ of ADAM iterations performed at each time step, and the time step
$\tts$ itself.  Neither is determined by the analysis, and the value of $J$
needed to reach a given accuracy varies significantly with the problem and with
the initial condition.  Since the whole training database of this paper is
generated with this scheme, we examine this point quantitatively.

In the case of a circle, for which the exact solution is known, the number of
inner iterations required to reach a satisfactory accuracy at $\tts = \ts/64$
is of the order of $2000$: as reported in \Table{tab:rerror_rumpf}, $J = 1000$
leaves a mean relative error of about $15\%$.  We also observe that increasing
$J$ beyond $2000$ brings no further gain, the error then being dominated by the
time discretization rather than by the inner optimization; reducing it further
would require decreasing $\tts$ as well.

The radius $R_{J}(n \tts)$ of the numerical solution is estimated from
$U_{n,J}$ through its volume, namely
\begin{equation}\label{eq:radius_estimate}
 R_{J}(n \tts) = \left( \frac{1}{\pi} \int_{\Omega}  U_{n,J}\,dx\right)^{1/2},
\end{equation}
which is exact in the sharp interface limit for a disk and, for a diffuse
profile, correct up to $\mathcal{O}(\varepsilon^2)$.

\begin{figure}[H]
    \centering
    \begin{subfigure}[b]{0.23\textwidth}
        \centering
        \includegraphics[width=\linewidth]{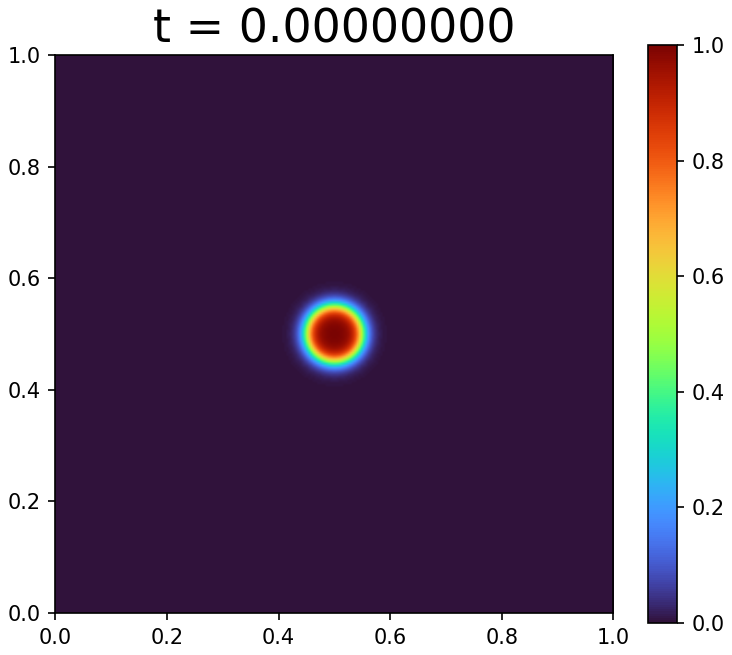}
    \end{subfigure}
    \hfill
    \begin{subfigure}[b]{0.23\textwidth}
        \centering
        \includegraphics[width=\linewidth]{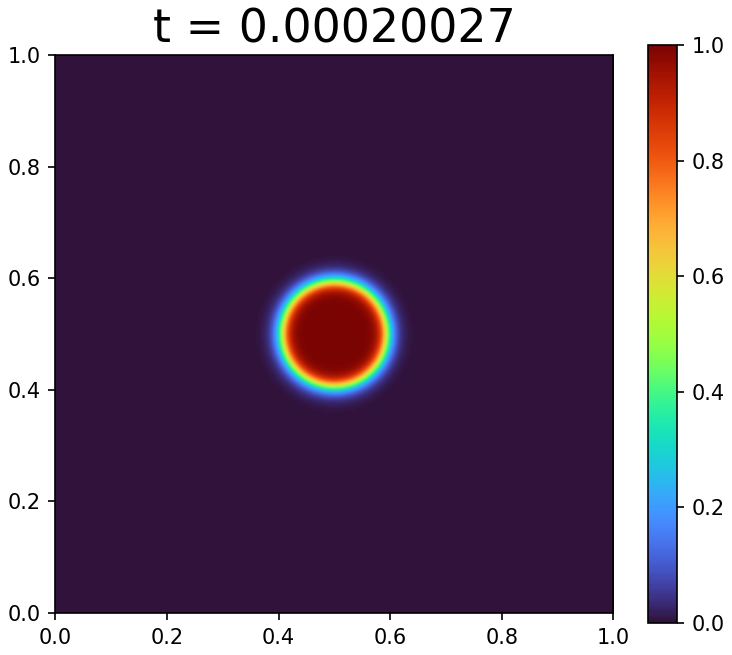}
    \end{subfigure}
    \hfill
    \begin{subfigure}[b]{0.23\textwidth}
        \centering
        \includegraphics[width=\linewidth]{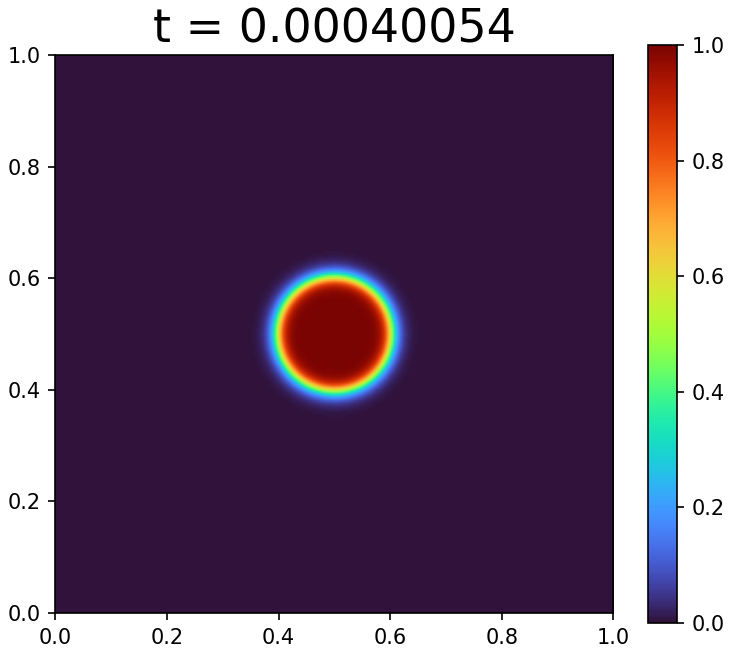}
    \end{subfigure}
    \hfill
    \begin{subfigure}[b]{0.23\textwidth}
        \centering
        \includegraphics[width=\linewidth]{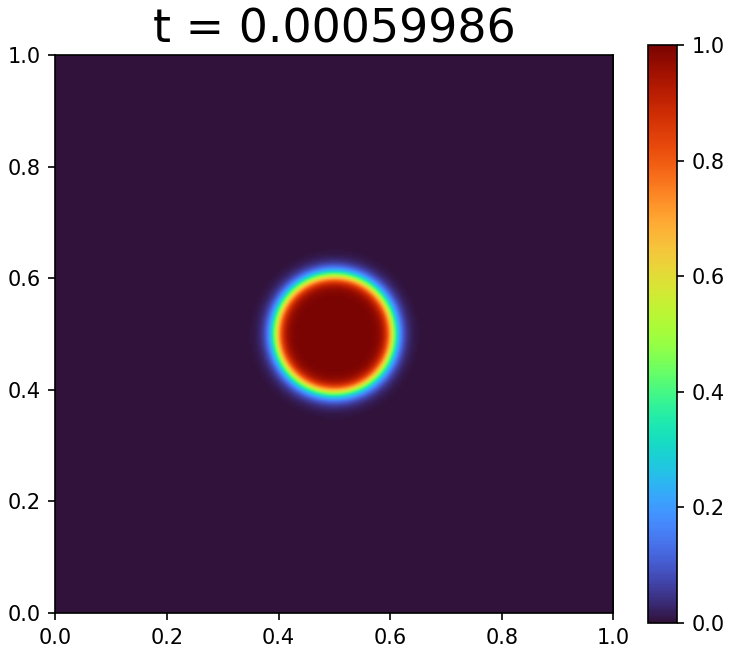}
    \end{subfigure}
    \caption{Evolution of a circle by the RSS scheme with $\tts = \ts / 64$ and $J = 1000$ by $T = 0.0006$.}
    \label{fig:circle_rumpf}
\end{figure}

\begin{figure}[H]
    \centering
    \begin{subfigure}[b]{0.8\textwidth}
        \centering
        \includegraphics[width=\linewidth]{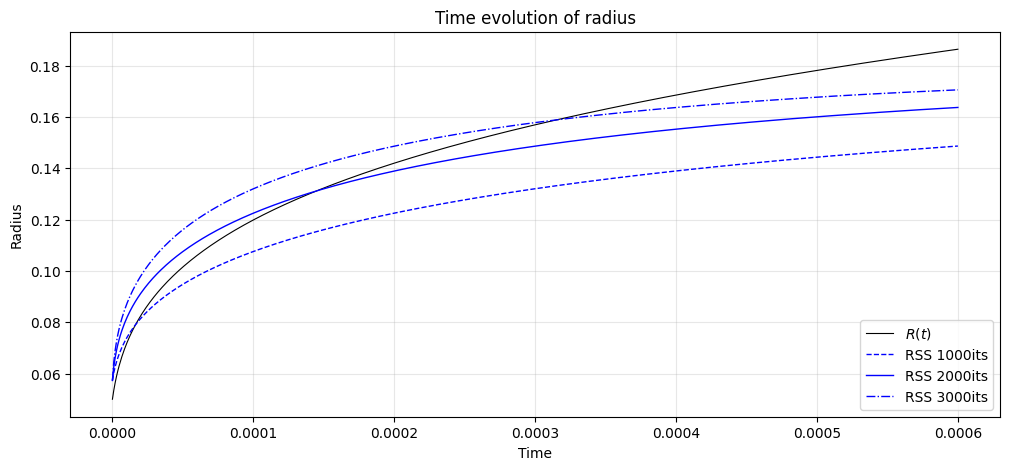}
    \end{subfigure}
    \caption{Evolution of the radius of a circle computed with the RSS scheme for
    $\tts = \ts / 64$ up to $T = 0.0003$, for $J = 1000$, $2000$ and $3000$
    inner iterations, compared with the exact law.
    The radius is estimated from the numerical solution through
    \eqref{eq:radius_estimate}.}
    \label{fig:radius_re_rumpf}
\end{figure}

\begin{table}[H]
    \centering
    \begin{tabular}{|c|c|c|c|}
        \toprule
        \textbf{Inner iterations $J$} &  \textbf{MRE} & \textbf{Time [s]} & \textbf{Time / inner iter.\ [s]} \\
        \midrule
        1000 & $0.1470$ & $1819$ & $0.0029$ \\
        2000 & $0.0647$ & $3718$ & $0.0030$ \\
        3000 & $0.0628$ & $5446$ & $0.0029$ \\
        \bottomrule
    \end{tabular}
    \caption{Mean relative error (MRE) of the RSS scheme with respect to the exact
    solution, and computational time, as functions of the number $J$ of inner
    ADAM iterations, for the evolution of a circle with $\tts = \ts/64$ and
    $T = 0.0003$.}
    \label{tab:rerror_rumpf}
\end{table}

In \Table{tab:rerror_rumpf}, the mean relative error (MRE) is defined by

\begin{equation}\label{eq:aerror}
    MRE := \frac{1}{T}\int_0^T RE(s) ds\qquad\mbox{ with }\quad RE(s) := \frac{\|U(s) - U_{\mathrm{exact}}(s)\|_{L^2(\Omega)}}{\|U_{\mathrm{exact}}(s)\|_{L^2(\Omega)}},
\end{equation}
where $U_{exact}$ denotes the phase-field function of the exact solution to the Willmore flow in 2D at time $s$;
the integration is calculated by the trapezoidal rule.

There is therefore a trade-off between accuracy and computational cost, the cost
being proportional to $J$.  The purpose of this paper is to remove that
trade-off at inference time by learning the RSS scheme with a neural network: if
the strategy succeeds, interfaces are evolved by a single forward pass of the
network instead of an inner ADAM optimization, the cost of the latter being paid
once and for all when the database is built.

\subsection{Dataset generation and training procedure}\label{sec:database_construction}

In this section we run the RSS scheme \eqref{eq:adam} to generate the database
used to train the neural networks $\wnnSimple{\theta}{1}$ and
$\wnnSimple{\theta}{2}$.

We build a database from the evolution of disks, ovals, unions of two or three
disks, and star-shaped curves.  Except the disk, for which the
Willmore flow is explicit, all these evolutions are computed numerically with
the RSS scheme.
Throughout this section the time step of the RSS scheme is $\tts = \ts / 64$,
where we recall that $\ts = 2^{-14}$ is the time step used in
\cite{BretinDenisMasnouTerii2022} to train the mean curvature networks
$\snnSimple{\theta}{1}$ and $\snnSimple{\theta}{2}$, and the number of inner
iterations is $J = 300$.  We work on the Cartesian grid $\mathcal{M}$ of mesh
size $h = 1/N$ with $N = 2^8$, and set $\varepsilon := 2/N = 2^{-7}$.

\begin{figure}[H]
    \centering
    \begin{subfigure}[b]{0.90\textwidth}
        \includegraphics[width=\linewidth]{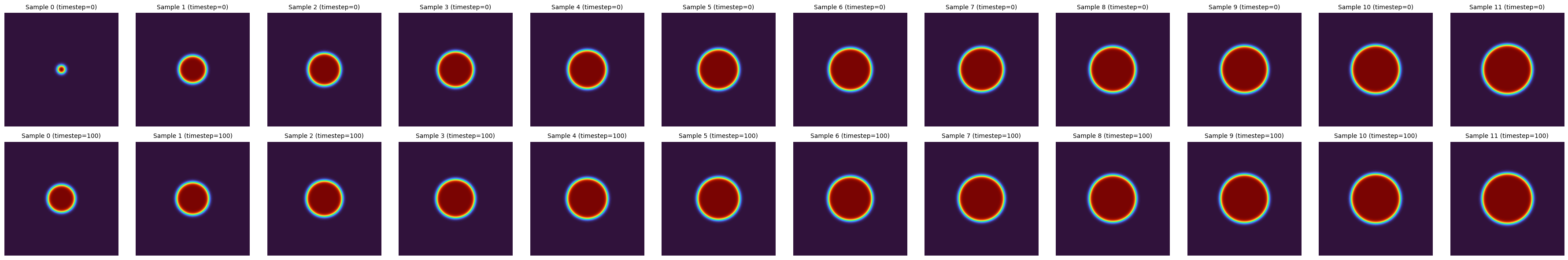}
        \caption{Circles}
        \label{fig:dataset_oriented_circles}
    \end{subfigure}
    \hfill
    \begin{subfigure}[b]{0.3\textwidth}
        \includegraphics[width=\linewidth]{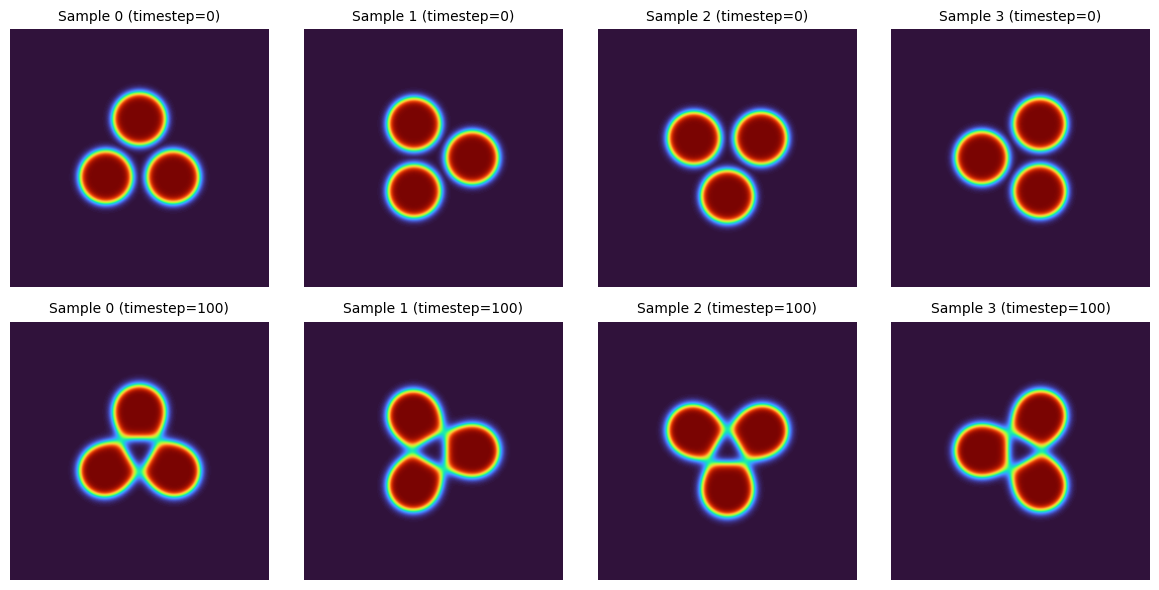}
        \caption{Three circles}
        \label{fig:dataset_oriented_three_circles}
    \end{subfigure}
    \hfill
    \begin{subfigure}[b]{0.3\textwidth}
        \includegraphics[width=\linewidth]{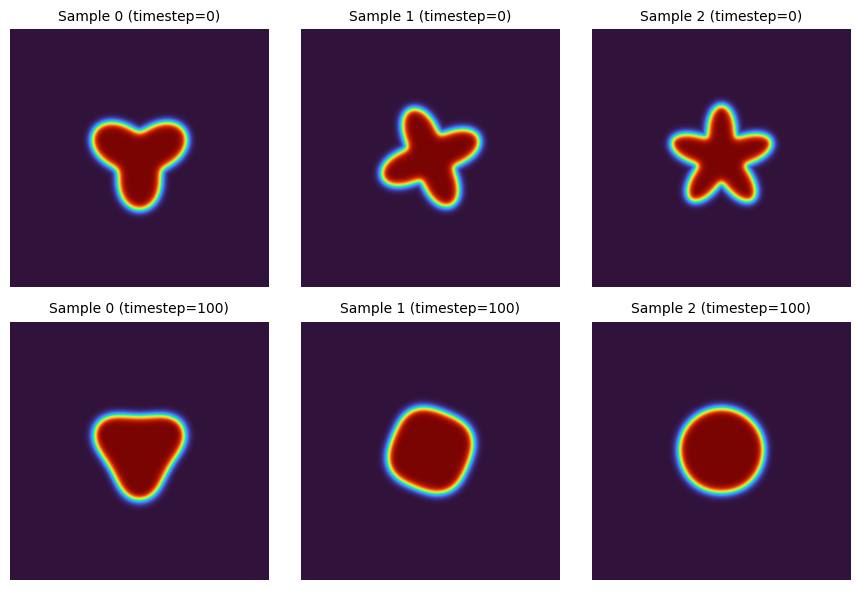}
        \caption{Stars}
        \label{fig:dataset_oriented_stars}
    \end{subfigure}
    \hfill
    \begin{subfigure}[b]{0.3\textwidth}
        \includegraphics[width=\linewidth]{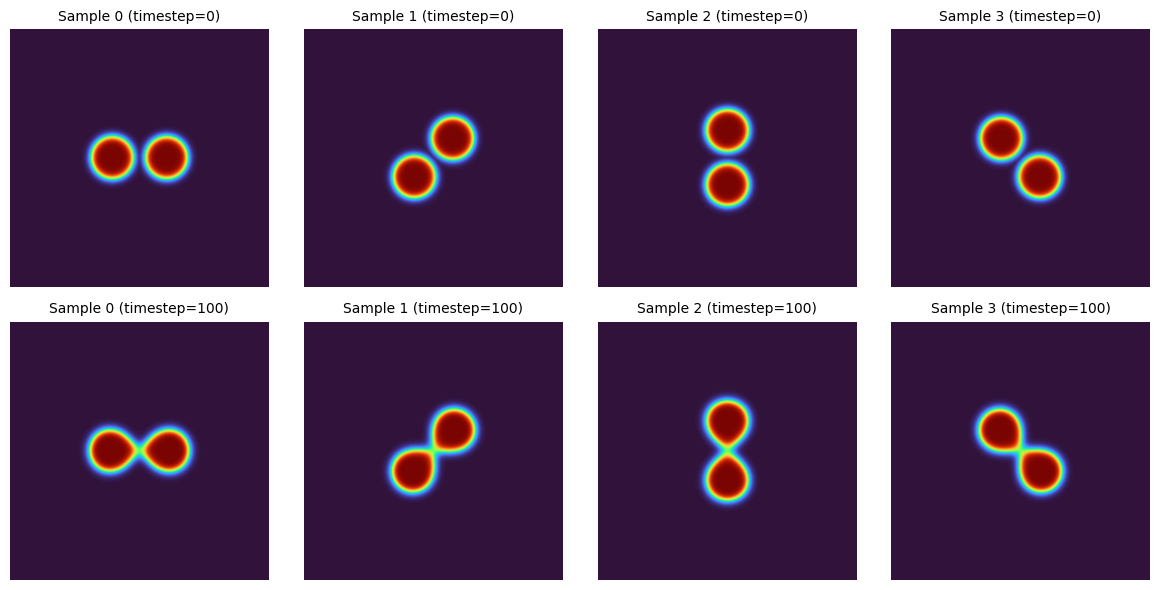}
        \caption{Horizon circles}
        \label{fig:dataset_oriented_horizon_circles}
    \end{subfigure}
    \caption{Datasets of oriented curves}
    \label{fig:datasets_oriented}
\end{figure}

\paragraph{The database.}Throughout the paper, a sample is a contiguous
subsequence of $101$ phase field functions extracted from one computed
evolution, each evolution being run for $100$ time steps of the RSS scheme
\eqref{eq:adam}.  \Table{tab:dataset_2d} collects the four families of initial
curves and their parameters, and \Figure{fig:datasets_oriented} shows a few of
the resulting evolutions.  From each sample we extract all the consecutive
triplets $(V_{k}, V_{k+1}, V_{k+2})$ for $k = 1,\dots,99$ representing one input and the two
successive targets used by the multistep loss \eqref{eq:loss_multistep} with
$S = 2$; it amounts to $5346$ triplets in total.

\begin{table}[H]
    \centering
    \small
    \begin{tabular}{l >{\raggedright\arraybackslash}p{0.50\textwidth} c}
        \toprule
        \textbf{Family} & \textbf{Initial curve and parameters} & \textbf{Samples}\\
        \midrule
        Disks & $B((a,b),R)$, several radii & $43$\\
        Two disks & union of $B((0.5-m/2,0.5),r)$ and $B((0.5+m/2,0.5),r)$
        rotated by $\theta$; $r = 0.08$, $m = 0.05$,
        $\theta \in \{0,\pi/4,\pi/2,3\pi/4\}$ & $4$\\
        Three disks & union of the $B((0.5 + R\sin\theta_i, 0.5 + R\cos\theta_i), r)$,
        $\theta_i = 2i\pi/3$, rotated by $\theta$; $r = 0.1$, $R = 1.5\,r$,
        $\theta$ as above & $4$\\
        Stars & polar curve $r(\theta) = s + \frac{s}{3}\sin(f\theta)$;
        $s = 0.15$, $f \in \{3,4,5\}$ & $3$\\
        \bottomrule
    \end{tabular}
    \caption{The four families of initial curves used to build the
    two-dimensional database.  Each sample is one evolution of $101$ phase
    field functions computed with the RSS scheme at $\tts = \ts/64$ and
    $J = 300$. Here, $B((a,b),R)$ denotes the disk with center
 $(a,b)\in\mathbb{R}^2$ and radius $R > 0$.}
    \label{tab:dataset_2d}
\end{table}

Training is implemented with the nnpf \footnote{\url{https://github.com/PhaseFieldICJ/nnpf}} package \cite{BretinDenisMasnouTerii2022}
which is based on \emph{PyTorch Lightning} \cite{Falcon_PyTorch_Lightning_2019}
to handle automatic logging and check-pointing.
We split the dataset into $80\%$ for training and $20\%$ for validation, the
split being performed by evolution: all the triplets extracted from a
given computed trajectory lie on the same side of the split.  This point
matters here.  Two consecutive triplets of one trajectory are strongly
correlated, since the phase field function changes very little from one time
step to the next, so a split performed at the level of individual triplets
would leave in the validation set states almost identical to states seen during
training, and the validation loss would measure interpolation along known
trajectories rather than generalization.

For tensors $X$ and $Y$,
the loss functions are defined by
\begin{equation}\label{eq:loss_multistep}
    \loss{\mathcal{N}}(X, \bv Y)
    := \frac{1}{BSN^2}\sum_{b=1}^B\sum_{s=1}^{S}\sum_{i=1}^N\sum_{j=1}^N\left({\mathcal{N}}^{(s)}(X^{(b)}(i,j)) - Y_s^{(b)}(i,j)\right)^2,
\end{equation}
with $\mathcal{N}$ being either $\wnnSimple{\theta}{1}$, $\wnnSimple{\theta}{2}$,
and where $B$ denotes the batch size of the training data;
$X^{(b)}(i,j)$ and $Y^{(b)}_{s}(i,j)$ denote the $i$-th and $j$-th components of the $b$-th batch of the training data, respectively;
$(\wnnSimple{\theta}{1})^{(s)}(X)$ and $(\wnnSimple{\theta}{2})^{(s)}(X)$ denote the $s$ times applications of the neural networks $\wnnSimple{\theta}{1}$ and $\wnnSimple{\theta}{2}$ to the input $X$.
In this experiment, we set $B = 32$ and $S = 2$;
each epoch ends when all the training data are used at least once.
The averaging over the $S$ steps in \eqref{eq:loss_multistep} is what makes the
training loss, evaluated at $S = 2$, and the validation loss, evaluated at
$S = 5$, comparable with each other; the same normalization is used in 3D, see
\eqref{eq:loss_multistep_3d}.
We optimize the trainable parameters $\theta$ of $\wnnSimple{\theta}{1}$ and $\wnnSimple{\theta}{2}$ using the ADAM optimizer with the learning rate $10^{-3}$ to be fixed throughout the training.
This outer use of ADAM, which acts on the parameters $\theta$, should not be
confused with the inner ADAM iterations of \eqref{eq:adam}, which act on the
phase field function itself when the database is generated.
Selection of the models $\wnnSimple{\theta}{1}$ and $\wnnSimple{\theta}{2}$ is based on the validation loss which is defined by \eqref{eq:loss_multistep} with $S = 5$.

\begin{rem}\label{rem:multistep}
    Since the evolution by the Willmore flow is quite slow, we have considered the multistep loss function in \eqref{eq:loss_multistep} with $S = 2$ rather than the single-step loss functions.
    We will also adopt the multistep loss function in \eqref{eq:loss_multistep} in applications for non-oriented curves and surfaces in 3D.
\end{rem}

\begin{figure}[H]
    \centering
    \begin{subfigure}[b]{0.48\textwidth}
        \centering
        \includegraphics[width=\linewidth]{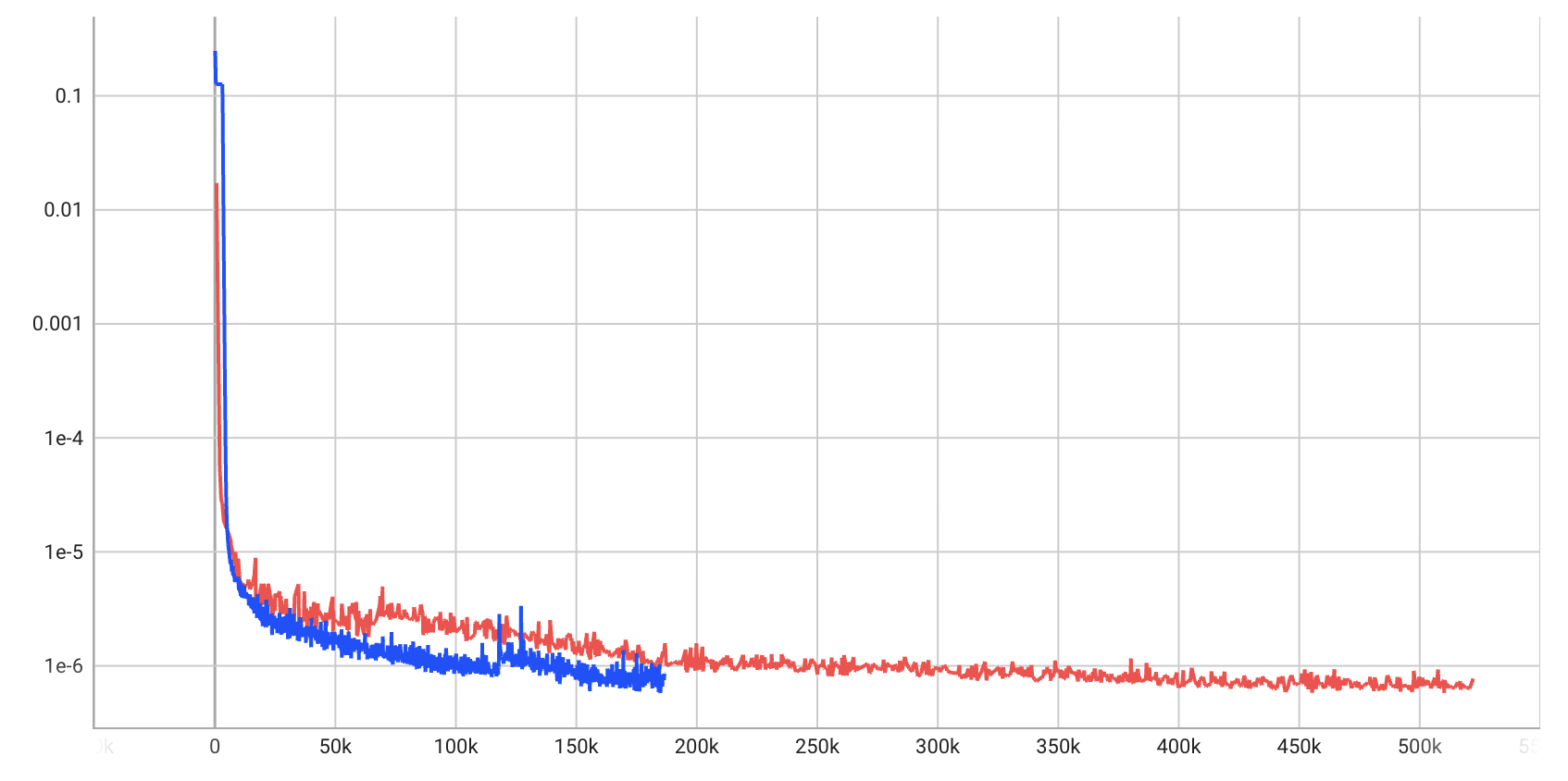}
    \end{subfigure}
    \hfill
    \begin{subfigure}[b]{0.48\textwidth}
        \centering
        \includegraphics[width=\linewidth]{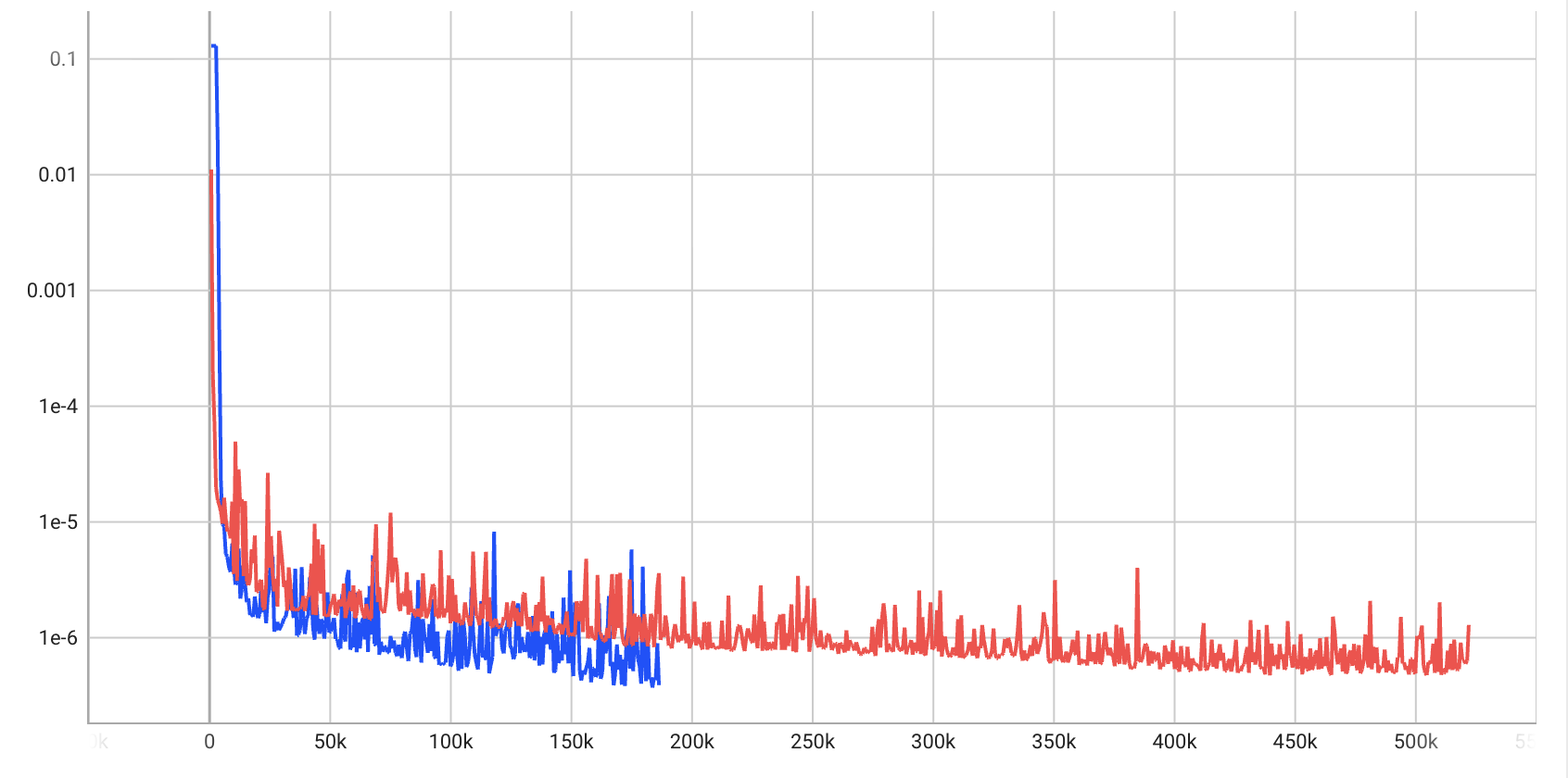}
    \end{subfigure}
    \caption{Evolution of the training (left) and validation (right) losses at each training step for the oriented curve models. Red: $\wnnSimple{\theta}{1}$, blue: $\wnnSimple{\theta}{2}$.}
    \label{fig:training_loss_oriented}
\end{figure}
\Figure{fig:training_loss_oriented} displays the training and validation losses,
and \Table{tab:training-nn} reports the validation losses of
$\wnnSimple{\theta}{1}$ and $\wnnSimple{\theta}{2}$.  The checkpoint retained
for each model is the one achieving the lowest validation loss once the
training curve has reached a plateau; hereafter we denote the selected
checkpoints by $\wnn{\theta}{1}$ and $\wnn{\theta}{2}$.

\begin{rem}[Reading the loss tables]\label{rem:loss_tables}
The validation losses reported in \Table{tab:training-nn} and in the tables of
the following sections should be compared within a table and not across
tables: the datasets, the time steps and the representations of the interface
differ from one experiment to the next, so the corresponding losses are not
commensurable.  Within \Table{tab:training-nn}, the two architectures are in
fact very close, that is $4.73\times 10^{-7}$ against $3.71\times 10^{-7}$, and this
$21\%$ gap is by itself too small to justify a preference.  We shall see in
\Section{sec:no_willmore_nn_2d} a more striking instance of the same phenomenon:
in the non-oriented case, the model with the lower validation loss is the
one that produces a qualitatively wrong flow.
\end{rem}

\begin{table}[H]
    \centering
    \begin{tabular}{|c|c|c|c|}
        \toprule
        \textbf{Model} &  \textbf{Epochs} & \textbf{Steps} & \textbf{Validation loss} \\
        \midrule
        $\wnn{\theta}{1}$ & $3279$ & $439519$  & $4.73\times 10^{-7}$ \\
        $\wnn{\theta}{2}$ & $1369$ & $183579$  & $3.71\times 10^{-7}$ \\
        \bottomrule
    \end{tabular}
    \caption{Validation losses of the selected checkpoints for the oriented curve models.}
    \label{tab:training-nn}
\end{table}

\subsection{Application for the approximation of the Willmore flow}\label{sec:generalization-test}
In this section we assess the behaviour of the trained model on initial shapes
that do not appear in the training set.  We compare the repeated application of
$\wnn{\theta}{2}$ with the RSS scheme run with the Allen--Cahn model
$\snn{\theta}{1}$, the latter playing the role of the reference since it is
what generated the training data.  We use the Cartesian grid $\mathcal{M}$ with
$N = 256$; for the RSS scheme we take $\tts = \ts / 64$ and $J = 300$, as when
generating the training data.

\paragraph{Smooth flower.} The initial shape for the smooth flower is defined by
\begin{equation*}
    \left\{(r(n, r_0, a, k, \theta)\cos(\theta) + c_x, r(n, r_0, a, k, \theta)\sin(\theta) + c_y)^\top \biggm| \theta \in (0,2\pi]\right\}
\end{equation*}
with $r(n, r_0, a, k, \theta) := r_0 + a\sin^{2k}(n\theta/2)$,
where we set $c_x = c_y := 0.5$, $n := 5$, $r_0 := 0.15$, $a := 0.06$, and $k := 1$.

%

\begin{figure}[H]
  \centering
  \setlength{\tabcolsep}{4pt}   
  \renewcommand{\arraystretch}{1.0} 
  \begin{tabular}{m{0.14\textwidth} m{0.19\textwidth} m{0.19\textwidth} m{0.19\textwidth} m{0.19\textwidth}}
    \raggedleft $RSS$ &
      \includegraphics[width=\linewidth]{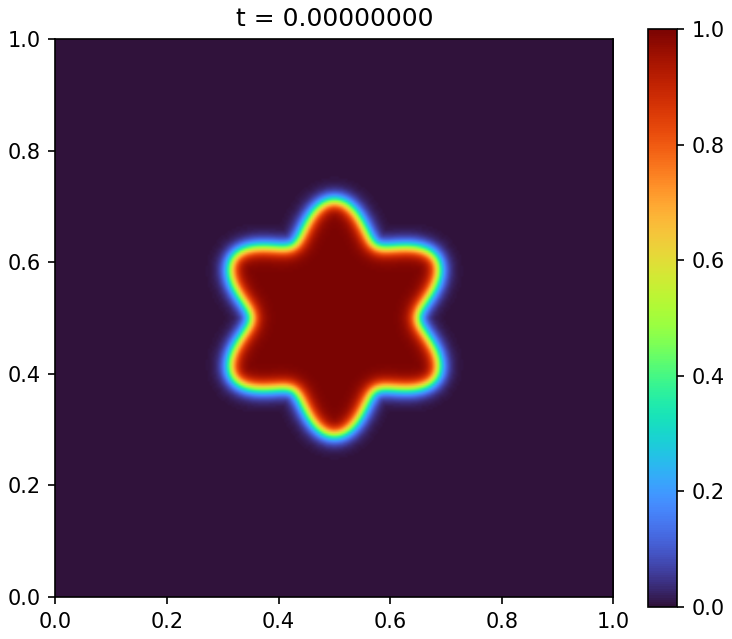} &
      \includegraphics[width=\linewidth]{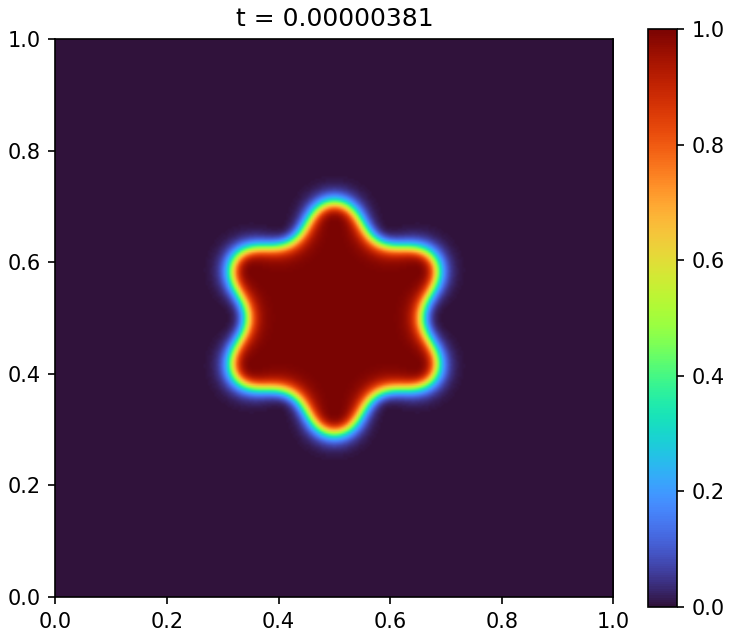} &
      \includegraphics[width=\linewidth]{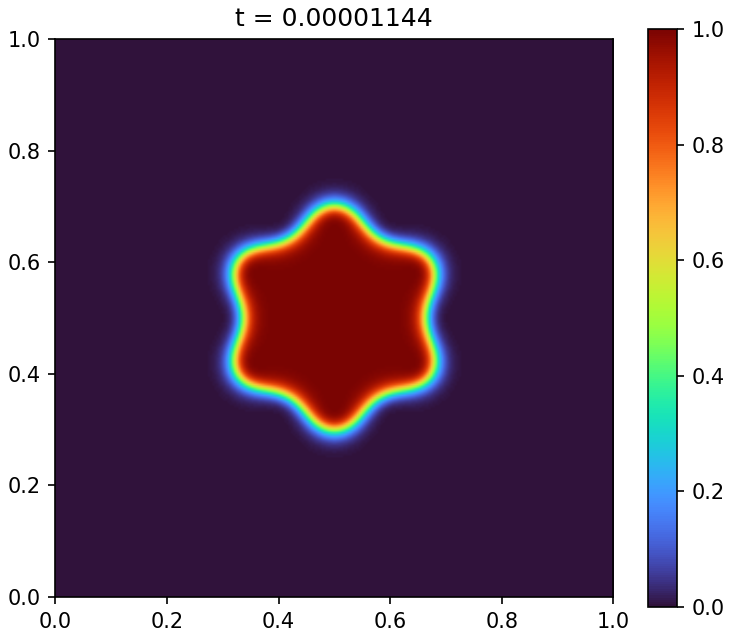} &
      \includegraphics[width=\linewidth]{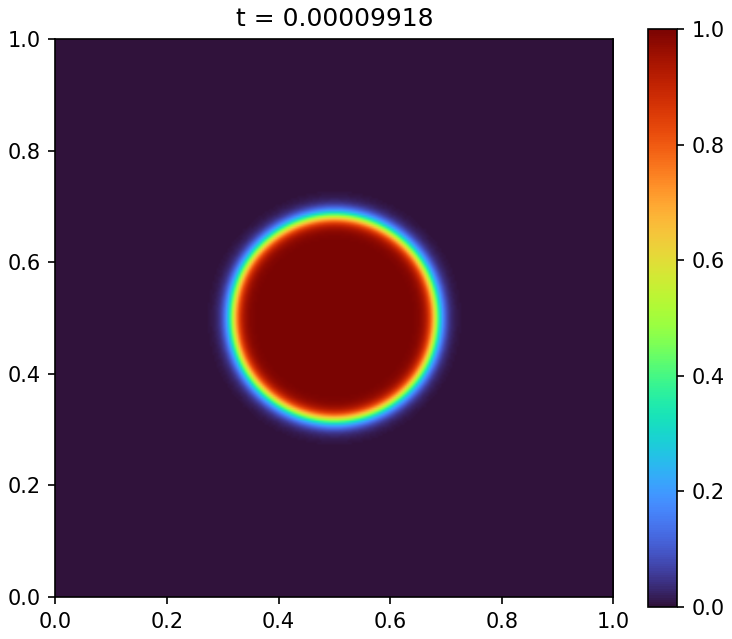} \\
    \raggedleft $\wnn{\theta}{2}$ &
      \includegraphics[width=\linewidth]{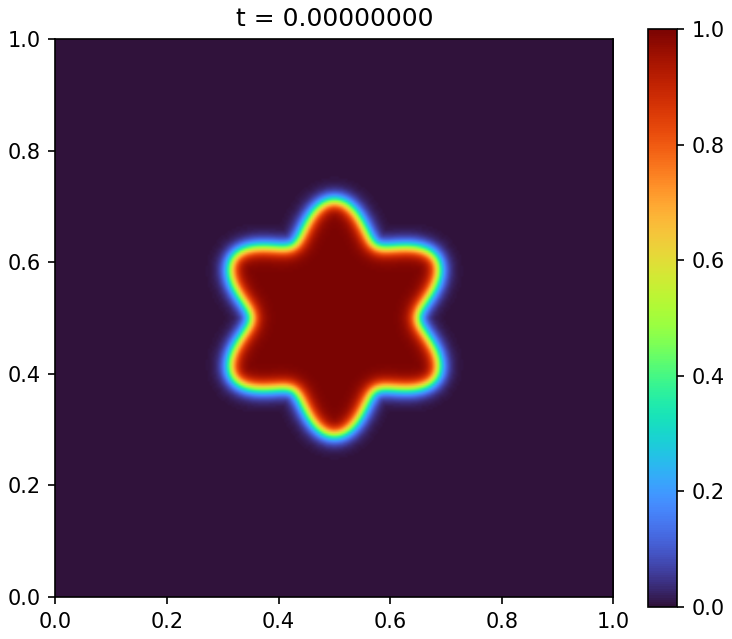} &
      \includegraphics[width=\linewidth]{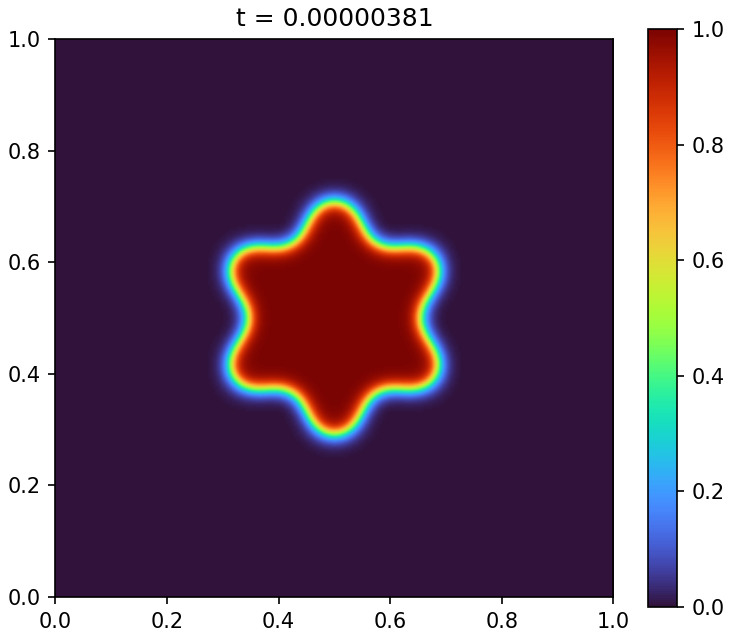} &
      \includegraphics[width=\linewidth]{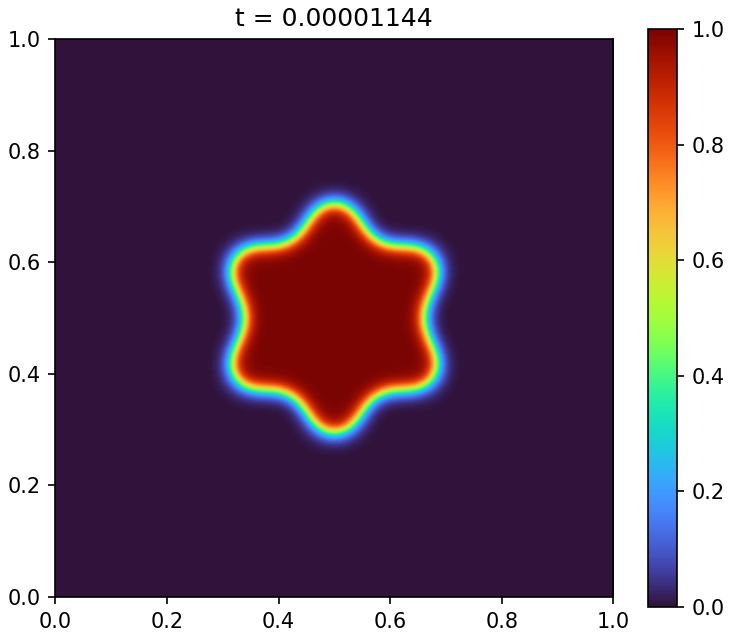} &
      \includegraphics[width=\linewidth]{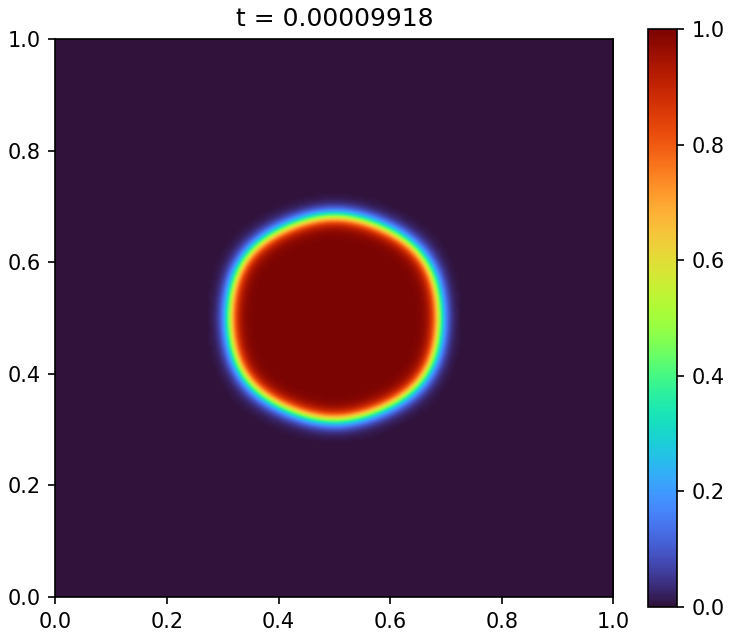} \\
  \end{tabular}
  \caption{Approximation of the Willmore flow at different time steps for the
  smooth flower curve.  First row: RSS scheme; second row: repeated application
  of $\wnn{\theta}{2}$.  }
  \label{fig:smooth_flower6_2d}
\end{figure}
\begin{figure}[H]
    \centering
    \begin{subfigure}[b]{0.70\textwidth}
        \includegraphics[width=\linewidth]{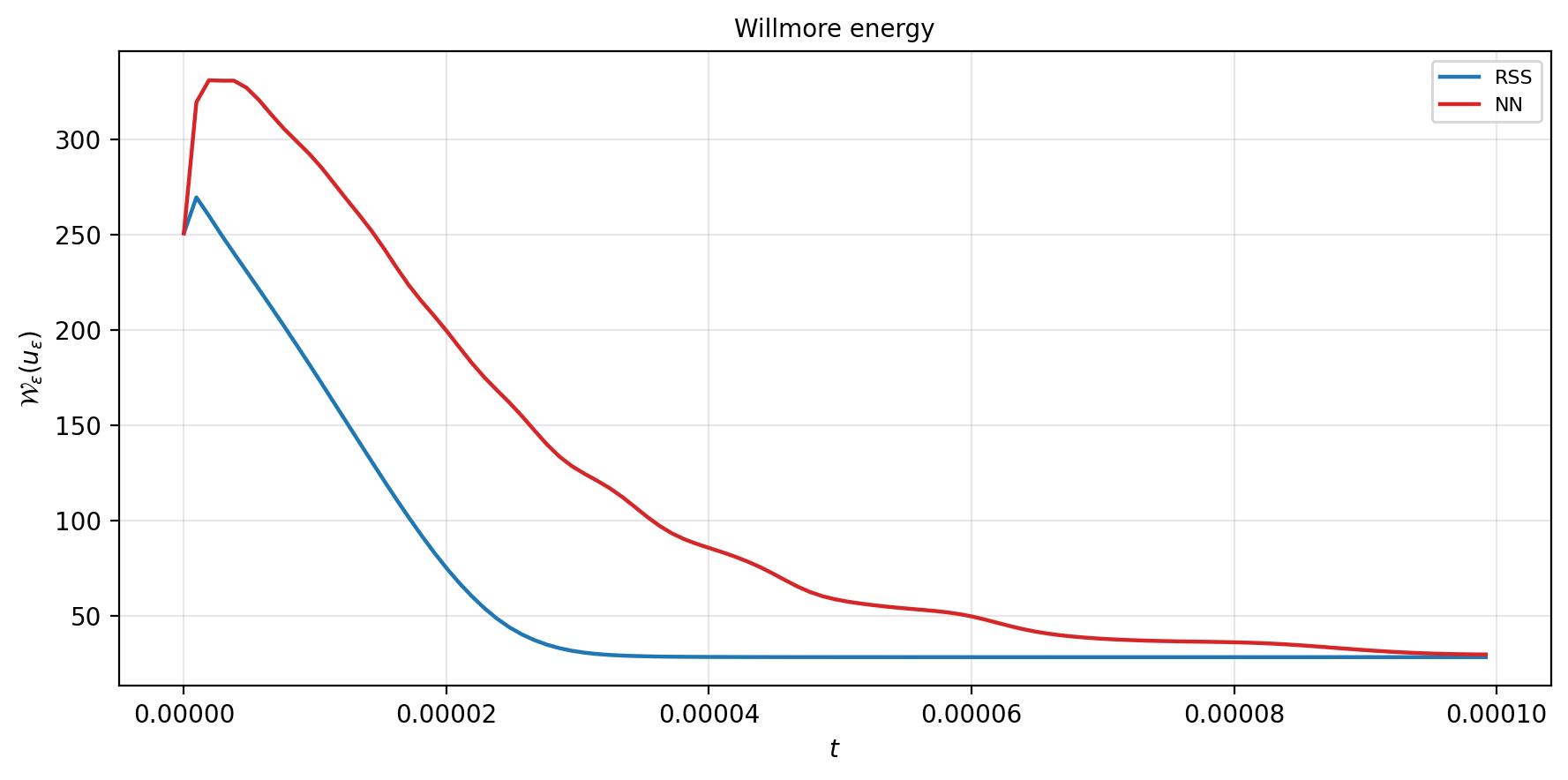}
    \end{subfigure}
    \caption{Plot of the Willmore energy defined by \eqref{eq:classical-willmore-energy} in the case of \Figure{fig:smooth_flower6_2d}.}
\end{figure}

\paragraph{Four circles in a square configuration.} The initial shape is the union of four disks of radius $0.12$ whose centers are located at the vertices of a square: $(0.35,0.35)$, $(0.65,0.35)$, $(0.35,0.65)$, and $(0.65,0.65)$.

\begin{figure}[H]
  \centering
  \setlength{\tabcolsep}{4pt}   
  \renewcommand{\arraystretch}{1.0} 
  \begin{tabular}{m{0.14\textwidth} m{0.19\textwidth} m{0.19\textwidth} m{0.19\textwidth} m{0.19\textwidth}}
    \raggedleft $RSS$ &
      \includegraphics[width=\linewidth]{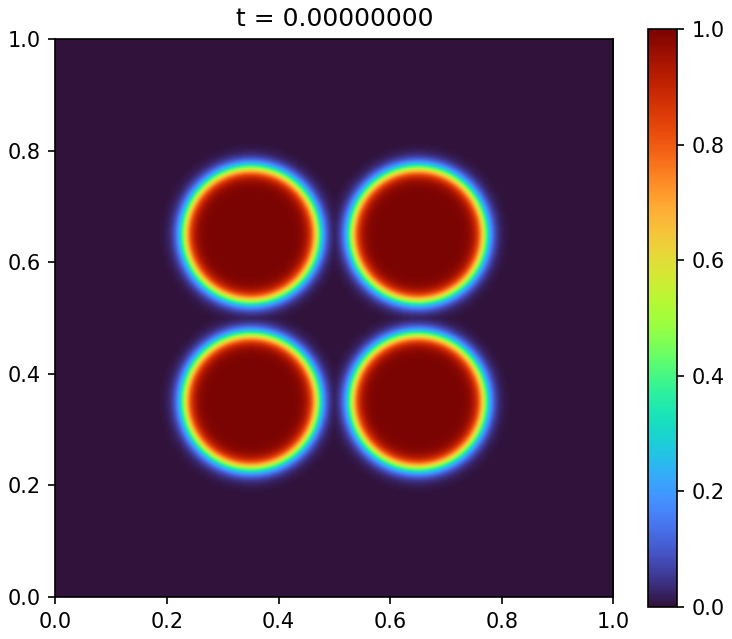} &
      \includegraphics[width=\linewidth]{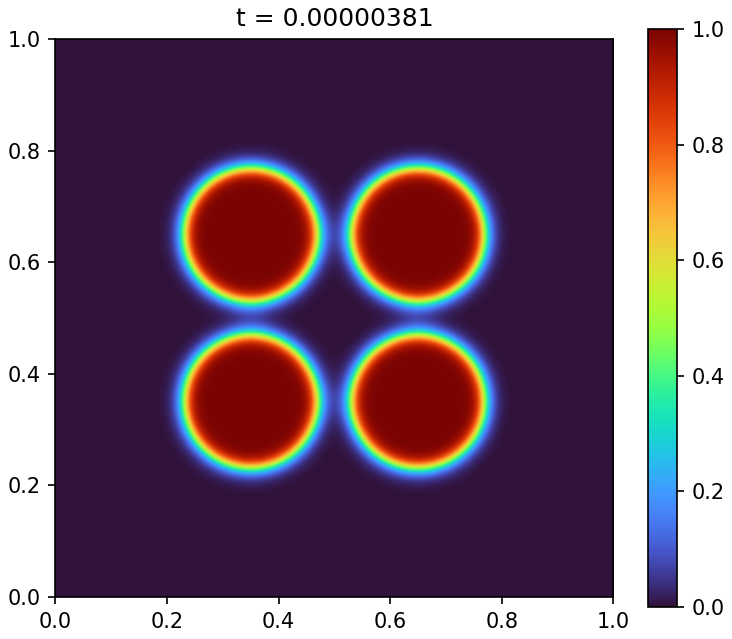} &
      \includegraphics[width=\linewidth]{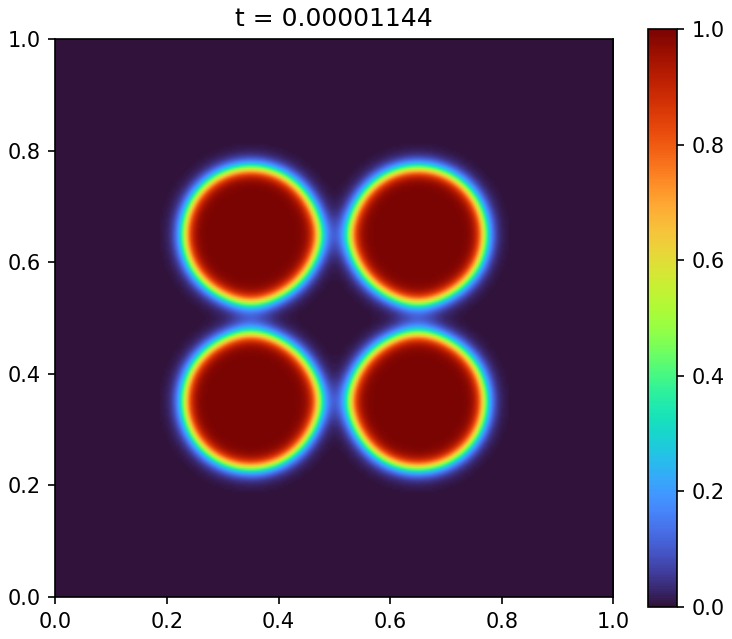} &
      \includegraphics[width=\linewidth]{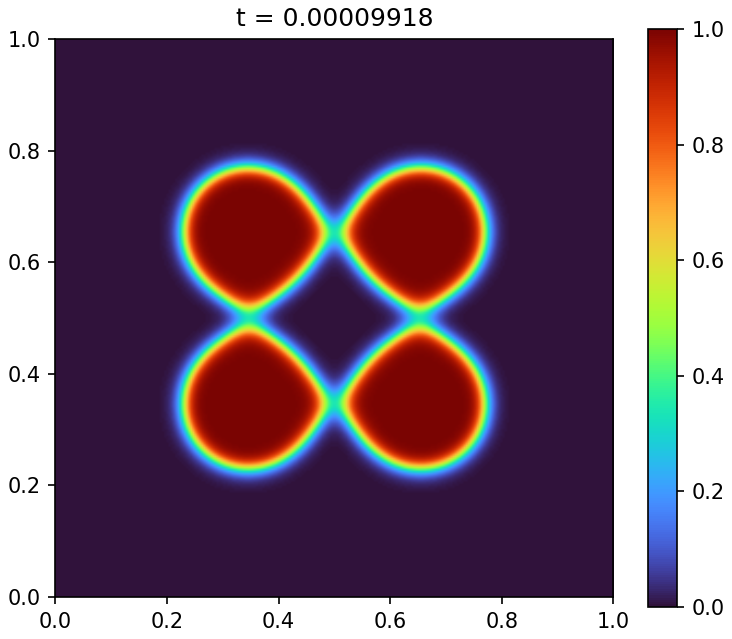} \\
    \raggedleft $\wnn{\theta}{2}$ &
      \includegraphics[width=\linewidth]{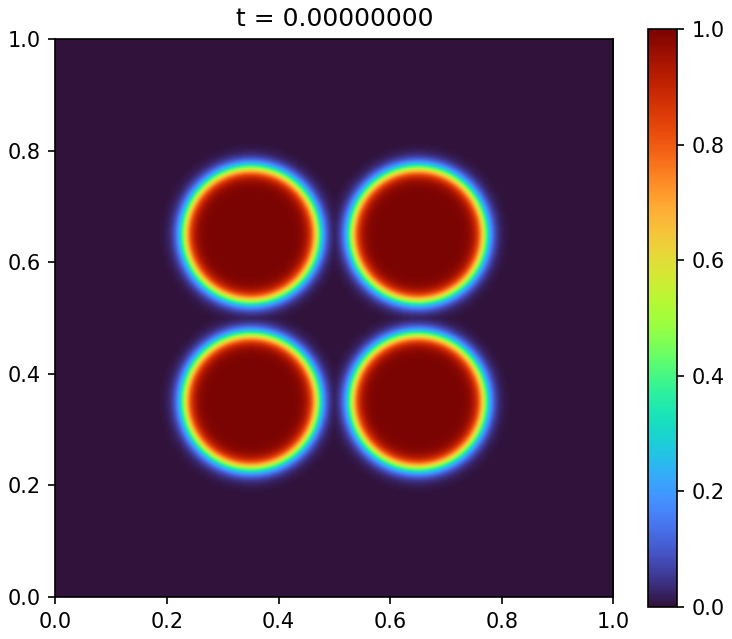} &
      \includegraphics[width=\linewidth]{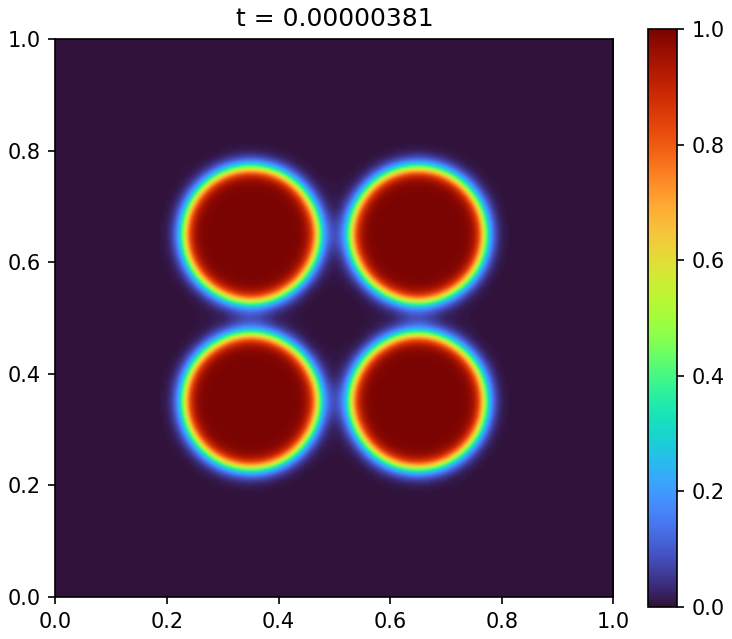} &
      \includegraphics[width=\linewidth]{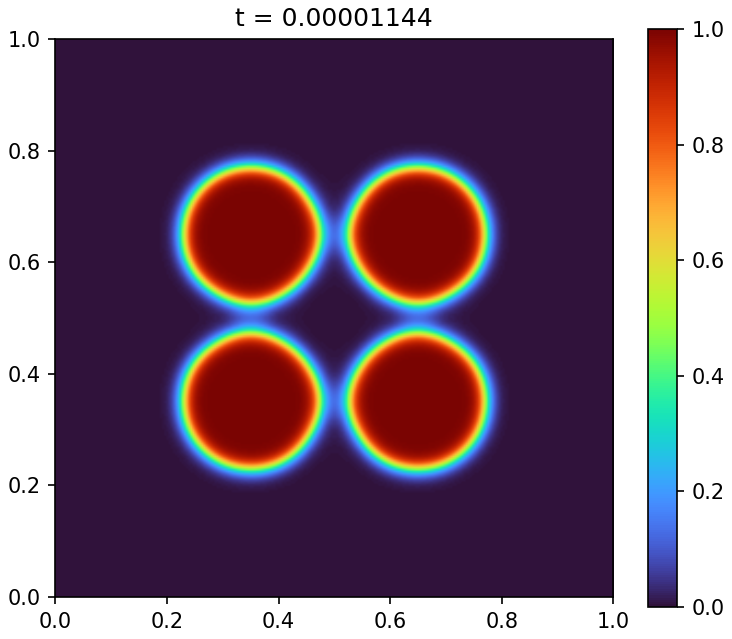} &
      \includegraphics[width=\linewidth]{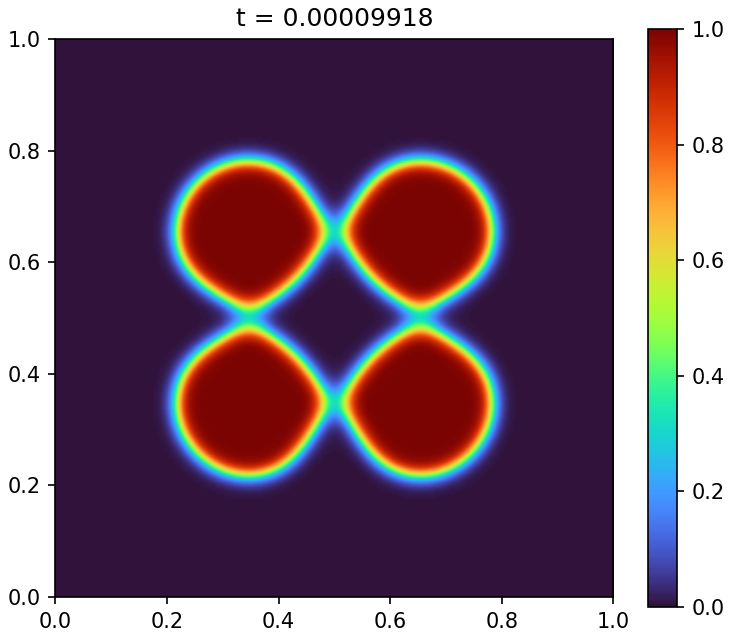} \\
  \end{tabular}
  \caption{Approximation of the Willmore flow at different time steps for four circles in a square configuration.
  First row: RSS scheme, second row: $\wnn{\theta}{2}$.}
  \label{fig:4circles_2d}
\end{figure}

These experiments indicate that the proposed network generalizes to smooth
interfaces as well as to more complex configurations absent from the training
set.  The evolutions are visually close to those produced by the RSS scheme, at
a significantly lower computational cost: the RSS scheme solves an optimization
problem at each time step, which requires $J$ evaluations of
$\snn{\theta}{2}$, whereas one time step of our scheme requires a single
evaluation of $\wnn{\theta}{2}$.  \Table{tab:computation-time-2d-q} compares the
computation times for the smooth flower of \Figure{fig:smooth_flower6_2d}.
\begin{table}[H]
    \centering
    \begin{tabular}{|c|c|c|}
        \toprule
        \textbf{Scheme} & \textbf{Total time (s)} & \textbf{Time / Step (s)}\\
        \midrule
        $RSS$ & $335.0$ & $3.22$ \\
        $\wnn{\theta}{2}$ & $1.03$ & $0.00994$ \\
        \bottomrule
    \end{tabular}
    \caption{Comparison of computation times of the numerical schemes. The times include only implementation of the time-stepping
    (e.g., excluding initialization of the initial phase field functions and other operations like plotting).
    First row: RSS scheme. Second row: Repetitions of $\wnn{\theta}{2}$.}
    \label{tab:computation-time-2d-q}
\end{table}

\begin{rem}
The speed-up reported in \Table{tab:computation-time-2d-q} is not an artefact of
an unnecessarily large value of $J$.  The cost of the RSS scheme is proportional
to $J$, and the timing above was obtained with $J = 300$, the very value used to
generate the training database; decreasing it further degrades the dynamics, as
for the smooth flower, where $J = 30$ produces a markedly slower evolution which
is no longer the Willmore flow.  If anything the comparison understates the
gain: for the circle, \Figure{fig:radius_re_rumpf} and \Table{tab:rerror_rumpf}
indicate that $J$ of the order of $2000$ is needed for the RSS scheme itself to
be accurate at this time step, so a fully converged reference would be several
times more expensive.
\end{rem}
\section{Approximation of  Willmore flow for non-oriented interfaces}
\label{sec:no_willmore_nn_2d}

A notable advantage of the proposed neural network framework is that it can be
applied to non-oriented interfaces.  As in \cite{BretinDenisMasnouTerii2022},
where this was carried out for the mean curvature flow, this is achieved by
replacing the profile $q$ associated with the signed distance function $\sdist{\Omega}$ with the profile $q'$ combined with the unsigned distance function. In the context of the Willmore flow, however, the dataset is unavailable in closed form as it is generated using a numerical scheme. Therefore, the proposed approach relies on exploiting the asymptotic profile of the solutions and performing suitable post-processing of the initial dataset to construct a non-oriented representation of the phase-field function.
More precisely, to generate dataset, we apply the differential operator $\operatorprim : u\mapsto -\varepsilon|\nabla u|$ to the output of the RSS scheme for oriented curves (see \Section{sec:database_construction}).
Indeed, if $u_0$ is the profile function of an oriented curve $\Omega_0$, say $u_0 = q(\sdist{\Omega_0} / \varepsilon)$,
then $\operatorprim[u_0]$ can be explicitly computed as

\begin{equation}\label{eq:initial_data_qprim}
    \operatorprim[u_0] = -\varepsilon \left|\nabla u_0\right|
    = -\varepsilon \left|\frac{1}{\varepsilon}q'\left(\frac{\sdist{\Omega_0}}{\varepsilon}\right)\nabla \sdist{\Omega_0}\right| = q'\left(\frac{\operatorname{dist}(\cdot,\partial \Omega_0)}{\varepsilon}\right).
\end{equation}

This section will show that adopting a non-oriented representation of the interface substantially increases the difficulty of the learning problem. In the experiments below, unconstrained kernels lead to anisotropic dynamics. We introduce radial diffusion kernels to address this behavior and compare the resulting evolutions. Furthermore, we demonstrate that this approach facilitates the treatment of non-oriented configurations involving triple junctions. Finally, we explain how the proposed schemes can be used to create interface approximations from point cloud data.

\subsection{Dataset, training and fast approximation of the Willmore flow}
\begin{figure}[H]
    \centering
    \begin{subfigure}[b]{0.90\textwidth}
        \includegraphics[width=\linewidth]{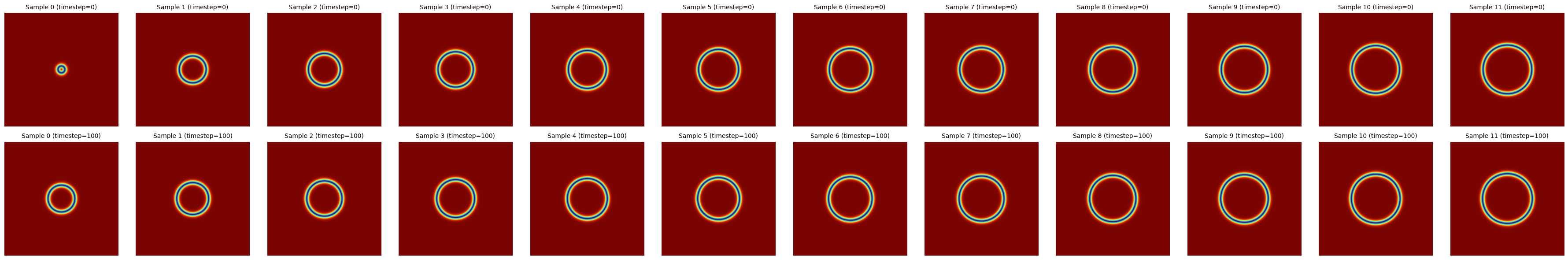}
        \caption{Circles}
    \end{subfigure}
    \hfill
    \begin{subfigure}[b]{0.3\textwidth}
        \includegraphics[width=\linewidth]{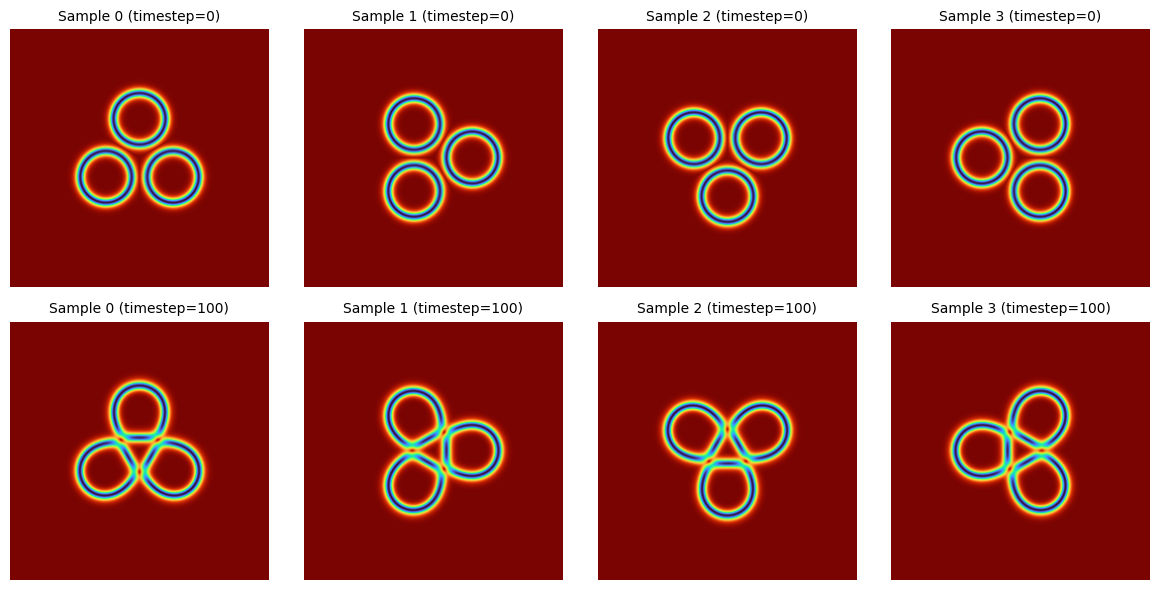}
        \caption{Three circles}
    \end{subfigure}
    \hfill
    \begin{subfigure}[b]{0.3\textwidth}
        \includegraphics[width=\linewidth]{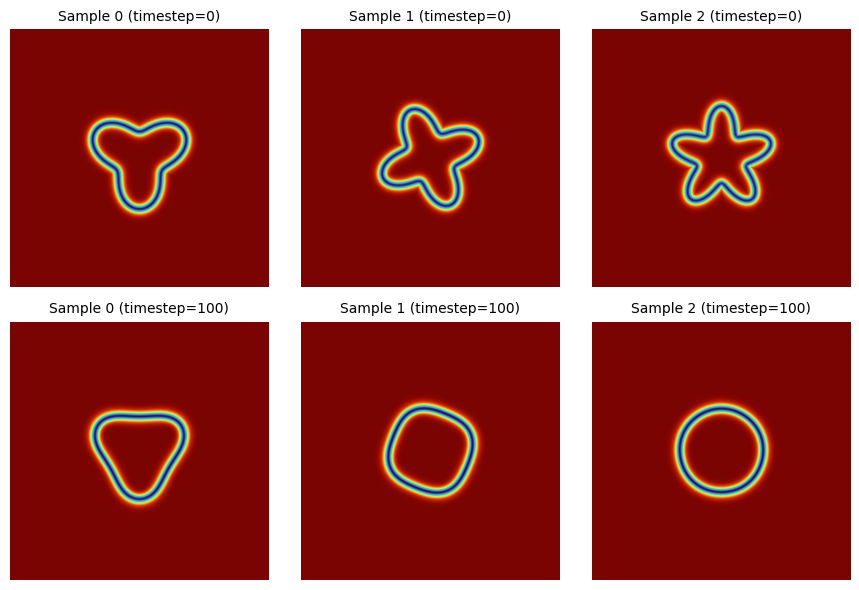}
        \caption{Stars}
    \end{subfigure}
    \hfill
    \begin{subfigure}[b]{0.3\textwidth}
        \includegraphics[width=\linewidth]{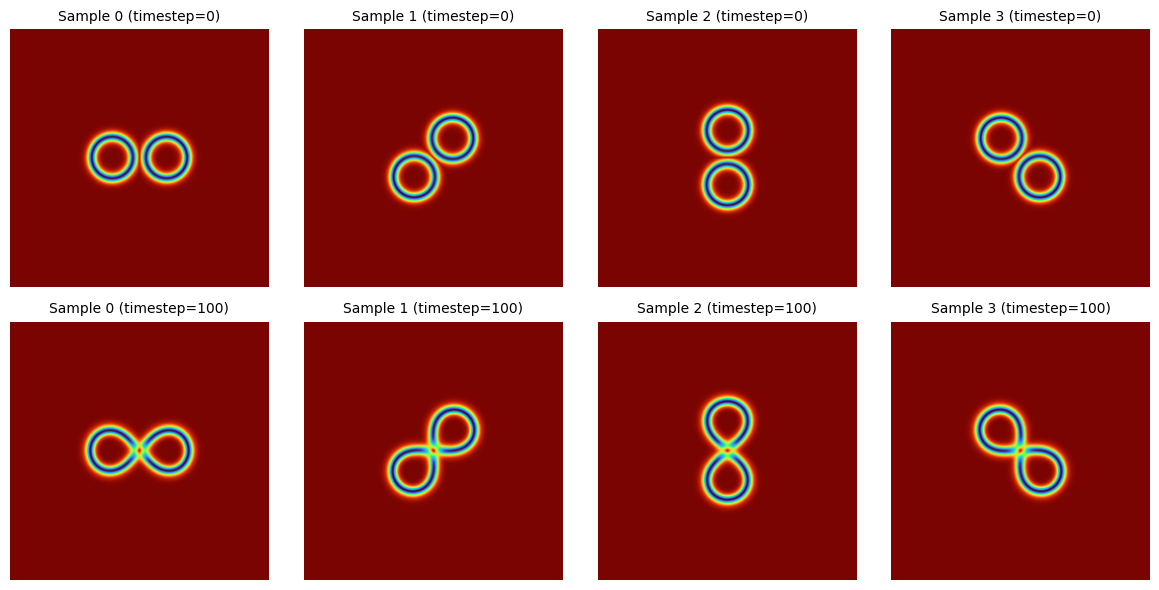}
        \caption{Horizon circles}
    \end{subfigure}
    \caption{Datasets of non-oriented curves}
    \label{fig:datasets_non_oriented_small}
\end{figure}

In the sequel, we set the resolution of the computational domain to $N = 2^8 = 256$,
the interface parameter to $\varepsilon = 2.0/N = 2^{-7}$, and the time step to $\tts = \varepsilon^2/ 64 = 2^{-20}$.  This is the same time step as in the oriented 2D dataset. The non-oriented
profile is concentrated on the interface and decays on both sides.  In the formulas above the
gradient $\nabla u$ of a discrete function $u:\mathcal{M}\to\mathbb{R}$ is
computed by central finite differences.
See \Figure{fig:datasets_non_oriented_small} for the datasets of non-oriented curves.

We now train the networks $\wnnSimple{\theta}{1}$ and $\wnnSimple{\theta}{2}$,
whose structure is recalled in \Figure{fig:nn_structure}, on the dataset of
non-oriented curves.  The training strategy is the one of
\Section{sec:database_construction}, with the multistep loss
\eqref{eq:loss_multistep} with $S = 2$ and a learning rate
$\eta = 5.0\times 10^{-4}$.

\begin{figure}[H]
    \centering
    \begin{subfigure}[b]{0.45\textwidth}
        \includegraphics[width=\linewidth]{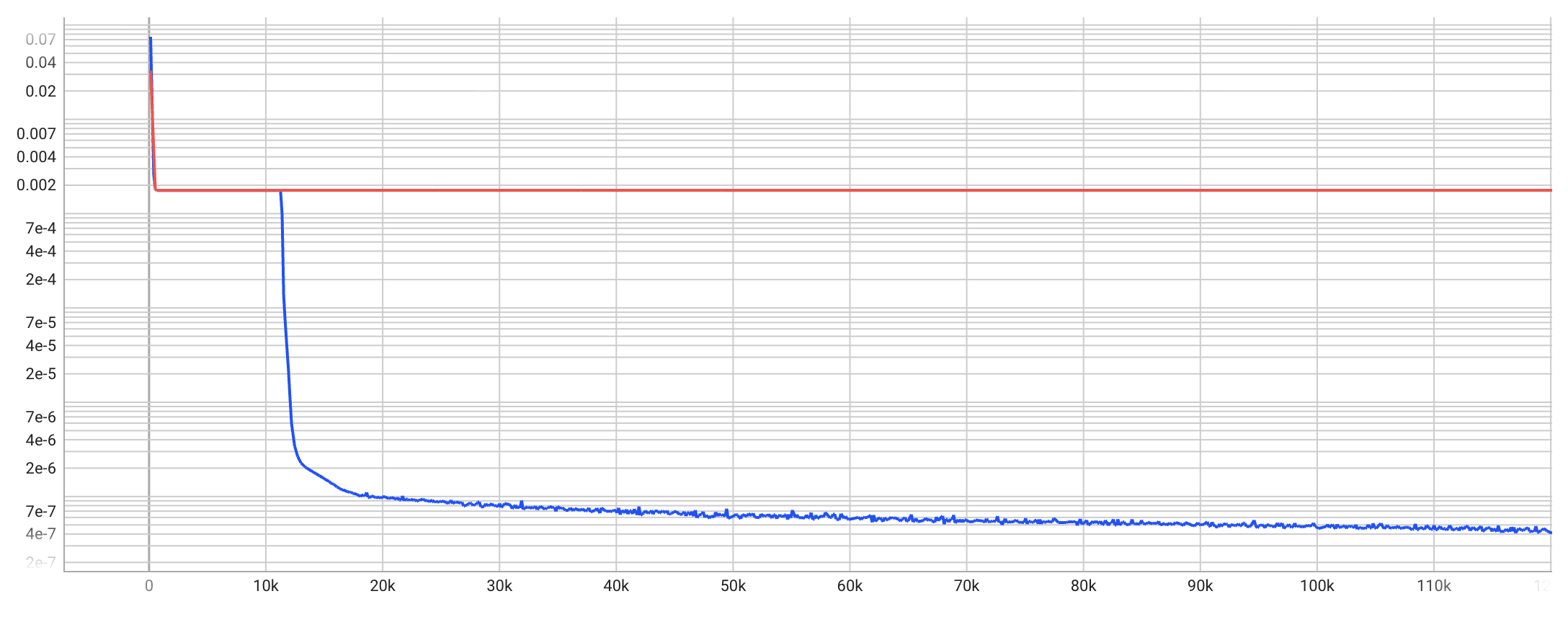}
        \caption{Training losses}
    \end{subfigure}
    \hfill
    \begin{subfigure}[b]{0.45\textwidth}
        \includegraphics[width=\linewidth]{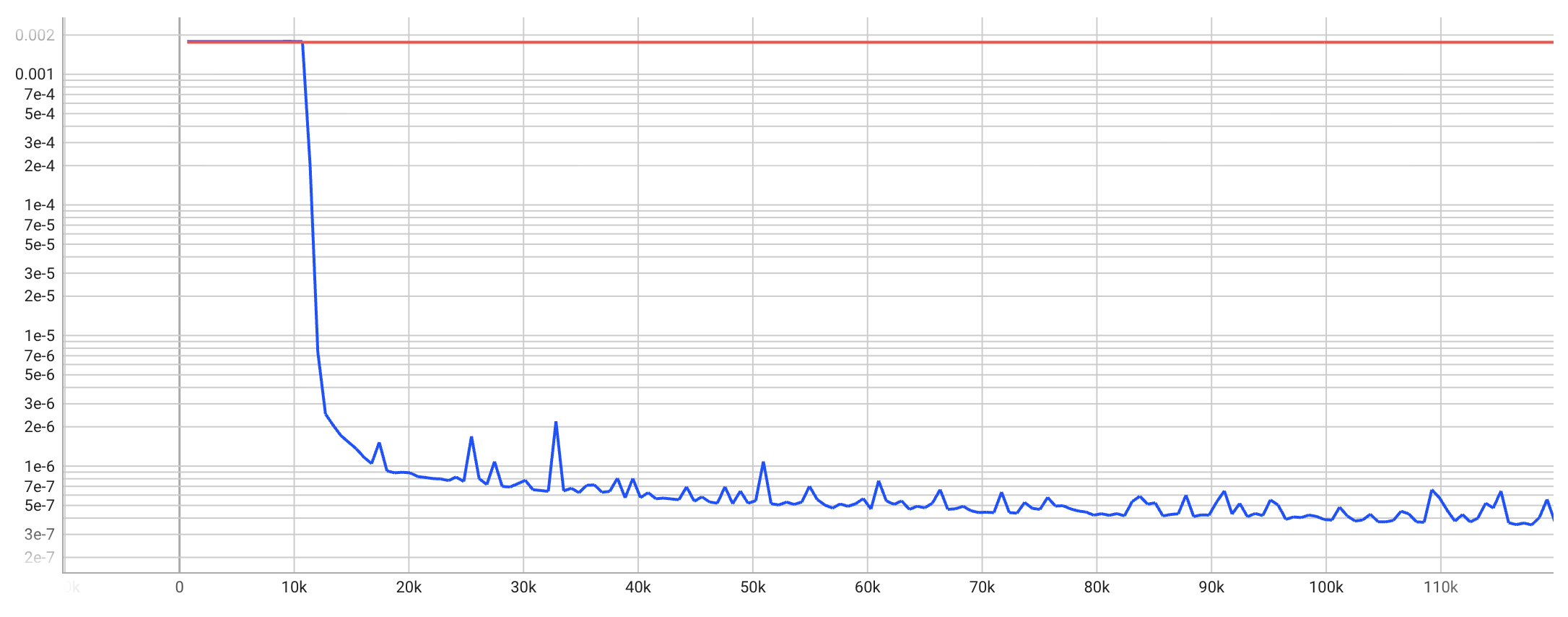}
        \caption{Validation losses}
    \end{subfigure}
    \caption{Evolution of the training (left) and validation (right) losses at each training step for the non-oriented curve models. Red: $\wnnSimple{\theta}{1}$, blue: $\wnnSimple{\theta}{2}$.}
    \label{fig:training_loss_nonoriented_prim}
\end{figure}

\Figure{fig:training_loss_nonoriented_prim} shows that the two architectures
behave very differently here.  The $2$-parallel model $\wnnSimple{\theta}{1}$,
which has two branches and therefore only two of the time parameters $\alpha_i$ of
\Section{sec:one-step-NN} at its disposal, stalls almost immediately: its
training and validation losses remain flat at about $1.8\times 10^{-3}$.  The
$4$-parallel model $\wnnSimple{\theta}{2}$ decreases by four orders of
magnitude.  This is consistent with the expansion of \Section{sec:one-step-NN}:
with two time parameters only, the cancellation of the Cahn--Hilliard term leaves an
error of order $t^{3/2}$, whereas four time parameters allow the next order to be
cancelled as well.  In the oriented case, that difference was not visible
(\Table{tab:training-nn}); in the non-oriented case it is decisive, and we
therefore retain only the $4$-parallel model in what follows.

To distinguish the selected models for the oriented and non-oriented cases,
we write $\wnnprim{\theta}{1}$ and $\wnnprim{\theta}{2}$ for the selected models for the non-oriented case.

\begin{table}[H]
    \centering
    \begin{tabular}{|c|c|c|c|}
        \toprule
        \textbf{Model} &  \textbf{Epochs} & \textbf{Steps} & \textbf{Validation loss} \\
        \midrule
        $\wnnprim{\theta}{1}$ & $864$ & $115909$  & $1.76\times 10^{-3}$ \\
        $\wnnprim{\theta}{2}$ & $389$ & $52259$  & $5.04\times 10^{-7}$ \\
        \bottomrule
    \end{tabular}
    \caption{Validation loss of the selected checkpoints for the non-oriented curve models.}
    \label{tab:training-nn-nonoriented}
\end{table}

\begin{figure}[H]
  \centering
  \setlength{\tabcolsep}{4pt}   
  \renewcommand{\arraystretch}{1.0} 
  \begin{tabular}{m{0.14\textwidth} m{0.19\textwidth} m{0.19\textwidth} m{0.19\textwidth} m{0.19\textwidth}}
    \raggedleft $\operatorprim\circ RSS$ &
      \includegraphics[width=\linewidth]{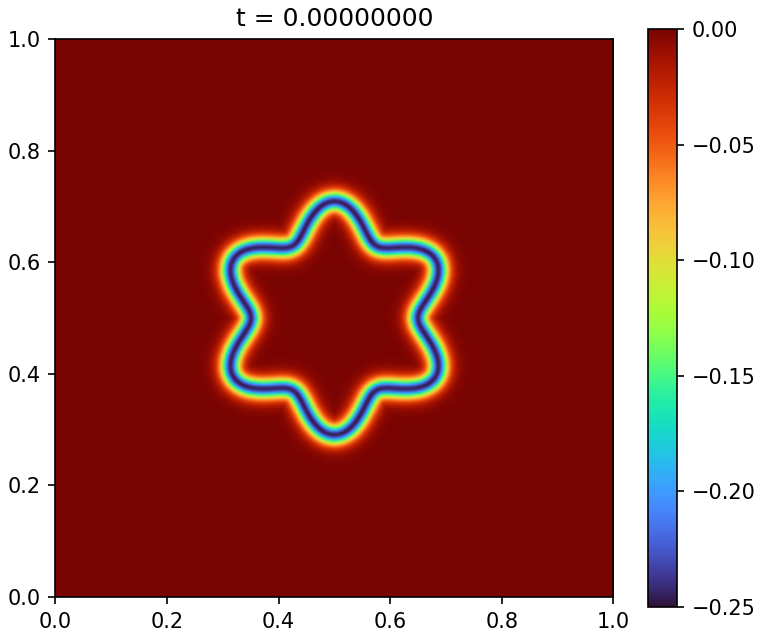} &
      \includegraphics[width=\linewidth]{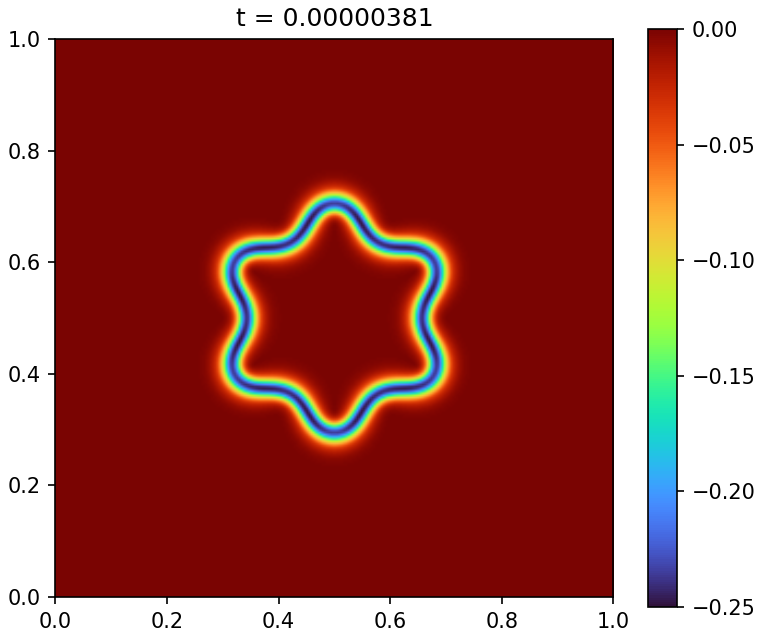} &
      \includegraphics[width=\linewidth]{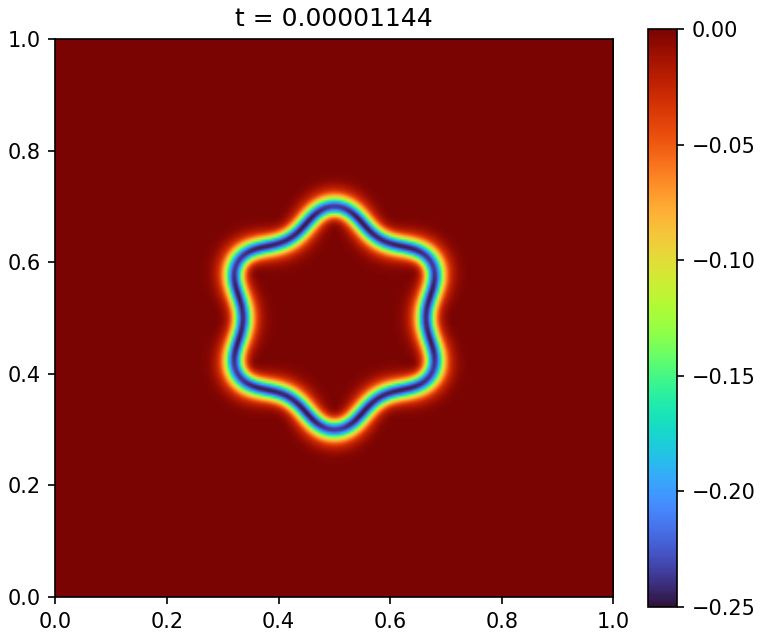} &
      \includegraphics[width=\linewidth]{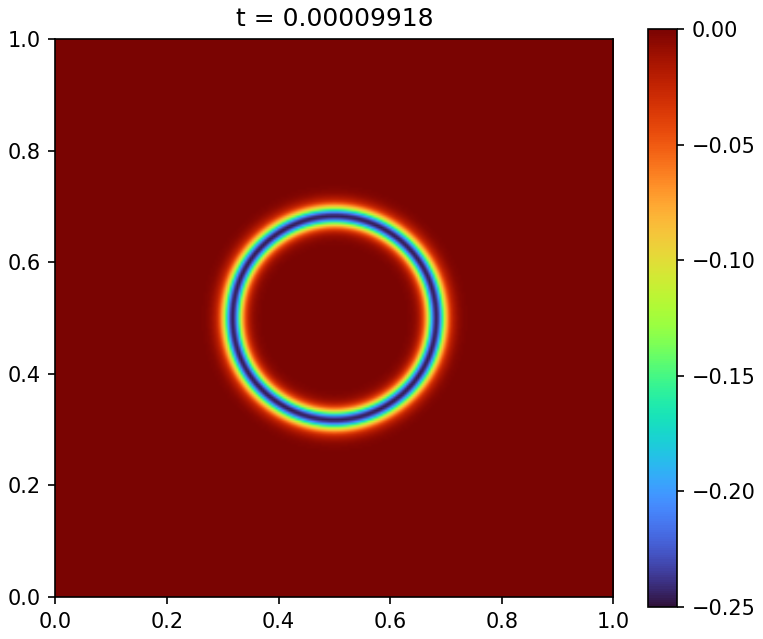} \\
    \raggedleft $\wnnprim{\theta}{2}$ &
      \includegraphics[width=\linewidth]{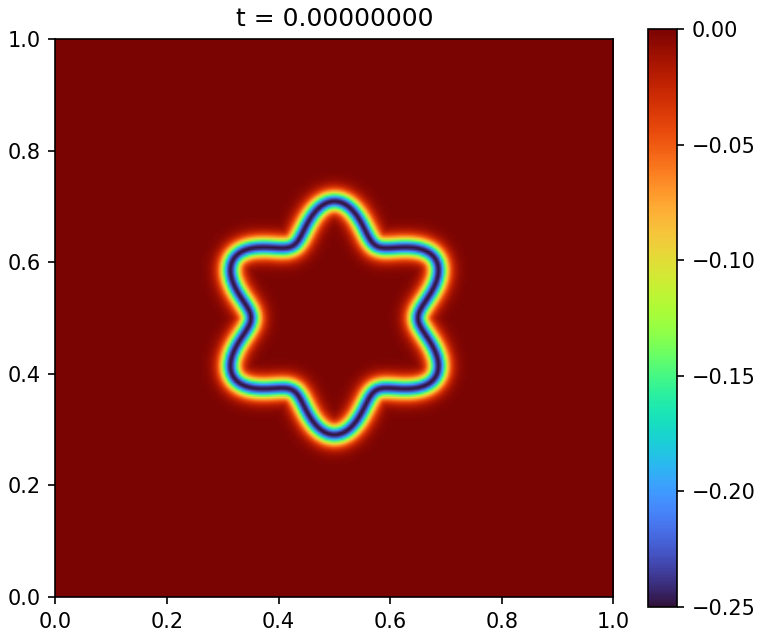} &
      \includegraphics[width=\linewidth]{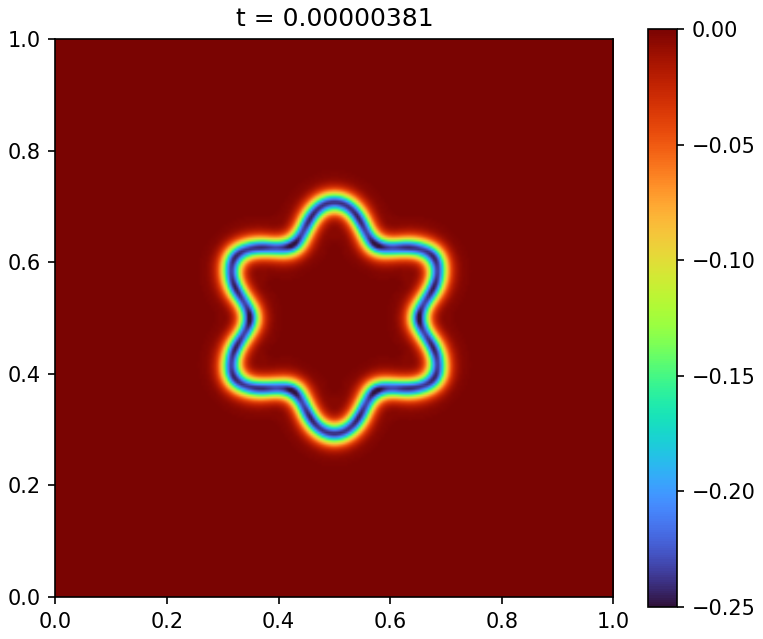} &
      \includegraphics[width=\linewidth]{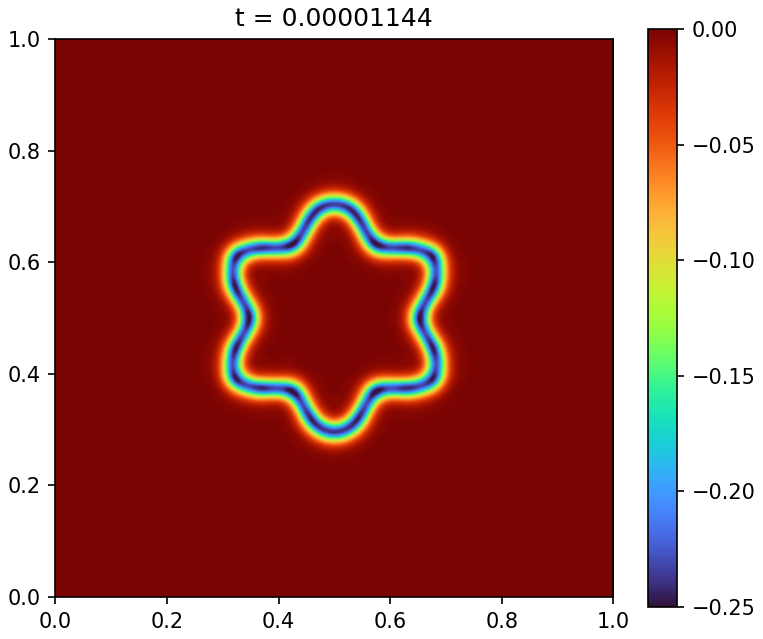} &
      \includegraphics[width=\linewidth]{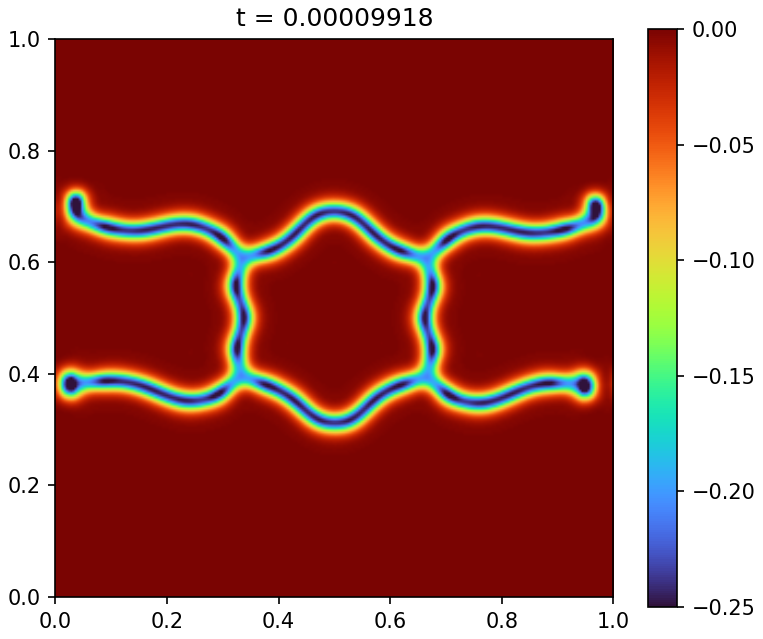} \\
  \end{tabular}
  \caption{Approximation of the non-oriented Willmore flow at different time steps for a smooth flower-shaped initial curve.
First row: $\operatorprim[u^n]$ computed using the RSS scheme.
Second row: $u^n$ computed using the neural network $\wnnprim{\theta}{2}$.}
  \label{fig:smooth_flower6_prim_2d}
\end{figure}

We present now the result of a generalization test for $\wnnprim{\theta}{2}$ on a non-oriented curve.
We compare the prediction of the kernel model $\wnnprim{\theta}{2}$
with the reference given by the RSS scheme \eqref{eq:adam} for oriented curves
followed by the discrete operator $\operatorprim$, a combination that we denote
by $\operatorprim\circ RSS$ in the sequel.  As explained above, the narrow
model $\wnnprim{\theta}{1}$ is not tested, since it does not reach a sufficient
accuracy during training.

\Figure{fig:smooth_flower6_prim_2d} shows the outcome, the initial shape being
the one defined in \Section{sec:generalization-test}.  The free kernel model
$\wnnprim{\theta}{2}$ fails, and it fails in a specific way that is worth
describing precisely.  Over the first steps the evolution is plausible, but the
curve does not contract isotropically: the horizontal and vertical directions
are favored over the diagonal ones, although the target model
\eqref{eq:willmore} is isotropic.  The defect then amplifies with the number of
applications of the operator, and in the last frame the curve has been completely destroyed and replaced by a spurious network of oscillating interfaces spanning the whole domain.  What we observe is therefore not merely
an anisotropic approximation of the Willmore flow but an instability, triggered
by the anisotropy and amplified by the repeated application of the operator.
We stress that this happens for a \emph{smooth} initial curve, and with a model
whose validation loss is low (see \Table{tab:training-nn-nonoriented}).

\begin{figure}[H]
    \centering
    \begin{subfigure}[b]{0.5\textwidth}
        \centering
        \includegraphics[width=\linewidth]{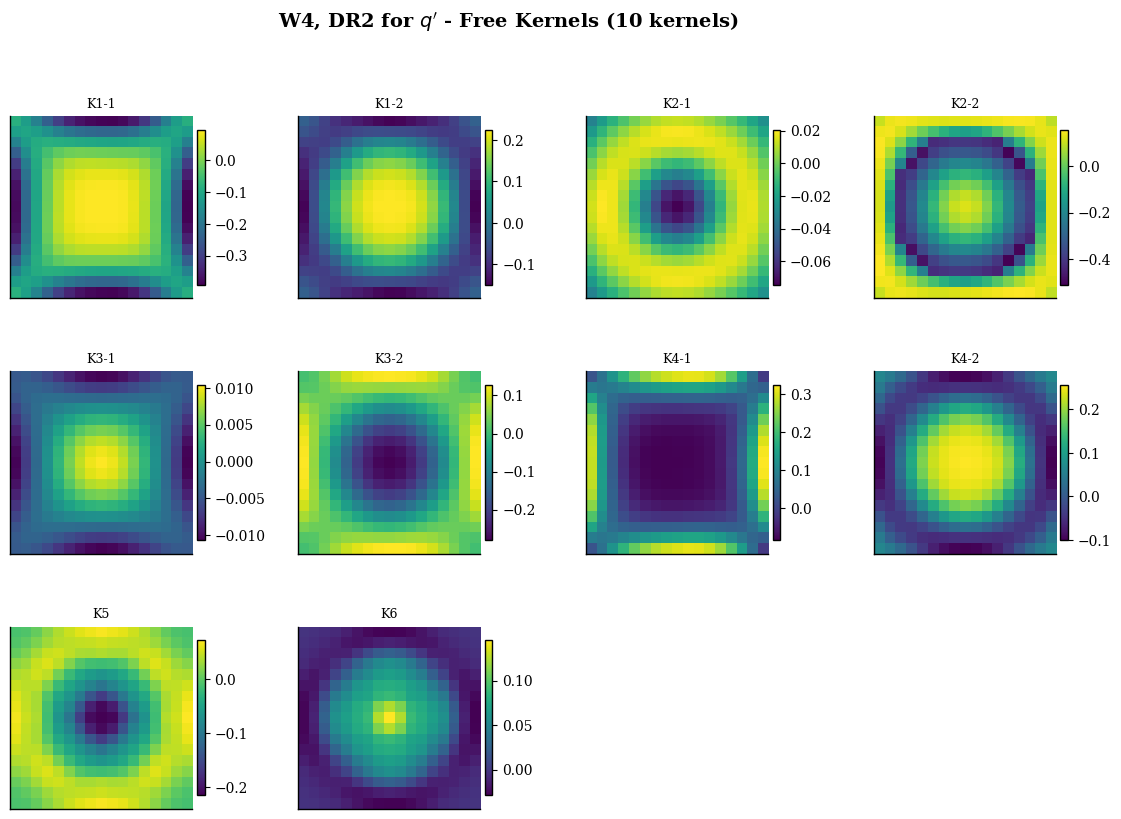}
    \end{subfigure}
    \hfill
    \caption{Representation of the $10$ kernels trained and used in our neural network  $\wnnprim{\theta}{2}$.}
    \label{fig:learned-kernels-w4dr2-qprim}
    \end{figure}

The origin of the defect is visible in the learned kernels.
\Figure{fig:learned-kernels-w4dr2-qprim} displays the ten diffusion kernels of
$\wnnprim{\theta}{2}$: five of them are clearly anisotropic, with a marked
preference for the axes of the grid.  Nothing in the loss penalizes this, since
the training shapes -- disks, unions of disks, stars -- are themselves aligned
with the grid often enough for an axis-aligned kernel to fit them well.  To
prevent the network from learning such a direction-dependent dynamics, we
constrain the architecture in the next section.

\subsection{Enforcing isotropy: rotationally invariant diffusion kernels}\label{sec:IsotropicNN_MCF}
We return to the continuous convolution \eqref{eq:convolution_continuous}
introduced in \Section{sec:intro_nn}, whose discrete counterpart
\eqref{eq:convolution_2d} is used in the diffusion layers of $S_\theta$.
We first impose rotational invariance on the continuous kernel and then
transpose this construction to the discrete setting.
Since the desired flow \eqref{eq:willmore} is isotropic, all diffusion kernels in $S_\theta$
should be independent of the direction (in other words, rotationally invariant).
Following this observation in the continuous setting, we shall enforce the invariance of direction of mean curvature flow by replacing the diffusion layers in $\snnSimple{\theta}{1}$ and $\snnSimple{\theta}{2}$ with radially invariant kernels.
To this end, we suppose that for each diffusion layer $\mathcal{D}$, there exists a radial function $\phi$ such that

\begin{equation*}
    \mathcal{D}(x) = \phi\left(\frac{\|x\|}{r}\right)\qquad \mbox{for}\quad x\in\R^2,
\end{equation*}
where $r > 0$ denotes the radius of the kernel.
Namely, the value of the convolutional filter $\mathcal{D}$ at $x$ depends only on the normalized distance from the origin.
Let us write such a diffusion layer as $\mathcal{D}_{iso}$.
We now transpose this to the discrete setting.
Each $\mathcal{D}_{iso}$ can be regarded as a discrete 2D function which is characterized by a radial function
$\phi : \R_{\geq 0}\to\mathbb{R}$ instead of a convolutional filter $\phi : \mathbb{R}^2\to\mathbb{R}$.
In other words, the function $\phi$ is expected to be a function of the normalized distance from the center of the kernel,
and the convolutional filter is recovered by applying $\phi$ radially.
Precisely speaking, for each index $(i,j)\in\{-\Krad,\ldots,\Krad\}^2$ the component $\mathcal{D}_{ij}$ in $\mathcal{D}_{iso}$ is defined by

\begin{equation}\label{eq:diffusion_kernel-2D}
    \mathcal{D}_{ij} := \phi\left(\frac{\sqrt{i^2 + j^2}}{\Krad\sqrt{2}}\right).
\end{equation}
The kernel has size $(2\Krad+1)\times(2\Krad+1)$, and
$\Krad\sqrt{2}$ is the distance, in grid units, from its center to a corner.
As in \Section{sec:intro_nn}, we use $\Krad=8$ in practice.

Thus, we now aim to learn the radial function $\phi$ instead of the full
convolution filter $\mathcal{D}$, and $\phi$ is approximated by a multilayer
perceptron (MLP) with a single input and a single output.

\begin{figure}[H]
    \centering
    \begin{tikzpicture}[font=\small,
        neuron/.style={circle,draw=black!75,fill=white,minimum size=3.2mm,inner sep=0pt},
        connection/.style={draw=black!30,line width=0.35pt}]
        \coordinate (p0-0) at (0.000,0.000);
        \coordinate (p1-0) at (2.650,0.600);
        \coordinate (p1-1) at (2.650,0.200);
        \coordinate (p1-2) at (2.650,-0.200);
        \coordinate (p1-3) at (2.650,-0.600);
        \coordinate (p2-0) at (5.300,1.400);
        \coordinate (p2-1) at (5.300,1.000);
        \coordinate (p2-2) at (5.300,0.600);
        \coordinate (p2-3) at (5.300,0.200);
        \coordinate (p2-4) at (5.300,-0.200);
        \coordinate (p2-5) at (5.300,-0.600);
        \coordinate (p2-6) at (5.300,-1.000);
        \coordinate (p2-7) at (5.300,-1.400);
        \coordinate (p3-0) at (7.950,0.000);
        \draw[connection] (p0-0) -- (p1-0);
        \draw[connection] (p0-0) -- (p1-1);
        \draw[connection] (p0-0) -- (p1-2);
        \draw[connection] (p0-0) -- (p1-3);
        \draw[connection] (p1-0) -- (p2-0);
        \draw[connection] (p1-0) -- (p2-1);
        \draw[connection] (p1-0) -- (p2-2);
        \draw[connection] (p1-0) -- (p2-3);
        \draw[connection] (p1-0) -- (p2-4);
        \draw[connection] (p1-0) -- (p2-5);
        \draw[connection] (p1-0) -- (p2-6);
        \draw[connection] (p1-0) -- (p2-7);
        \draw[connection] (p1-1) -- (p2-0);
        \draw[connection] (p1-1) -- (p2-1);
        \draw[connection] (p1-1) -- (p2-2);
        \draw[connection] (p1-1) -- (p2-3);
        \draw[connection] (p1-1) -- (p2-4);
        \draw[connection] (p1-1) -- (p2-5);
        \draw[connection] (p1-1) -- (p2-6);
        \draw[connection] (p1-1) -- (p2-7);
        \draw[connection] (p1-2) -- (p2-0);
        \draw[connection] (p1-2) -- (p2-1);
        \draw[connection] (p1-2) -- (p2-2);
        \draw[connection] (p1-2) -- (p2-3);
        \draw[connection] (p1-2) -- (p2-4);
        \draw[connection] (p1-2) -- (p2-5);
        \draw[connection] (p1-2) -- (p2-6);
        \draw[connection] (p1-2) -- (p2-7);
        \draw[connection] (p1-3) -- (p2-0);
        \draw[connection] (p1-3) -- (p2-1);
        \draw[connection] (p1-3) -- (p2-2);
        \draw[connection] (p1-3) -- (p2-3);
        \draw[connection] (p1-3) -- (p2-4);
        \draw[connection] (p1-3) -- (p2-5);
        \draw[connection] (p1-3) -- (p2-6);
        \draw[connection] (p1-3) -- (p2-7);
        \draw[connection] (p2-0) -- (p3-0);
        \draw[connection] (p2-1) -- (p3-0);
        \draw[connection] (p2-2) -- (p3-0);
        \draw[connection] (p2-3) -- (p3-0);
        \draw[connection] (p2-4) -- (p3-0);
        \draw[connection] (p2-5) -- (p3-0);
        \draw[connection] (p2-6) -- (p3-0);
        \draw[connection] (p2-7) -- (p3-0);
        \node[neuron] at (p0-0) {};
        \node[align=center] at (0.000,2.10) {Input\\1 value};
        \node[align=center,anchor=north] at (0.000,-1.85) {$\rho$};
        \node[neuron] at (p1-0) {};
        \node[neuron] at (p1-1) {};
        \node[neuron] at (p1-2) {};
        \node[neuron] at (p1-3) {};
        \node[align=center] at (2.650,2.10) {Layer 1\\4 neurons};
        \node[align=center,anchor=north] at (2.650,-1.85) {ReLU};
        \node[neuron] at (p2-0) {};
        \node[neuron] at (p2-1) {};
        \node[neuron] at (p2-2) {};
        \node[neuron] at (p2-3) {};
        \node[neuron] at (p2-4) {};
        \node[neuron] at (p2-5) {};
        \node[neuron] at (p2-6) {};
        \node[neuron] at (p2-7) {};
        \node[align=center] at (5.300,2.10) {Layer 2\\8 neurons};
        \node[align=center,anchor=north] at (5.300,-1.85) {ReLU};
        \node[neuron] at (p3-0) {};
        \node[align=center] at (7.950,2.10) {Layer 3\\1 neuron};
        \node[align=center,anchor=north] at (7.950,-1.85) {Affine\\$\phi_\theta(\rho)$};
        \node at (3.975,-2.85) {ReLU activation: $\relu(s)=\max\{0,s\}$};
    \end{tikzpicture}
    \caption{Architecture of the MLP $\phi_\theta$ used to approximate the radial profile $\phi$: three affine layers with ReLU
    activation after the first two layers ($57$ learnable parameters).
    Within each diffusion block, the same network is applied independently
    to the normalized radius of every kernel position, with shared parameters.}
    \label{fig:radial_profile_mlp}
\end{figure}
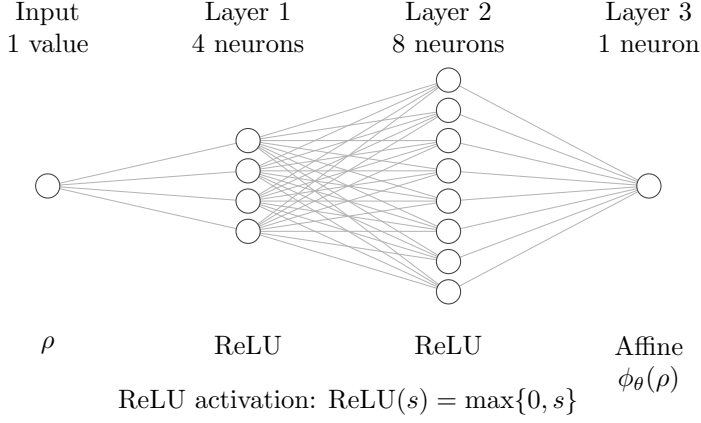

We note that constraining a convolution kernel in order to make the layer
equivariant under a group of symmetries is a well-established idea in the
machine learning literature, where it is known as group equivariant or
steerable convolution \cite{cohen2016group,weiler2019general};
radial kernels are the simplest instance of it for the rotation group.  What is specific here
is not the construction of the kernels but the role the constraint plays. In our non-oriented experiments, the constraint corrects substantial
anisotropic distortions of the learned evolution, rather than merely
improving the training loss.  We
also note that, since the kernel is ultimately sampled on a Cartesian
stencil, the resulting discrete operator is equivariant only up to the
discretization error, a point we come back to at the end of this section.
The radial profile $\phi$ is represented by a small MLP $\phi_\theta$ whose
layer widths and ReLU activations are shown in
\Figure{fig:radial_profile_mlp}. Within each diffusion block, the same
scalar map is evaluated pointwise at the normalized radius of each
kernel position, with shared parameters across the stencil. Positions
with the same radius therefore receive the same coefficient. This MLP
constructs the kernel, the resulting diffusion layer acts on the phase
field through convolution.

\begin{figure}[H]
    \centering
    \begin{subfigure}[b]{0.45\textwidth}
        \centering
        \includegraphics[width=\linewidth]{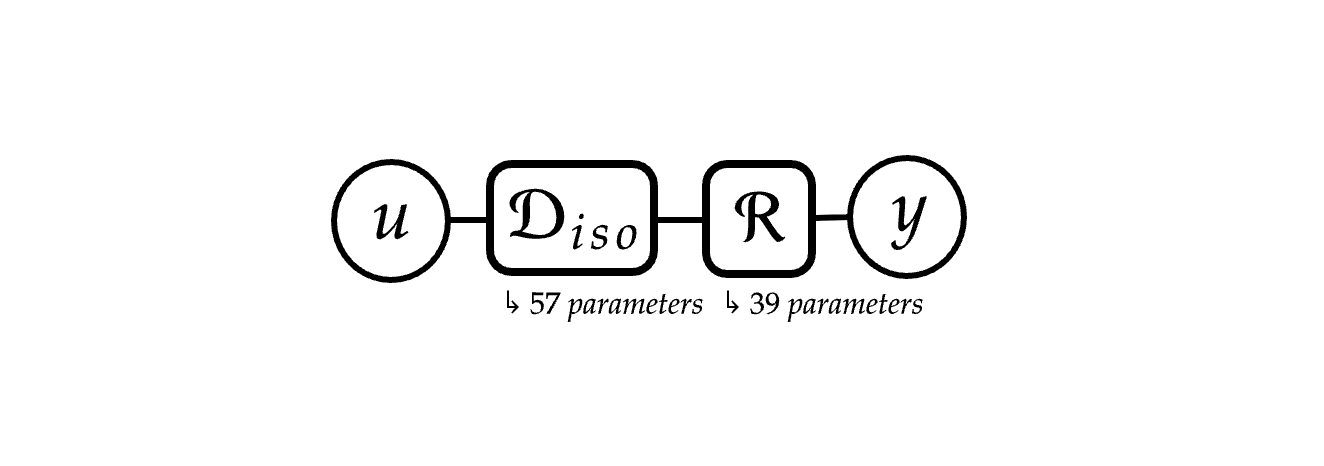}
        \caption{$\spnnSimple{\theta}{1}$ ($96$ learnable parameters)}
    \end{subfigure}
    \begin{subfigure}[b]{0.45\textwidth}
        \centering
        \includegraphics[width=\linewidth]{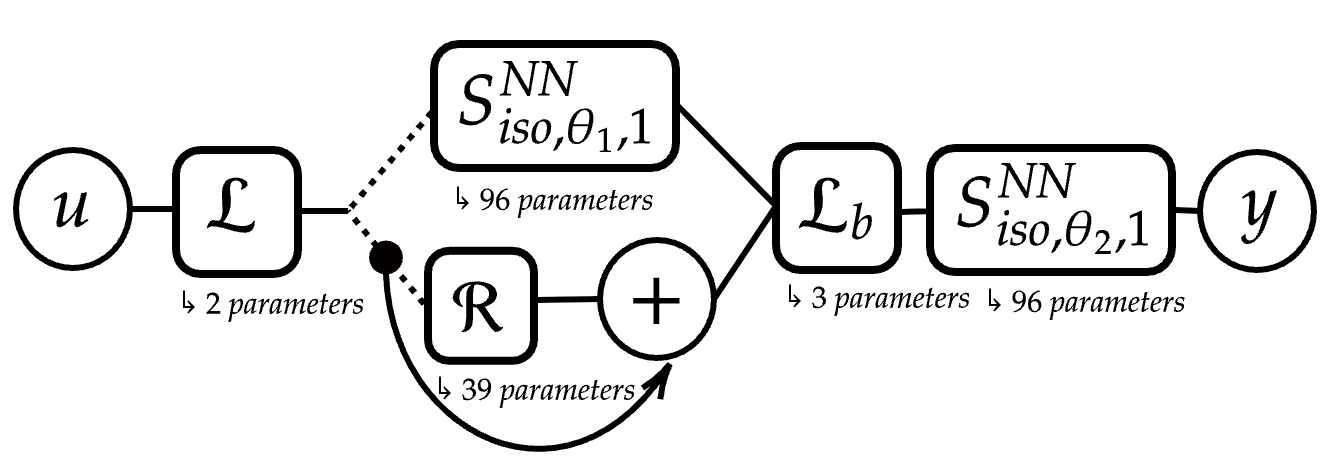}
        \caption{$\spnnSimple{\theta}{2}$ ($236$ learnable parameters)}
    \end{subfigure}
    \caption{
        The network structures of first (left) and second (right) order models with rotationally invariant kernels to approximate the phase-field representation of the mean curvature flow.
        The number of the learnable parameters in each component is displayed below the corresponding block.
        }
    \label{fig:nn_structure_iso}
\end{figure}

We let $\spnnSimple{\theta}{1}$ and $\spnnSimple{\theta}{2}$ denote the neural networks $\snnSimple{\theta}{1}$ and $\snnSimple{\theta}{2}$ with the diffusion layer $\mathcal{D}_{iso}$, respectively. Their structures are shown in~\Figure{fig:nn_structure_iso} (compare with \Figure{fig:phase_field_neural_network_structure}).

\begin{figure}[H]
    \centering
    \begin{subfigure}[b]{0.48\textwidth}
        \centering
        \includegraphics[width=\linewidth]{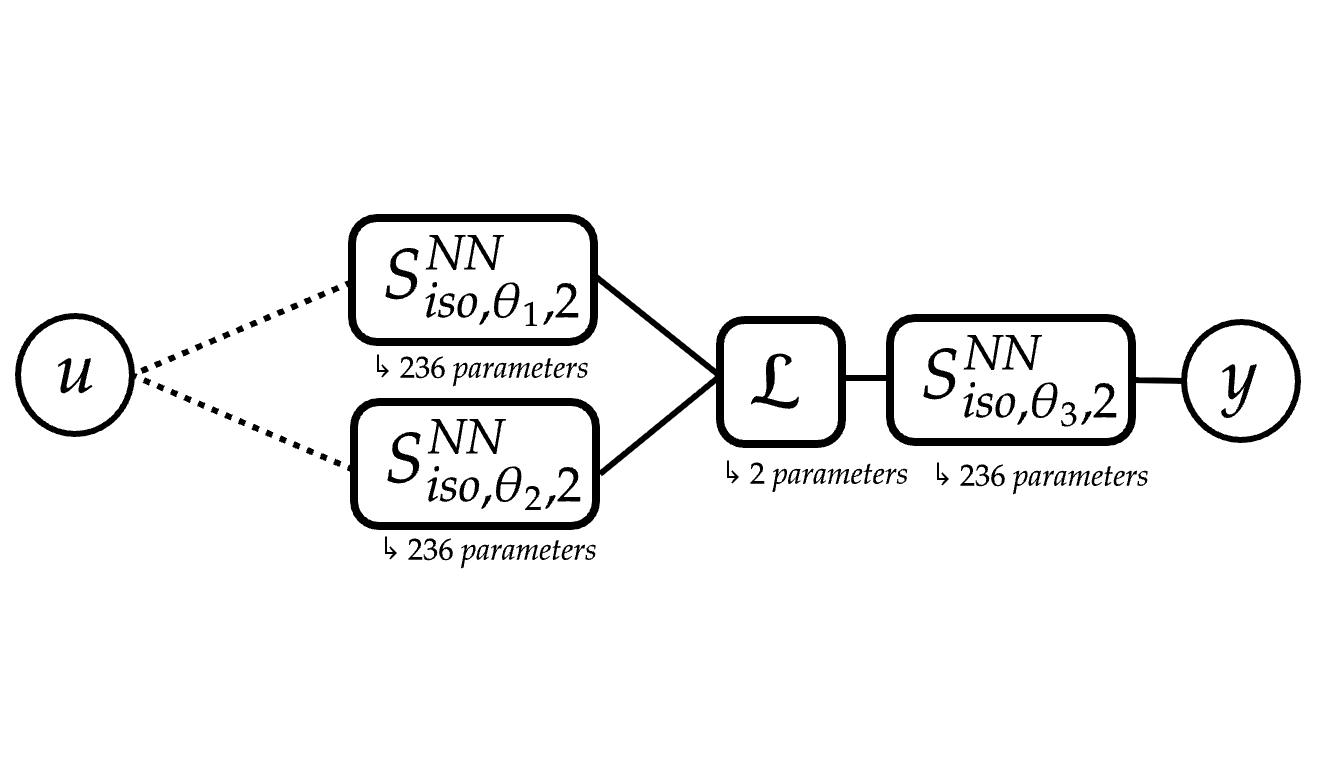}
        \caption{$\wpnnSimple{\theta}{1}$ ($710$ learnable parameters)}
    \end{subfigure}
    \hfill
    \begin{subfigure}[b]{0.48\textwidth}
        \centering
        \includegraphics[width=\linewidth]{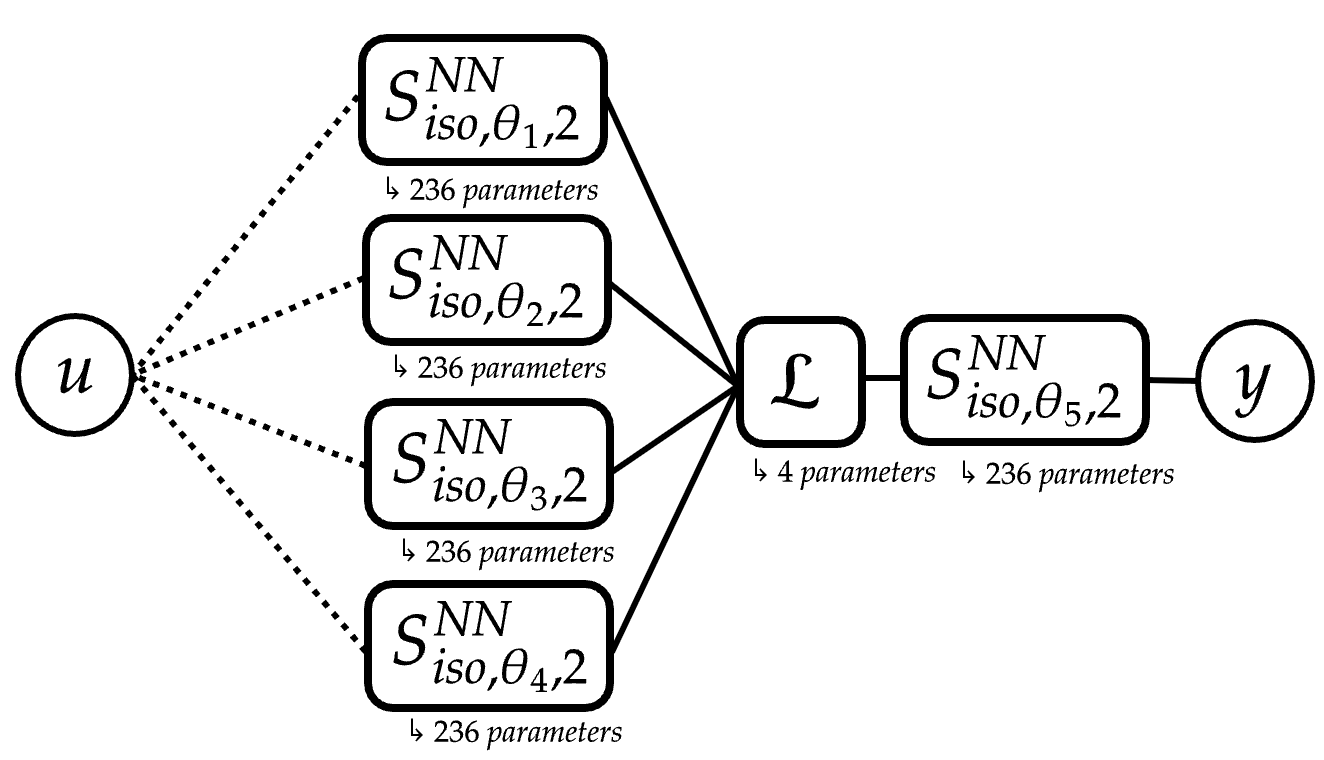}
        \caption{$\wpnnSimple{\theta}{2}$ ($1184$ learnable parameters)}
    \end{subfigure}
    \caption{
        The network structures of the two (left) and four (right) parallel models with rotationally invariant kernels used to approximate the Willmore flow.
        The number of learnable parameters in each component is displayed below the corresponding block.
        }
    \label{fig:nn_structure_iso_willmore}
\end{figure}
We then replace $\snnSimple{\theta}{2}$ with the isotropic-kernel model $\spnnSimple{\theta}{2}$ in the free-kernel models $\wnnSimple{\theta}{1}$ and $\wnnSimple{\theta}{2}$, thereby defining the isotropic-kernel models $\wpnnSimple{\theta}{1}$ and $\wpnnSimple{\theta}{2}$, respectively. The explicit formulas of $\wpnnSimple{\theta}{1}$ and $\wpnnSimple{\theta}{2}$ are given in \eqref{eq:formula_WillmoreNN}, and the resulting structures are shown 
\Figure{fig:nn_structure_iso_willmore} (compare with \Figure{fig:nn_structure}).

We train $\wpnnSimple{\theta}{1}$ and $\wpnnSimple{\theta}{2}$ on the same
dataset and with the same strategy as in \Section{sec:database_construction}.
As in the free kernel case, the $2$-parallel model stalls, and only the
$4$-parallel model is trained; we therefore keep the latter for the generalization
test and call the selected model $\wpnnprim{\theta}{2}$.

\begin{table}[H]
    \centering
    \begin{tabular}{|c|c|c|c|}
        \toprule
        \textbf{Model} &  \textbf{Epochs} & \textbf{Steps} & \textbf{Validation loss} \\
        \midrule
        $\wpnnprim{\theta}{1}$ & $444$ & $59629$  & $1.78\times 10^{-3}$ \\
        $\wpnnprim{\theta}{2}$ & $229$ & $30819$  & $4.26\times 10^{-6}$ \\
        \bottomrule
    \end{tabular}
    \caption{Validation loss of the checkpoints selected for the non-oriented
    curve models with rotationally invariant kernels; the selection rule is the
    one described in \Section{sec:database_construction}.
    }
    \label{tab:training-nn-nonoriented-rot}
\end{table}

\begin{figure}[H]
    \centering
    \begin{subfigure}[b]{0.95\textwidth}
        \centering
        \includegraphics[width=\linewidth]{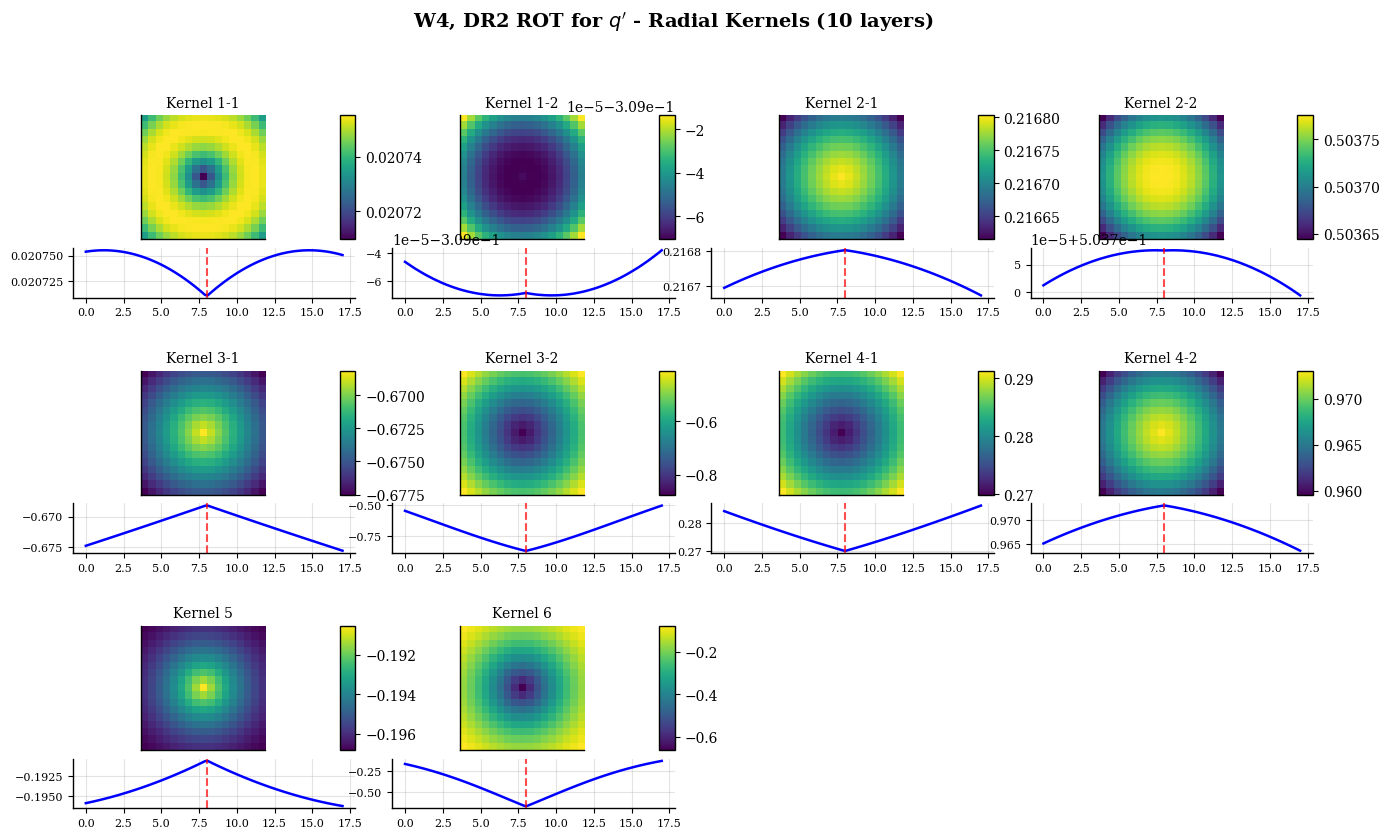}
    \end{subfigure}
    \caption{
        Representation of the $10$ kernels trained and used in our neural network $\wpnnprim{\theta}{2}$ (upper row: heatmaps of the kernels, lower row: radial profiles).
        Each kernel is recovered using the MLP $\radial{\theta}$ in the radial formula \eqref{eq:diffusion_kernel-2D}.
        }
    \label{fig:learned-kernels-w4dr2-qprim-rot}
\end{figure}

In \Figure{fig:learned-kernels-w4dr2-qprim-rot}, we plot the diffusion kernels which are obtained by the above training.
Thanks to the definition \eqref{eq:diffusion_kernel-2D}, each kernel is enforced to be rotationally invariant.

\begin{figure}[H]
  \centering
  \setlength{\tabcolsep}{4pt}   
  \renewcommand{\arraystretch}{1.0} 
  \begin{tabular}{m{0.14\textwidth} m{0.19\textwidth} m{0.19\textwidth} m{0.19\textwidth} m{0.19\textwidth}}
    \raggedleft $\operatorprim\circ RSS$ &
      \includegraphics[width=\linewidth]{figs/general_qprim_smooth_flower6_rssrot1_1.png} &
      \includegraphics[width=\linewidth]{figs/general_qprim_smooth_flower6_rssrot1_2_2.png} &
      \includegraphics[width=\linewidth]{figs/general_qprim_smooth_flower6_rssrot1_3_2.png} &
      \includegraphics[width=\linewidth]{figs/general_qprim_smooth_flower6_rssrot1_4.png} \\
    \raggedleft $\wpnnprim{\theta}{2}$ &
      \includegraphics[width=\linewidth]{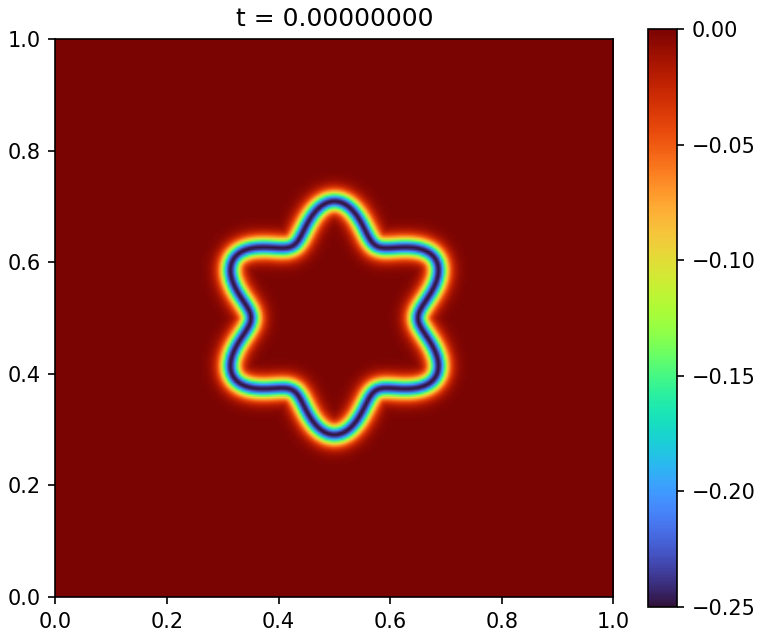} &
      \includegraphics[width=\linewidth]{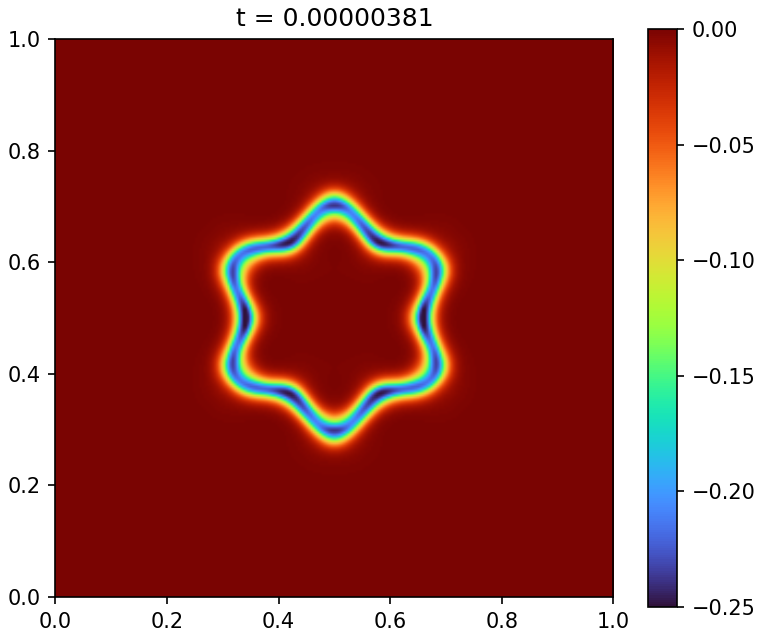} &
      \includegraphics[width=\linewidth]{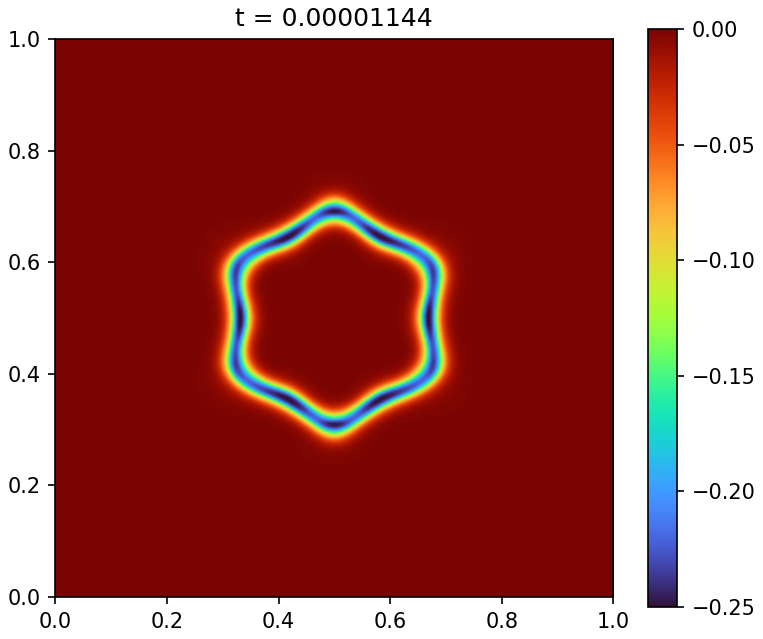} &
      \includegraphics[width=\linewidth]{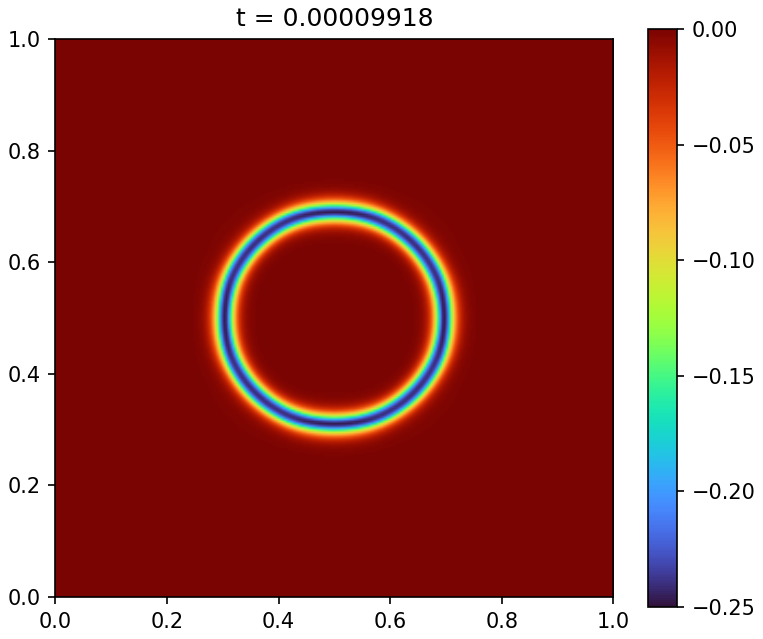} \\
  \end{tabular}
  \caption{
    Approximation of the non-oriented Willmore flow at different time steps for a smooth flower-shaped initial curve.
    First row: $\operatorprim[u^n]$ computed using the RSS scheme. Second row: $u^n$ computed using the isotropic kernel neural network $\wpnnprim{\theta}{2}$.
  }
  \label{fig:smooth_flower6_prim_2d_fixed}
\end{figure}

\Figure{fig:smooth_flower6_prim_2d_fixed} shows that, in contrast with the free
kernel model $\wnnprim{\theta}{2}$, the isotropic model $\wpnnprim{\theta}{2}$ produces a stable
evolution that follows the reference: the curve contracts towards a circle
without developing the spurious structures of \Figure{fig:smooth_flower6_prim_2d}.
A residual grid effect is nevertheless visible at intermediate times, where the computed curve is noticeably
hexagonal, whereas the reference is not.  This is expected: the kernel is
radial as a function of the continuous variable, but it is sampled on a
$(2\Krad+1)\times(2\Krad+1)$ Cartesian stencil, so that the resulting discrete operator is only
approximately equivariant under rotations, and the pointwise nonlinearity
together with the repeated application of the operator amplifies the residual
directional bias.  Enforcing isotropy at the level of the radial profile
therefore removes the instability but does not completely eliminate grid anisotropy.

\begin{figure}[H]
  \centering
  \setlength{\tabcolsep}{4pt}   
  \renewcommand{\arraystretch}{1.0} 
  \begin{tabular}{m{0.14\textwidth} m{0.19\textwidth} m{0.19\textwidth} m{0.19\textwidth} m{0.19\textwidth}}
    \raggedleft $\operatorprim\circ RSS$ &
      \includegraphics[width=\linewidth]{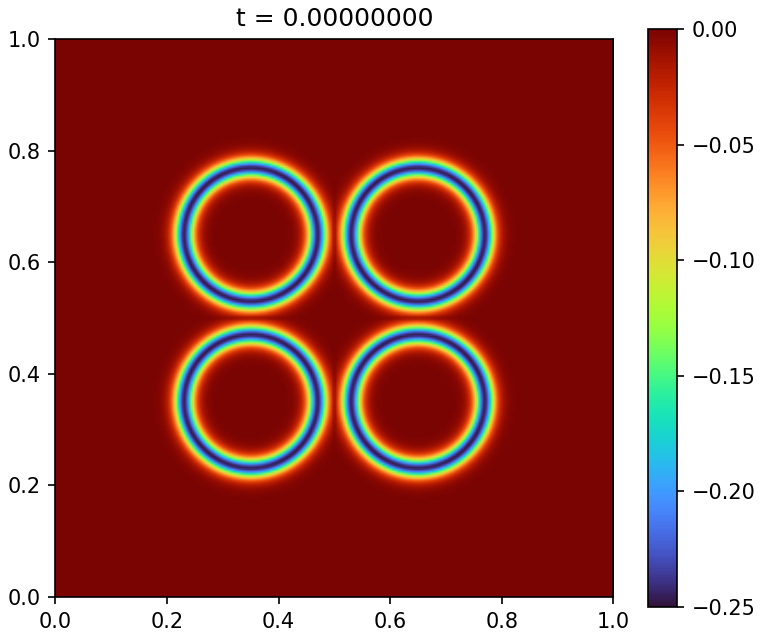} &
      \includegraphics[width=\linewidth]{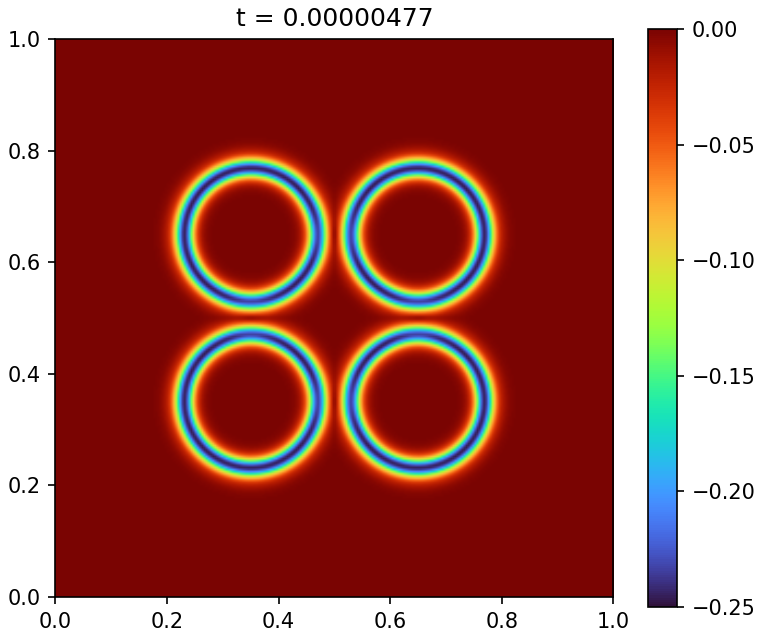} &
      \includegraphics[width=\linewidth]{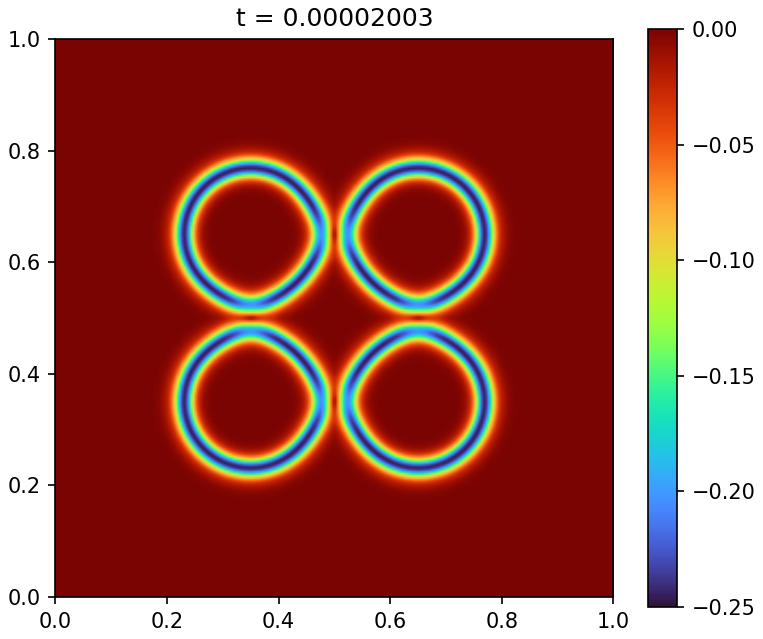} &
      \includegraphics[width=\linewidth]{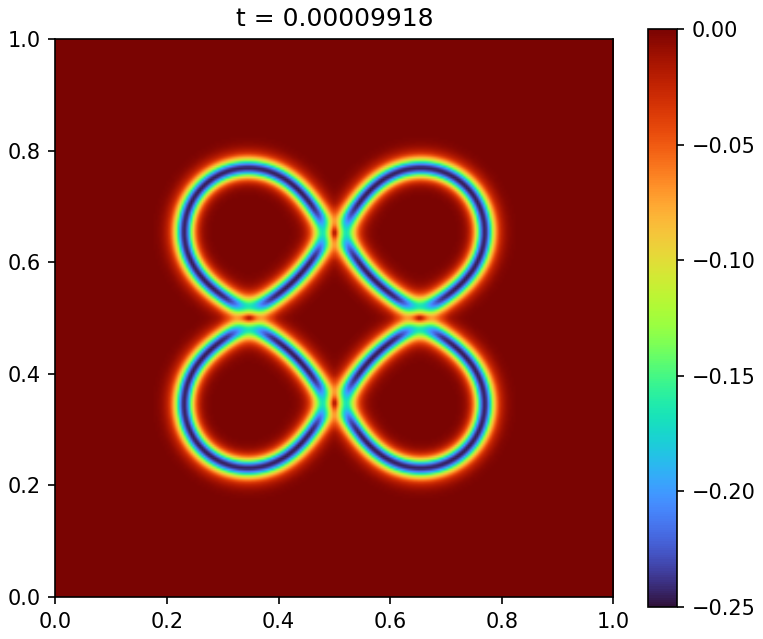} \\
    \raggedleft $\wpnnprim{\theta}{2}$ &
      \includegraphics[width=\linewidth]{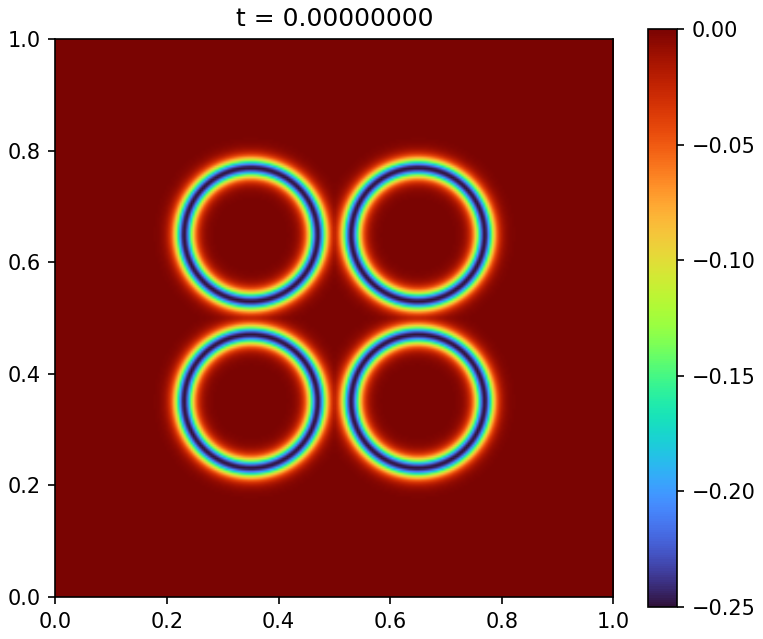} &
      \includegraphics[width=\linewidth]{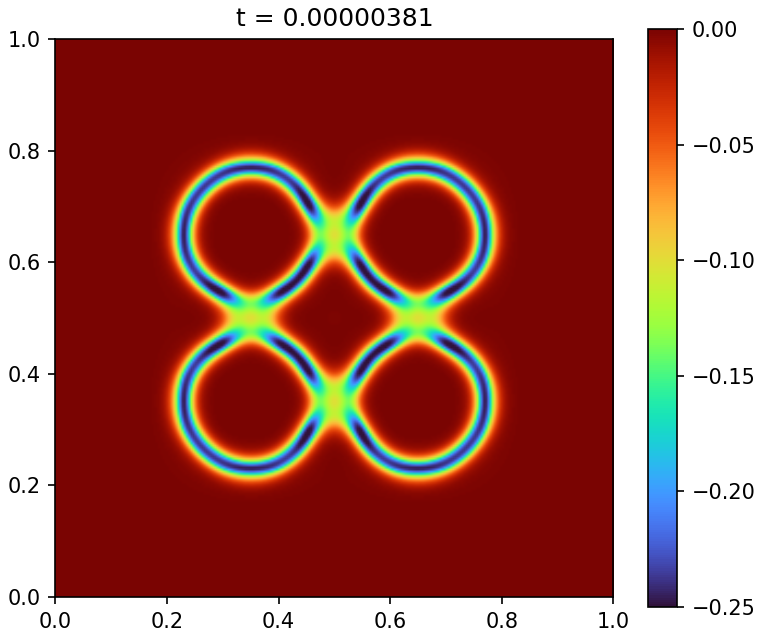} &
      \includegraphics[width=\linewidth]{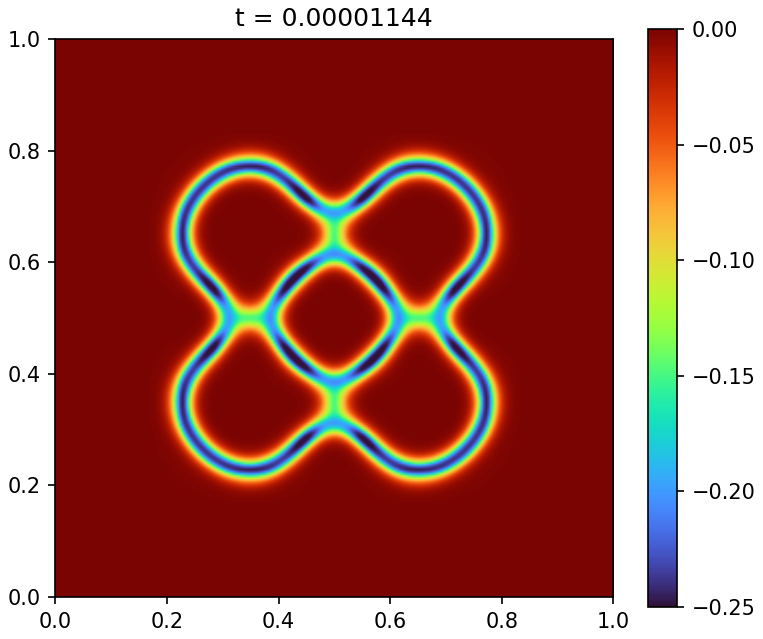} &
      \includegraphics[width=\linewidth]{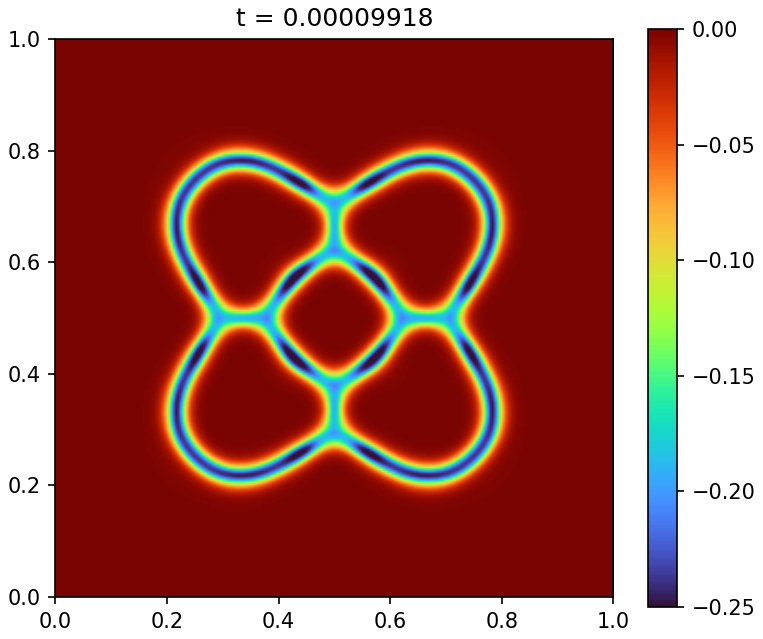} \\
  \end{tabular}
  \caption{
    Approximation of the non-oriented Willmore flow at different time steps for an initial shape consisting of four circles in a square configuration.
    First row: $\operatorprim[u^n]$ computed using the RSS scheme. Second row: $u^n$ computed using the isotropic kernel neural network $\wpnnprim{\theta}{2}$.
  }
  \label{fig:four_square_circles_prim_2d_rot}
\end{figure}
We now carry out one more experiment, starting from four disjoint circles placed at the
vertices of a square (see \Figure{fig:four_square_circles_prim_2d_rot}). This initial interface is the boundary of a disconnected domain. When two neighboring circles touch, the interface develops a singularity at the contact point, where four branches meet. After this topological change, the interface evolves as a curve network that is no longer the boundary of a domain. We note that the training data, namely the evolutions generated by $\operatorprim\circ RSS$, are obtained from the oriented scheme, so that the network $\wpnnprim{\theta}{2}$ is exposed only to domain's boundaries during training. Nevertheless, during the evolution it produces a triple junction, that is, a configuration outside the class represented in the training data. This is possible precisely because the representation based on $q'$ does not distinguish between inside and outside. It also suggests that the learned operator has captured some of the local structure of the flow, rather than merely interpolating between configurations seen during training. We use this property in the next section.

\subsection{Application to reconstruction of curves from point clouds}\label{sec:qprim_application}

Reconstructing a curve or surface from an unorganized set of sample points is a classical inverse problem.
We refer to \cite{berger2016survey} for a review and focus here only on the few works that directly inform our construction.
Implicit variational formulations can be traced back to the work of Zhao, Osher, Merriman and Kang \cite{zhao2000implicit}, who evolve a level set function in order to minimize the mean distance to the data.
However, such data-fidelity models are ill-posed on sparse or noisy samples, and subsequent methods differ mainly in the regularization they add.  Methods such as Poisson reconstruction \cite{kazhdan2006poisson} rely on normals at the data points and are accurate when reliable normals are available. However, the normals must first be estimated and consistently oriented, which requires the target object to be orientable.  The family of methods that regularize by curvature rather than area is closest to the present work. For example, He, Kang and Liu ~\cite{he2020curvature} minimize a combination of the distance to the data and a power of the mean curvature. This prevents the reconstruction from collapsing onto the data and closes the regions not covered by the samples.  However, this approach leads to a fourth-order, non-convex energy whose minimization is computationally challenging. To alleviate this difficulty, threshold dynamics inspired by the Merriman--Bence--Osher scheme \cite{merriman1994motion} have been proposed in \cite{wang2021efficient}, although in the context of the distance-based model.

Our aim is not to compete with these methods on their own terms, but rather to point out that the network introduced in \Section{sec:IsotropicNN_MCF}  provides a regularizer of this type at the cost of a single forward pass. The learned Willmore time-step map provides the curvature regularization, while its training cost is incurred once, prior to reconstruction, rather than at each iteration. Moreover, the non-oriented representation based on the profile $q'$ requires neither normal estimation nor a consistent orientation. In particular, the reconstructed object need not be the boundary of a set. We illustrate below configurations with triple junctions, and present surface reconstructions in \Section{sec:applications_3d}. Handling interfaces that do not bound a set requires an appropriate implicit representation, see also \cite{wan2012reconstructing}. In this respect, our approach also extends the slice-based reconstruction method of Bretin, Dayrens and Masnou~\cite{Bretin_Dayrens_Masnou17} to unoriented data.

Given a point cloud $\widetilde{\Gamma} := \{c_i\}_{i=1}^M$, we set
$U_0 := q'\left(\dist(\cdot,\widetilde{\Gamma})/\varepsilon\right)$ and iterate
\begin{equation}\label{eq:reconstruction_iteration}
    U^{(n+1)} = \min\left\{\frac{U_0}{\alpha}, \wpnnprim{\theta}{2}\left(U^{(n)}\right)\right\},
    \qquad n = 0,1,2,\ldots,
\end{equation}
from $U^{(0)} = U_0$, with $\alpha \geq 1$.

Each iteration combines two effects. The network provides a learned curvature regularization of the current curve, while the pointwise minimum enforces attachment to the data. Since $U_0$ is close to its minimum value $-1/4$ near the prescribed points, taking the minimum forces $U^{(n+1)}$ to remain low in their vicinity. The parameter $\alpha$ controls how strictly this obstacle constraint is enforced. Taking $\alpha=1$ imposes the full pointwise constraint, whereas $\alpha>1$ shifts $U_0/\alpha$ towards $0$, allowing the curve to pass at some distance from the data points, as is desirable for noisy point clouds.

We now present several numerical reconstruction examples. In the following figures, the white dots indicate the prescribed points on $\widetilde{\Gamma}$.

\begin{figure}[H]
  \centering
  \setlength{\tabcolsep}{4pt}   
  \renewcommand{\arraystretch}{1.0} 
  \begin{tabular}{m{0.12\textwidth} m{0.26\textwidth} m{0.26\textwidth} m{0.26\textwidth}}
    \raggedleft $\alpha = 1.0$ &
      \includegraphics[width=\linewidth]{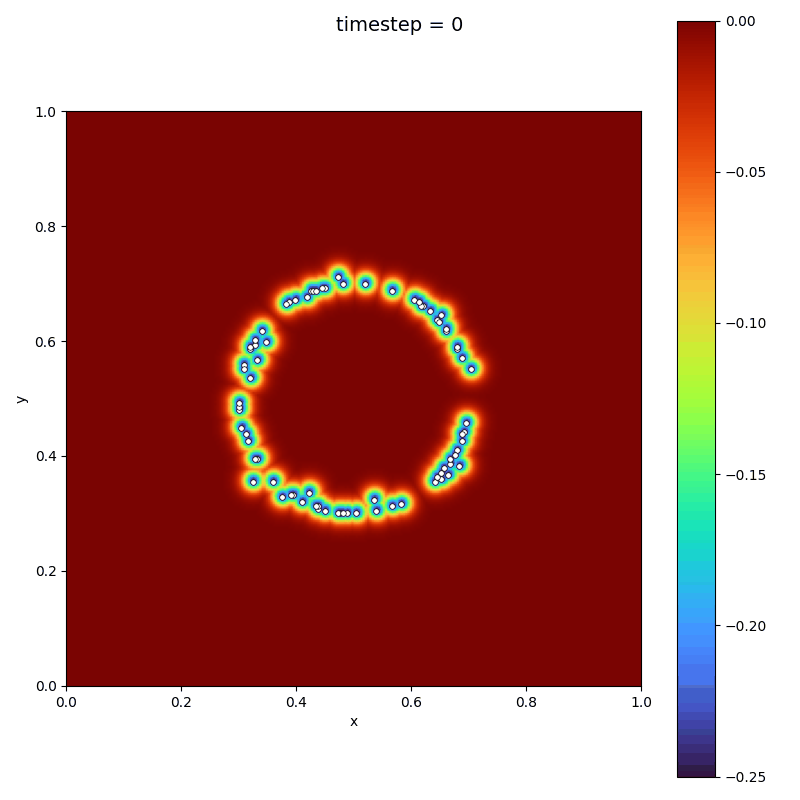} &
      \includegraphics[width=\linewidth]{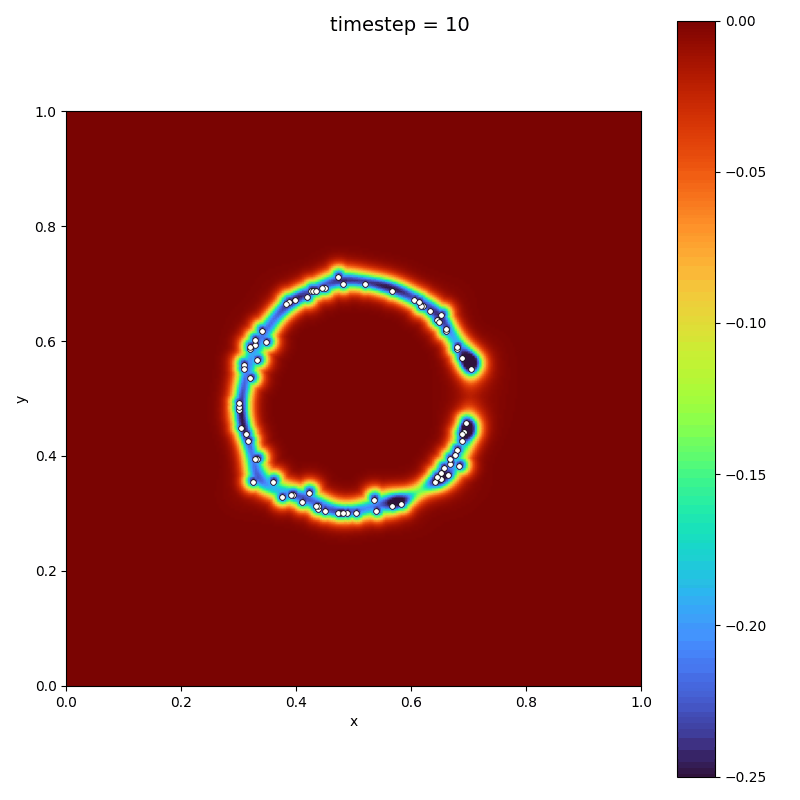} &
      \includegraphics[width=\linewidth]{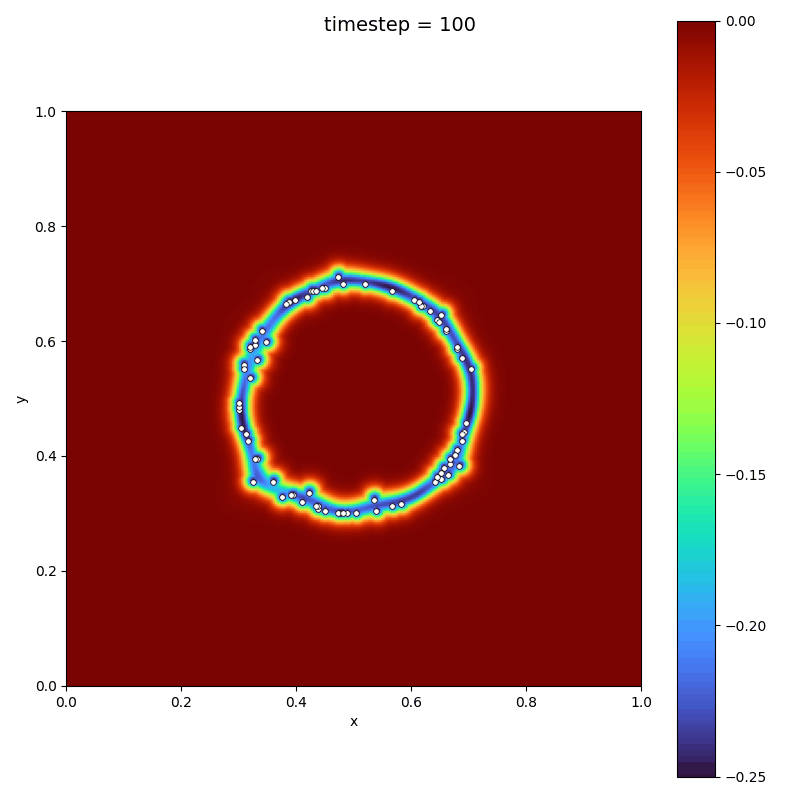} \\
    \raggedleft $\alpha = 1.5$ &
      \includegraphics[width=\linewidth]{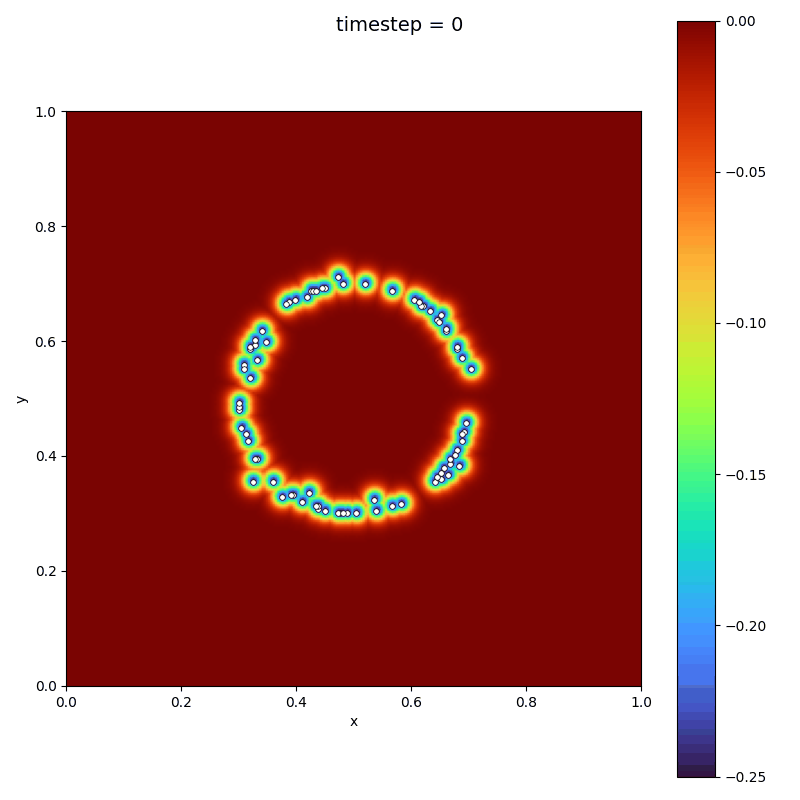} &
      \includegraphics[width=\linewidth]{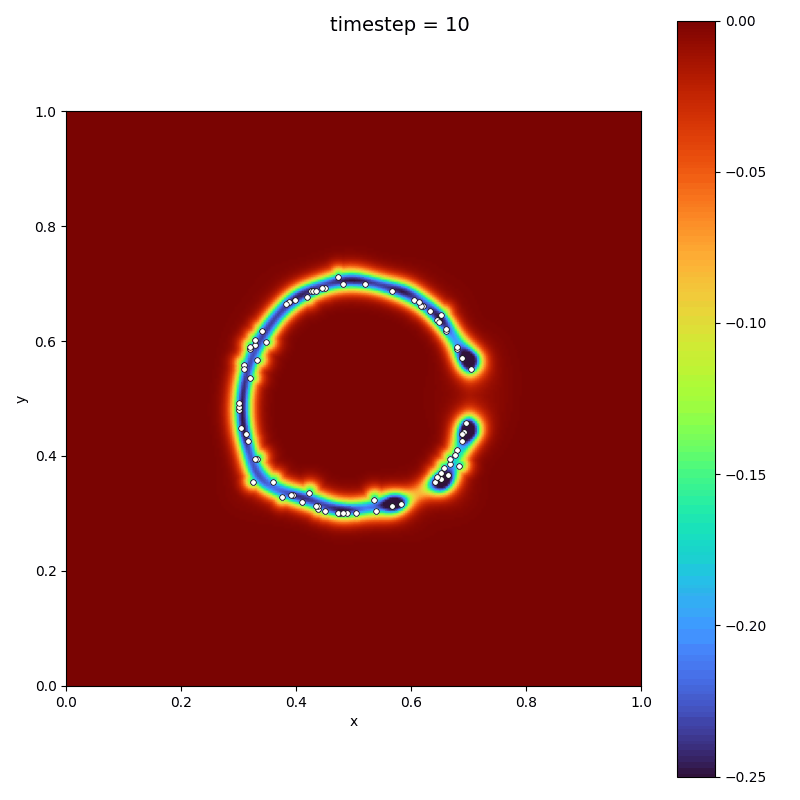} &
      \includegraphics[width=\linewidth]{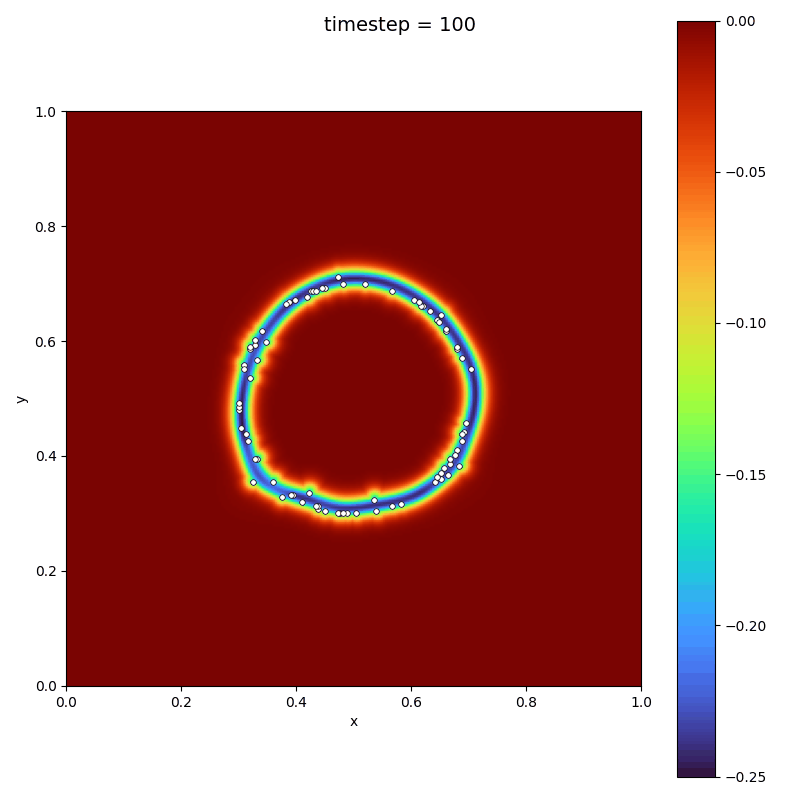} \\
    \raggedleft $\alpha = 2.0$ &
      \includegraphics[width=\linewidth]{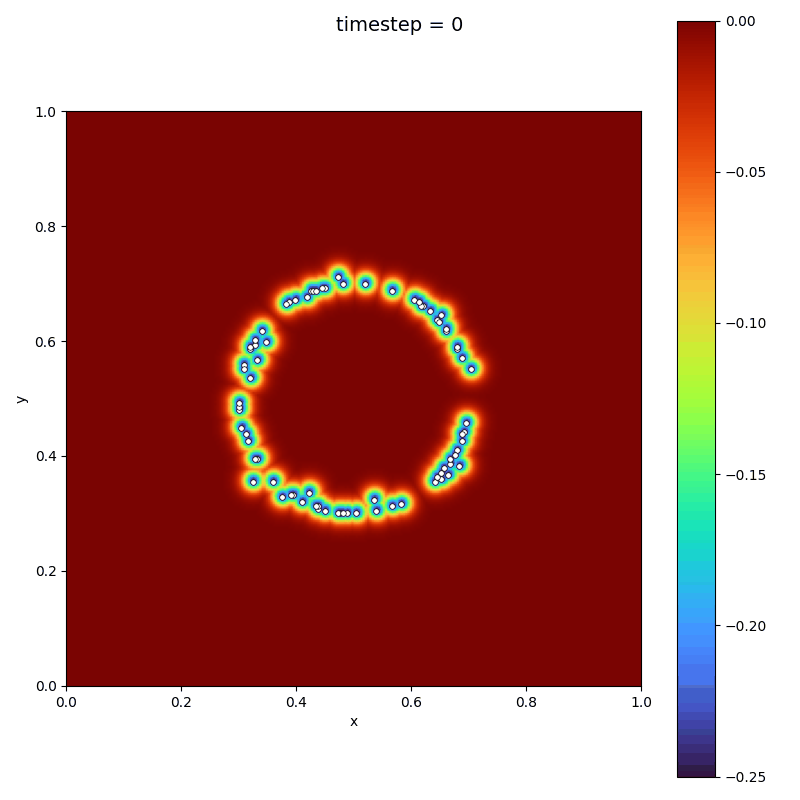} &
      \includegraphics[width=\linewidth]{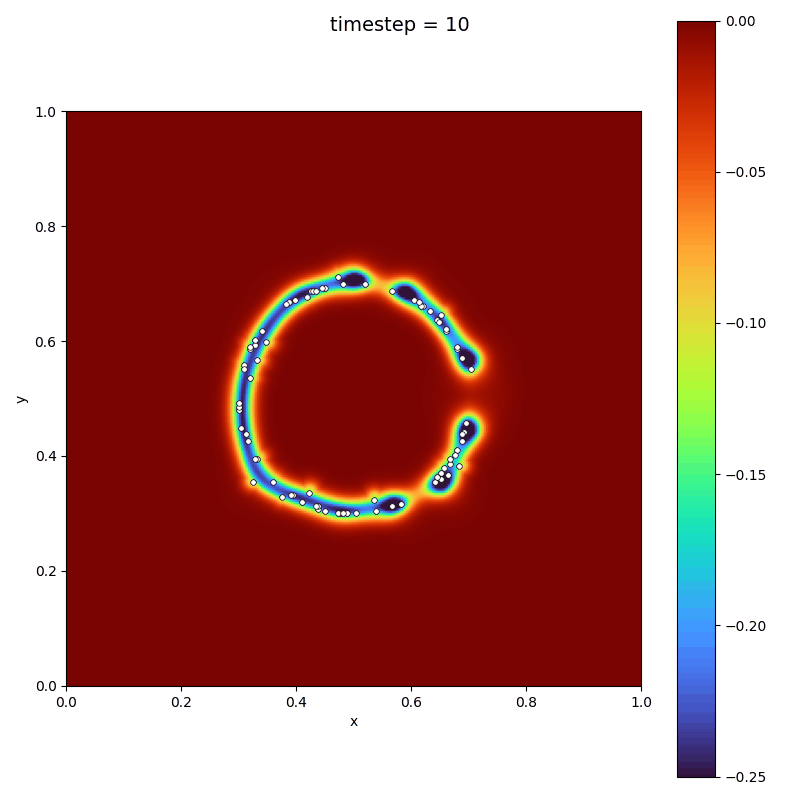} &
      \includegraphics[width=\linewidth]{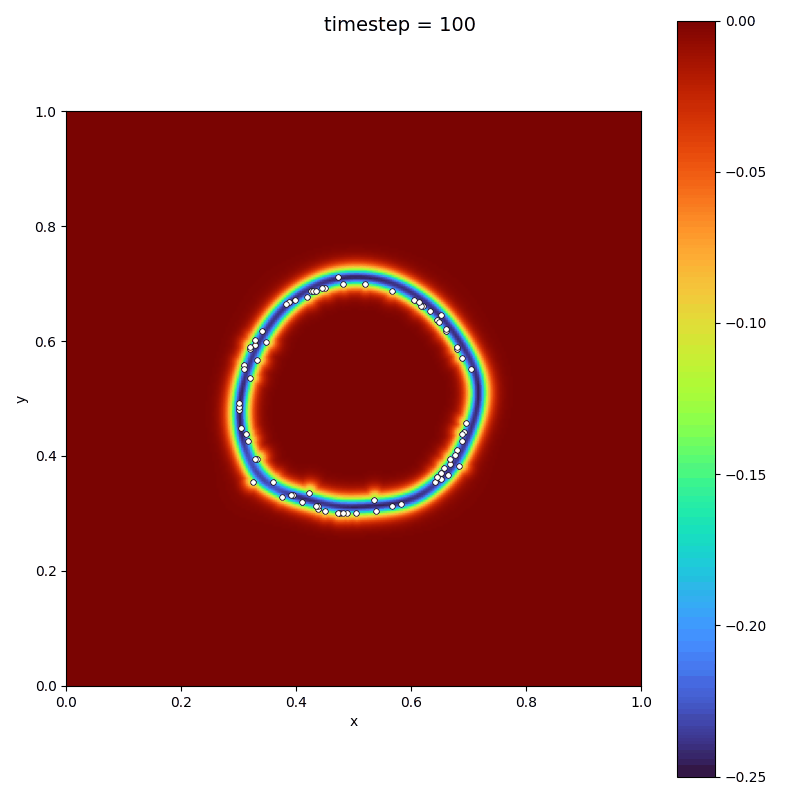} \\
  \end{tabular}
  \caption{Reconstruction of a curve from points sampled on a circle and perturbed by noise; Each row shows the phase-field produced by the scheme for a different value of $\alpha$, namely
 $\alpha =1.0,\,1.5,\,2.0$ at iterations $0$ (initial), $10$, and $100$ (final).}
  \label{fig:reconstruction_circle_noise_2d}
\end{figure}
\paragraph{Circle with noise.}
We randomly sample $80$ points on a circle and perturb every tenth point with noise. We then apply iteration  \eqref{eq:reconstruction_iteration} for $100$ steps using three different values of $\alpha$. For smaller values of $\alpha$, the reconstructed curves tend to pass closer to the prescribed points, whereas for the largest value of $\alpha$, the reconstruction is more flexible and can deviate further from them.

\paragraph{Double bubble.} We consider the reconstruction of a double-bubble configuration, which is a non-oriented curve network with triple junctions. We observe that the endpoint of the open curve separating the two bubbles moves toward the intersection of the two outer curves, where a triple junction eventually forms. This behavior is consistent with the decrease of the Willmore energy along the approximate flow. Remarkably, the steady state of the evolution is itself a non-oriented curve network with triple junctions, even though the model $\wpnnprim{\theta}{2}$ was trained only on oriented configurations (see again \Figure{fig:four_square_circles_prim_2d_rot}).

\begin{figure}[H]
  \centering
  \setlength{\tabcolsep}{4pt}   
  \renewcommand{\arraystretch}{1.0} 
  \begin{tabular}{m{0.25\textwidth} m{0.25\textwidth} m{0.25\textwidth}}
      \includegraphics[width=\linewidth]{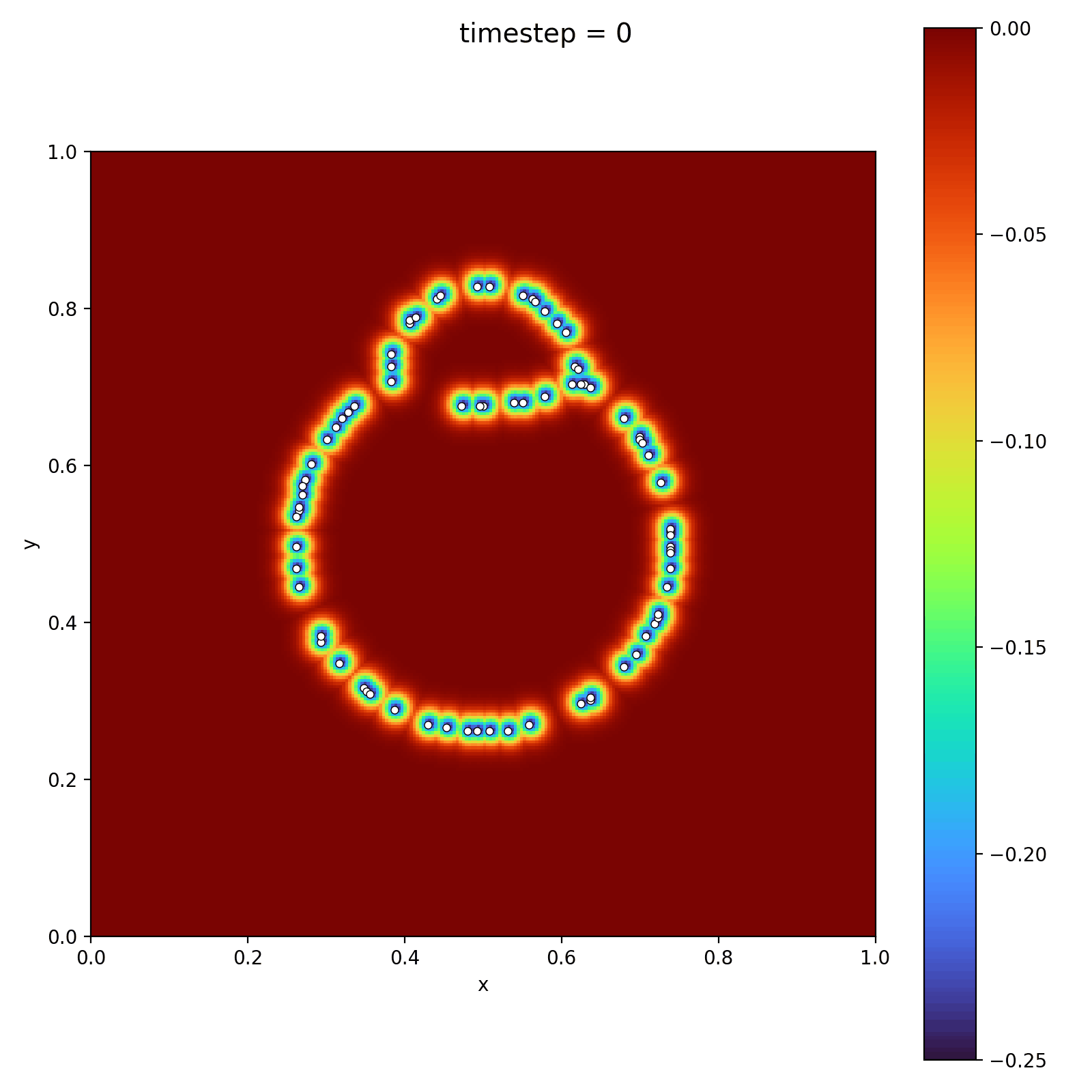} &
      \includegraphics[width=\linewidth]{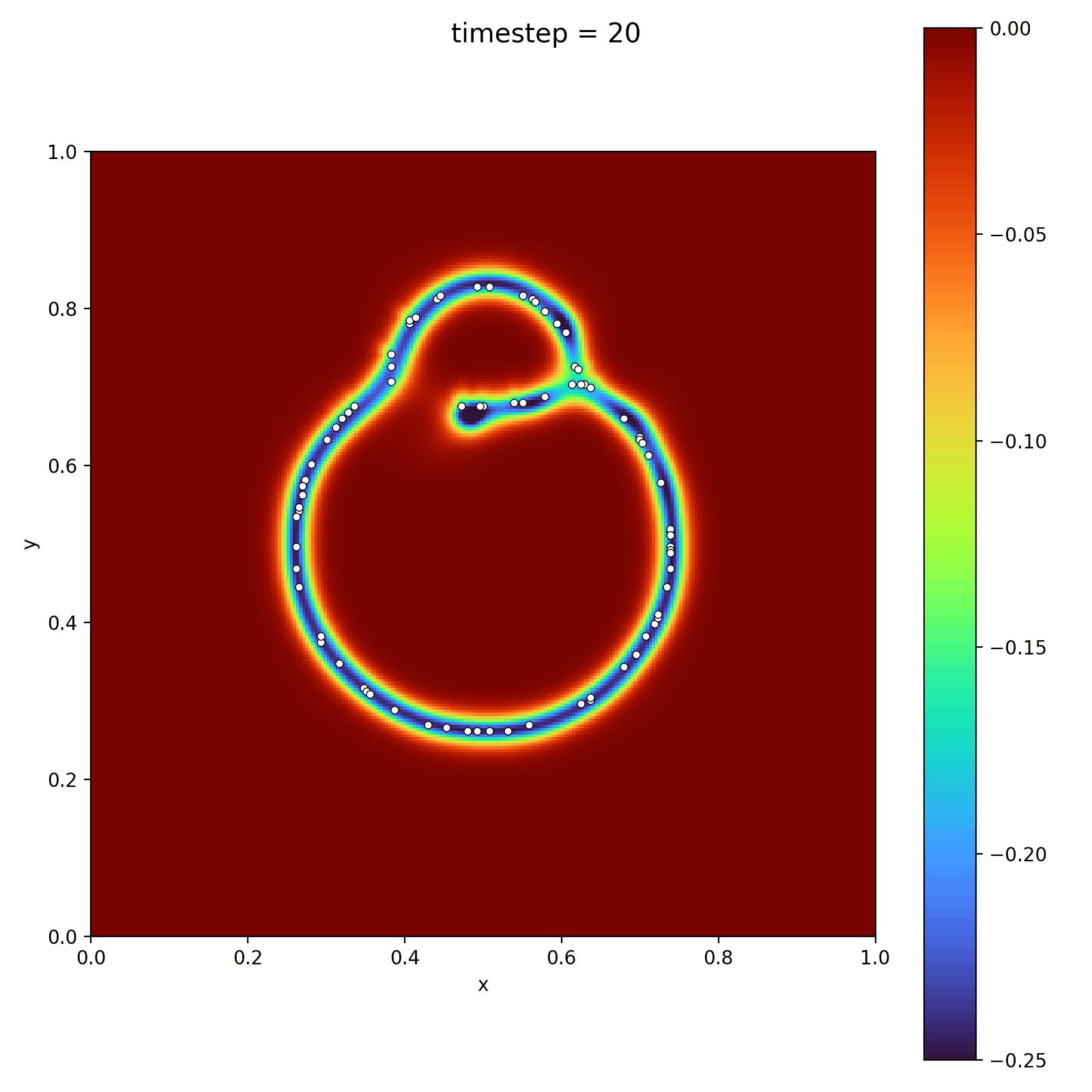} &
      \includegraphics[width=\linewidth]{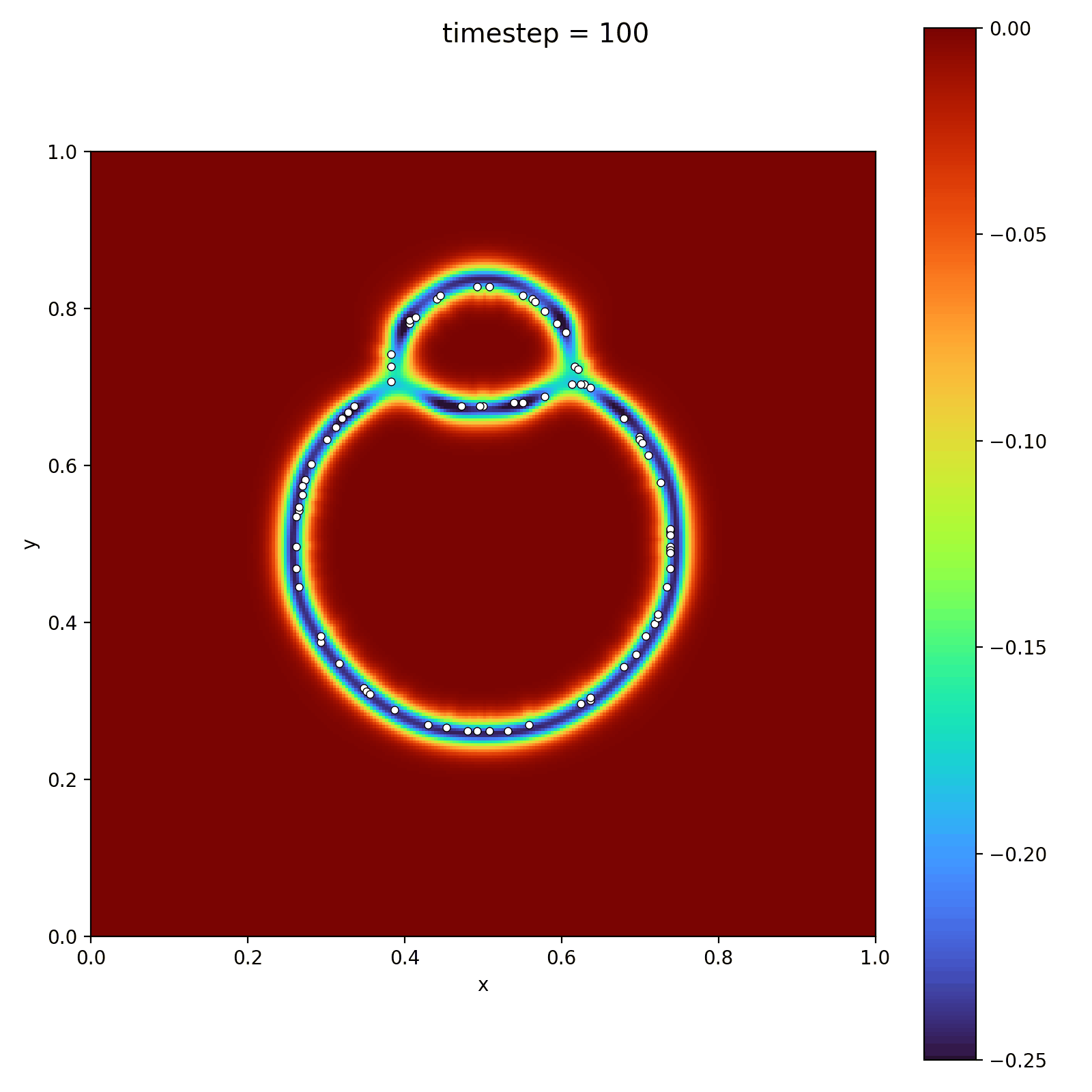} \\
  \end{tabular}
  \caption{Reconstruction of a double bubble curve from $80$ points after $100$
  iterations of \eqref{eq:reconstruction_iteration}, using the rotationally
  invariant neural network $\wpnnprim{\theta}{2}$.}
  \label{fig:reconstruction_double_bubble_2d}
\end{figure}

\section{Approximation of the Willmore flow in 3D}\label{sec:o_willmore_nn_3d}

This section investigates the extension of the two-dimensional results to three dimensions. The first challenge is to construct a sufficiently rich dataset to accurately represent three-dimensional Willmore flows. We then assess how well the networks generalize to smooth, oriented interfaces.

As in the two-dimensional setting, non-oriented interfaces are substantially more challenging. Training the networks in this case requires an initialization based on previously trained two-dimensional models. Assuming isotropic diffusion kernels, the rotational symmetry of the flow allows the three-dimensional kernels to be constructed from their two-dimensional counterparts. With this initialization, the networks produce reasonable approximations of the three-dimensional Willmore flow, which can subsequently be used for interface reconstruction from point clouds.

\subsection{Construction of a dataset in 3D}\label{sec:rss_3d}

\begin{figure}[H]
    \centering
    \begin{subfigure}[b]{0.24\textwidth}
        \centering
        \includegraphics[width=\linewidth]{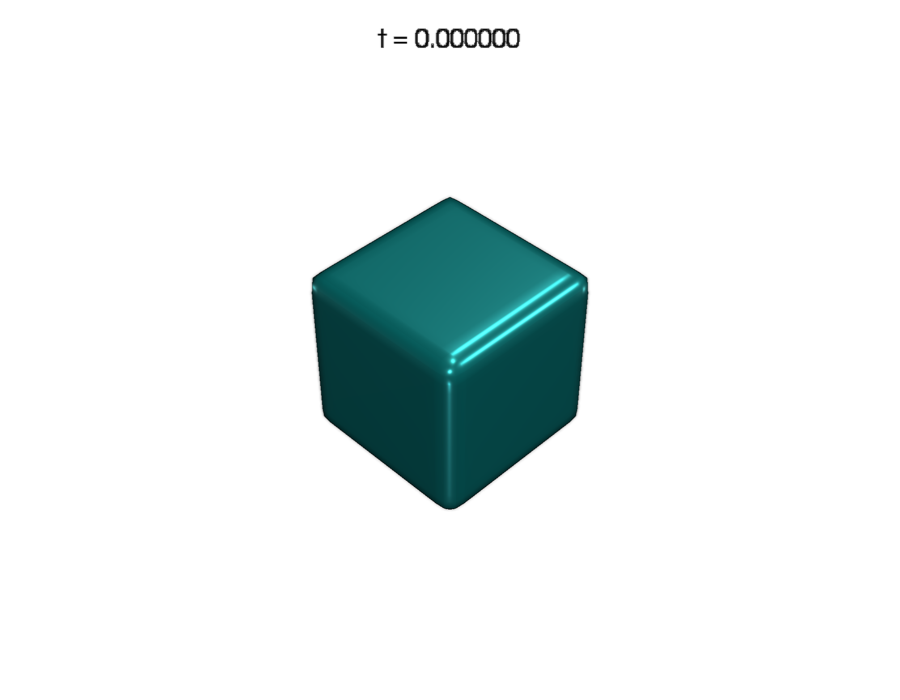}
    \end{subfigure}
    \hfill
    \begin{subfigure}[b]{0.24\textwidth}
        \centering
        \includegraphics[width=\linewidth]{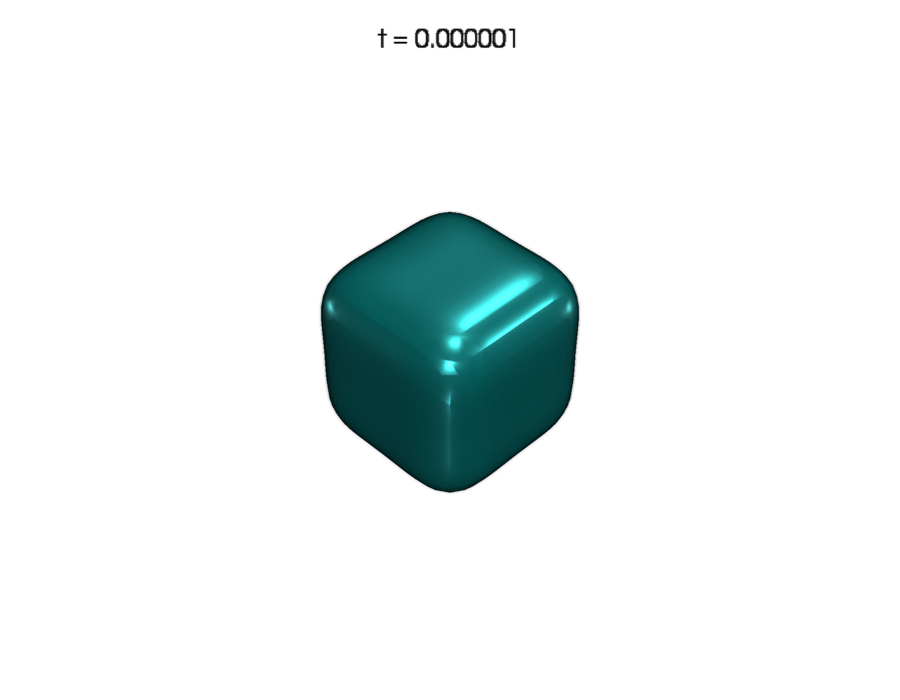}
    \end{subfigure}
    \hfill
    \begin{subfigure}[b]{0.24\textwidth}
        \centering
        \includegraphics[width=\linewidth]{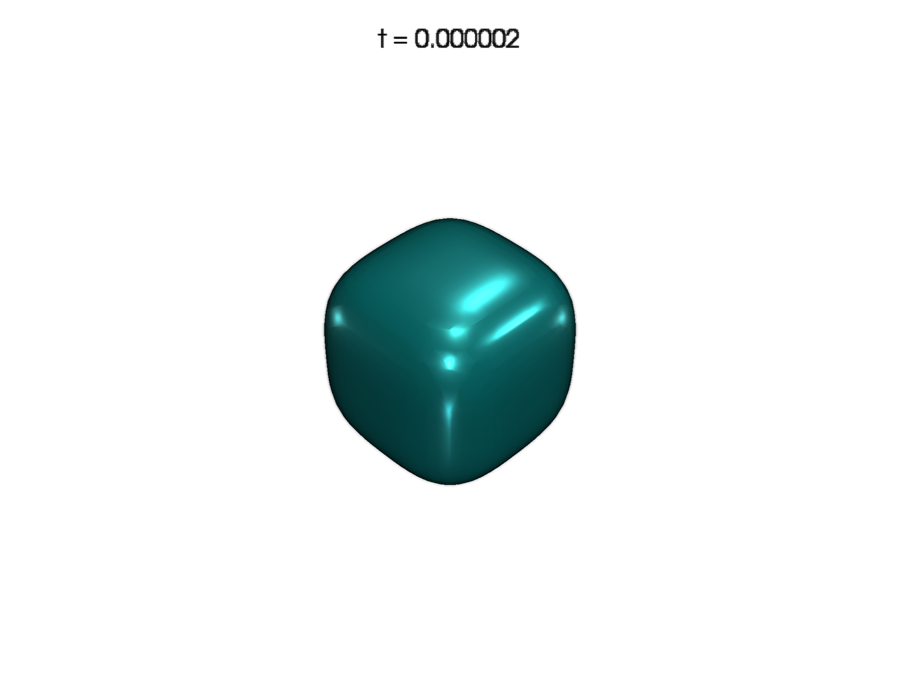}
    \end{subfigure}
    \hfill
    \begin{subfigure}[b]{0.24\textwidth}
        \centering
        \includegraphics[width=\linewidth]{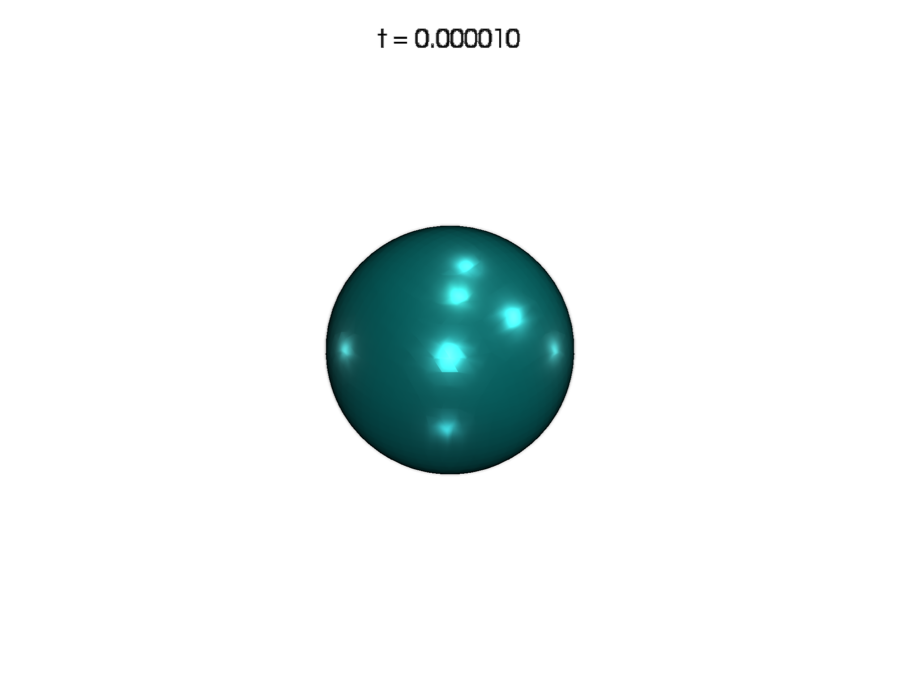}
    \end{subfigure}
    \begin{subfigure}[b]{0.24\textwidth}
        \centering
        \includegraphics[width=\linewidth]{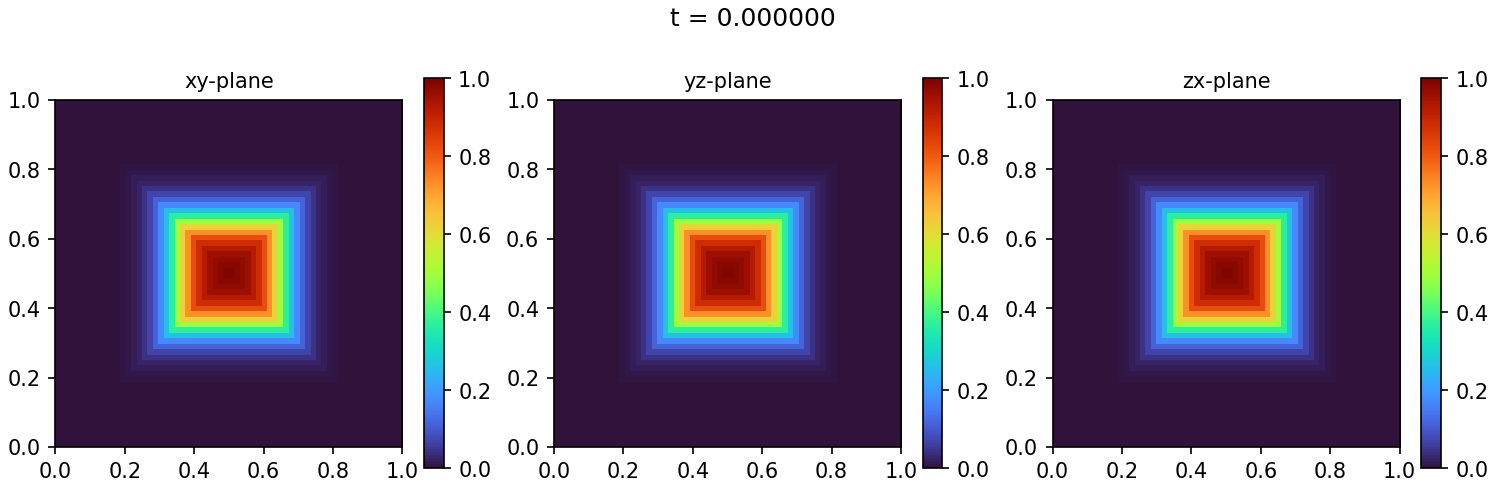}
    \end{subfigure}
    \hfill
    \begin{subfigure}[b]{0.24\textwidth}
        \centering
        \includegraphics[width=\linewidth]{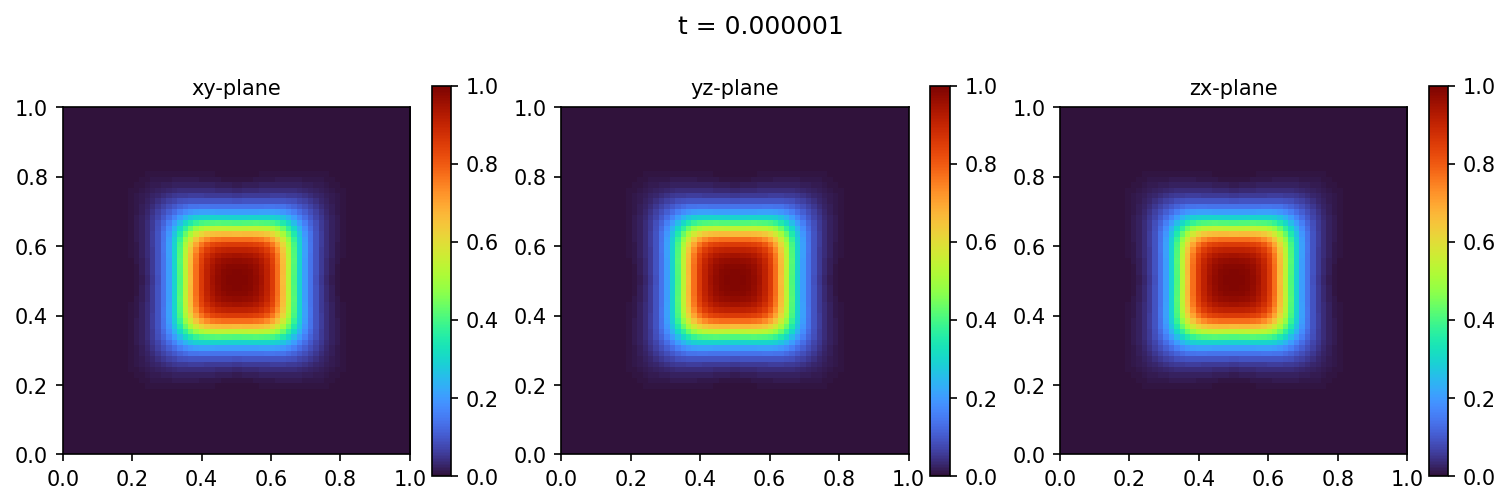}
    \end{subfigure}
    \hfill
    \begin{subfigure}[b]{0.24\textwidth}
        \centering
        \includegraphics[width=\linewidth]{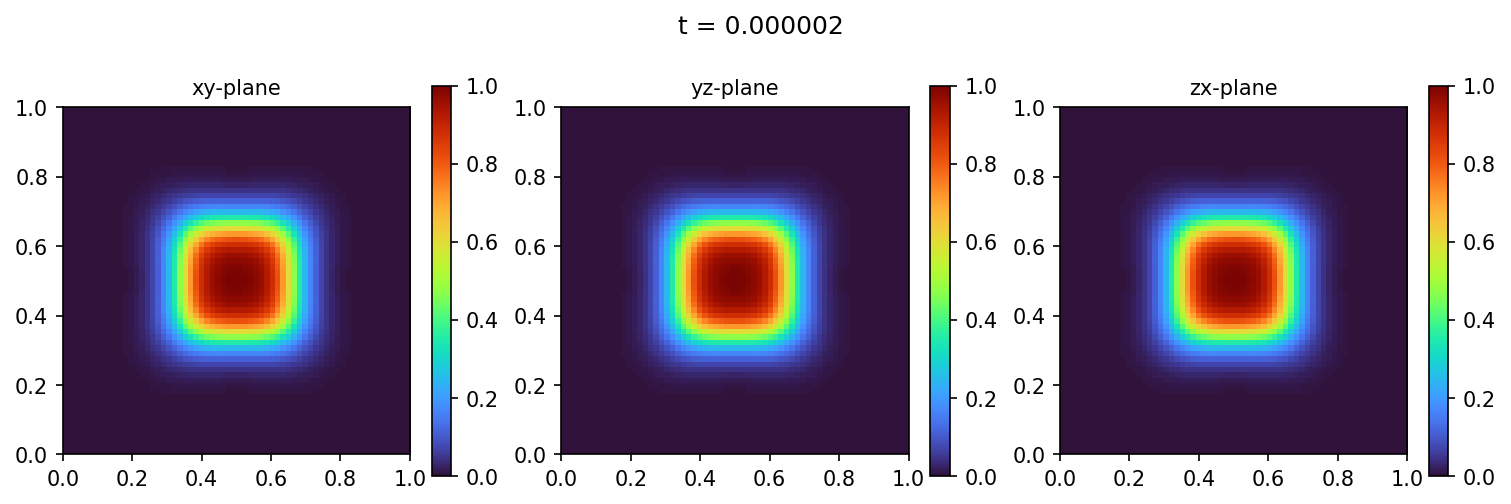}
    \end{subfigure}
    \hfill
    \begin{subfigure}[b]{0.24\textwidth}
        \centering
        \includegraphics[width=\linewidth]{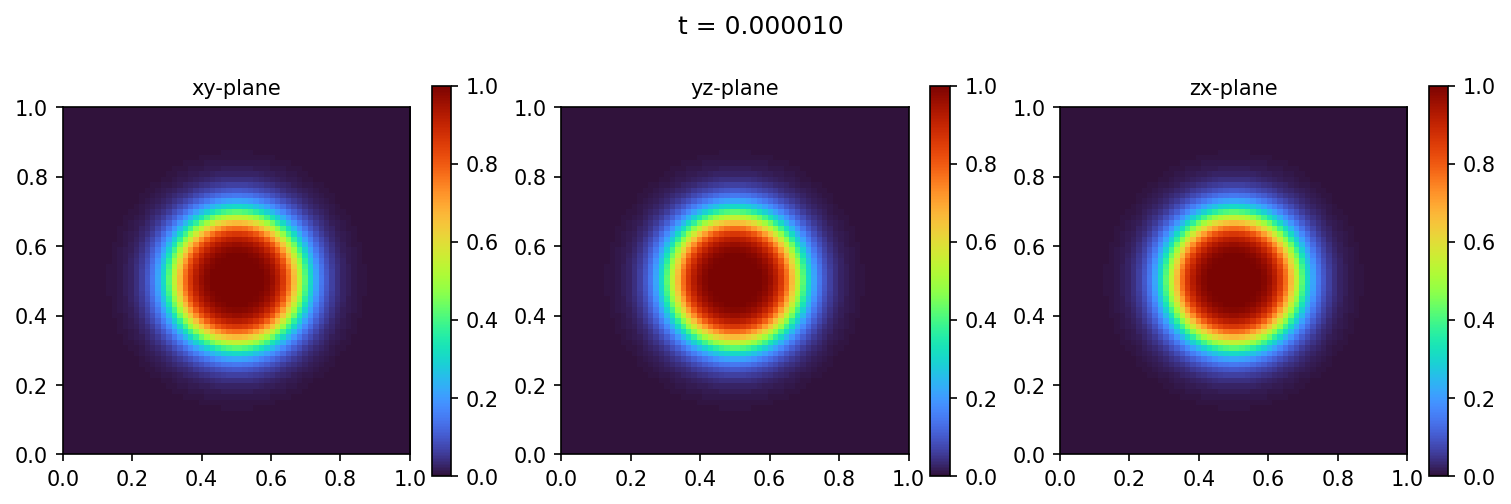}
    \end{subfigure}
    \caption{Evolution of a cube (a genus $0$ shape) under the Willmore flow
    computed with the RSS scheme in 3D, with $\tts = \ts / 64$ and
    $T = 10^{-5}$.  First row: the $0.5$-level set of the phase field function;
    second row: its projections onto the $XY$-, $YZ$- and $ZX$-planes.}
    \label{fig:rss_3d_cube}
\end{figure}

As in the 2D case, we adopt the RSS scheme in 3D to simulate the Willmore flow of oriented surfaces and regard its outcome as the ground truth to train our neural network.
We perform the RSS scheme in 3D starting from several shapes. We set the resolution to $N = 64$ and implement the RSS scheme in 3D.
For each time step, the learning rate is set to $\eta = 1.0$ and the number of
inner ADAM iterations to $J = 50$ for the illustrative evolution of
\Figure{fig:rss_3d_cube}, and to $J = 300$ for the database itself.
\Figure{fig:rss_3d_cube} shows one such evolution, starting from a cube.
Each surface represents the $0.5$-level set of the phase field function $u$ at
the corresponding time.
The time step $\tts$ and the time horizon $T$ are different for each shape and are specified in the following sections.
The first and second rows of the following figures show the evolution of the contour sets of the induced phase field functions
and the projections of the phase field functions onto the $XY$-, $YZ$-, and $ZX$-planes, respectively.

\begin{figure}[H]
    \centering
    \begin{subfigure}[b]{0.9\textwidth}
        \centering
        \includegraphics[width=0.45\linewidth]{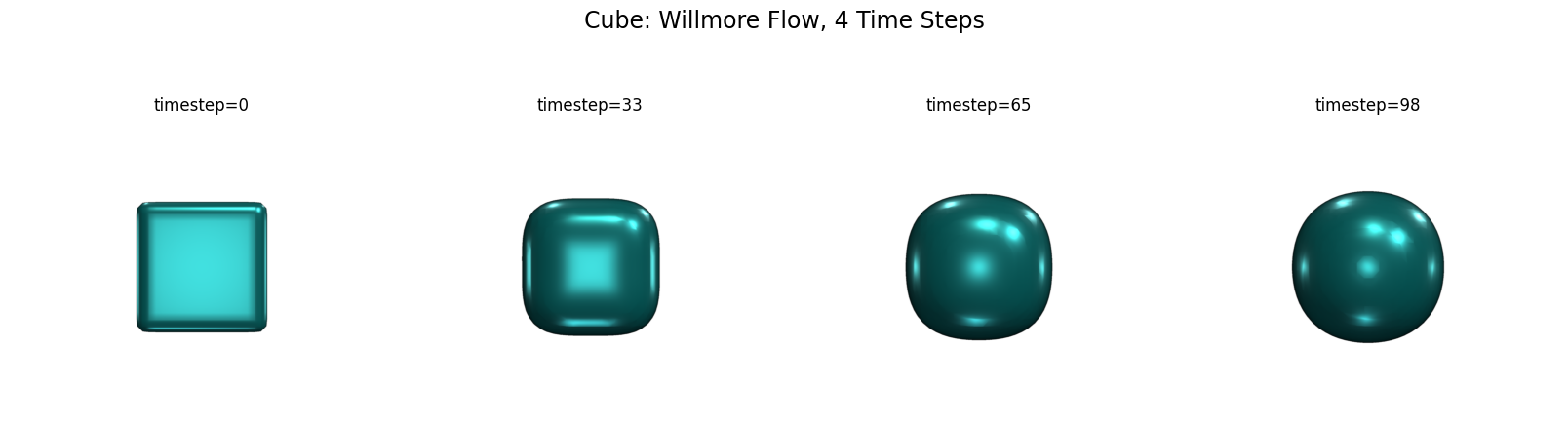} $\quad \quad $
           \includegraphics[width=0.45\linewidth]{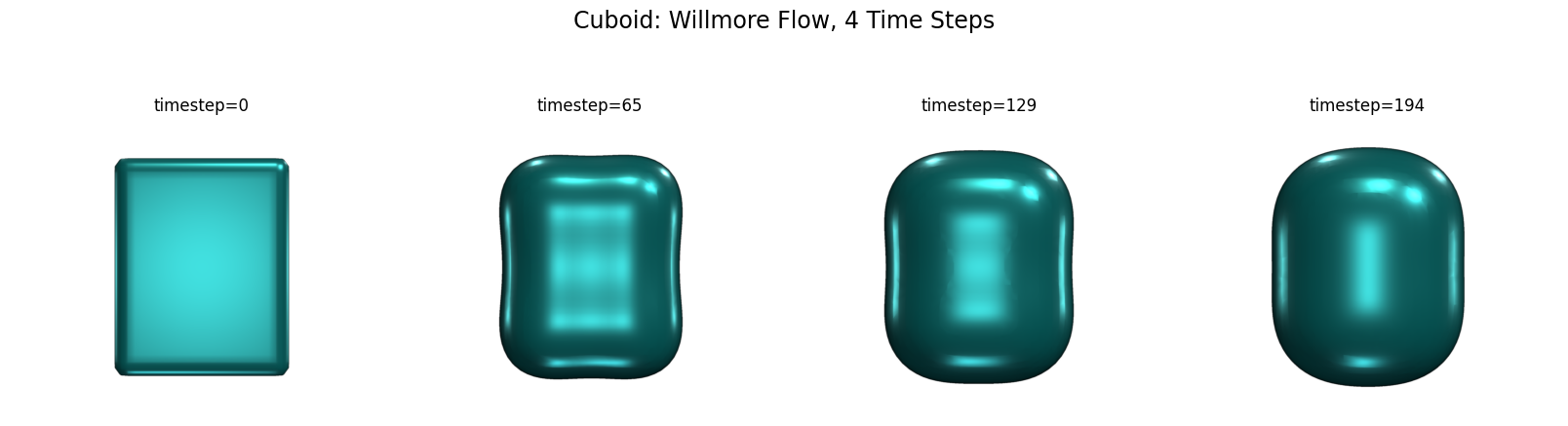} \\
              \includegraphics[width=0.45\linewidth]{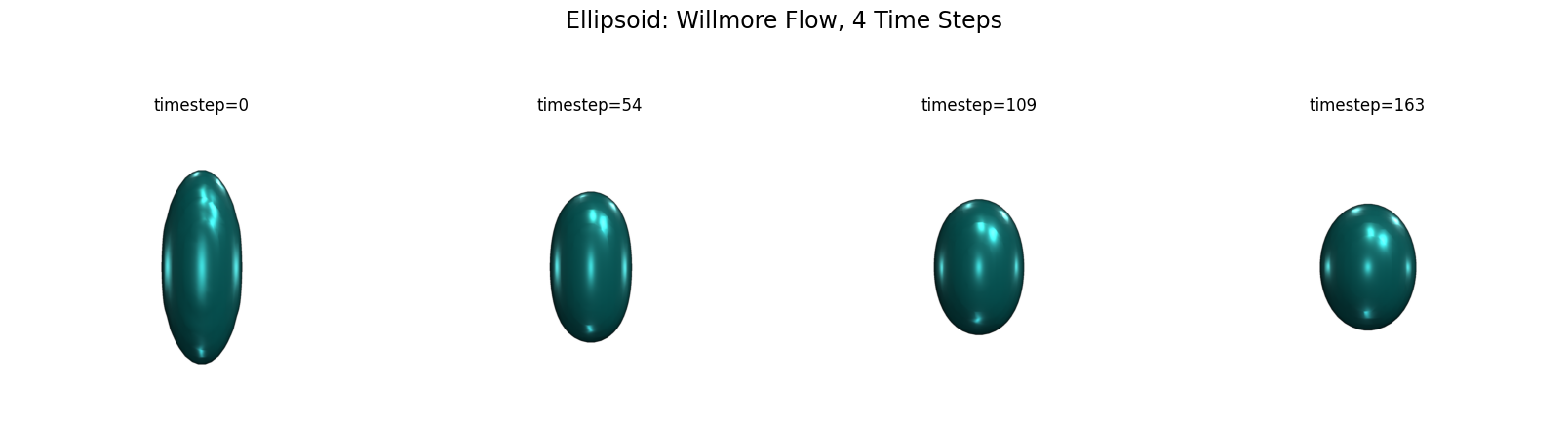} $\quad \quad $
                \includegraphics[width=0.45\linewidth]{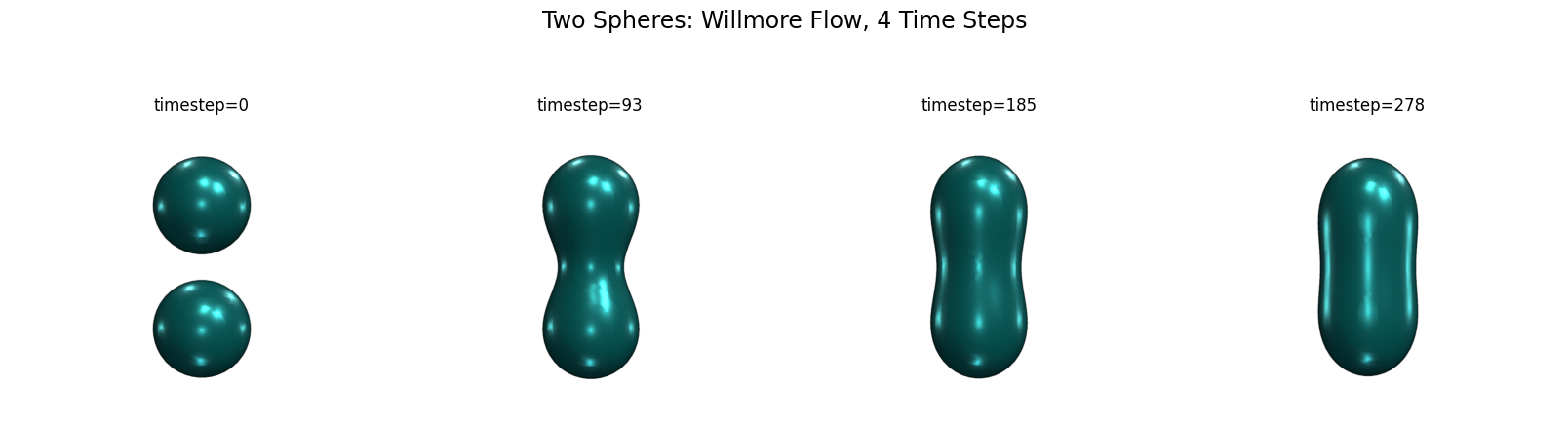}\\
                   \includegraphics[width=0.45\linewidth]{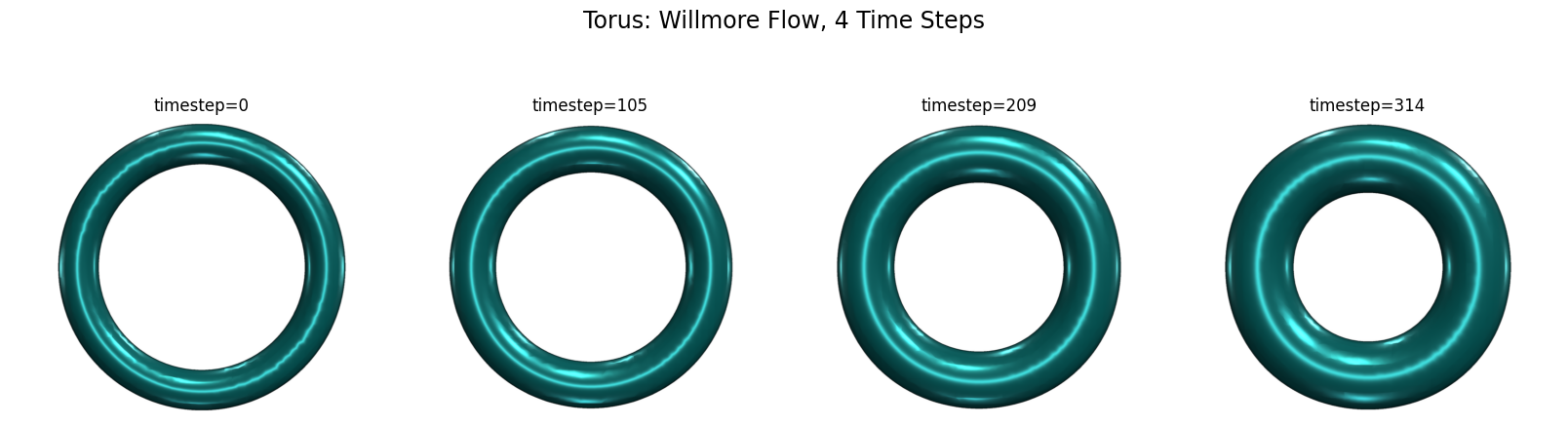}
    \end{subfigure}
    \caption{Dataset of the Willmore flow generated by the RSS scheme in 3D. Each surface represents the $0.5$ level set of the phase-field function at each time step. Each sequence of four pictures corresponds to a different initial configuration.}
    \label{fig:dataset_oriented_3d_set1}
\end{figure}

More specifically, the database is built from the Willmore flow of initial sets
consisting of cubes, cuboids, ellipsoids, unions of two spheres and tori. \Figure{fig:dataset_oriented_3d_set1} shows examples of the corresponding
evolutions computed with the RSS scheme.  These families are the
three-dimensional counterparts of the ones used in 2D, and were chosen for the
same reason: they provide surfaces with non-constant and strongly varying mean
curvature, non-trivial topology and, for the union of two spheres, a
singularity.  Other families -- cylinders, star-shaped surfaces, or the
stationary sphere itself -- could be added at no conceptual cost.

The database to train $\wnn{\theta}{2}$ consists of time series of the phase field functions of the Willmore flow generated by the RSS scheme in 3D.
Each phase field function is defined on the cubic domain $\Omega = [0,1]^3$ with the resolution $N = 64$.
The hyperparameters are set to $\varepsilon := 2 / N = 2^{-5}$, $\tts := \ts / 64$ and $J = 300$.
We have selected several representative shapes to enrich the training dataset.

\begin{rem}
We do not include the sphere in the training dataset. It is a stationary solution of the three-dimensional Willmore flow, and its stability can actually be learned from the datasets described above, which contain shapes that are close to spherical.
\end{rem}

\subsection{Training and application to the approximation of oriented Willmore flow}

As in the non-oriented 2D case, we first propose to force the radial invariance of the diffusion kernels $\diffusion$.
In the 3D case, each $\mathcal{D}_{iso}$ is a discrete three-dimensional function determined by a one-dimensional radial profile $\radial{}:\R{\geq 0}\to\mathbb{R}$. Rather than parameterizing the convolution kernel directly as a function on $\mathbb{R}^3$, we define it through $\radial{}$, evaluated at the normalized distance from the kernel center. The resulting three-dimensional kernel is thus obtained by radial extension of $\radial{}$.

Precisely speaking, for each index $(i,j,k)\in\{-\Krad,\ldots,\Krad\}^3$,
the component $\mathcal{D}_{i,j,k}$ of $\mathcal{D}_{iso}$ is defined as follows (compare with \eqref{eq:diffusion_kernel-2D}):

\begin{equation}\label{eq:diffusion_kernel-3D}
    \mathcal{D}_{i,j,k} := \radial{}\left(\frac{\sqrt{i^2 + j^2 + k^2}}{\Krad\sqrt{3}}\right).
\end{equation}

The kernel has size $(2\Krad+1)\times(2\Krad+1)\times(2\Krad+1)$,
and $\Krad\sqrt{3}$ is the distance, in grid units, from its center to a
corner. We use the same value $\Krad=8$ as in two dimensions.
As in the 2D case (see \eqref{eq:convolution_2d}), we define the discrete convolutional operation in 3D as follows:

\begin{equation}\label{eq:discrete_convolution_3d}
    (\mathcal{D}*u)_{i,j,k} := \sum_{i'=-\Krad}^{\Krad}\sum_{j'=-\Krad}^{\Krad}\sum_{k'=-\Krad}^{\Krad}\mathcal{D}_{i',j',k'}u_{i-i',j-j',k-k'}.
\end{equation}

In order to prevent the training from stalling on a plateau, all the linear
channels are taken without bias.

\begin{figure}[H]
    \centering
    \includegraphics[width=0.6\textwidth]{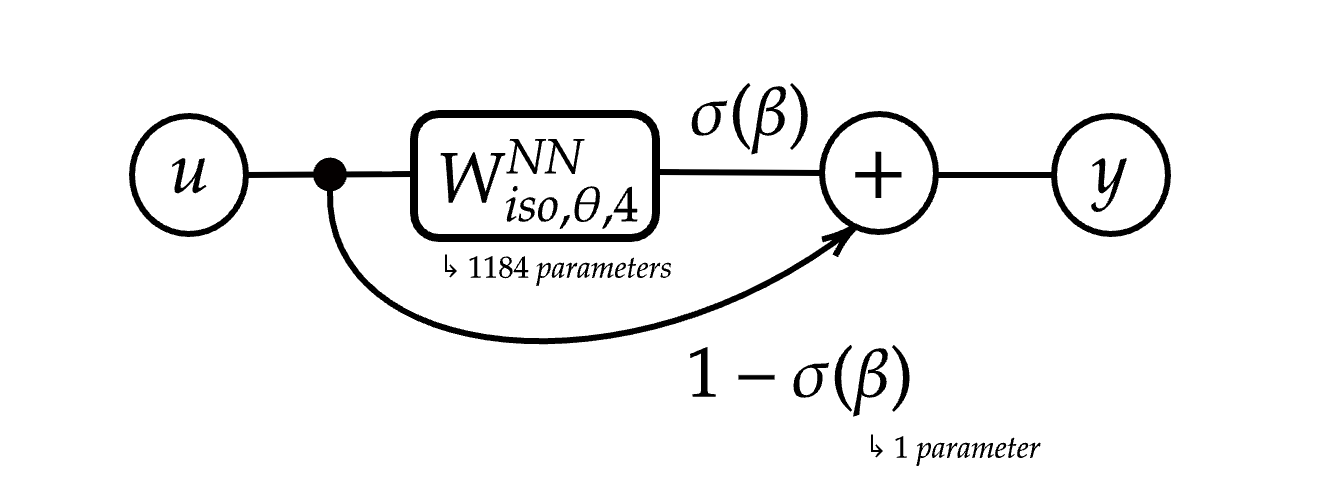}
    \caption{Structure of $\wpnn{\theta}{2}$ in 3D ($1185$ learnable parameters)}
    \label{fig:structure_wnn2_3d}
\end{figure}

Moreover, we introduce a learnable scaling parameter $\beta\in\mathbb{R}$ at the output of the neural network, providing additional flexibility in fitting the training dynamics. More precisely, we redefine the network output as
\begin{equation}\label{eq:scaling_3d}
    \mathcal{\tilde{N}}(\bv u) :=  \sigma(\beta)\mathcal{N}(\bv u) + (1 - \sigma(\beta))\bv u
\end{equation}
for $\bv u\in\mathbb{R}^{N\times N\times N}$. Here, $\mathcal{N}$ equals either $\wnnSimple{\theta}{2}$ or $\wpnnSimple{\theta}{2}$, and $\sigma:\mathbb{R}\to(0,1)$ denotes the sigmoid function $\sigma(s)=\frac{1}{1+e^{-s}}$. For simplicity, we continue to denote the scaled versions $\wnnSimple{\theta}{2}$ and $\wpnnSimple{\theta}{2}$ (see \Figure{fig:structure_wnn2_3d}).
We note that the number of learnable parameters in $\wpnnSimple{\theta}{2}$ is almost the same as in the 2D case, since only the radial MLPs $\radial{\theta}$ are trained rather than the full convolution kernels.

To learn the correct dynamics, we incorporate first-order differences in the loss function.
More precisely, we modify the definition of the loss function as follows (typically $S = 2$):

\begin{equation}\label{eq:loss_multistep_3d}
    \begin{aligned}
    \loss{\mathcal{N}}(X, \bv Y)
    := \frac{1}{BSN^3}\sum_{b=1}^B\sum_{s=1}^{S}\sum_{i=1}^N\sum_{j=1}^N\sum_{k=1}^N
    \Biggl[\left({\mathcal{N}}^{(s)}(X^{(b)}(i,j,k)) - Y_s^{(b)}(i,j,k)\right)^2\\
    + \left|\left(\mathcal{N}^{(s)}(X^{(b)}(i,j,k)) - Y_{s-1}^{(b)}(i,j,k)\right)^2
    - \left(Y_s^{(b)}(i,j,k) - Y_{s-1}^{(b)}(i,j,k)\right)^2\right|^2
    \Biggr].
    \end{aligned}
\end{equation}

\begin{rem}\label{rem:practical_3d}
We observe that when the output of equation \eqref{eq:scaling_3d} is not scaled and a residual connection is assumed,
the dynamics yielded by the learned model are much faster than the reference one.
The convex combination of equation \eqref{eq:scaling_3d} enables the network to learn the magnitude of the update, regardless of its form, thereby reducing uncertainty in the flow dynamics.
Moreover, the choice of loss function \eqref{eq:loss_multistep_3d} compares the squared increment produced by the network to the squared increment of the data rather than to the states themselves.
This penalizes networks that predict the correct direction but with an incorrect amplitude. In our numerical experiment, using only the first term causes the scale parameter  $\beta$ in equation \eqref{eq:scaling_3d} to converge to a very small value, leaving the phase field virtually unchanged. Adding the second term restores dynamics similar to those of the reference case.
\end{rem}

\begin{figure}[H]
    \centering
    \begin{subfigure}[b]{0.45\textwidth}
        \centering
        \includegraphics[width=\linewidth]{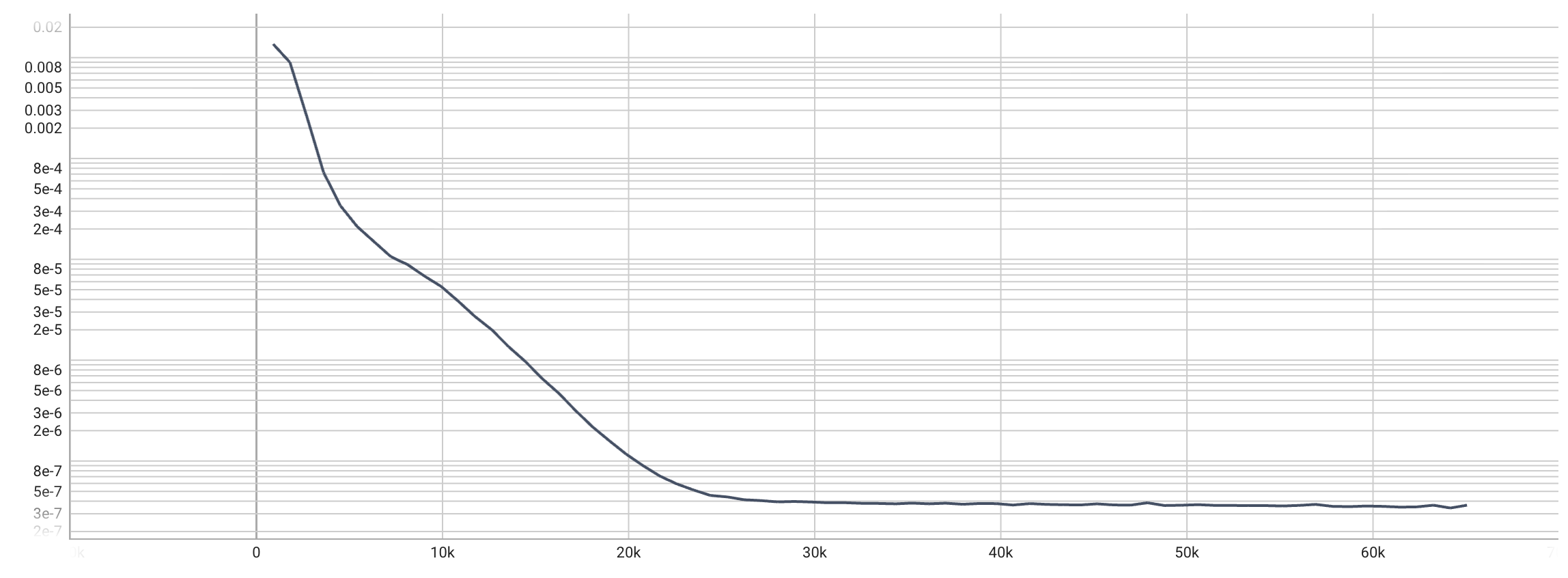}
        \caption{Training loss}
    \end{subfigure}
    \hfill
    \begin{subfigure}[b]{0.45\textwidth}
        \centering
        \includegraphics[width=\linewidth]{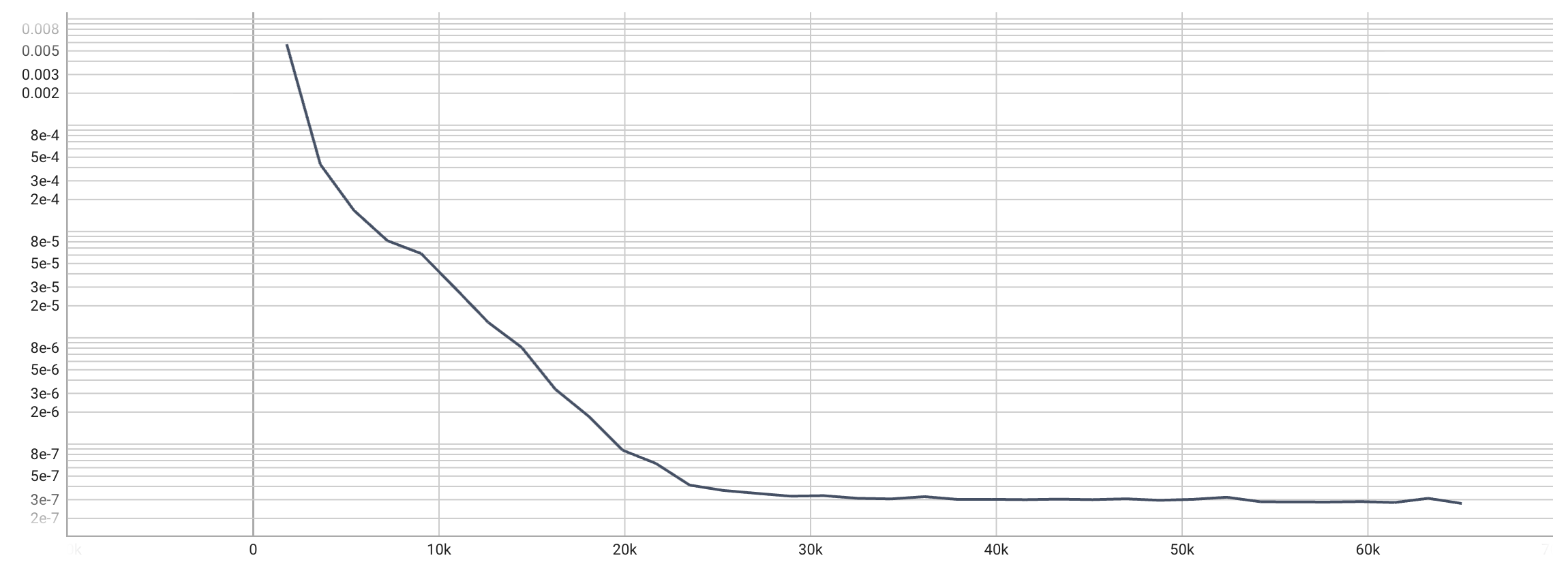}
        \caption{Validation loss}
    \end{subfigure}
    \caption{Evolution of the training (left) and validation (right) losses at each training step for the rotationally invariant model $\wpnnSimple{\theta}{2}$ in 3D.}
    \label{fig:train_loss_wnn4_3d_enriched}
\end{figure}

\begin{figure}[H]
  \centering
  \setlength{\tabcolsep}{4pt}   
  \renewcommand{\arraystretch}{1.0} 
  \begin{tabular}{m{0.14\textwidth} m{0.19\textwidth} m{0.19\textwidth} m{0.19\textwidth} m{0.19\textwidth}}
    \raggedleft \small $RSS$ &
      \includegraphics[width=\linewidth]{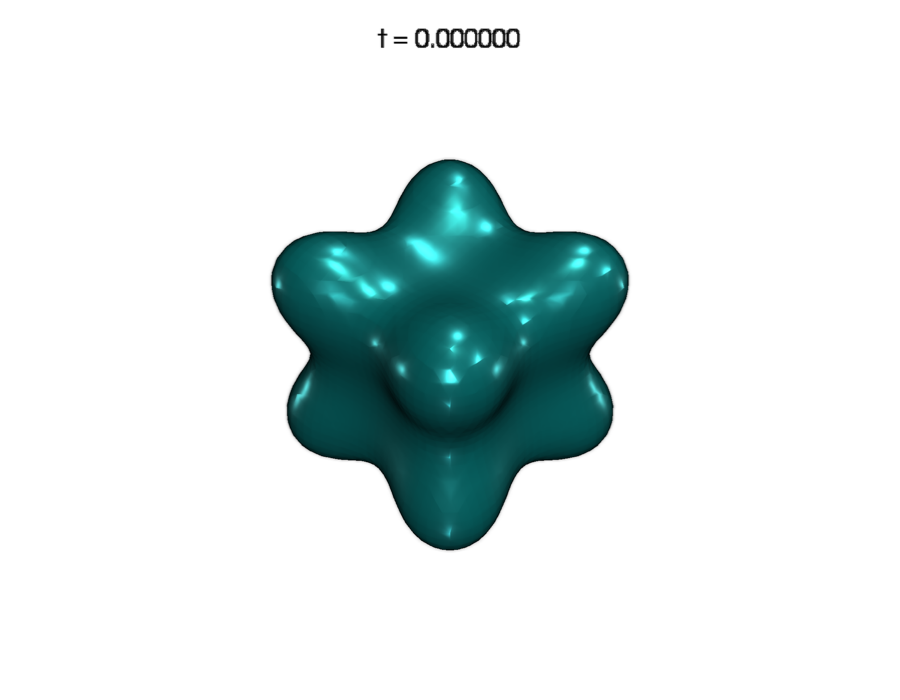} &
      \includegraphics[width=\linewidth]{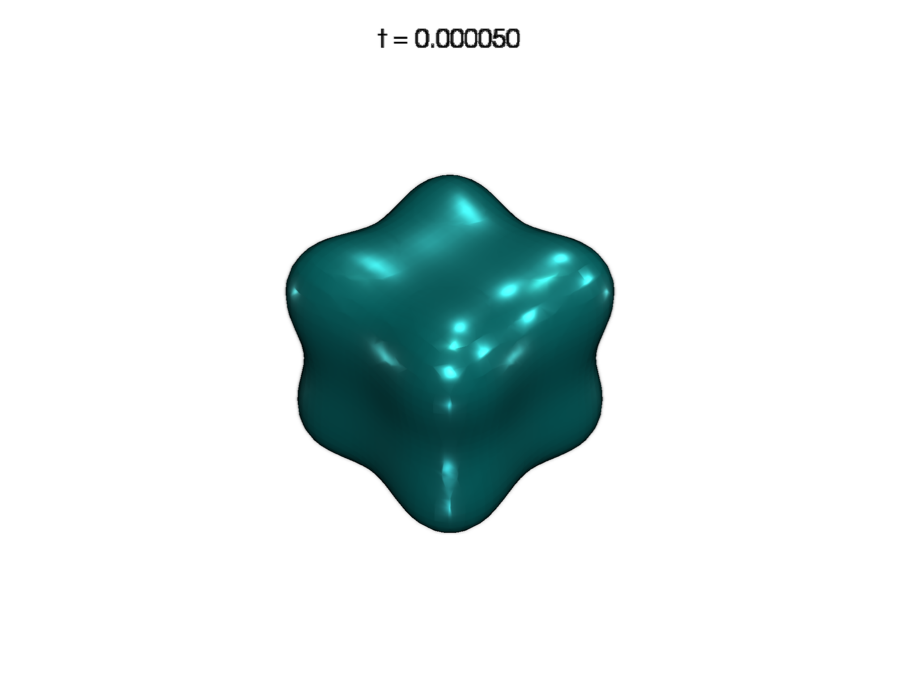} &
      \includegraphics[width=\linewidth]{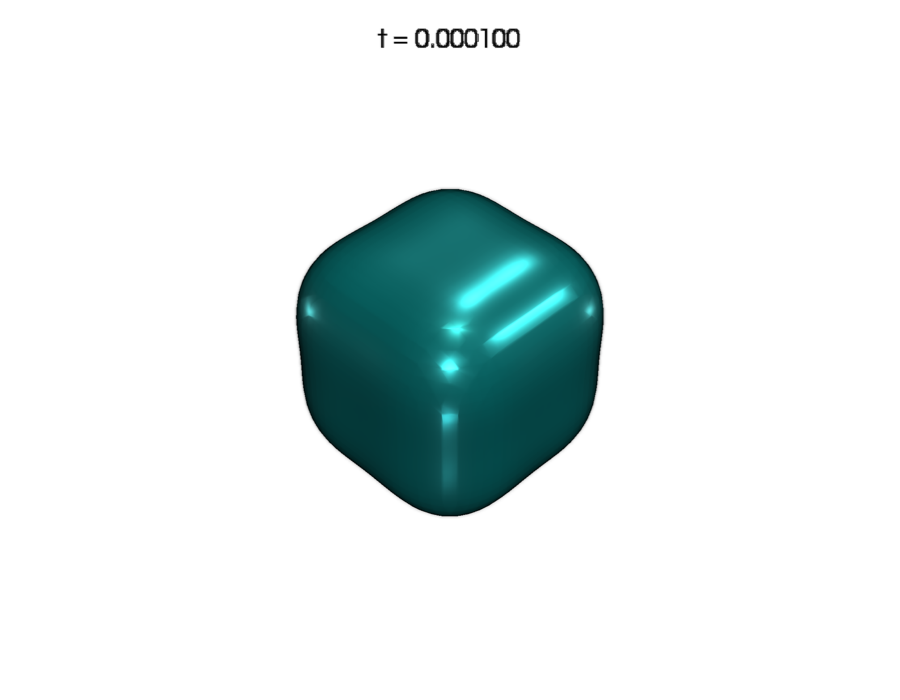} &
      \includegraphics[width=\linewidth]{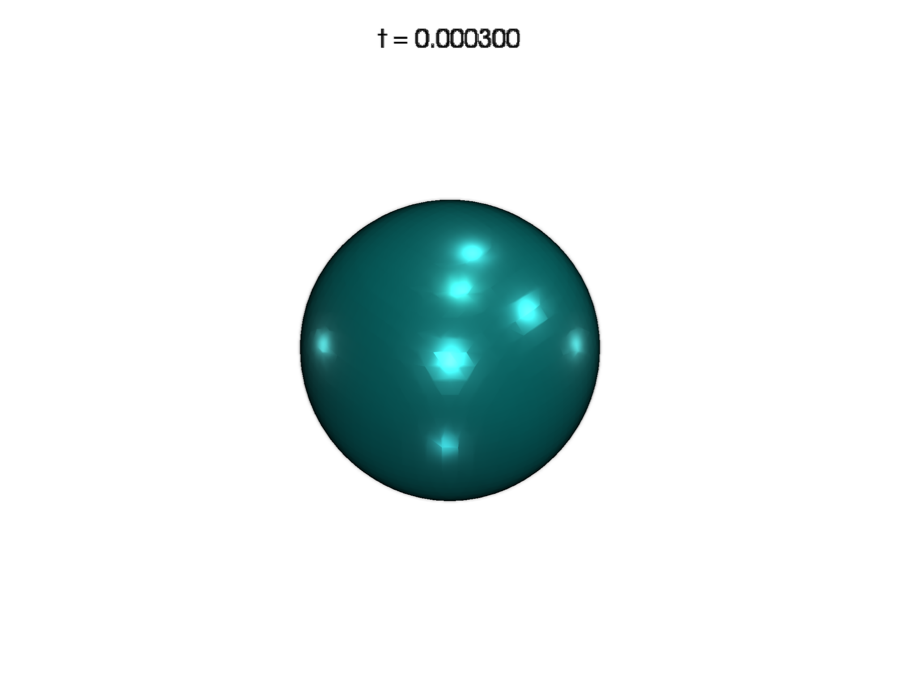} \\
    \raggedleft $\wpnn{\theta}{2}$ &
      \includegraphics[width=\linewidth]{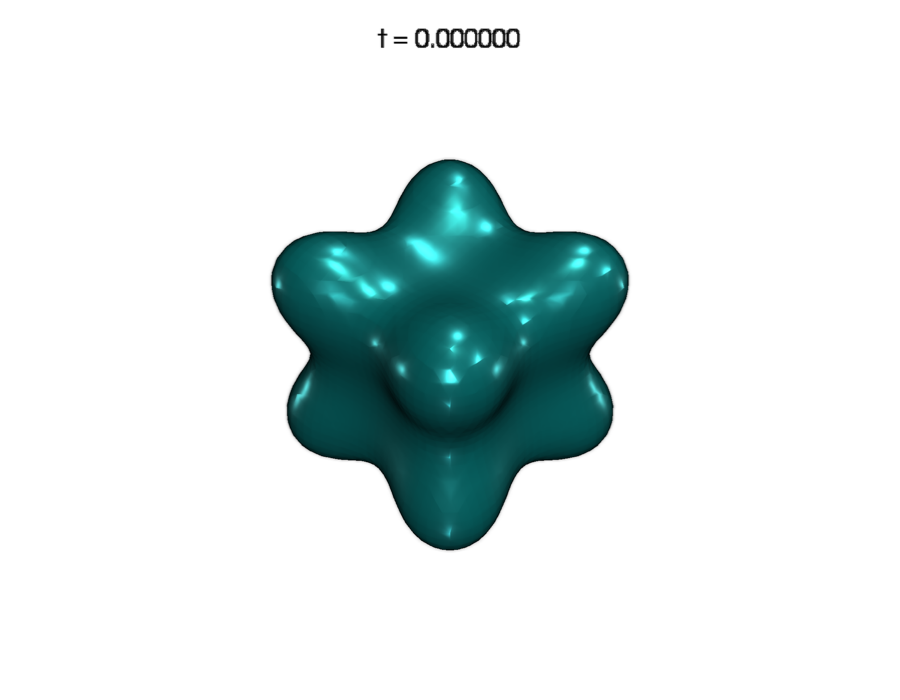} &
      \includegraphics[width=\linewidth]{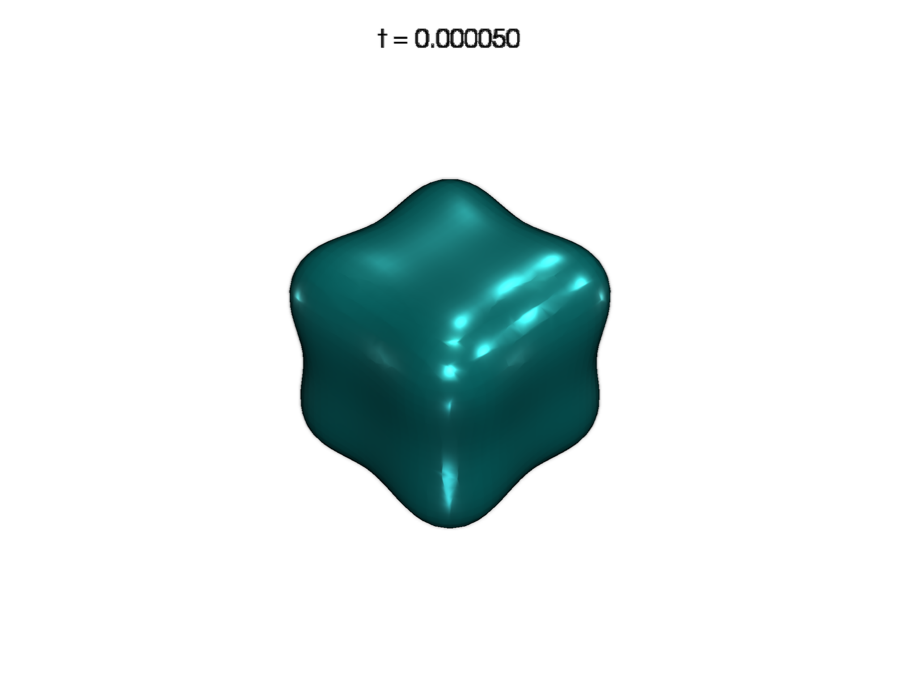} &
      \includegraphics[width=\linewidth]{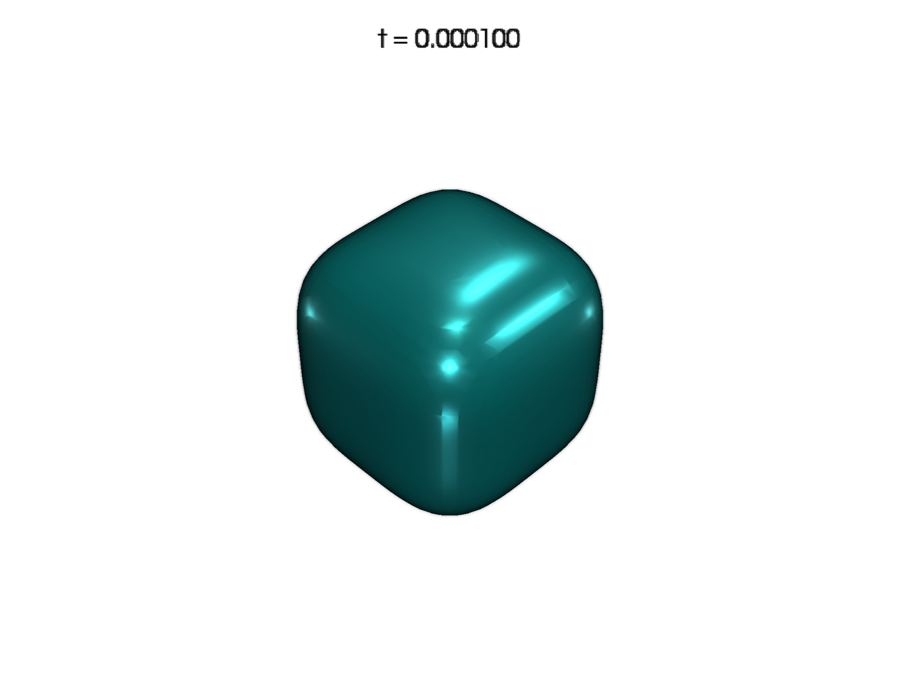} &
      \includegraphics[width=\linewidth]{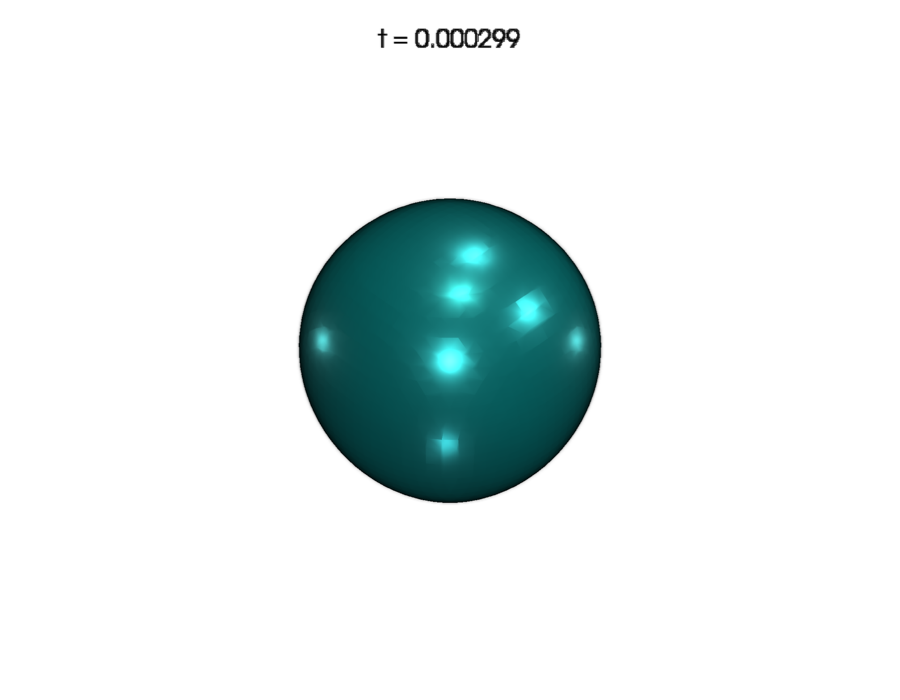} \\
  \end{tabular}
  \caption{Approximation of the oriented Willmore flow at different time steps for a smooth spiky initial surface.
First row: RSS scheme.
Second row: Repetition of neural network $\wpnn{\theta}{2}$.}
  \label{fig:flower_3d}
\end{figure}

\begin{figure}[H]
    \centering
    \begin{subfigure}[b]{0.70\textwidth}
        \centering
        \includegraphics[width=\linewidth]{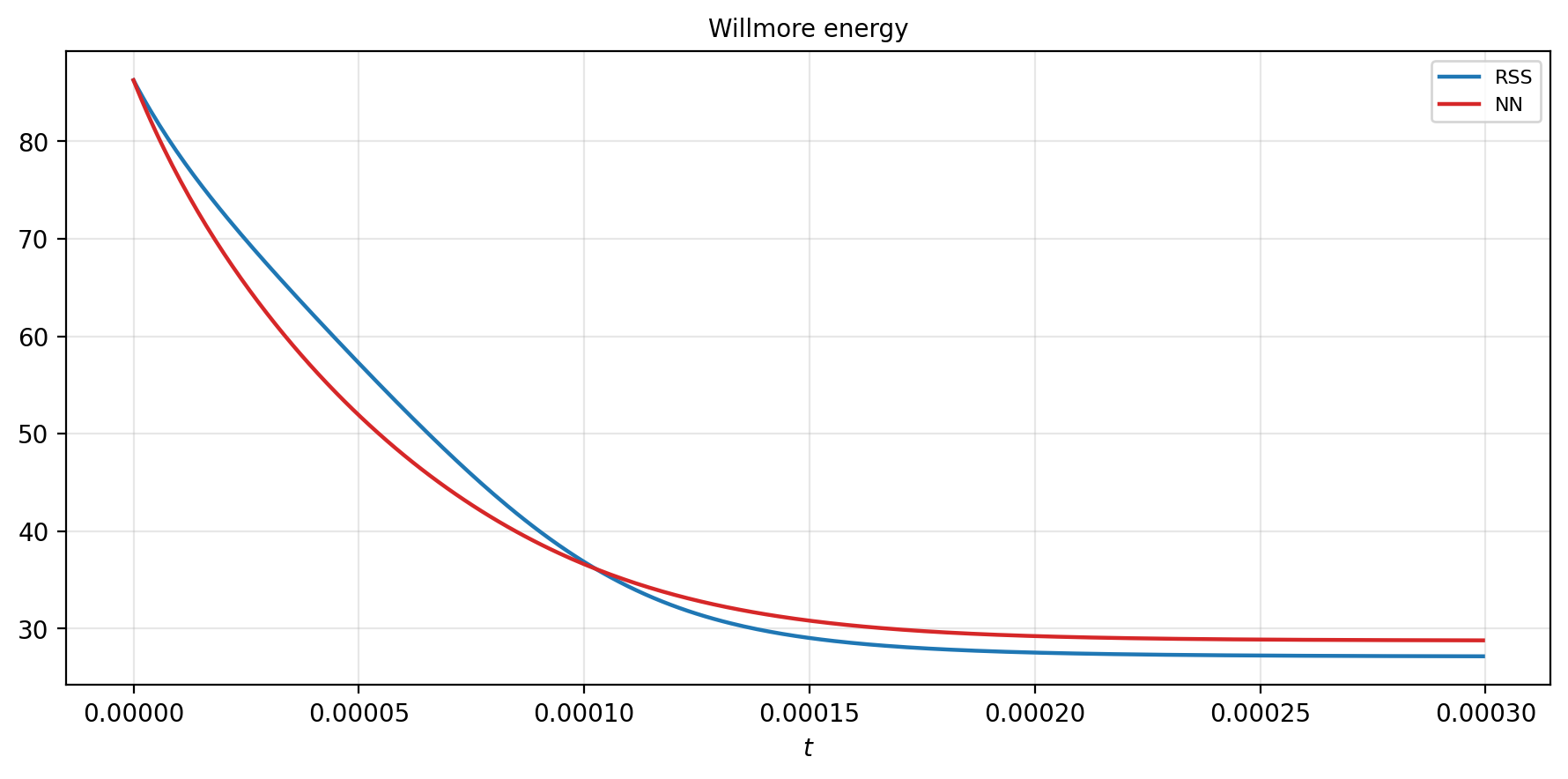}
    \end{subfigure}
    \caption{Plot of the Willmore energy corresponding to \Figure{fig:flower_3d} defined by \eqref{eq:classical-willmore-energy}.}
    \label{fig:energy_flower_3d}
\end{figure}

\Figure{fig:train_loss_wnn4_3d_enriched} illustrates the loss function values during training,
particularly highlighting successful learning and stabilization on both the training and validation data.
\Table{tab:training-nn-3d-enriched} shows the loss function value at the selected checkpoint.
\Figure{fig:flower_3d} shows a comparison between the numerical flow obtained using the neural network $\wpnn{\theta}{2}$
and the 3D RSS scheme, with a smooth flower as the initial set $\Omega_0$.
We also plot the phase-field approximation of the Willmore flow during iteration in \Figure{fig:energy_flower_3d}.
Both schemes reduce the energy and appear very similar to one another.

\begin{table}[H]
    \centering
    \begin{tabular}{|c|c|c|c|}
        \toprule
        \textbf{Model} &  \textbf{Epochs} & \textbf{Steps} & \textbf{Validation loss} \\
        \midrule
        $\wpnn{\theta}{2}$ & $67$ & $61471$ &  $2.82\times 10^{-7}$ \\
        \bottomrule
    \end{tabular}
    \caption{Validation loss of the checkpoint selected for the oriented surface
    model in 3D; the selection rule is the one described in
    \Section{sec:database_construction}.}
    \label{tab:training-nn-3d-enriched}
\end{table}


As in the two-dimensional case, this test shows that our networks can closely
match the flows obtained with the RSS scheme at a significantly lower
computational cost (see \Table{tab:computation-time-3d-q}).

\begin{table}[H]
    \centering
    \begin{tabular}{|c|c|c|}
        \toprule
        \textbf{Scheme} & \textbf{Total time (s)} & \textbf{Time / Step (s)} \\
        \midrule
        $RSS$ & $1040$ & $10.1$ \\
        $\wpnn{\theta}{2}$ & $4.10$ & $0.0394$ \\
        \bottomrule
    \end{tabular}
    \caption{Comparison of computation times of the numerical schemes.
    First row: RSS scheme.
    Second row: Repetitions of $\wpnn{\theta}{2}$.}
    \label{tab:computation-time-3d-q}
\end{table}

\section{Case of the non-oriented Willmore flow in 3D}\label{sec:nonoriented_3d}

\begin{figure}[H]
    \centering
    \begin{subfigure}[b]{0.8\textwidth}
        \centering
        \includegraphics[width=\linewidth]{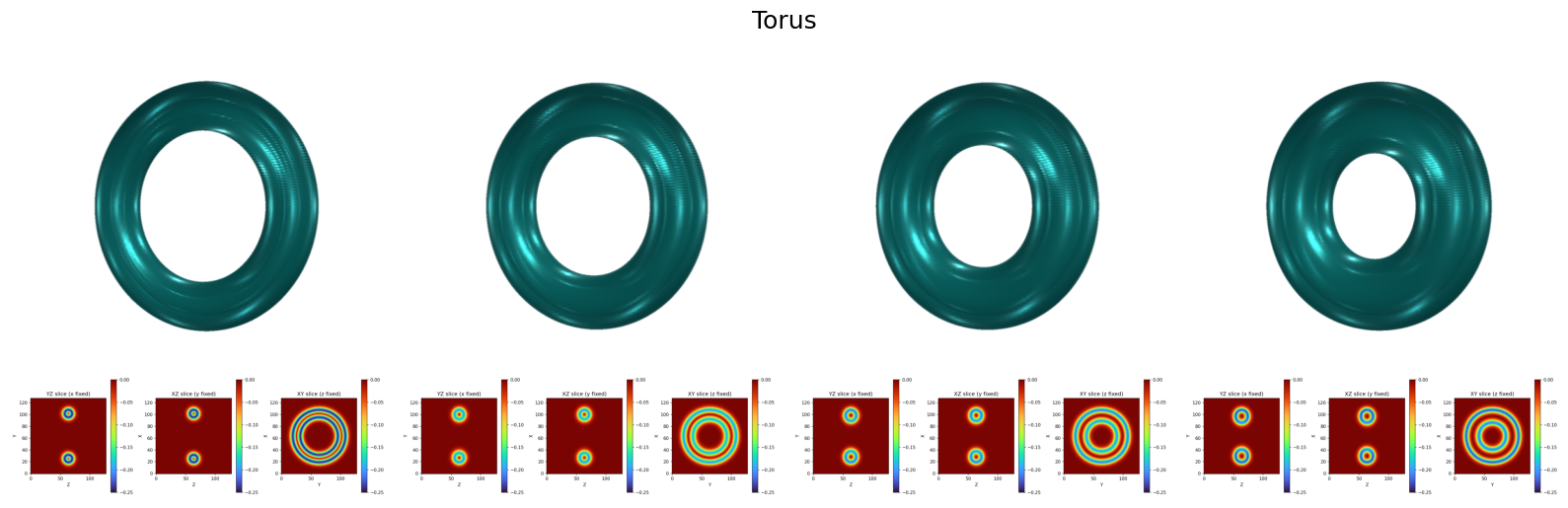}
    \end{subfigure}
    \caption{Dataset of the Willmore flow of a torus generated by the oriented 3D RSS scheme and post-processed using the operator $\operatorprim$. Each surface represents the $(-0.25,-0.15)$-level set of the phase-field function at each time step.}
    \label{fig:dataset_nonoriented_3d_torus1}
\end{figure}

We now turn to the non-oriented case, thereby illustrating the flexibility of our approach in a setting that the RSS method does not currently seem able to handle. As in two dimensions, we first construct a non-oriented training dataset by post-processing the oriented data with the operator
$
\operatorprim(U) = -\varepsilon\,|\nabla U|
$, 
which transforms the oriented phase-field representation into a non-oriented one.
The dataset consists of triplets extracted from the numerical approximation of the flows of cubes, rectangular prisms, ellipsoids, unions of two spheres and tori (see \Figure{fig:dataset_oriented_3d_set1}).
\Figure{fig:dataset_nonoriented_3d_torus1} shows the non-oriented representation of a torus obtained after applying the operator $\operatorprim $.

In practice, however, we were unable to train these networks from scratch as in the 2D case. Despite experimenting with several learning rates, optimization on the non-oriented three-dimensional data remained difficult. We therefore exploit the isotropy of the Willmore flow to initialize the 3D network from its already trained 2D counterpart. More precisely, we transfer the radial Fourier profiles of the learned convolution kernels, keep the pointwise parameters unchanged, and then fine-tune the resulting model on the 3D data.

\subsection{Transferring the trained 2D kernels to three dimensions}
\label{sec:transfer_2d_3d}

An isotropic convolution operator has a radial Fourier symbol. This does
not by itself imply that the symbol is independent of dimension. Our
initialization strategy is to preserve the learned radial multiplier when
passing from two to three dimensions, as one can do for diffusion operators
specified by a function of $|\xi|^2$. The resulting transfer is a modeling
choice, whose usefulness is assessed by the experiments below, rather than
a consequence of isotropy alone.

Formally, for $d = 2$ and $3$, let $\ft{d}{\cdot}$ and $\ift{d}{\cdot}$ denote
the $d$-dimensional Fourier transform and its inverse.  Given a radially
invariant kernel $\radial{2D}:\mathbb{R}^2\to\mathbb{R}$, we define its radial
Fourier profile
\begin{equation*}
    \phi_{\text{Fourier}}(\rho) := \frac{1}{2\pi}\int_{\mathbb{S}^1} \ft{2}{\radial{2D}}(\rho \theta)\,d\theta,\qquad \rho \geq 0,
\end{equation*}
and the three-dimensional kernel it generates,
$\radial{3D} := \ift{3}{\phi_{\text{Fourier}}(|\cdot|)}$.

In practice, the transfer is performed directly on the learned kernels, discretely and layer by layer. The trained 2D $4$-parallel model contains five diffusion layers, each involving two kernels, for a total of ten rotationally invariant kernels (see again \Figure{fig:nn_structure}). For each kernel, we consider the learned 2D stencil
$K^{(\ell)}_{2\mathrm{D}} \in \mathbb{R}^{(2\Krad+1)\times(2\Krad+1)}$, pad it
with zeros to the size of the computational grid and apply a centered FFT.
Averaging the result over the frequency shells, that is, over the sets of
frequency indices whose radius $|\xi|$ falls in a common bin, yields samples of
the radial profile $\phi_{\text{Fourier}}^{(\ell)}$, which we extend to $\mathbb{R}_+$ by linear
interpolation.  Setting $\widehat{K}^{(\ell)}_{3\mathrm{D}}(\xi) :=
\phi_{\text{Fourier}}^{(\ell)}(|\xi|)$ on the three-dimensional frequency grid and inverting the
transform produces a spatial kernel $K^{(\ell)}_{3\mathrm{D}}$ of size
$(2\Krad+1)\times(2\Krad+1)\times(2\Krad+1)$, which we assign to the corresponding 3D
diffusion module.  All the remaining parameters of the 2D model, which do not
depend on the space dimension, are copied unchanged. \\

Before going further, it is worth seeing what this transfer alone produces.
\Figure{fig:spikyball_onlyinit_3d} shows the flow obtained, without additional
3D training, starting from a non-oriented spiky ball with $8$ spikes.
The transferred network maintains a coherent interface in this experiment without additional 3D training. Secondly, the interface appears to converge towards a sphere, the critical point of the Willmore energy. However, upon closer inspection and by comparing this evolution with that shown in \Figure{fig:flower_3d}, it becomes clear that, although similar, the two evolutions do not exactly correspond, suggesting the need to retrain the network.

\begin{figure}[H]
  \centering
  \setlength{\tabcolsep}{4pt}   
  \renewcommand{\arraystretch}{1.0} 
  \begin{tabular}{m{0.23\textwidth} m{0.23\textwidth} m{0.23\textwidth} m{0.23\textwidth}}
      \includegraphics[width=\linewidth]{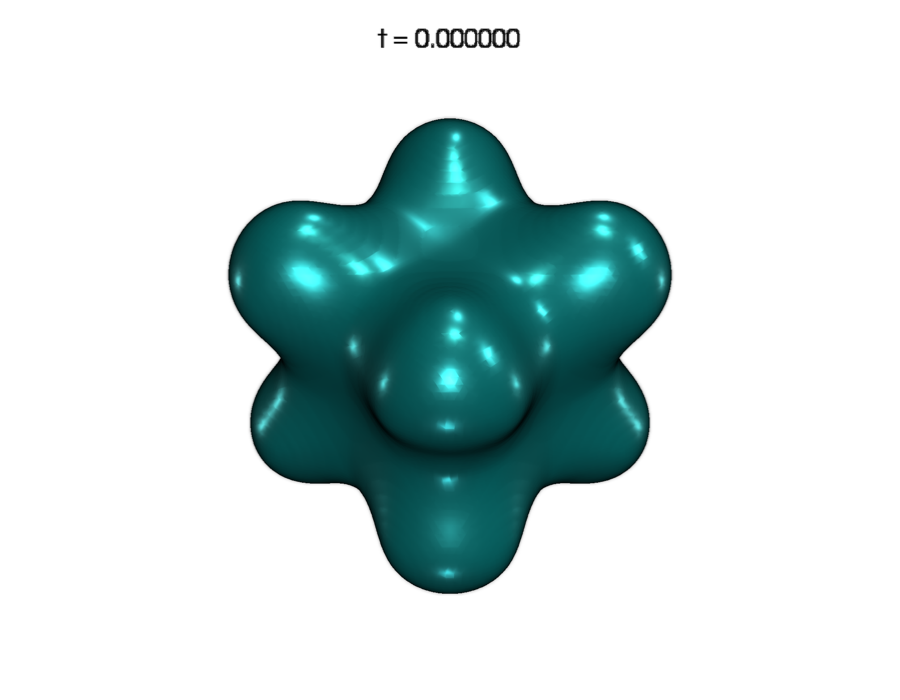} &
      \includegraphics[width=\linewidth]{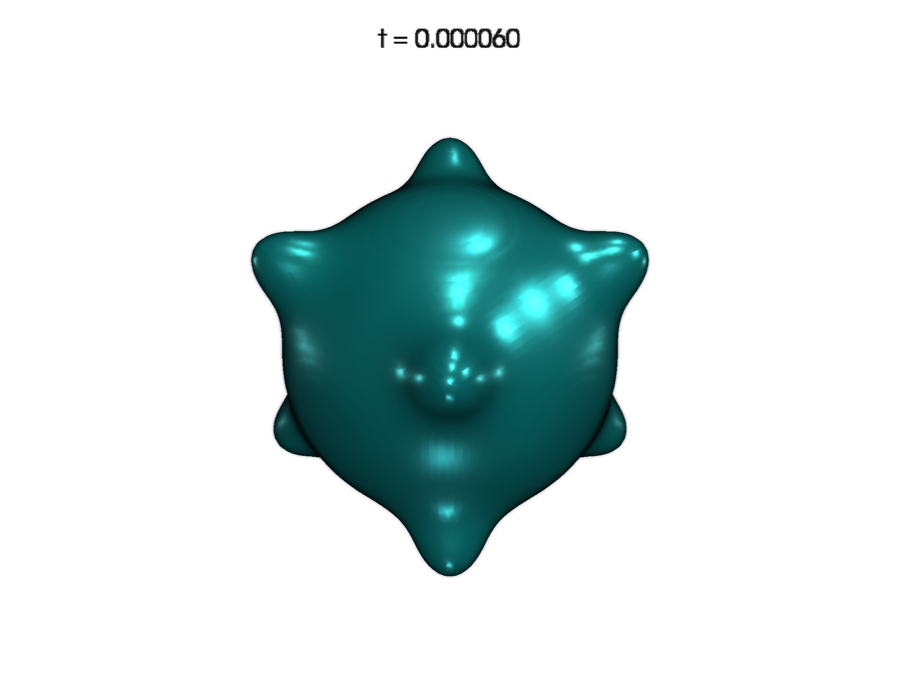} &
      \includegraphics[width=\linewidth]{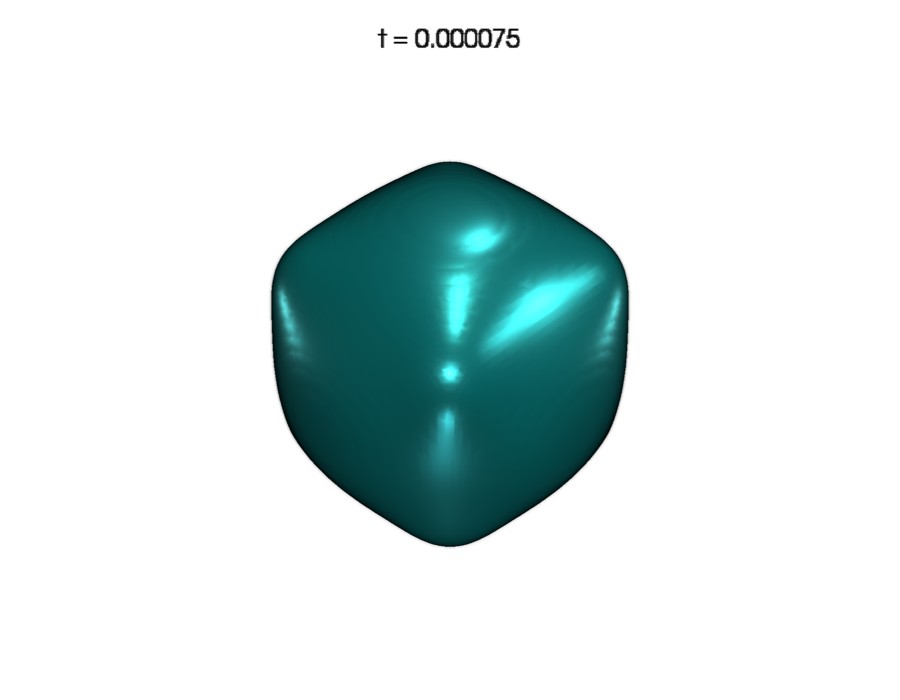} &
      \includegraphics[width=\linewidth]{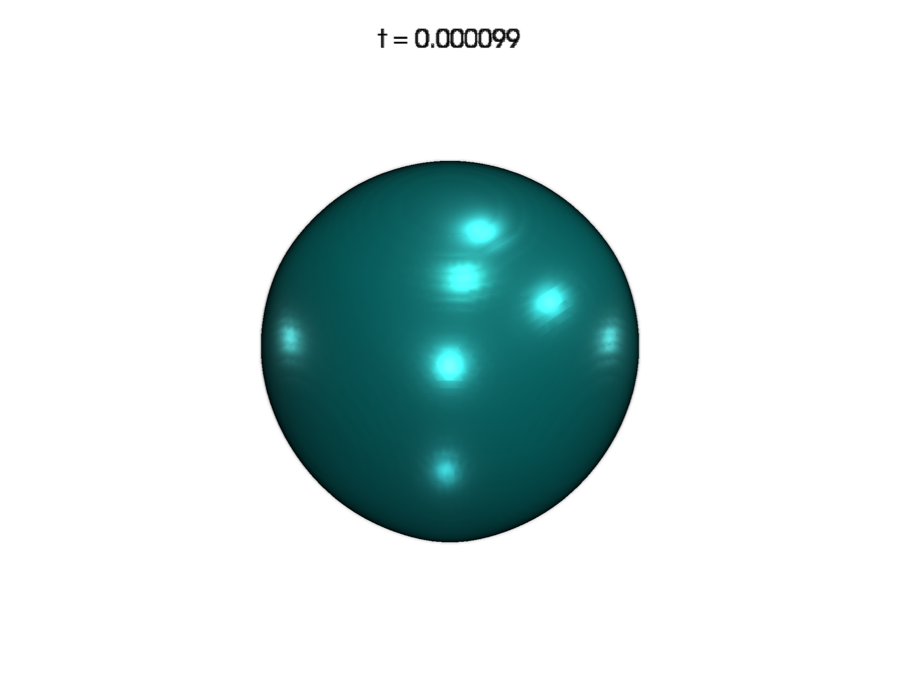} \\
  \end{tabular}
  \caption{Evolution of a non-oriented spiky ball in 3D using the free kernel model $\wnnprim{\theta}{2}$ initialized with the rotationally invariant kernel model $\wpnnprim{\theta}{2}$ in 2D.
    The surfaces represent the isosurfaces inferred by $\wnnprim{\theta}{2}$ with $\delta = 0.15$.}
  \label{fig:spikyball_onlyinit_3d}
\end{figure}

The construction of the above-described kernels results in free kernels with a non-isotropic representation. Without special treatment, allowing these kernels to remain free during training should lead to the anisotropic behavior observed in this 2D case, thereby degrading the overall approximation of the Willmore flow. To avoid this, we must express the transferred kernels in the radial parameterization from \Section{sec:IsotropicNN_MCF}.
Each 3D diffusion block uses the same pointwise radial parameterization
with ReLU activation as in \Figure{fig:radial_profile_mlp}. Its scalar map
$r\mapsto w_\theta(r)$ is evaluated independently at each kernel position,
with parameters shared across the stencil, using the normalized radii
\[
    \rho_m := \frac{\sqrt{i^2 + j^2 + k^2}}{\Krad\sqrt{3}},
\]
where $(i,j,k)$ is the offset of the $m$-th voxel from the center of the
stencil, with $-\Krad\leq i,j,k\leq\Krad$ and $m=1,\dots,(2\Krad+1)^3$.  We group the voxels into spherical shells
of equal discretized radius, average the transferred kernel
$K^{(\ell)}_{3\mathrm{D}}$ over each shell to obtain targets
$(\bar\rho_p, \bar w^{(\ell)}_p)_{p = 1,\dots,P}$, and fit the MLP to them in
the least-squares sense,
\[
    \min_{\theta^{(\ell)}}\ \frac{1}{P}\sum_{p=1}^{P}
    \Bigl(w_{\theta^{(\ell)}}(\bar\rho_p)-\bar w^{(\ell)}_p\Bigr)^2,
\]
with a few hundred ADAM steps on $\theta^{(\ell)}$ alone.  The parameters that
do not concern the diffusion blocks are copied from the transferred model.  The
rotationally invariant 3D model obtained in this way is the one we train in the next subsection.

\subsection{Training and results}

We optimize the parameters of $\wpnnSimple{\theta}{2}$ by minimizing the loss
function \eqref{eq:loss_multistep_3d} on the training dataset described in \Figure{fig:dataset_oriented_3d_set1}.
This training was run for $2$ epochs only, that is, the model saw the whole
training set twice.  This short training illustrates the benefit of the construction: the transfer described above places the network in a regime where only a mild correction is needed, whereas training from a random initialization
did not converge at all.

After that, we again implement the inference of the trained model $\wpnnprim{\theta}{2}$ in 3D with the same setting as in \Figure{fig:spikyball_onlyinit_3d}.
In \Figure{fig:spikyball_init_et_training_3d}, we observe a similar flow of surfaces to that obtained by the model for the oriented case, which is shown in \Figure{fig:flower_3d}.

\begin{figure}[H]
  \centering
  \setlength{\tabcolsep}{4pt}   
  \renewcommand{\arraystretch}{1.0} 
  \begin{tabular}{m{0.23\textwidth} m{0.23\textwidth} m{0.23\textwidth} m{0.23\textwidth}}
      \includegraphics[width=\linewidth]{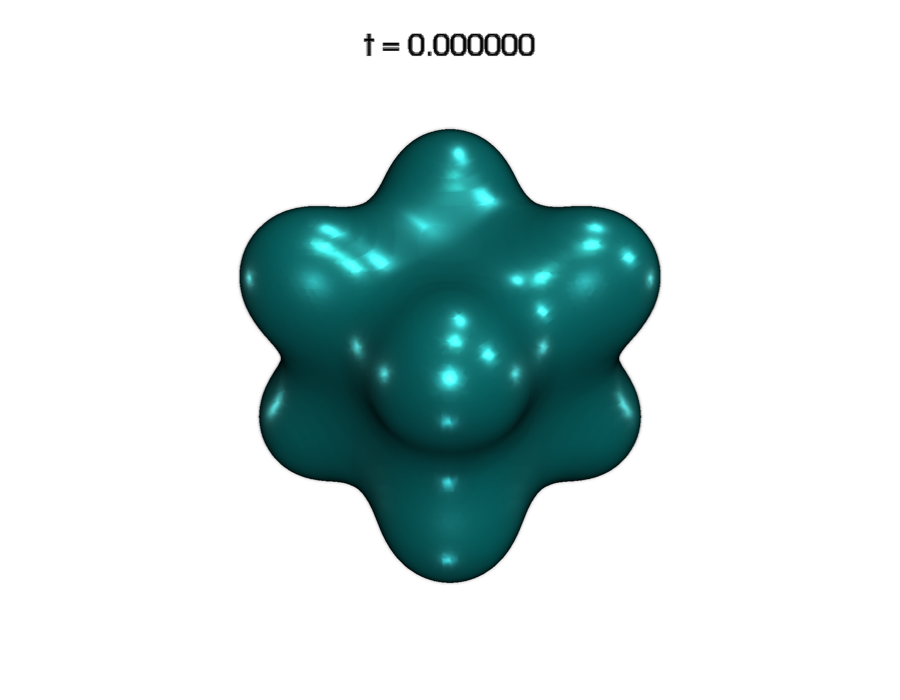} &
      \includegraphics[width=\linewidth]{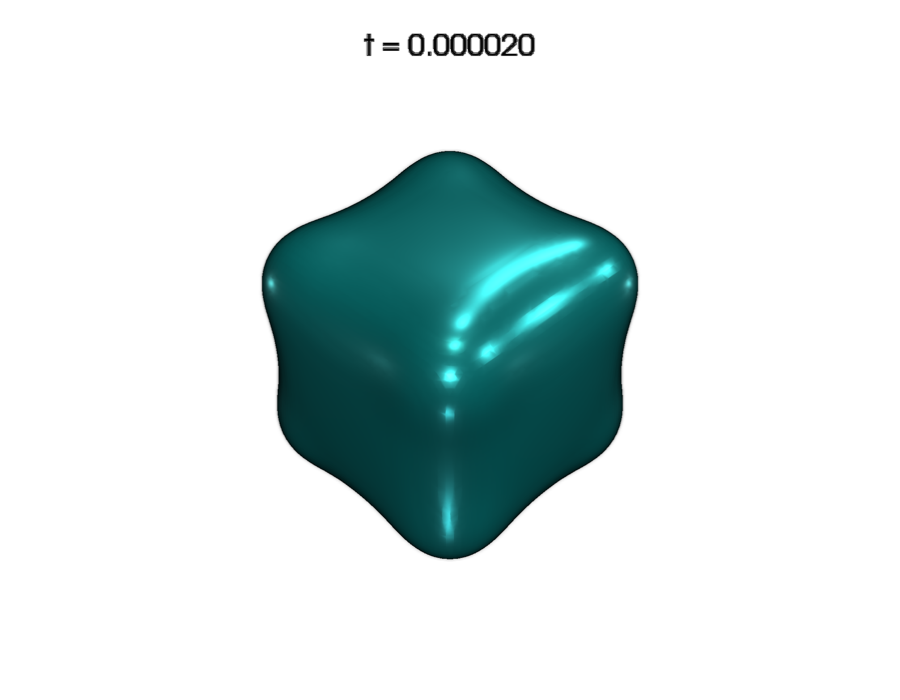} &
      \includegraphics[width=\linewidth]{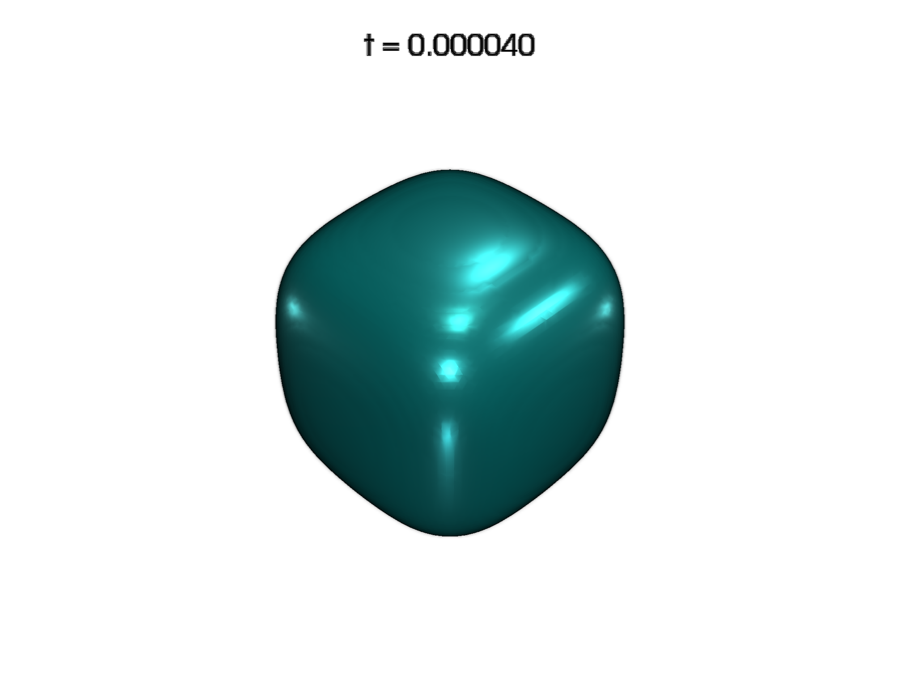} &
      \includegraphics[width=\linewidth]{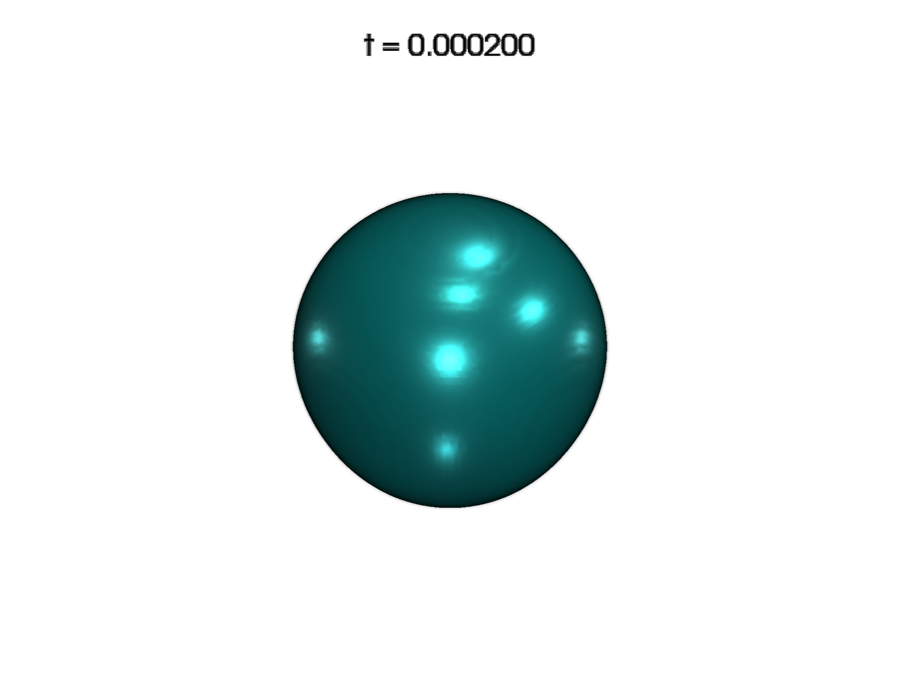} \\
      \includegraphics[width=\linewidth]{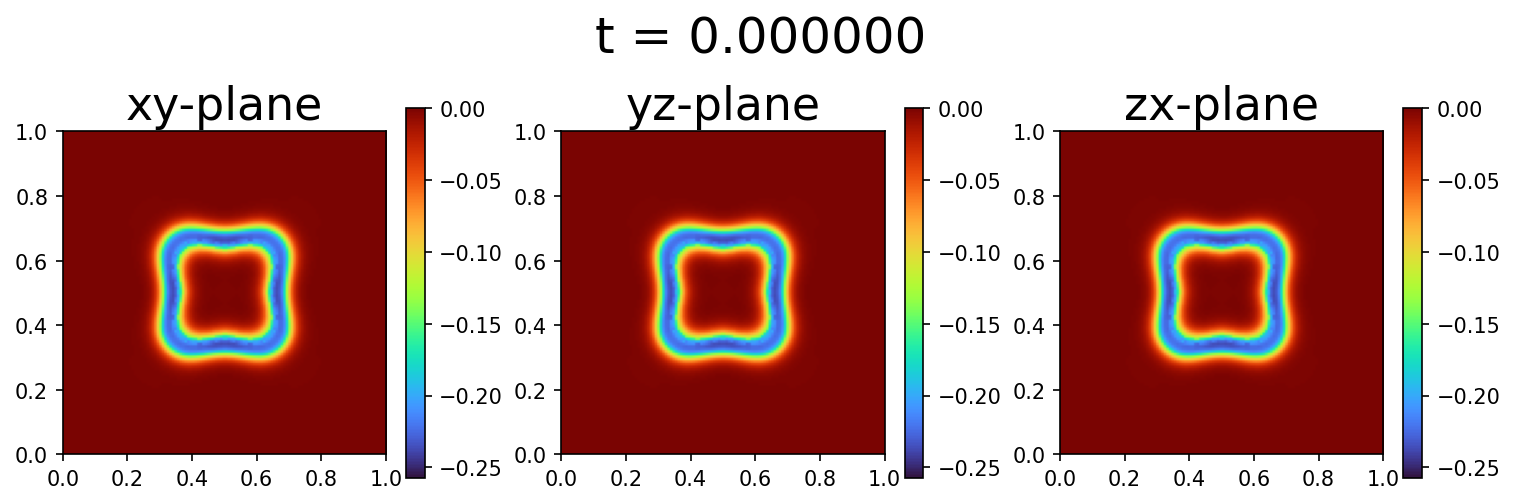} &
      \includegraphics[width=\linewidth]{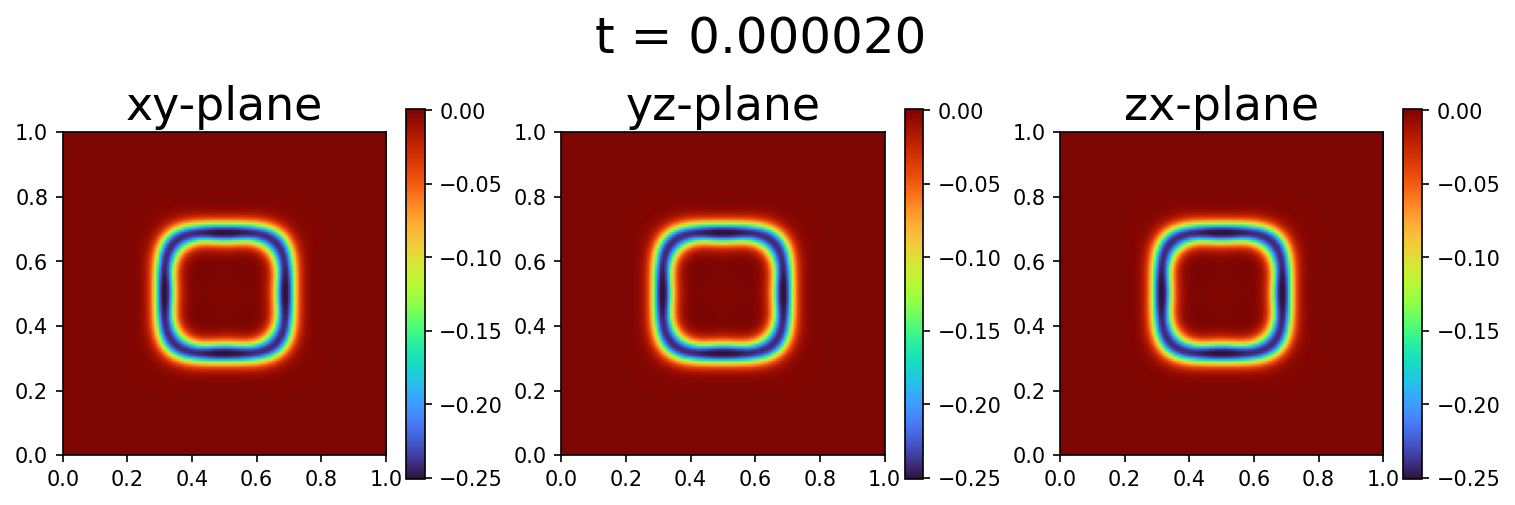} &
      \includegraphics[width=\linewidth]{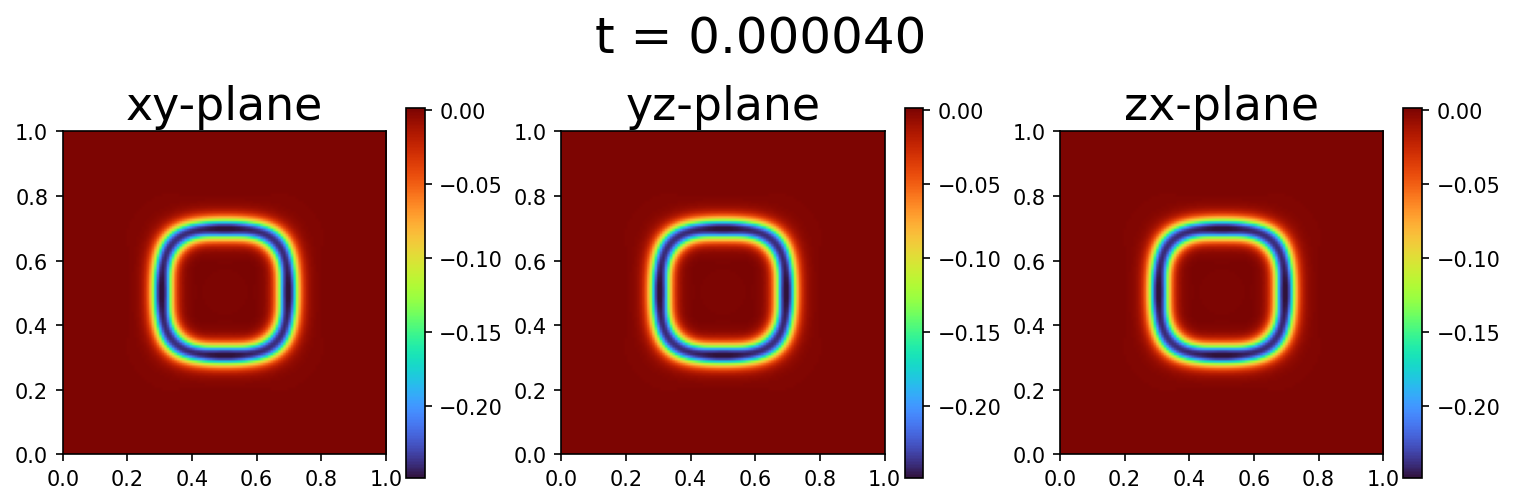} &
      \includegraphics[width=\linewidth]{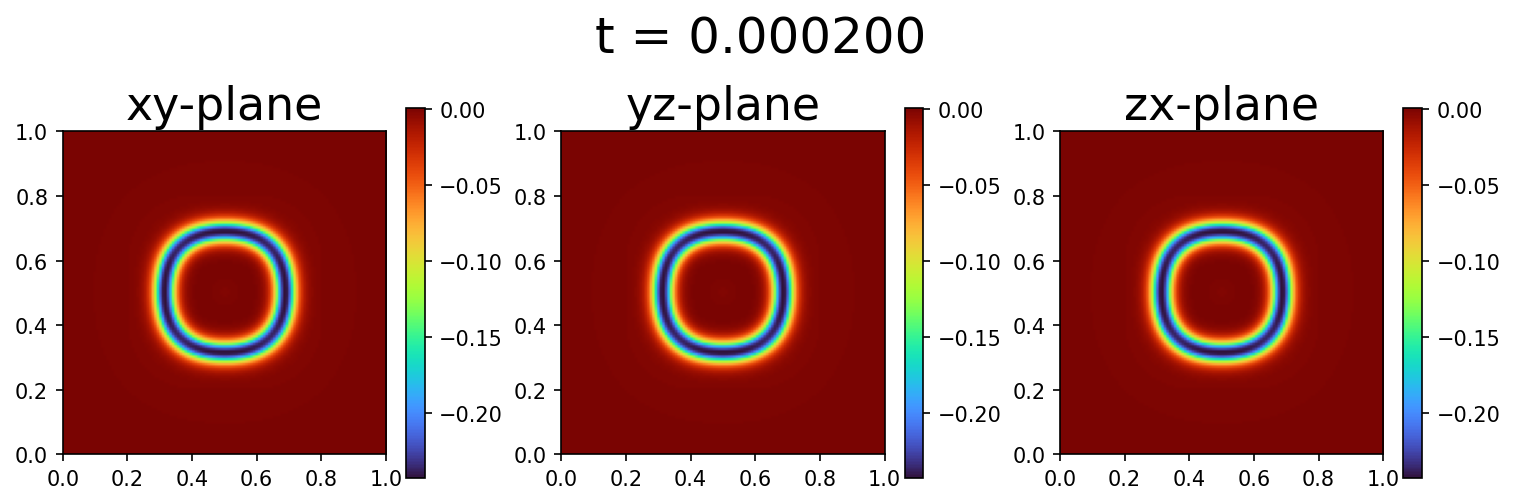} \\
  \end{tabular}
  \caption{
    Approximation of the Willmore flow for a non-oriented spiky ball in 3D using the rotationally invariant model $\wpnnprim{\theta}{2}$
    initialized with the rotationally invariant kernel model $\wpnnprim{\theta}{2}$ in 2D and trained afterwards.
    First row: isosurfaces of $\wpnnprim{\theta}{2}$, second row: slices of $\wpnnprim{\theta}{2}$ at the same time steps.}
  \label{fig:spikyball_init_et_training_3d}
\end{figure}

We conclude that the rotation-invariant model $\wpnnprim{\theta}{2}$ provides a good approximation of the Willmore flow, at least in this example involving unoriented 3D data. It should be noted that this is also thanks to our kernel initialization procedure, which reduced the training time to just a few epochs.

\subsection{Application to surface reconstruction from point clouds}\label{sec:applications_3d}

We will now focus on using our networks for surface reconstruction from point clouds.
This extends the curve reconstruction procedure of \Section{sec:qprim_application} to the case of three dimension.

Another point specific to this section is that we use the simply initialized model $\wnnprim{\theta}{2}$ from \Section{sec:nonoriented_3d};
although it does not accurately approximate the Willmore flow, it proves to be an effective curvature regulator.
The advantage is that reconstruction requires no additional 3D training database; the transferred parameters still originate from the trained 2D model.
Iteration \eqref{eq:reconstruction_iteration} can be interpreted as a decomposition scheme for a curvature-regularized reconstruction model,
as defined in ~\cite{he2020curvature}, where the regularization sub problem is solved by a single network evaluation rather than an internal optimization.

Given a point cloud, we first construct an initial phase field function $U_0$ as in the 2D case.
We then perform several iterations of \eqref{eq:reconstruction_iteration}, followed by several
unconstrained applications of $\wnnprim{\theta}{2}$ to obtain the final smooth
surface.  This experiment extends the reconstruction of smooth surfaces from
point clouds of Bretin, Dayrens and Masnou \cite{Bretin_Dayrens_Masnou17} to the
case where the point cloud carries no information for orientation.

Throughout this section, we set the resolution of the grid to $N = 256$ and the interface
thickness to $\varepsilon = 2/N = 2^{-7}$; the value of the relaxation
parameter $\alpha$ is given in each caption.  The computational grid is
therefore finer than the one used during training, while the ratio
$\varepsilon/h$ is preserved, in accordance with the discussion of \Section{sec:intro_nn}. In each figure below,
we display the isosurface of the phase field function $u$ produced at successive iterations,
the first panel showing the initial point cloud.

First, we consider a point cloud of $2000$ points sampled on a torus. The reconstructed surface, shown in \Figure{fig:reconstruction_torus_qprim_3d_fulltransfer_nn}, is smooth and preserves the hole, indicating that the method respects the shape's global topology.

\begin{figure}[H]
    \centering
    \begin{subfigure}[b]{0.24\textwidth}
        \centering
        \includegraphics[width=\linewidth]{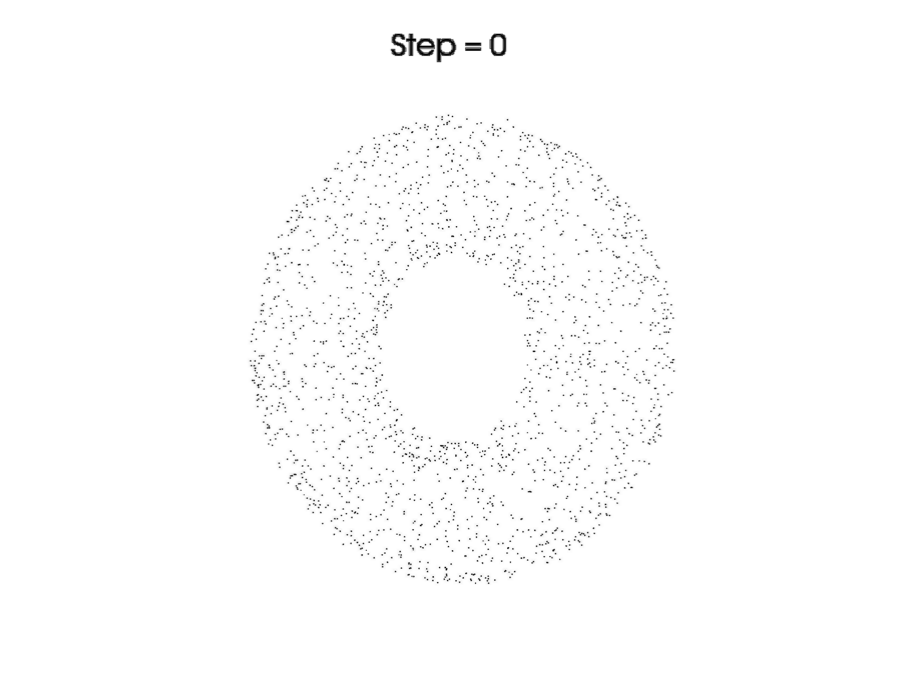}
    \end{subfigure}
    \hfill
    \begin{subfigure}[b]{0.24\textwidth}
        \centering
        \includegraphics[width=\linewidth]{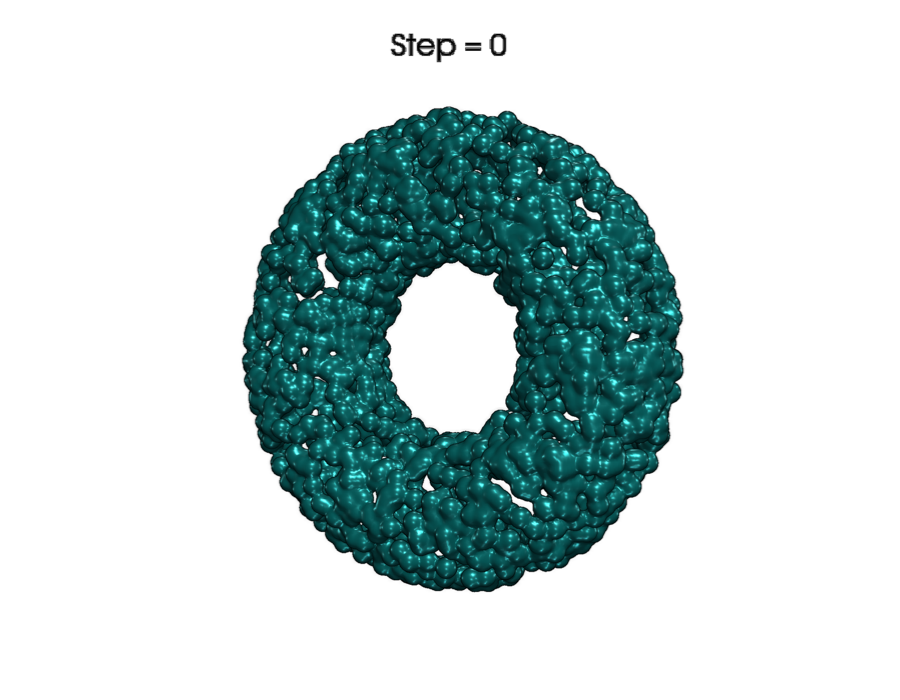}
    \end{subfigure}
    \hfill
    \begin{subfigure}[b]{0.24\textwidth}
        \centering
        \includegraphics[width=\linewidth]{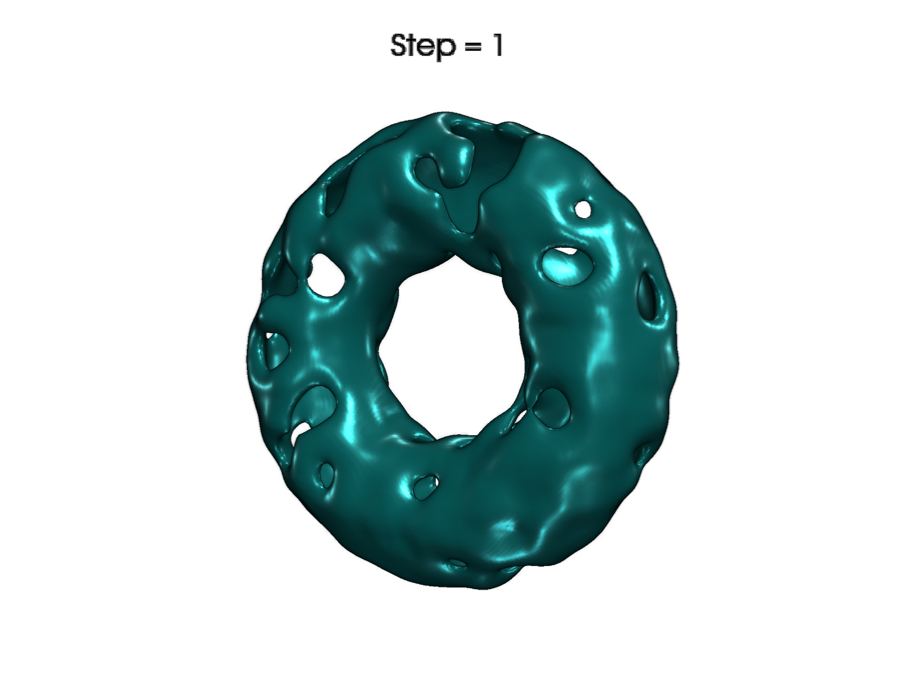}
    \end{subfigure}
    \hfill
    \begin{subfigure}[b]{0.24\textwidth}
        \centering
        \includegraphics[width=\linewidth]{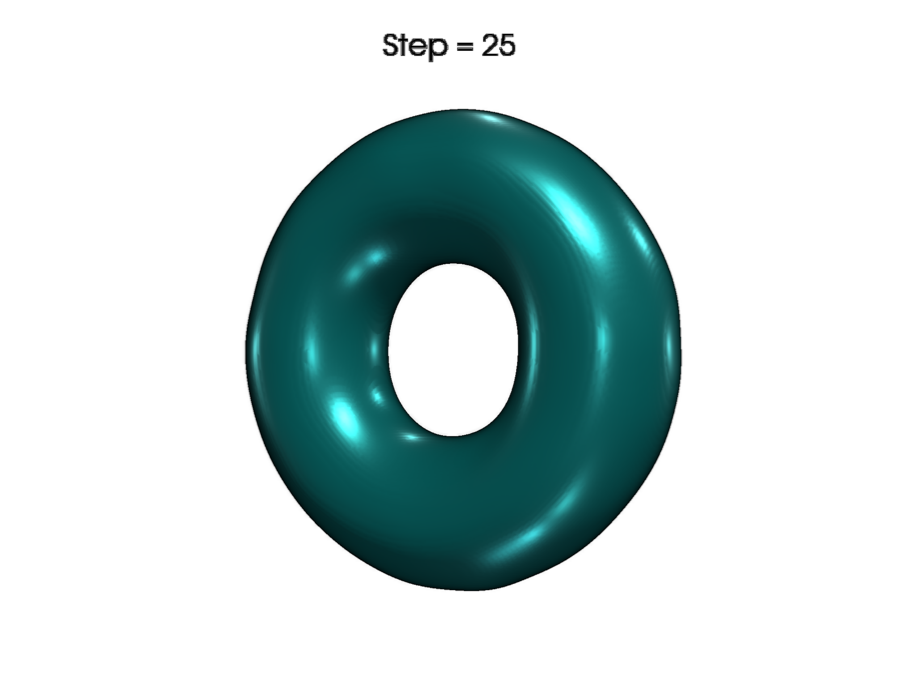}
    \end{subfigure}
        \begin{subfigure}[b]{0.24\textwidth}
        \rule{\textwidth}{0pt}
    \end{subfigure}
    \hfill
    \begin{subfigure}[b]{0.24\textwidth}
        \centering
        \includegraphics[width=\linewidth]{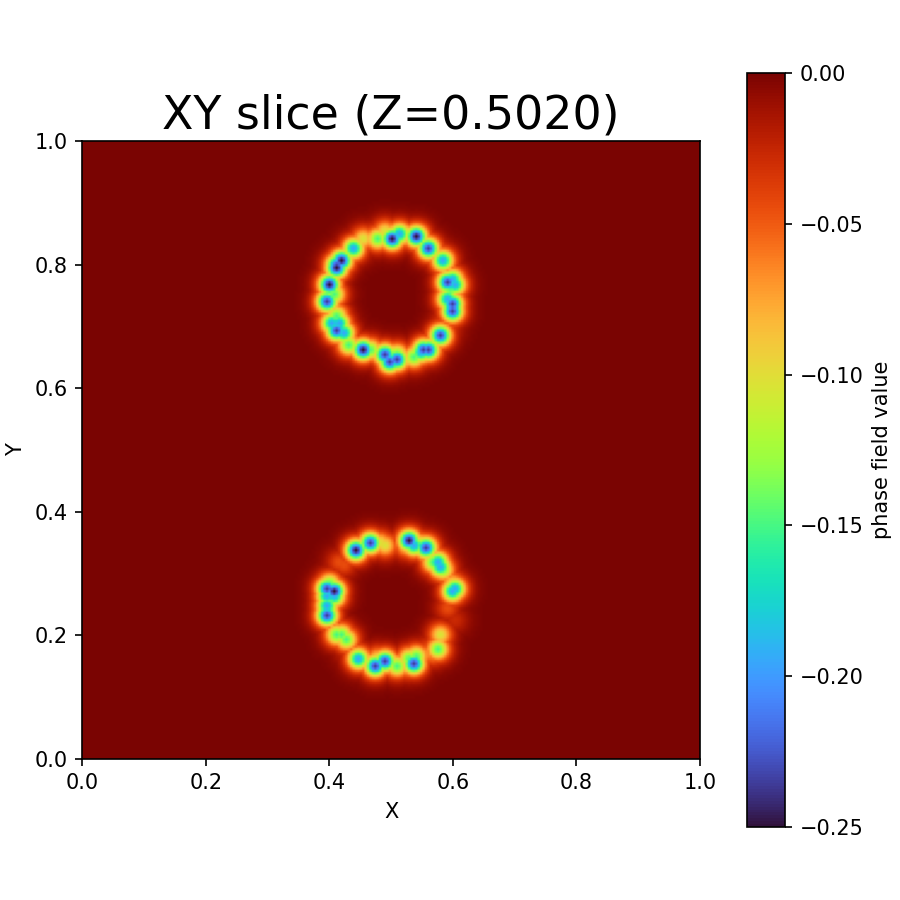}
    \end{subfigure}
    \hfill
    \begin{subfigure}[b]{0.24\textwidth}
        \centering
        \includegraphics[width=\linewidth]{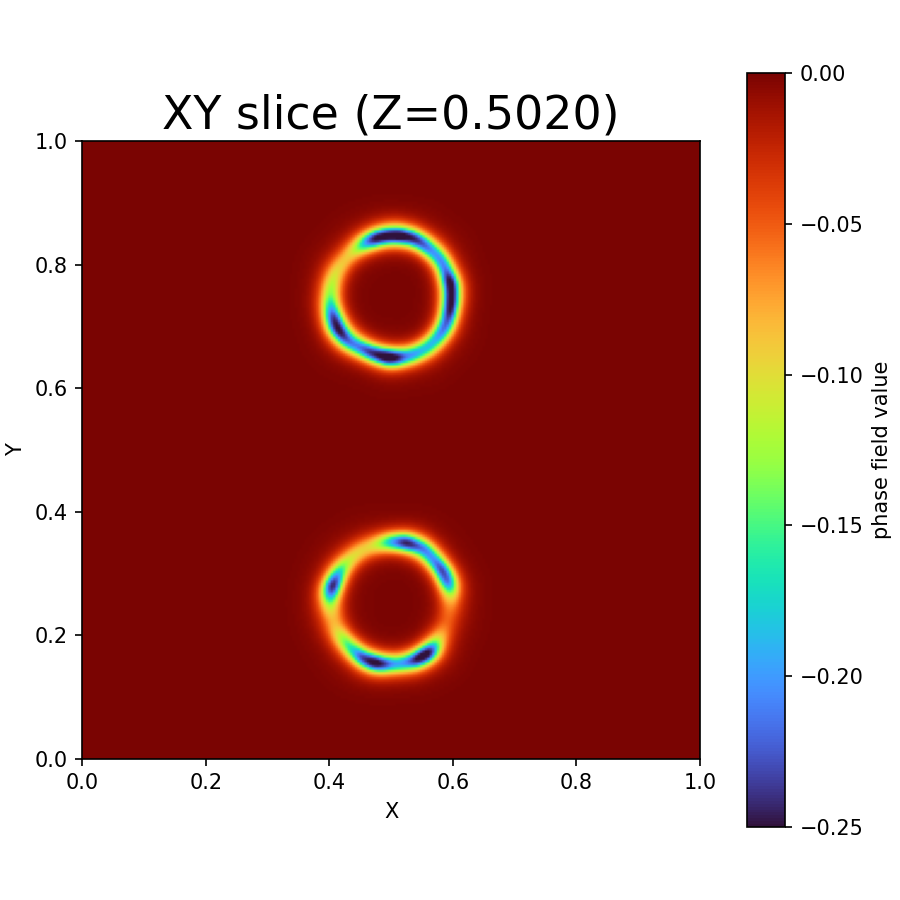}
    \end{subfigure}
    \hfill
    \begin{subfigure}[b]{0.24\textwidth}
        \centering
        \includegraphics[width=\linewidth]{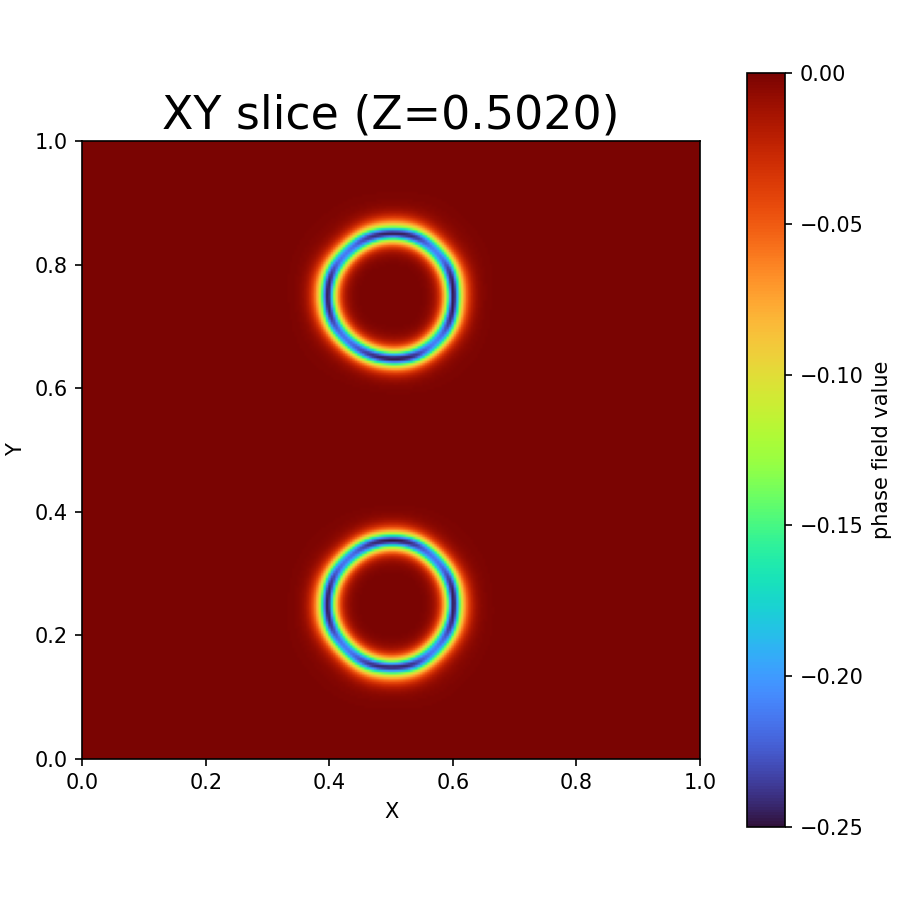}
    \end{subfigure}
    \caption{Reconstruction of a torus from $2000$ points using the free 4-parallel model $\wnnprim{\theta}{2}$ generated from the rotationally invariant model $\wpnnprim{\theta}{2}$ in 2D.
    The number of iterations of \eqref{eq:reconstruction_iteration} is $20$.
    Relaxation and visualization parameters are $\alpha = 1.0$ and $\delta = 0.15$.
    First row: isosurface of the phase-field function at each time step; second row: the heatmap of the phase-field function on the plane $Z = 0.5$.
    First column: initial point cloud; second column: iteration $0$; third column: iteration $1$; and fourth column: iteration $25$.}
    \label{fig:reconstruction_torus_qprim_3d_fulltransfer_nn}
\end{figure}

\begin{figure}[H]
    \centering
    \begin{subfigure}[b]{0.24\textwidth}
        \centering
        \includegraphics[width=\linewidth]{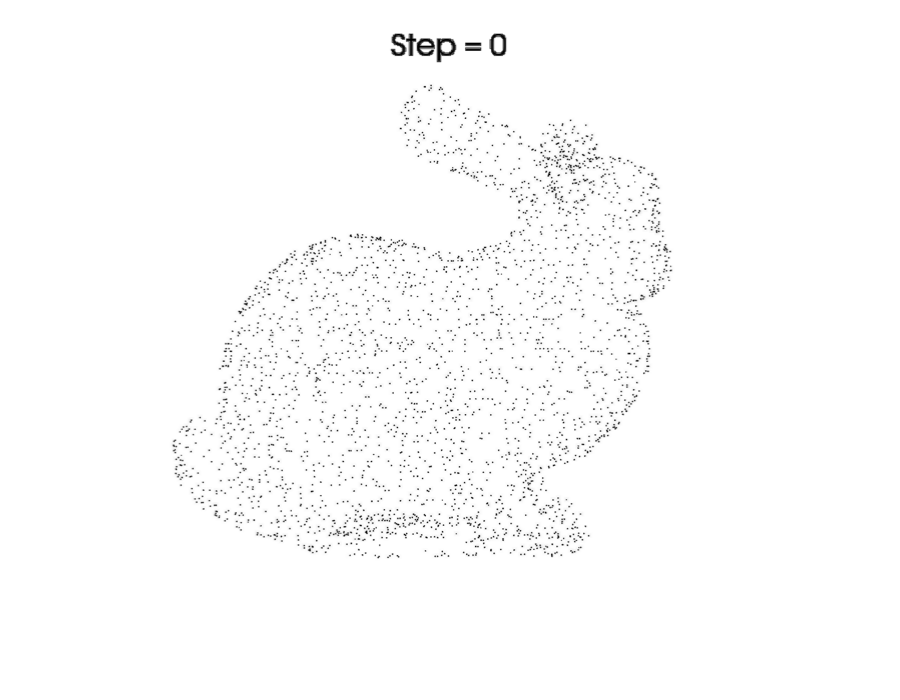}
    \end{subfigure}
    \hfill
    \begin{subfigure}[b]{0.24\textwidth}
        \centering
        \includegraphics[width=\linewidth]{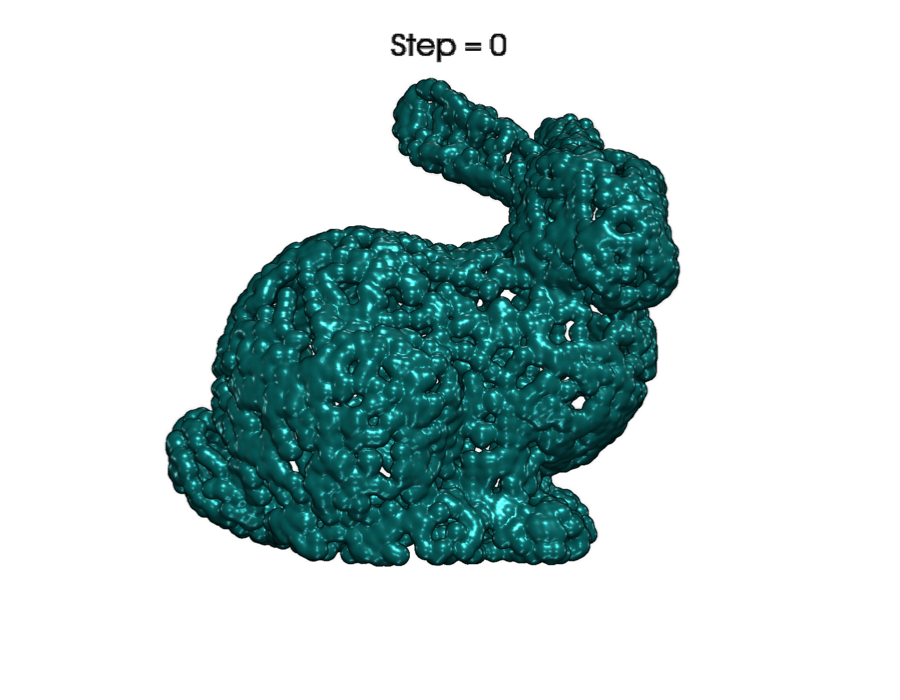}
    \end{subfigure}
    \hfill
    \begin{subfigure}[b]{0.24\textwidth}
        \centering
        \includegraphics[width=\linewidth]{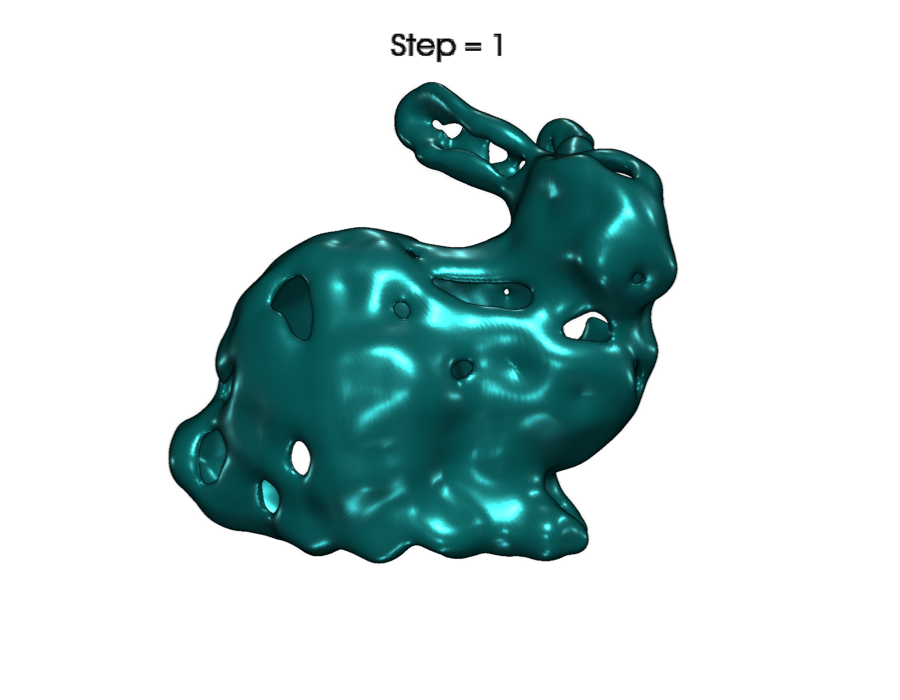}
    \end{subfigure}
    \hfill
    \begin{subfigure}[b]{0.24\textwidth}
        \centering
        \includegraphics[width=\linewidth]{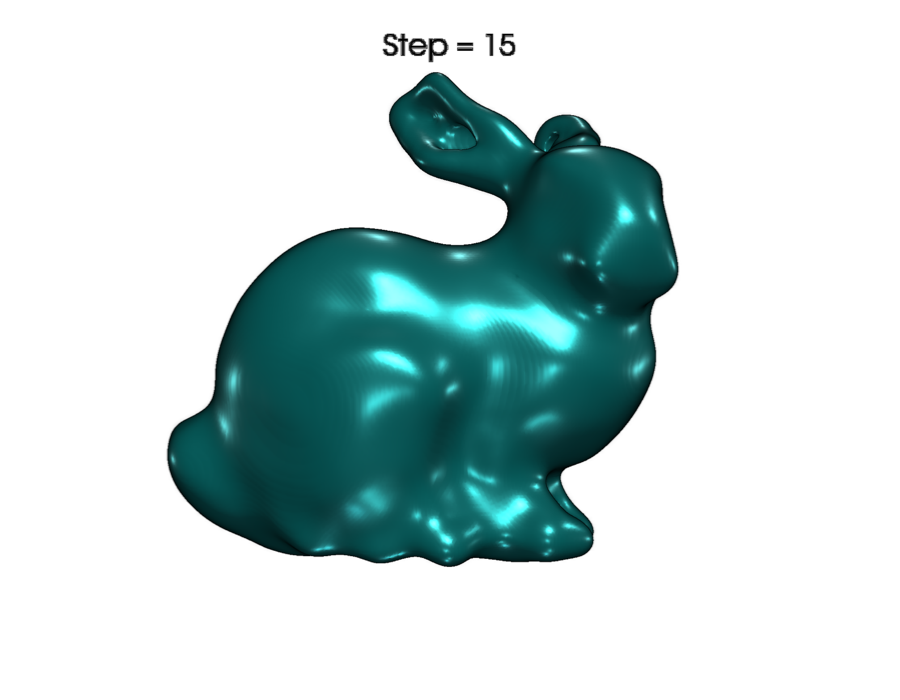}
    \end{subfigure}
        \begin{subfigure}[b]{0.24\textwidth}
        \rule{\textwidth}{0pt}
    \end{subfigure}
    \hfill
    \begin{subfigure}[b]{0.24\textwidth}
        \centering
        \includegraphics[width=\linewidth]{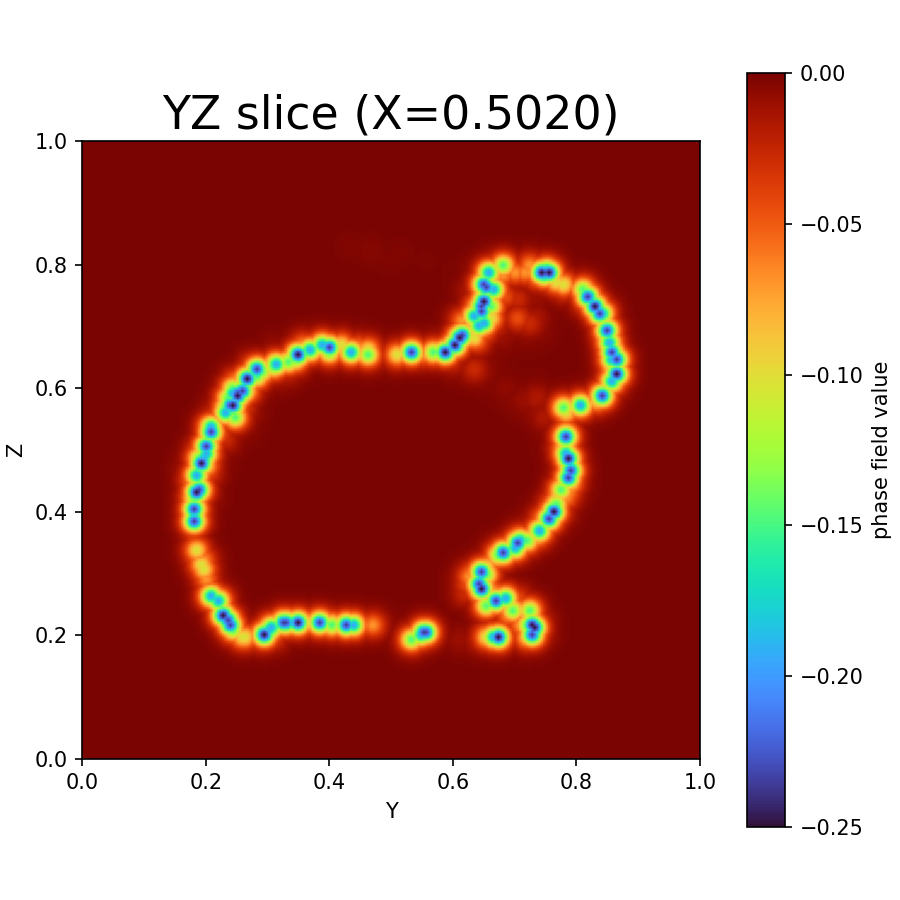}
    \end{subfigure}
    \hfill
    \begin{subfigure}[b]{0.24\textwidth}
        \centering
        \includegraphics[width=\linewidth]{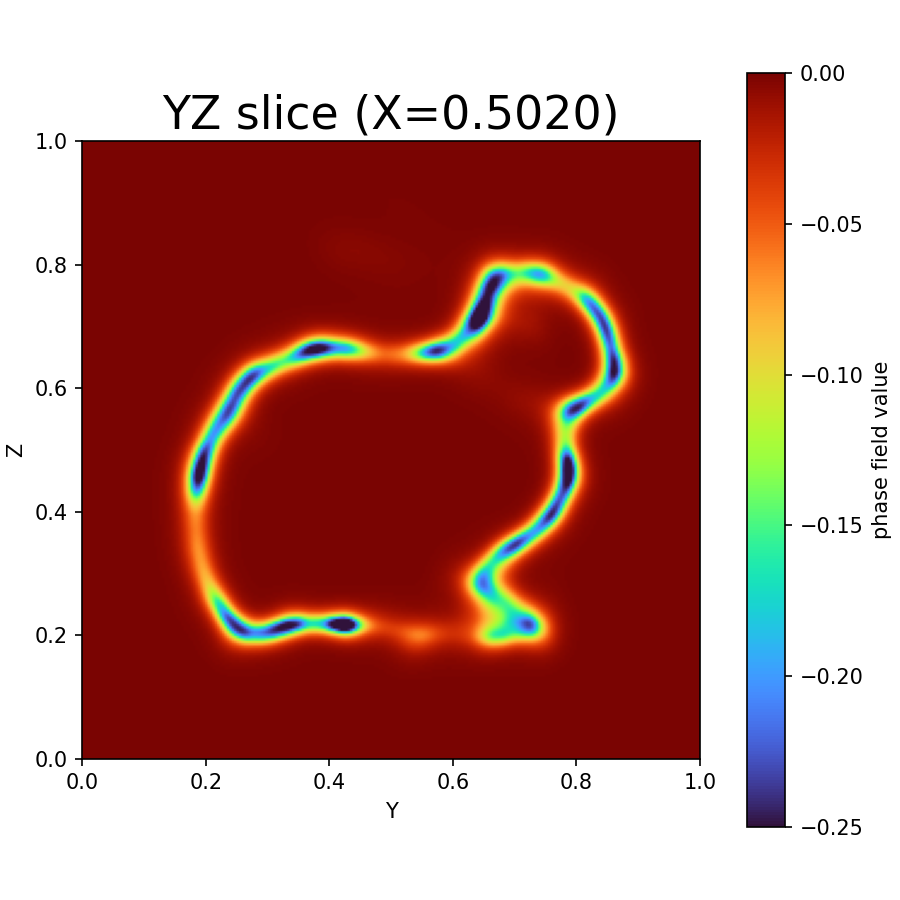}
    \end{subfigure}
    \hfill
    \begin{subfigure}[b]{0.24\textwidth}
        \centering
        \includegraphics[width=\linewidth]{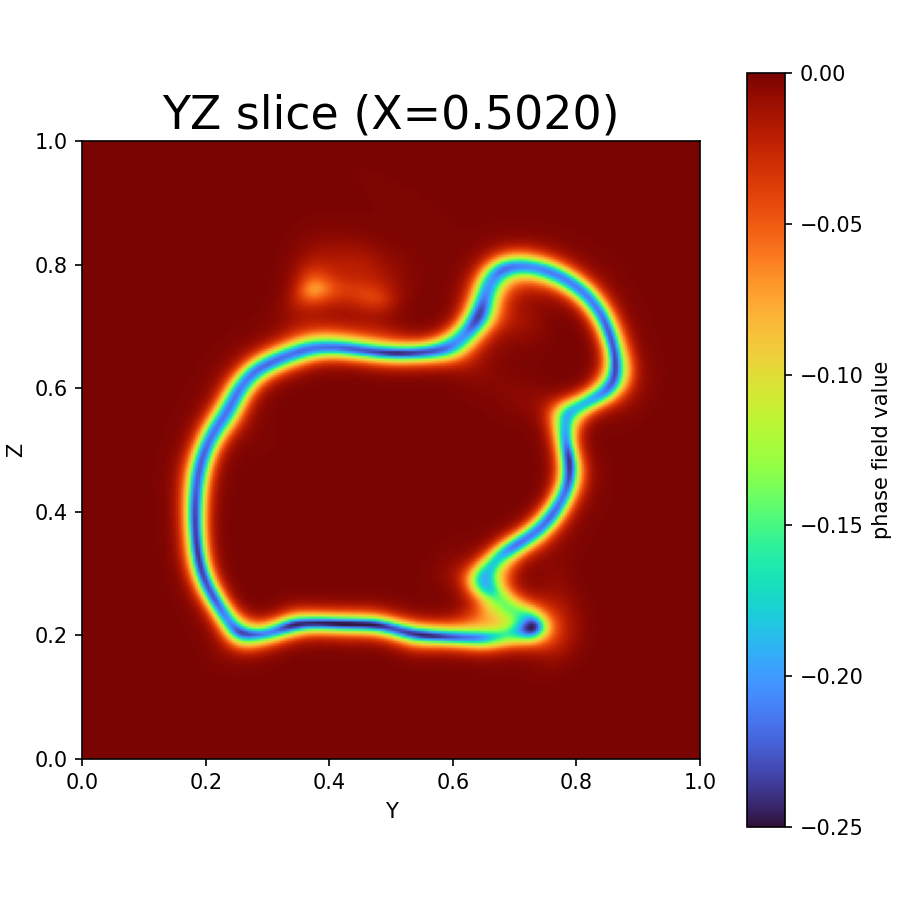}
    \end{subfigure}
    \caption{
    Reconstruction of a bunny from a point cloud ($3000$ points) with initialization by the free model initialized by $\wpnnprim{\theta}{2}$ in 2D.
    The number of iterations of \eqref{eq:reconstruction_iteration} is $14$.
    $\alpha = 1.0$ and $\delta = 0.15$.
    First row: isosurface of the phase-field function at each time step; second row: the heatmap of the phase-field function on the plane $X = 0.5$.}
    \label{fig:reconstruction_bunny_qprim_3d}
\end{figure}

To confirm the robustness of the method,
we now use it to reconstruct a more complex surface from a point cloud taken from the \textit{Stanford 3D Scanning Repository}~\cite{stanford3dscanrep},
a standard benchmark for this task. The sample points are chosen at random from the provided points and rescaled with respect to the centroid and the largest side of the bounding box to fit in the unit cube $(0,1)^3$.
The reconstruction, shown in  \Figure{fig:reconstruction_bunny_qprim_3d}, is smooth while preserving the global geometry of the original shape.
This is consistent with the variational nature of the Willmore flow, which penalizes the concentration of curvature and promotes globally regular surfaces.

\begin{figure}[H]
    \centering
    \begin{subfigure}[b]{0.24\textwidth}
        \centering
        \includegraphics[width=\linewidth]{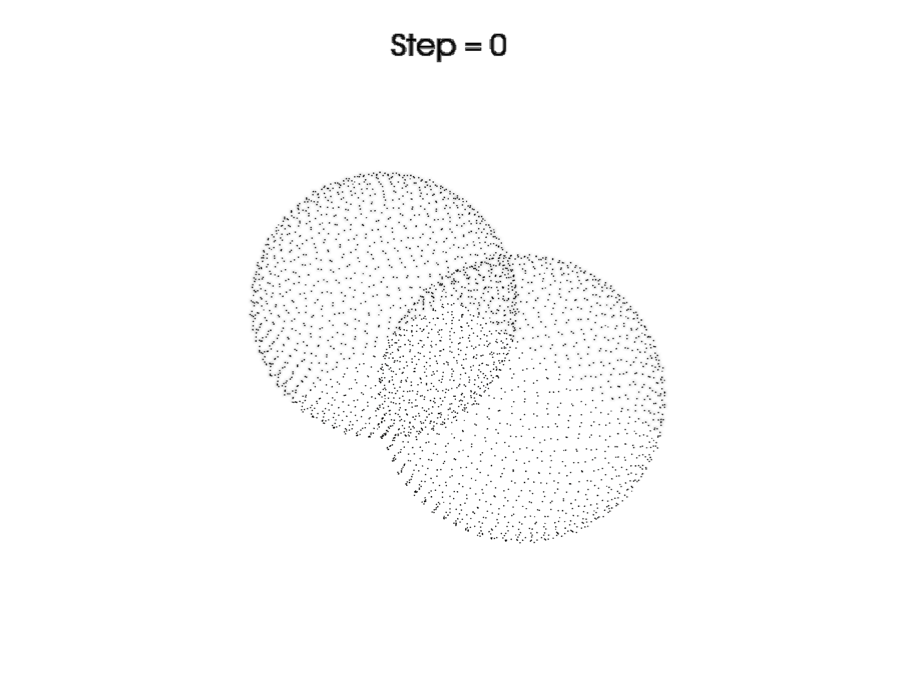}
    \end{subfigure}
    \hfill
    \begin{subfigure}[b]{0.24\textwidth}
        \centering
        \includegraphics[width=\linewidth]{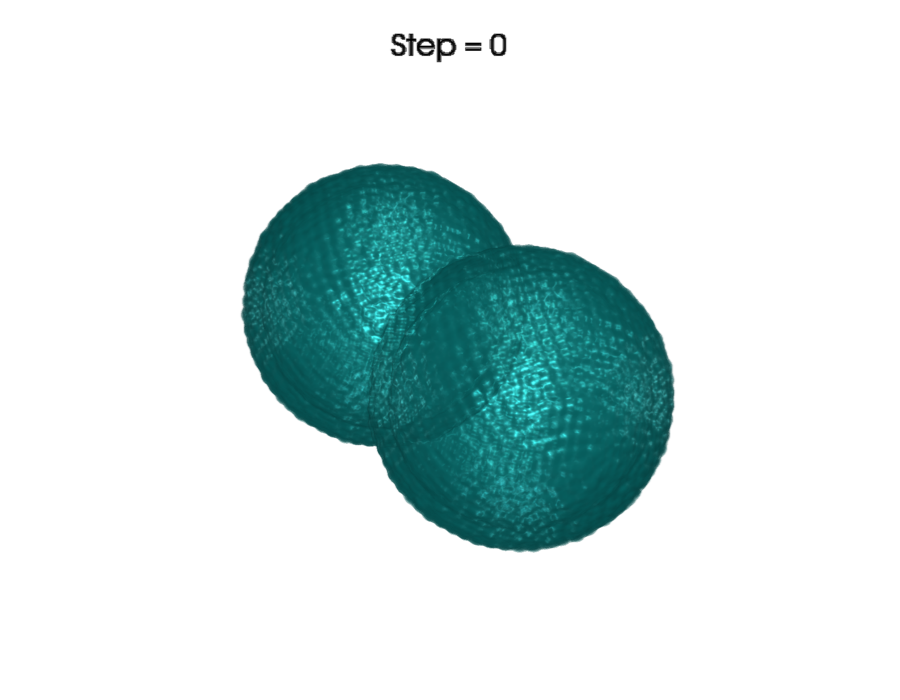}
    \end{subfigure}
    \hfill
    \begin{subfigure}[b]{0.24\textwidth}
        \centering
        \includegraphics[width=\linewidth]{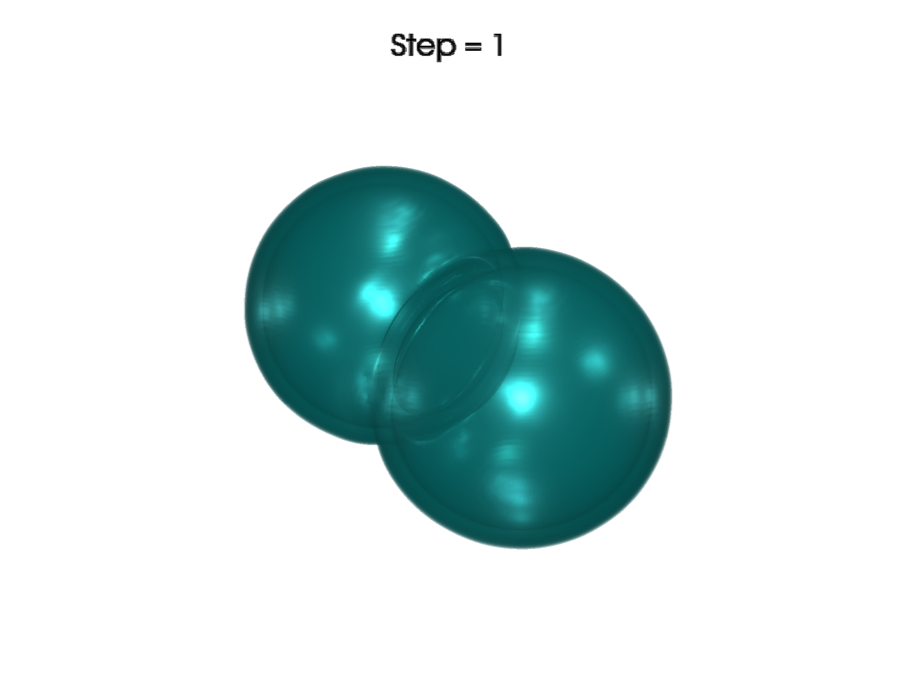}
    \end{subfigure}
    \hfill
    \begin{subfigure}[b]{0.24\textwidth}
        \centering
        \includegraphics[width=\linewidth]{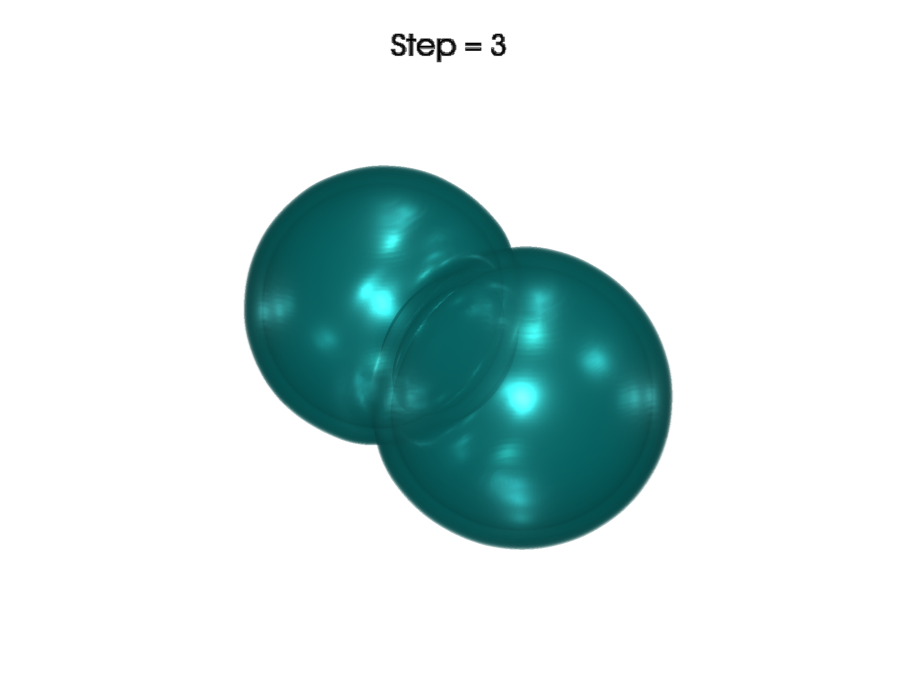}
    \end{subfigure}
    \begin{subfigure}[b]{0.24\textwidth}
    \end{subfigure}
    \hfill

     \hspace{4cm}
    \begin{subfigure}[b]{0.2\textwidth}
        \centering
        \includegraphics[width=\linewidth]{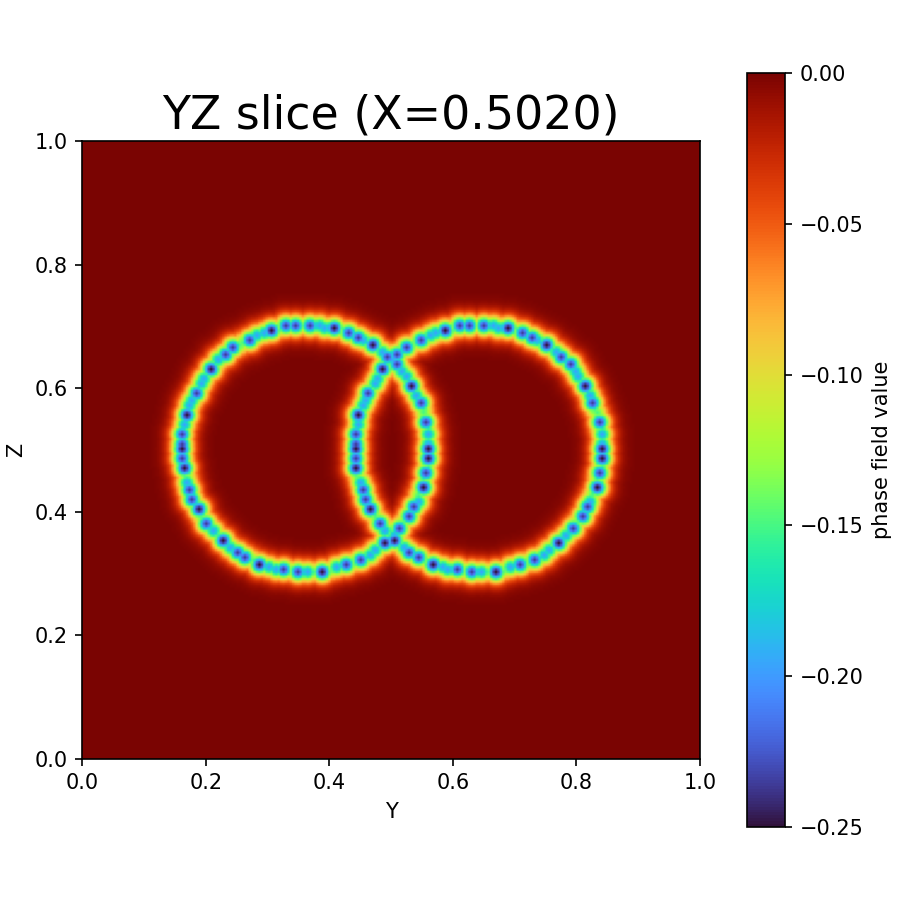}
    \end{subfigure}
    \hfill
    \begin{subfigure}[b]{0.2\textwidth}
        \centering
        \includegraphics[width=\linewidth]{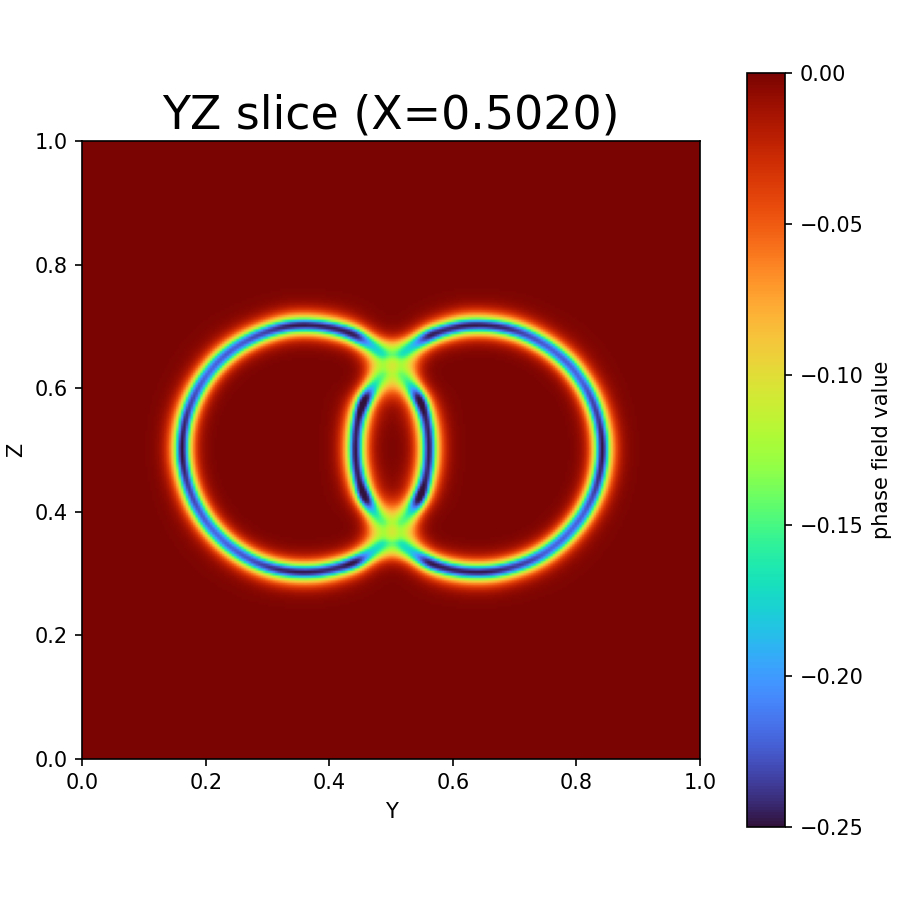}
    \end{subfigure}
    \hfill
    \begin{subfigure}[b]{0.2\textwidth}
        \centering
        \includegraphics[width=\linewidth]{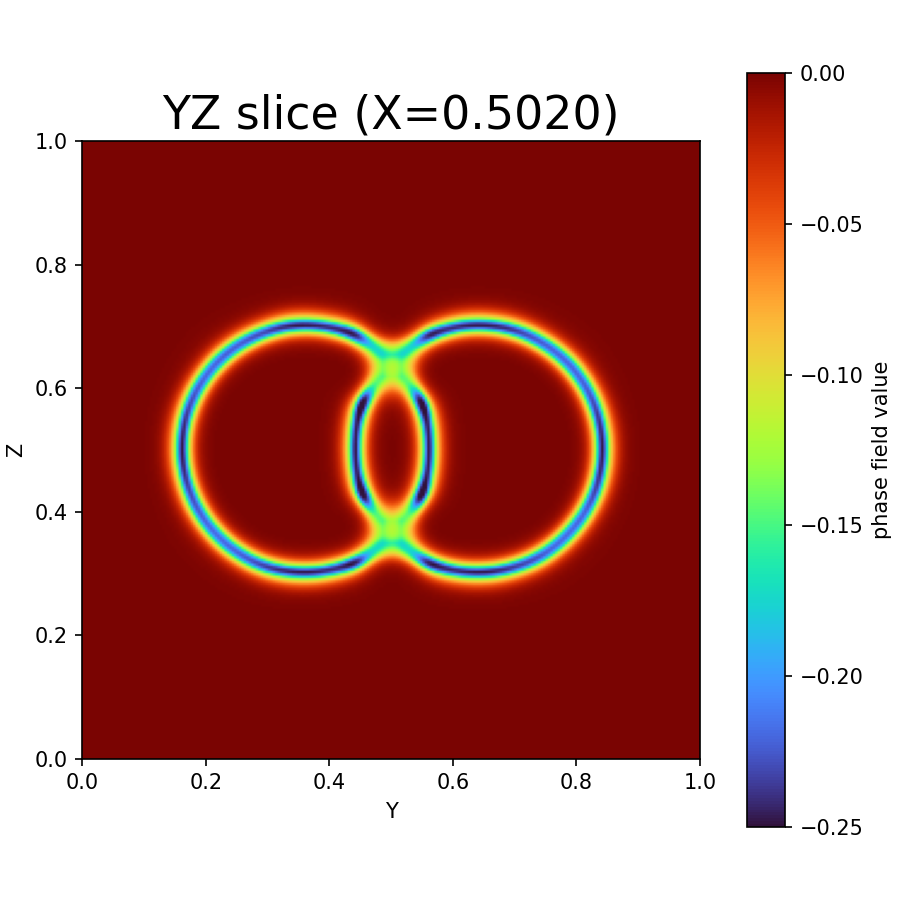}
    \end{subfigure}
    \caption{Reconstruction of a double bubble from $3000$ points using the free 4-parallel model $\wnnprim{\theta}{2}$ generated from the rotationally invariant model $\wpnnprim{\theta}{2}$ in 2D.
    The number of iterations of \eqref{eq:reconstruction_iteration} is $2$.
    Relaxation and visualization parameters are $\alpha = 1.0$ and $\delta = 0.15$.
    First row: isosurface of the phase field function; second row: slices of the phase field function at iterations $0$, $1$, and $3$ on the plane $X = 0.5$.
    }
    \label{fig:reconstruction_double_bubble_qprim_3d}
\end{figure}

\begin{figure}[H]
    \centering
    \begin{subfigure}[b]{0.24\textwidth}
        \centering
        \includegraphics[width=\linewidth]{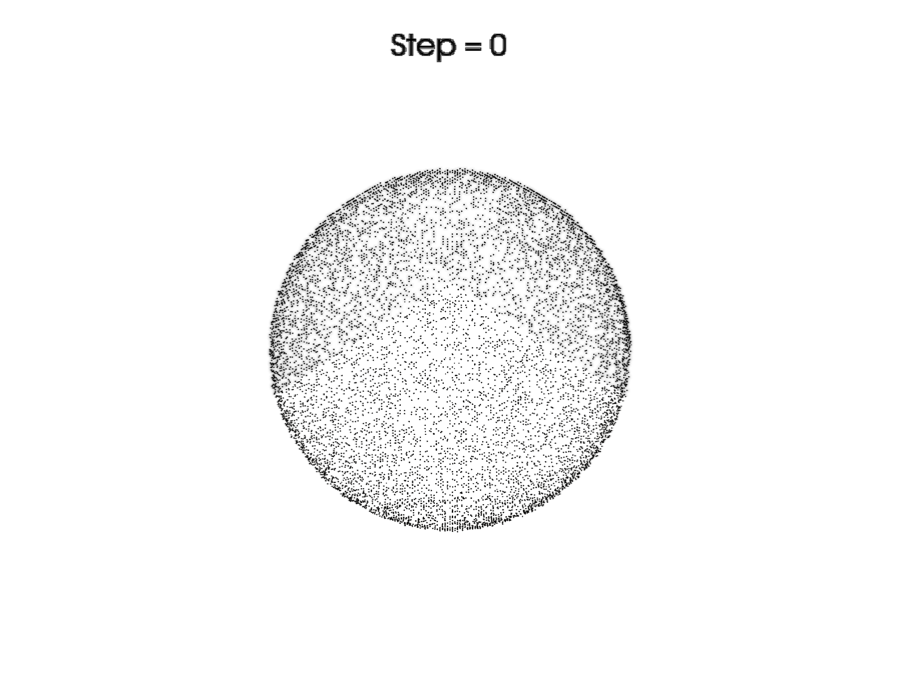}
    \end{subfigure}
    \hfill
    \begin{subfigure}[b]{0.24\textwidth}
        \centering
        \includegraphics[width=\linewidth]{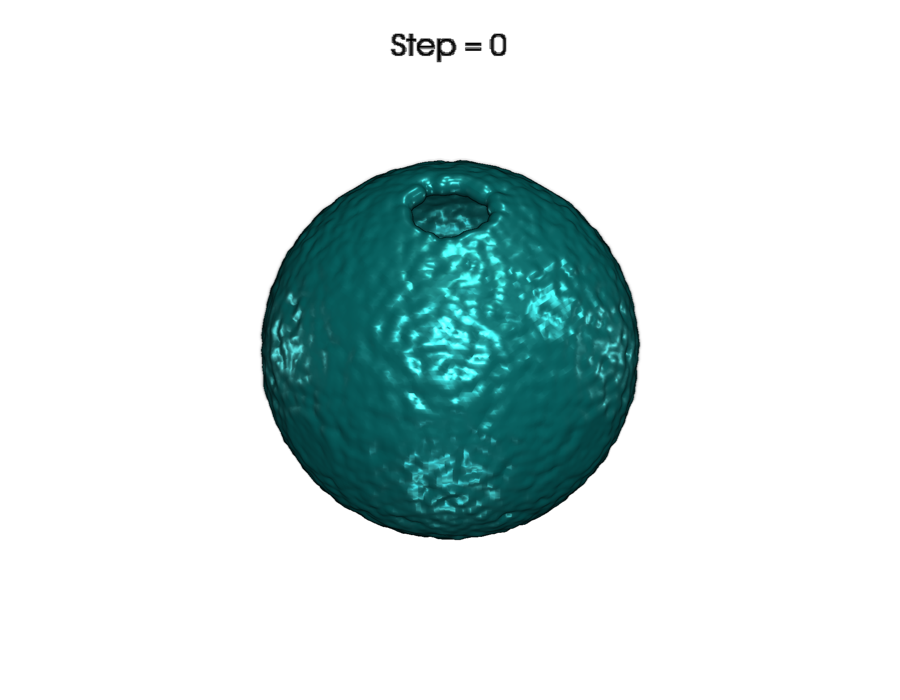}
    \end{subfigure}
    \hfill
    \begin{subfigure}[b]{0.24\textwidth}
        \centering
        \includegraphics[width=\linewidth]{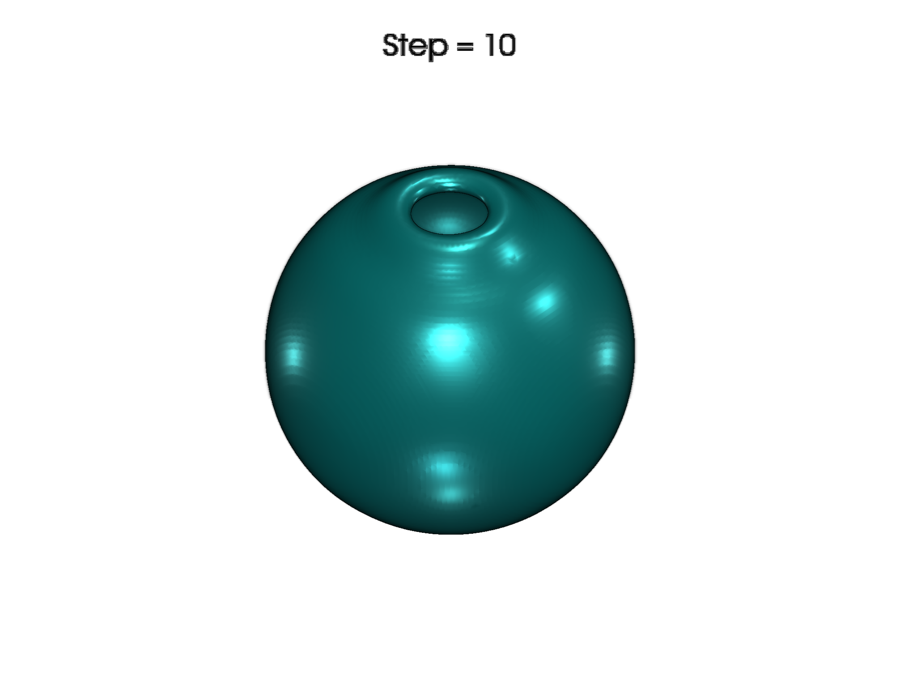}
    \end{subfigure}
    \hfill
    \begin{subfigure}[b]{0.24\textwidth}
        \centering
        \includegraphics[width=\linewidth]{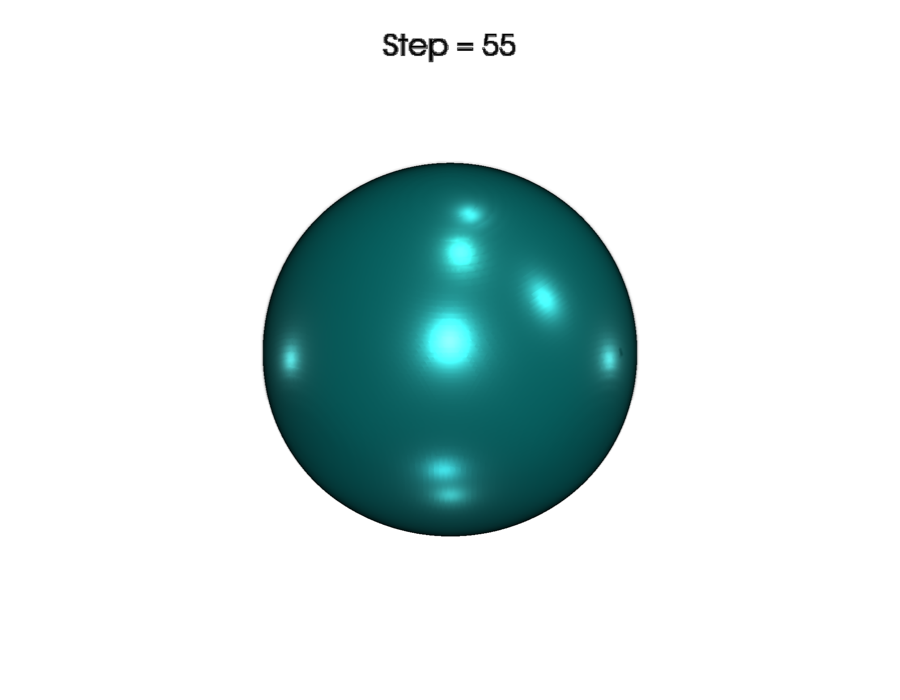}
    \end{subfigure}
       \begin{subfigure}[b]{0.24\textwidth}
        \rule{\textwidth}{0pt}
    \end{subfigure}
    \hfill
    \begin{subfigure}[b]{0.24\textwidth}
        \centering
        \includegraphics[width=\linewidth]{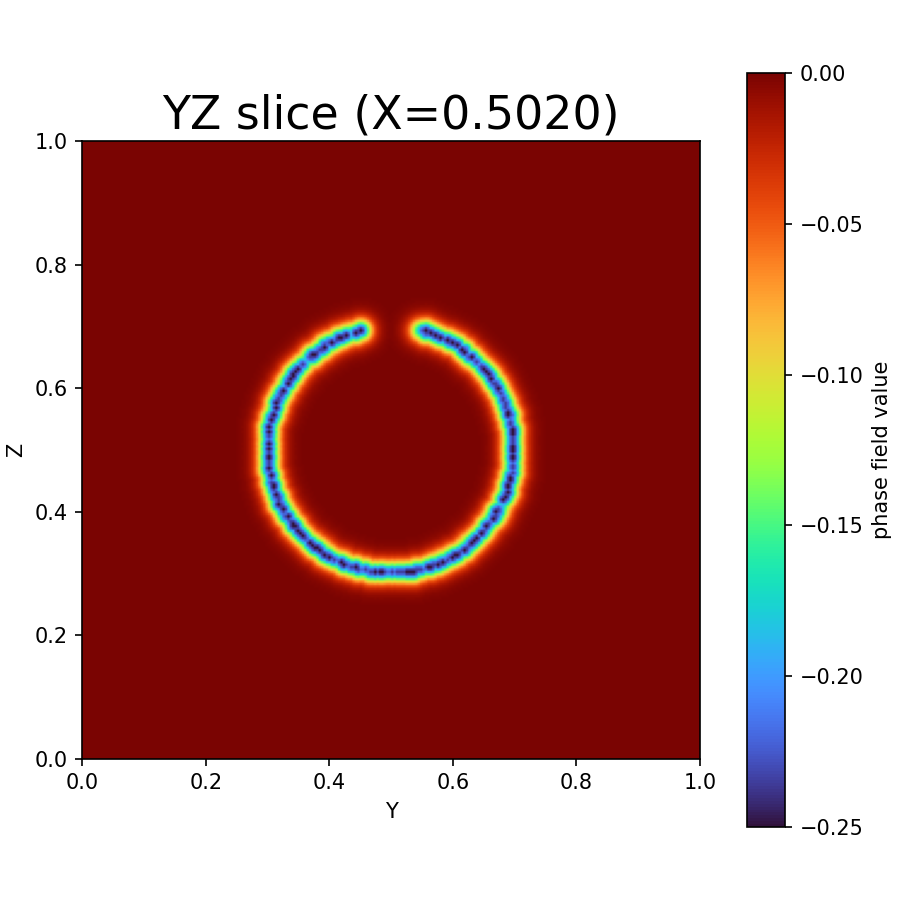}
    \end{subfigure}
    \hfill
    \begin{subfigure}[b]{0.24\textwidth}
        \centering
        \includegraphics[width=\linewidth]{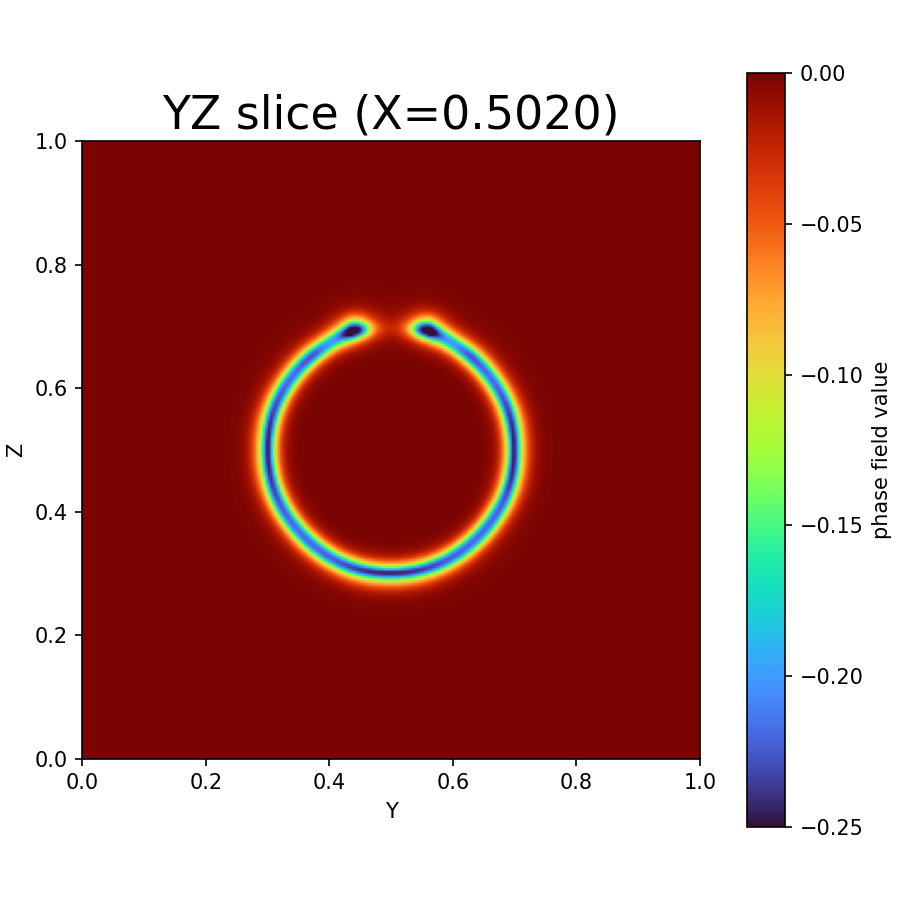}
    \end{subfigure}
    \hfill
    \begin{subfigure}[b]{0.24\textwidth}
        \centering
        \includegraphics[width=\linewidth]{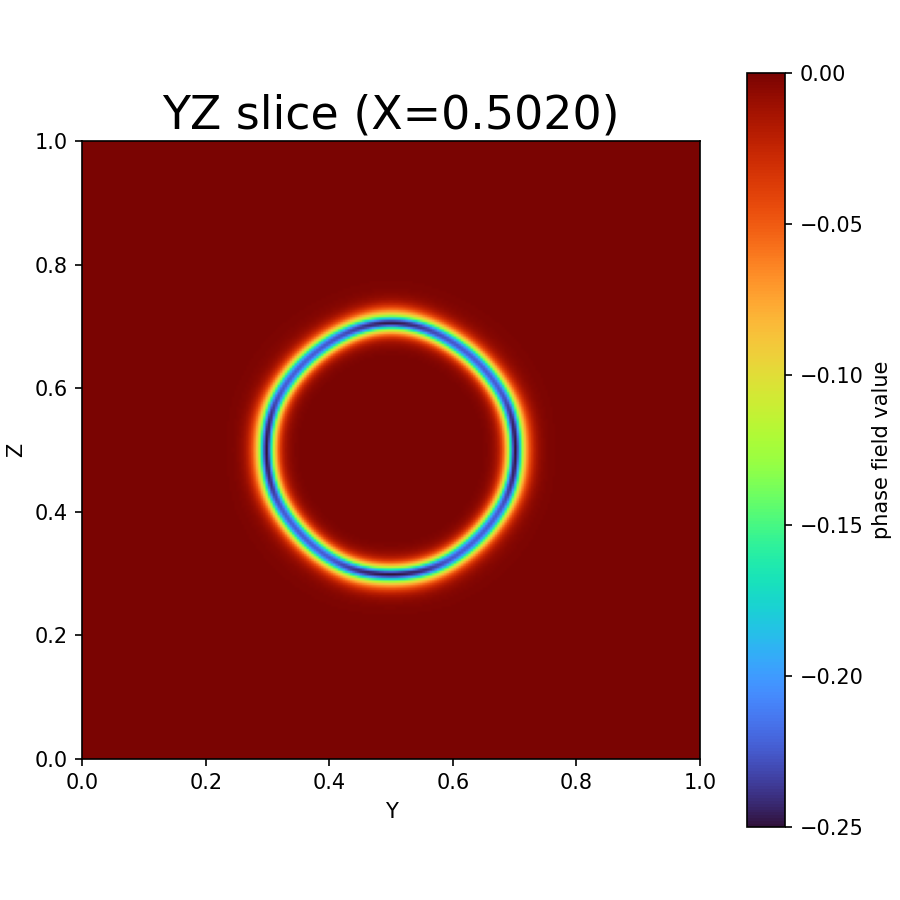}
    \end{subfigure}
    \caption{Reconstruction of a sphere with a large hole from $10000$ points using the free 4-parallel model $\wnnprim{\theta}{2}$ generated from the rotationally invariant model $\wpnnprim{\theta}{2}$ in 2D.
    The number of iterations of \eqref{eq:reconstruction_iteration} is $49$.
    Relaxation and visualization parameters are $\alpha = 1.2$ and $\delta = 0.1$.
    First row: isosurface of the phase-field function at each time step; second row: the heatmap of the phase-field function on the plane $X = 0.5$.
    First column: initial point cloud; second column: iteration $0$; third column: iteration $10$; and fourth column: iteration $55$.}
    \label{fig:reconstruction_sphere_with_hole_qprim_3d_large}
\end{figure}

The third example is a double-bubble configuration with a triple junction line, represented without an orientation. The reconstruction is shown in   \Figure{fig:reconstruction_double_bubble_qprim_3d} and demonstrates that our approach can handle this type of configuration, without requiring a prescribed normal field at the junction.

A further advantage of using the Willmore energy for reconstruction is its ability to fill holes, addressing the problem of surface inpainting. In this experiment, the learned regularization closes the missing region. This is an observed behavior of the reconstruction procedure; no boundary-energy asymptotic or general hole-filling guarantee is established here. The fourth example, shown in \Figure{fig:reconstruction_sphere_with_hole_qprim_3d_large}, illustrates this with a point cloud covering only part of a sphere.

\section*{Acknowledgments}
The authors acknowledge support from the "France 2030" program under grants ANR-23-PEIA-0004 (PEPR IA, PDE-AI) and from the French National Research Agency (ANR) under grant ANR-24-CE40-2216 (STOIQUES). 
\par Part of this work was also supported by the LABEX MILYON (ANR-10-LABX-0070) of Université de Lyon,
within the program "Investissements d’Avenir" (ANR-11-IDEX- 0007) operated by the French National Research Agency (ANR).

\section*{Declaration of generative AI use}
All scientific ideas, methodological developments, analyses, and interpretations presented in this work were developed by the authors. AI tools were used only to assist with language editing and presentation of the manuscript, as well as with the generation and revision of parts of the computational code. 

\section*{Code and data availability}
Our implementation builds upon the \texttt{nnpf} package
\cite{BretinDenisMasnouTerii2022}, built on \emph{PyTorch Lightning} \cite{Falcon_PyTorch_Lightning_2019},
extending its existing model classes and incorporating custom data-generation procedures for the numerical experiments presented in this paper.
The \texttt{nnpf} is available at \url{https://github.com/PhaseFieldICJ/nnpf}.

\bibliographystyle{abbrv}
\bibliography{cite}

\end{document}